\RequirePackage{fix-cm}
\documentclass{article}
\usepackage{arxiv}
\usepackage{lipsum}
\usepackage{placeins}
\usepackage[italicComments=true,indLines=true,noEnd=false]{algpseudocodex}
\usepackage{float}
\usepackage[english]{babel}
\usepackage[utf8]{inputenc} 
\usepackage[T1]{fontenc}    
\usepackage{hyperref}       
\usepackage{url}            
\usepackage{booktabs}       
\usepackage{stackrel,amssymb}
\usepackage{amsfonts,amsmath,amstext,amsbsy,amsthm}
\usepackage{nicefrac}       
\usepackage[version=4]{mhchem}
\usepackage[nopatch=eqnum]{microtype}      
\usepackage{graphicx}
\usepackage[numbers,sort&compress]{natbib}
\usepackage{doi}

\usepackage{xcolor}
\usepackage{enumitem}
\usepackage{caption}
\usepackage[title]{appendix}
\usepackage[capitalise,nameinlink]{cleveref}

\usepackage[bottom]{footmisc}
\usepackage{multirow}
\usepackage{subfigure}

\usepackage{algorithm}
\usepackage{algpseudocodex}

\title{A physics-informed neural network  method for the approximation of slow invariant manifolds for the general class of stiff systems of ODEs}

\author{
\textbf{Dimitrios G. Patsatzis\textcolor{blue}{$^{1}$}, Lucia Russo\textcolor{blue}{$^{2}$}, Constantinos Siettos\textcolor{blue}{$^{3,}$}\thanks{Corresponding author, email: \texttt{constantinos.siettos@unina.it}}}
{}\\
\textcolor{blue}{$^{(1)}$}Modelling Engineering Risk and Complexity, \emph{Scuola Superiore Meridionale}, Naples 80138, Italy \\
\textcolor{blue}{$^{(2)}$}Institute of Science and Technology for Energy and Sustainable Mobility, \emph{Consiglio Nazionale}\\ \emph{ delle Ricerche}, Naples 80125, Italy\\
\textcolor{blue}{$^{(3)}$}Dipartimento di Matematica e Applicazioni ‘‘Renato Caccioppoli", \emph{Universit\`a degli Studi di Napoli}\\ \emph{Federico II}, Naples 80126, Italy\\
}

\begin{document}
\maketitle
\begin{abstract}
We present a physics-informed neural network (PINN) approach  for the discovery of slow invariant manifolds (SIMs), for the most general class of fast/slow dynamical systems of ODEs. 
In contrast to other machine learning (ML) approaches that construct reduced order black box surrogate models using simple regression, and/or require a priori knowledge of the fast and slow variables per se, our approach, simultaneously decomposes the vector field into fast and slow components and provides a functional of the underlying SIM in a closed form. The decomposition is achieved by finding a transformation of the state variables to the fast and slow ones, which enables the derivation of an explicit, in terms of fast variables, SIM functional.~The latter is obtained by solving a PDE corresponding to the invariance equation within the Geometric Singular Perturbation Theory (GSPT) using a single-layer feedforward neural network with symbolic differentiation. The performance of the proposed physics-informed ML framework is assessed via three benchmark problems, namely the Michaelis-Menten, the target mediated drug disposition (TMDD) reaction mechanism  and, a fully competitive substrate-inhibitor (fCSI) mechanism. We also provide a comparison with other GPST methods, namely the quasi steady state approximation (QSSA), the partial equilibrium approximation (PEA) and CSP with one and two iterations. We show that the proposed PINN scheme provides SIM approximations, of equivalent or even higher accuracy, than those provided by QSSA, PEA and CSP, especially close to the boundaries of the underlying SIMs. Finally, we note that for the TMDD and fCSI problems, the CSP with two iterations (but also QSSA and PEA for the latter problem) does not provide the SIM approximations in an explicit closed form. In these cases, the computational time required for training the proposed PINN is comparable to the one required for the CSP to numerically approximate the SIM with Newton-iterations.
\end{abstract}


\section{Introduction}

The dynamics of complex systems, often deploy across multiple scales in time and space.~To enable comprehension and systematic numerical analysis, design of controllers and facilitate uncertainty quantification, simplification is required, provided by reduced order models (ROMs).~The construction of ROMs assumes that the effective dynamics evolves on low-dimensional topological spaces, (differentiable) manifolds.~For systems of ODEs/DAEs/PDEs with multiple timescales, the dynamics evolve on a fast timescale towards the neighborhood of an (attractive) manifold, and then evolve on it with a much slower timescale \citep{fenichel1979geometric,kaper1999introduction,jones1995geometric,gear2005projecting}.~Due to the latter slower dynamics, the emergent low-dimensional topological space is referred as Slow Invariant Manifold (SIM).~The approximation of SIMs is crucial for providing physical understanding, and for constructing ``parsimonious'' ROMs that are free of fast dynamics \citep{kuehn2015multiple,zagaris2004analysis}; the latter enabling accurate and robust to uncertainties reconstruction \citep{valorani2001explicit,goussis2012quasi}.

The approximation of SIMs is traditionally obtained in the context of Geometric Singular Perturbation Theory (GSPT) \citep{levinson1950small,tikhonov1952systems,vasil1963asymptotic,fenichel1979geometric,jones1995geometric}, and its predecessor, Singular Perturbation Theory (SPT); for a detailed review see \citep{kuehn2015multiple,verhulst2005methods}.~These theories were originally developed for singularly perturbed systems of ODEs, i.e., stiff systems of ODEs exhibiting fast-slow dynamics with an explicit timescale separation. In these systems the gap between fast and slow dynamics is explicitly expressed by a ``sufficiently'' small, constant in time, perturbation parameter, usually denoted by $\epsilon$.~The determination of $\epsilon$ naturally deems the system's variables as fast and slow (the variables that are associated with the fast/slow dynamics).~According to GSPT, the SIM is then locally described as a functional of the slow variables, $\epsilon$ and the fast variables.~For the approximation of the SIM functional, a variety of analytical and numerical model reduction methods have been developed by exploiting the invariance equation (IE), a PDE arising from the \textit{invariance} property of the SIM \citep{gorban2003method,gorban2004constructive,roussel1991geometry,fraser1988steady,gear2005projecting,zagaris2009analysis,ginoux2008slow}. 
The employment of the above methods is, however, restrictive, since usually there is no explicit knowledge about the timescale splitting, nor a-priori knowledge of the fast/slow variables. Discovering, such slow variables and computing SIMs for the construction of ROMs based on them, is a challenging and open research field. 

In this work, a more general class of stiff dynamical systems is considered, in which the fast/slow timescale decomposition exists, but it is not explicitly known; nor the decomposition of fast/slow variables and their number.~For these systems, some of the above methods \citep{gorban2003method,gorban2004constructive,gear2005projecting,zagaris2009analysis,ginoux2008slow} may provide a SIM approximation under intuition-based assumptions about the fast/slow variables.~However, more sophisticated methods have been developed to tackle particularly this general class of dynamical systems, such as the Computational Singular Perturbation (CSP) method \citep{lam1989understanding,goussis1992study,hadjinicolaou1998asymptotic,valorani2005higher,goussis2012quasi}, the Invariant Low-Dimensional Method (ILDM) \citep{maas1992simplifying,maas1992implementation,maas1994laminar} and the Tangential Stretching Rate (TSR) method  \citep{valorani2015dynamical,valorani2017tangential}.~These methods differ from the former in that they exploit the properties of the \textit{tangent} space, along which the fast/slow dynamics act, without requiring a priori knowledge of the fast/slow variables \citep{kuehn2015multiple,goussis2011model,kaper1999introduction,zagaris2004analysis,kaper2015geometry}.~In fact, they decompose the vector field into fast/slow components and compute the SIM approximations in terms of the original state variables.~Discovering the SIM as a functional of the state variables per se, is crucial for revealing the physical processes that lead the dynamical system onto the emergent SIM \citep{lam1994csp,goussis1992study,goussis2006model,patsatzis2019new}.~On the other hand, this frequently results in SIM approximations (either analytical or numerical) of implicit form \citep{goussis2006efficient,kuehn2015multiple}.~This feature introduces significant computational cost when the SIM expressions are utilized for the reconstruction of the ambient high-dimensional space (e.g., after the time integration of the ROM), since  additional numerical root-finding algorithms for obtaining explicit expressions of the SIM are required \citep{massias1999algorithm}.~In spite of the above drawback, CSP provides a systematic framework to discover SIM approximations resulting from traditional chemical kinetics model reduction techniques, such as the Quasi Steady-State Approximation (QSSA) and Partial Equilibrium Approximation (PEA) \citep{roussel1991geometry,bowen1963singular,rein1992partial,patsatzis2019new,zagaris2009analysis,goussis2012quasi}.~In addition, due to its numerical nature, CSP has been employed for approximating SIMs of high-dimensional stiff systems \citep{manias2016mechanism,manias2019investigation,tingas2021computational}.

In the framework of data-driven methods, a plethora of Machine Learning (ML) methodologies has been developed for  the construction of black or gray-box ROM surrogates of complex systems bypassing the explicit computation of SIMs.~This is usually performed in two steps. The first step involves the detection of the low-dimensional space using manifold learning algorithms \citep{balasubramanian2002isomap,bollt2007attractor,nadler2006diffusion,dsilva2016data,gallos2023data}, including the use of Proper orthogonal Decomposition (POD) \cite{sirisup2005equation,chen2021physics,koronaki2023partial} 
and the Koopman operator \cite{williams2015data,bollt2018matching,lusch2018deep,mauroy2020koopman,gallos2023data}. The second step is the construction of ROMs, using for example interpretable ML \cite{taprantzis1997fuzzy,siettos2002semiglobal,siettos2002truncated,tatar2020reduced}, Artificial Neural Networks (ANNs) \citep{lee2020coarse,chen2021physics,galaris2022numerical,fabiani2023tasks}, Gaussian Processes (GPs) \cite{wan2017reduced,lee2020coarse,sapsis2022optimal}, long short-term memory networks (LSTMs) \cite{wan2018data}, physics-informed  deep learning \cite{raissi2018deep,raissi2019physics,chen2021physics,karniadakis2021physics}, and DeepONet \cite{lu2019deeponet,demo2023deeponet}.~Autoencoders are also extensively used for performing both steps with a single training process \citep{chen2018molecular,vlachas2022multiscale,vijayarangan2024data,solera2024beta}.
Physics-informed ML methods have been also used to construct gray-box ROM surrogates from high-fidelity simulations, by finding parameters/closures of an effective system of PDEs or ODEs \cite{lee2022learning,galassi2022adaptive,kim2022fast,daryakenari2023ai}.
Recently, ``next generation" Equation-free (EF) \cite{kevrekidis2003equation} data-driven methodologies  have been proposed to build ROMs bypassing the need of constructing surrogate models, based on manifold learning and  k-nearest neighbors algorithm \cite{chin2022enabling} or geometric harmonics \cite{patsatzis2023data,evangelou2023double} for the construction of the map between the high and low dimensional spaces. 
The above ML methodologies, either the ones based on the construction of surrogate ML models or the ones within the EF framework,
construct ROMs, bypassing/avoiding the explicit computation of the underlying SIMs. A comparison of various ML approaches for the construction of ROMs in a data-driven way including deep convolutional autoencoders can be found in \cite{gruber2022comparison}; a comparison between ROM surrogates and ROMs based on the EF-approach can be found in \cite{fabiani2023tasks}.\par 
Most recently, in \cite{galassi2022adaptive}, the authors used an ANN to learn from data the projection basis
vectors produced by the implementation of the CSP method to define numerically
the local slow invariant manifold, thus constructing a ROM by projecting the vector field of a problem of stiff chemical kintetics onto the slow CSP basis vectors. In \citep{patsatzis2023slow}, we presented a physics-informed ML approach based on the GSPT framework for analytical derivation of SIMs of singularly perturbed systems in an explicit closed form. In particular, we used single-layer ANNs and random projection ANNs using symbolic differentiation to solve the IE PDE, that constraints the system dynamics to lie on the SIM. We demonstrated that the proposed physics-informed ML approach outperforms other classical GSPT-based methods, especially for relatively
large values of the perturbation parameter.
In \citep{serino2024intelligent}, the authors presented a data-driven ML approach for learning surrogate ROMs of singularly perturbed systems, by enforcing the existence of a trainable, attracting SIM as a hard constraint in a fast-slow ANN architecture. In both of the latter two approaches, it is assumed that there is an explicit time-splitting in the dynamics, i.e., that the system is a singularly perturbed system of ODEs, and, that the dimension of the attracting invariant manifold is known a-priori.\par  
Here, we make a step advancement, relaxing the two above restrictions, by proposing a physics-informed neural network (PINN) approach for learning explicit functionals of SIMs for the most general class of stiff dynamical systems,  for which the fast/slow timescale splitting as well as the dimension of the SIM are not explicitly known a-priori. 
Furthermore, unlike the other ML methods proposed so far, which create surrogate models via regression, our approach offers a functional expression of the SIM that solves the PDE corresponding to the IE via PINNs.~However, since the fast and slow variables are not \textit{a priori} known, the proposed PINN approach, based on the CSP diagonstic criterion of the dimension of the SIM, it (i) finds the transformations that decompose the state variables into fast and slow components, and (ii) solves the IE for computing the SIM functional in an explicit closed form. For the determination of the dimension of the SIM that is to be approximated, we use the CSP criterion introduced in \cite{valorani2001explicit} which involves the eigendecomposition of the linearized system which, can be computed analytically at any point, as the vector field is assumed to be known within the physics-informed modelling framework.

The performance of the proposed PINN scheme was assessed via three benchmark problems, namely, the Michaelis-Menten (MM) enzyme reaction scheme, the pharmacokinetic/pharmacodynamic Target Mediated Drug Disposition (TMDD) mechanism, and, a fully competitive substrate-inhibitor (fCSI) mechanism.~For a straightforward comparison with the GPST methods, we also derived the SIM approximations provided on the basis of QSSA, PEA and CSP with one and two iterations for all the above problems; the second CSP iteration is required in the case where the assumption on the fast variables is erroneous \cite{goussis2012quasi}.

Before delving into details about the methodology, we highlight that, except from the interpolation advantage over other ML methods using regression, the proposed PINN method provides explicit SIM functionals that can be evaluated at any set of values of the state variables.~This is a significant advantage also over the analytical or numerical GSPT methods, since when the latter result to implicit SIM approximations, additional numerical root-finding algorithms are required to estimate point-by-point the values of the fast variables; see e.g., the QSSA, PEA and CSP implicit expressions of the fCSI mechanism in \cref{app:Inh_SIMs}, \cref{tb:InhSt_convSIMs,tb:InhTr_convSIMs}.

The paper is structured as follows:~Firstly, we present the proposed PINN methodology, providing the analytical expressions required for the learning process, alongside discussions on its implementation and numerical accuracy.~Following this, we present the three benchmark problems - the MM, TMDD, and fCSI mechanisms - and provide the various SIM approximations constructed on the basis of GSPT.~Then, we present the numerical results, assessing and comparing the performance of the proposed PINN approach with respect to the other GSPT approximations.~Finally, we summarize the advantages and limitations of the proposed PINN approach in comparison to GSPT methods, and suggest future research directions.

\section{Methodology}
\label{sec:Meth}

For the completeness of the presentation, we begin by very briefly introducing fundamental concepts of Perturbation Theory that are pertinent to the proposed ML methodology.\par
Let us consider the autonomous system of ODEs:
\begin{equation}
    \dfrac{d \mathbf{z}}{dt} = \mathbf{F}(\mathbf{z}), 
    \qquad \mathbf{z}(t_0) = \mathbf{z}_0,
    \label{eq:gen}
\end{equation}
where $\mathbf{z}\in \mathbb{R}^N$ denotes the $N$-dim. state vector of variables, $\mathbf{F}:\mathbb{R}^N \rightarrow \mathbb{R}^N$ is a C$^r$ smooth vector field and $\mathbf{z}_0\in \mathbb{R}^N$ is the initial condition of the initial value problem (IVP) under study.~We assume that the components of the vector field $\mathbf{F}$ result in stiff dynamics.\par
Traditionally, such systems are tackled by  the Geometric Singular Perturbation Theory (GSPT) \citep{levinson1950small,tikhonov1952systems,vasil1963asymptotic,fenichel1979geometric}; 
a detailed summary can be found in \citep{kaper1999introduction,o1991singular,kuehn2015multiple,verhulst2005methods,kevorkian2013perturbation,jones1995geometric,hek2010geometric}.~GSPT was originally developed for multi-scale dynamical systems with explicit timescale splitting; i.e., systems where the gap between the characteristic fast, say $\tau_f$, and slow, say $\tau_s$, timescales, can be determined explicitly by a small perturbation parameter $\epsilon=\tau_f/\tau_s\ll1$.~In such cases, the fast (slow) subsystems can be extracted for describing the fast (slow) dynamics of the original system in the rescaled time $t_f=t/\tau_f$ ($t_s=t/\tau_s$).~In particular, it is usually assumed that $M$ linearly transformed 
variables of $\mathbf{z}\in\mathbb{R}^N$, say $\mathbf{\hat{x}}\in \mathbb{R}^M$, are associated with the fast dynamics, while the remaining $N-M$ ones, say $\mathbf{\hat{y}}\in \mathbb{R}^{N-M}$, are associated with the slow dynamics.~Following this decomposition, the original system in Eq.~\eqref{eq:gen} can be transformed to its slow subsystem form (a singularly perturbed system):
\begin{equation}
    \epsilon \dfrac{d \mathbf{\hat{x}}}{dt_s} = \mathbf{\hat{f}}(\mathbf{\hat{x}},\mathbf{\hat{y}},\epsilon), \qquad \dfrac{d \mathbf{\hat{y}}}{dt_s} = \mathbf{\hat{g}}(\mathbf{\hat{x}},\mathbf{\hat{y}},\epsilon),
    \label{eq:FSvecD}
\end{equation}
where $\mathbf{\hat{f}}: \mathbb{R}^N\times I\rightarrow  \mathbb{R}^M$ and $\mathbf{\hat{g}}: \mathbb{R}^N\times I\rightarrow  \mathbb{R}^{N-M}$ are the vector fields of the fast $\mathbf{\hat{x}}$ and slow $\mathbf{\hat{y}}$ variables and $I\subset \mathbb{R}$ is an interval containing $\epsilon=0$.~Similarly to the slow (singularly perturbed) subsystem in Eq.~\eqref{eq:FSvecD}, a fast subsystem can be also extracted for describing the fast dynamics of the original system.

In the view of GSPT, the system in Eq.~\eqref{eq:FSvecD} admits a \emph{critical manifold} $\mathcal{M}_c=\{ (\mathbf{\hat{x}},\mathbf{\hat{y}})\in \mathbb{R}^N: \mathbf{\hat{f}}(\mathbf{\hat{x}},\mathbf{\hat{y}},0)=0 \}$ in the limit $\epsilon\rightarrow 0$.~All the subsets $\mathcal{M}_0\subset \mathcal{M}_c$ are \emph{normally hyperbolic} if the matrix $\nabla_{\hat{x}} \mathbf{\hat{f}}(\mathbf{\hat{x}},\mathbf{\hat{y}},0)$ has no eigenvalues with zero real part for every $(\mathbf{\hat{x}},\mathbf{\hat{y}}) \in \mathcal{M}_0$.~In the cases of small, yet non-zero, $0<\epsilon \ll1$, according to the Fenichel-Tikhonov theorem \citep{tikhonov1952systems,fenichel1979geometric}, there exists a \emph{locally invariant} and \emph{normally hyperbolic manifold} $\mathcal{M}_\epsilon$ that (i) is diffeomorphic to $\mathcal{M}_0 \subset \mathcal{M}_c$, (ii) has the same stability properties with $\mathcal{M}_0$ w.r.t. the fast variables and (iii) the flow of the system on it converges to the slow flow of $\mathcal{M}_c$ as $\epsilon \rightarrow 0$.~The manifold $\mathcal{M}_\epsilon$, due to the latter property, is frequently called \emph{Slow Invariant Manifold} (SIM) and can be locally described as \citep{fenichel1979geometric,jones1995geometric}:
\begin{equation}
\mathcal{M}_\epsilon = \{ (\mathbf{\hat{x}},\mathbf{\hat{y}})\in \mathbb{R}^N: \mathbf{\hat{x}} = \mathbf{\hat{h}}(\mathbf{\hat{y}},\epsilon) \}.
\label{eq:SIM00}
\end{equation}
Substitution of the map $\mathbf{\hat{x}}=\mathbf{\hat{h}}(\mathbf{\hat{y}},\epsilon)$ in the differential equations of the slow variables $\mathbf{\hat{y}}$ in Eq.~\eqref{eq:FSvecD} yields the reduced order model (ROM) for the slow dynamics of the system; that is $d\mathbf{\hat{y}}/dt_s=\mathbf{\hat{g}}(\mathbf{\hat{h}}(\mathbf{\hat{y}},\epsilon),\mathbf{\hat{y}},\epsilon)$.

However, for the general class of stiff ODEs described by Eq.~\eqref{eq:gen} an explicit knowledge of the timescale splitting is usually not available.~Stiffness implies just the existence of $M$ fast timescales that are much faster than the characteristic slow ones (i.e., $\tau_1,<\ldots<\tau_M\leq \tau_f \ll \tau_s$) \cite{shampine1979user}. In the general case, $M$ can change over time; here we will assume that $M$ does not change during the period of interest.~Under normal hyperbolicity, the above assumption implies the emergence of an $(N-M)$-dim. SIM in the phase space, attracting or repelling all neighboring trajectories of the IVP.~On the SIM, the evolution of $M$ transformed variables of $\mathbf{z} \in \mathbb{R}^{N}$ are slaved by the remaining $N-M$ ones; thus, naturally deeming the former as fast variables and the latter as slow.~Denoting the fast variables $\mathbf{x}\in \mathbb{R}^{M}$ and the slow ones $\mathbf{y}\in \mathbb{R}^{N-M}$, the SIM can be locally described as:
\begin{equation}
    \mathcal{M} = \{ (\mathbf{x},\mathbf{y}) \in \mathbb{R}^N: \mathbf{x} = \mathbf{h}(\mathbf{y})\},
    \label{eq:SIM}
\end{equation}
where $\mathbf{h}:\mathbb{R}^{N-M}\rightarrow\mathbb{R}^{M}$ is a smooth map.~Since we seek to determine the SIM from data, we naturally restrict ourselves to the cases of \emph{attractive} SIMs.~In such cases, for any initial condition $\mathbf{z}_0$ close to $\mathcal{M}$, the solution of the IVP is attracted rapidly towards $\mathcal{M}$ under the action of the $M$ fast timescales.~Once on $\mathcal{M}$, the solution evolves along it, with the slow flow being characterized by $\tau_s$, until reaching a steady state or exiting through its boundary $\partial \mathcal{M}$ \citep{fenichel1979geometric,kaper1999introduction,zagaris2004analysis}.

Naturally, the existence of a SIM as given by Eq.~\eqref{eq:SIM} implies that the IVP in Eq.~\eqref{eq:gen} can be in principle decomposed  into fast and slow dynamics:
\begin{equation}
    \dfrac{d \mathbf{x}}{dt} = \mathbf{f}(\mathbf{x},\mathbf{y}), \qquad \dfrac{d \mathbf{y}}{dt} = \mathbf{g}(\mathbf{x},\mathbf{y})
    \label{eq:genVFdec}
\end{equation}
where $\mathbf{f}:\mathbb{R}^N \rightarrow \mathbb{R}^M$ and $\mathbf{g}:\mathbb{R}^N \rightarrow \mathbb{R}^{N-M}$ are the vector fields corresponding to the fast and slow variables, respectively. Following the above decomposition, one may construct a ROM for the slow variables as $d\mathbf{y}/dt=\mathbf{g}(\mathbf{h}(\mathbf{y}),\mathbf{y})$, since $\mathbf{x}=\mathbf{h}(\mathbf{y})$ when the solution evolves on $\mathcal{M}$.\par
Hence, the objective is the construction of the transformation operator that decomposes the original system in Eq.~\eqref{eq:gen} to the slow-fast dynamics system in Eq.~\eqref{eq:genVFdec}, such that the underlying SIM is provided by the functional form in Eq.~\eqref{eq:SIM}.

\subsection{The invariance equation and alternative methods for computing SIMs}

For the derivation of SIM approximations in the form of Eq.~\eqref{eq:SIM}, many reduction methods exploit the local invariance of $\mathcal{M}$.~In particular, differentiation of the map $\mathbf{x}=\mathbf{h}(\mathbf{y})$ w.r.t. time, yields:
\begin{equation*}
\dfrac{d\mathbf{x}}{dt} = \dfrac{d\mathbf{h}(\mathbf{y})}{dt} \Rightarrow \dfrac{d\mathbf{x}}{dt} - \nabla_{\mathbf{y}}\mathbf{h}(\mathbf{y}) \dfrac{d\mathbf{y}}{dt} = \mathbf{0}.
\end{equation*}
By substituting the expressions in Eq.~\eqref{eq:genVFdec}, one gets a PDE corresponding to the so-called \emph{invariance equation}:
\begin{equation}
      \mathbf{f}(\mathbf{x},\mathbf{y}) - \nabla_{\mathbf{y}}\mathbf{h}(\mathbf{y})  \mathbf{g}(\mathbf{x},\mathbf{y}) = \mathbf{0}.
    \label{eq:Inv}
\end{equation}
The matrix $\nabla_{\mathbf{y}}\mathbf{h}(\mathbf{y})\in \mathbb{R}^{M \times N-M}$ involves the partial derivatives of the map $\mathbf{h}$ with respect to the slow variables $\mathbf{y}$.~Assuming that the fast and slow variables, $\mathbf{x}$ and $\mathbf{y}$, are purely some of the components of the original variables in $\mathbf{z}$, various methods have been developed to obtain the function $\mathbf{h}(\mathbf{y})$ from Eq.~\eqref{eq:Inv};  either analytically, such as the Invariance Equation (IE) method \cite{verhulst2005methods,kuehn2015multiple}, the Method of Invariant Manifolds (MIM) \citep{gorban2003method}, the Rousel-Fraser method \citep{roussel1991geometry,fraser1988steady}, or numerical ones, such as the Zero Derivative Principle (ZDP)  method \citep{gear2005projecting,zagaris2009analysis} and a variant method of the Computational Singular Perturbation (CSP)  \cite{goussis2006efficient}. 

The above methods based on the invariance equation have been shown to be cases of more sophisticated computational methods in the context of GSPT for the derivation of SIMs approximations \citep{goussis2006efficient}; notably among them the computational singular perturbation (CSP) method \citep{lam1989understanding,goussis1992study,hadjinicolaou1998asymptotic,valorani2005higher,goussis2012quasi}, the invariant low-dimensional method (ILDM) \citep{maas1992simplifying,maas1992implementation,maas1994laminar} and the tangential stretching rate (TSR) method  \citep{valorani2015dynamical,valorani2017tangential}.~These methods follow iterative procedures to locally approximate the \emph{fast} and \emph{slow} subspaces resolving the \emph{tangent} space, along which the fast and slow dynamics act, respectively \citep{kuehn2015multiple,goussis2011model,kaper1999introduction,zagaris2004analysis}.~Accurate estimation of the basis vectors spanning the subspaces is critical, since the projection of the vector field on the fast subspace provides the SIM approximation, while that on the slow subspace provides the ROM.~In fact, CSP, ILDM and TSR provide different iterative procedures for approximating the fast/slow basis vectors \citep{kuehn2015multiple,kaper1999introduction,goussis2011model,goussis2006efficient}.~Usually, a few iterations of the above methods result, especially for systems with a low number of variables, to analytic SIM approximations in an explicit form; i.e., the map $\mathbf{x}=\mathbf{h}(\mathbf{y})$.~In general, however,  the SIM approximations provided by the computational GSPT methods 
are of implicit form, i.e., $\mathbf{h}(\mathbf{x},\mathbf{y})=\mathbf{0}$.~Thus, additional numerical root-finding algorithms are required for obtaining explicit expressions, which are performed point-by-point as $\mathbf{y}$ changes along the trajectory of the system in Eq.~\eqref{eq:gen} \citep{goussis2006efficient,kuehn2015multiple}. Simpler SIM approximations can be constructed by the implementation of the traditional (in chemical kinetics) model reduction methodologies of Quasi Steady-State Approximation (QSSA) and Partial Equilibrium Approximation (PEA) \cite{bowen1963singular,rein1992partial,fraser1988steady}.~The QSSA arises when the rate of change of a variable in $\mathbf{z}$ is much smaller than its production and consumption rate, while the PEA arises when the forward and backward directions of a unidirectional reaction equilibrate.~Recently, it was shown that the SIM approximations resulting by QSSA and PEA can be recovered, along with conditions for their validity, as limiting cases of those derived on the basis of CSP  \citep{goussis2012quasi,goussis2015model,patsatzis2023algorithmic}.~In this work, for comparison purposes with proposed PINN approach, we also derived SIM approximations on the basis of QSSA, PEA and CSP with one and two iterations; additional details for their implementation are provided in Section~\ref{sbsb:CSP_PEA_QSSA}.

\subsection{The proposed Physics-Informed Neural Network (PINN)-based methodology}
\label{sub:MLgen}

Here, we present a  machine learning methodology based on PINNs for the derivation of explicit, in terms of fast variables (in the form of Eq.~\eqref{eq:SIM}), SIM approximations.~We begin by considering the general form of the IVP in Eq.~\eqref{eq:gen} and assume that it exhibits a constant, during the period of interest, number $M$ of dissipative fast timescales, thus implying the emergence of an attractive $(N-M)$-dim. SIM $\mathcal{M}$.~We further assume that $\mathcal{M}$ can be locally approximated by the map $\mathbf{x}=\mathbf{h}(\mathbf{y})$ in Eq.~\eqref{eq:SIM}, without, however, having any a priori knowledge about the map or the fast and slow variables. \par
The first step is to estimate $M$. As discussed in the introduction, for complex systems, such an estimate can be attempted in a data-driven way using manifold learning.~Here, we exploit the availability of the vector field in Eq.~\eqref{eq:gen} to find $M$ via the CSP iterative criterion, introduced in \cite{valorani2001explicit}.~This criterion guarantees that any data point of the solution of the IVP is sufficiently close to the SIM $\mathcal{M}$ without diverting from it, due to the slow dynamics evolution.~In practice, one calculates the Jacobian matrix $\mathbf{J}=\partial \mathbf{F}(\mathbf{z})/\partial \mathbf{z}$ at any data point of the IVP's solution, in order to determine the eigenvalues, say $\lambda_j$, in descending order, for $j=1,\ldots,N$, and the corresponding right and left eigenvectors, say $\mathbf{v}_j$ and $\mathbf{u}^j$ respectively, of the linearized system.~Then, the number $M$ is determined by the largest integer ($M<N$) satisfying: 
\begin{equation}
\bigg \vert \dfrac{1}{\lVert \lambda_{M+1} \rVert} \sum_{m=1}^M \mathbf{v}_m \mathbf{u}^m \mathbf{F}(\mathbf{z}) \bigg \vert < \epsilon_{rel} \mathbf{z} + \boldsymbol{\epsilon}_{abs},
\label{eq:CSPcrit}
\end{equation}
where the relative and absolute errors $\epsilon_{rel}$  and $\boldsymbol{\epsilon}_{abs}$ express the maximum timescale splitting allowed (typically from $10^{-3}$ to $10^{-1}$) and the values of $\mathbf{z}$ that are considered zero, respectively \cite{valorani2001explicit}.

Having estimated $M$, we seek to construct two differentiable transformations from the state vector variables $\mathbf{z}\in\mathbb{R}^N$ in Eq.~\eqref{eq:gen} into the vectors of fast and slow variables, $\mathbf{x}\in\mathbb{R}^M$ and $\mathbf{y}\in \mathbb{R}^{N-M}$, as:
\begin{equation}
    \mathbf{x} = \mathcal{X}(\mathbf{z}), \quad \mathcal{X}: \mathbb{R}^N\rightarrow \mathbb{R}^M, \quad \text{and} \quad \mathbf{y} = \mathcal{Y}(\mathbf{z}), \quad \mathcal{Y}: \mathbb{R}^N\rightarrow \mathbb{R}^{N-M}.
    \label{eq:trans}
\end{equation}
Then, the system in Eq.~\eqref{eq:genVFdec} can be cast in terms of the state variables, as:
\begin{equation}
    \dfrac{d\mathbf{x}}{dt} = \nabla_{\mathbf{z}} \mathcal{X}(\mathbf{z})  \dfrac{d\mathbf{z}}{dt} = \mathbf{f}(\mathcal{X}(\mathbf{z}), \mathcal{Y}(\mathbf{z})), \qquad  \dfrac{d\mathbf{y}}{dt} = \nabla_{\mathbf{z}} \mathcal{Y} (\mathbf{z})  \dfrac{d\mathbf{z}}{dt} = \mathbf{g}(\mathcal{X}(\mathbf{z}), \mathcal{Y}(\mathbf{z})),
    \label{eq:TransDer}
\end{equation}
where $\nabla_{\mathbf{z}} \mathcal{X} (\mathbf{z}) \in \mathbb{R}^{M\times N}$ and $\nabla_{\mathbf{z}} \mathcal{Y} (\mathbf{z}) \in \mathbb{R}^{(N-M)\times N}$.~Now, given the existence of the SIM, we seek a mapping $\mathbf{x}=\mathbf{h}(\mathbf{y})$ in Eq.~\eqref{eq:SIM}, so that the transformed variables satisfy $(\mathcal{X}(\mathbf{z}),\mathcal{Y}(\mathbf{z}))\in \mathcal{M}$; i.e.,
\begin{equation}
    \mathcal{X}(\mathbf{z})=\mathbf{h}(\mathcal{Y}(\mathbf{z})))
    \label{eq:SIMTrans}
\end{equation}  
Due to the invariance of the SIM, the map for the transformed variables also satisfies the invariance equation in Eq.~\eqref{eq:Inv}, yielding:
\begin{equation}
    \nabla_{\mathbf{z}} \mathcal{X} (\mathbf{z})  \dfrac{d\mathbf{z}}{dt} - \nabla_{\mathbf{y}} \mathbf{h}  (\mathbf{y})  \nabla_{\mathbf{z}} \mathcal{Y} (\mathbf{z})  \dfrac{d\mathbf{z}}{dt} = 
    \left( \nabla_{\mathbf{z}} \mathcal{X} (\mathbf{z}) - \nabla_{\mathbf{y}} \mathbf{h} (\mathbf{y})  \nabla_{\mathbf{z}} \mathcal{Y} (\mathbf{z}) \right) \mathbf{F}(\mathbf{z}) = \mathbf{0}.
    \label{eq:InvTrans}
\end{equation}
The above is a set of $M$ PDEs, where the only known term is the state vector field $\mathbf{F}(\mathbf{z})$.\par
Here, we seek for  linear transformations, of the state variables $\mathbf{z}$.~Although this is restricting in the general case, this selection allows us to compare the proposed PINN framework with the GSPT methods, which are also based on such linear transformations of the state variables $\mathbf{z}$.~According to this, from Eq.~\eqref{eq:trans}, we have:
\begin{equation}
\mathbf{x} = \nabla_{\mathbf{z}} \mathcal{X} (\mathbf{z})  ~ \mathbf{z} = \mathbf{C} ~ \mathbf{z}, \quad x_m = \sum_{j=1}^N C_m^j z_j, 
\qquad 
\mathbf{y} = \nabla_{\mathbf{z}} \mathcal{Y} (\mathbf{z}) ~ \mathbf{z} = \mathbf{D} ~ \mathbf{z}, 
\quad y_d = \sum_{j=1}^N D_d^j z_j, 
\label{eq:transLin}
\end{equation}
where $\mathbf{C}\in \mathbb{R}^{M\times N}$ and $\mathbf{D} \in \mathbb{R}^{(N-M)\times N}$ are row-wise matrices formed by  $\mathbf{C}_m=[C_m^1,\ldots,C_m^j,\ldots,C_m^N]\in\mathbb{R}^N$ and $\mathbf{D}_d=[D_d^1,\ldots,D_d^j,\ldots,D_d^N]\in\mathbb{R}^N$, respectively, for $m=1,\ldots,M$, $d=1,\ldots,N-M$ and $j=1,\ldots,N$.~The $m$-th, $d$-th and $j$-th components of the fast, slow and original variables are hereby denoted by $x_m$, $y_d$ and $z_j$, respectively.~For ensuring that each of the state variables $\mathbf{z}$ is related purely to either the fast or slow variables, we require the following conditions: 
\begin{equation}
    \sum_{m=1}^M C^j_m + \sum_{d=1}^{N-M} D^j_d = 1, \qquad \sum_{j=1}^N C^j_m=\sum_{j=1}^N D^j_d=1, \qquad  D_d^j \sum^M_{m=1} C_m^j =0, \qquad C_m^j \sum^{N-M}_{d=1} D_d^j =0, 
    \label{eq:PinTrans}
\end{equation}
for every $m=1,\ldots,M$, $d=1,\ldots,N-M$ and $j=1,\ldots,N$.~The first two conditions normalize the transformation, while the second two ensure that each component of $\mathbf{z}$ is attributed solely to the fast or slow variables, $\mathbf{x}$ and $\mathbf{y}$.

Let's now assume a set of $n$ input points $\mathbf{z}^{(i)}\in \Omega \subset \mathbb{R}^N$ for $i=1,\ldots,n$, which, for all practical purposes, are sampled in the domain where the SIM $\mathcal{M}$ emerges.~Then, the numerical approximation of the SIM can be obtained via PINNs by 
\begin{enumerate}[label=(\roman*)]
    \item requiring $(\mathcal{X}(\mathbf{z}^{(i)}),\mathcal{Y}(\mathbf{z}^{(i)})) \in \mathcal{M}$, which according to Eq.~\eqref{eq:SIMTrans} implies the minimization of the loss function:
\begin{equation}
\mathcal{L}_1 (\mathbf{z}^{(i)},\mathbf{C},\mathbf{D},\mathbf{P};\mathbf{Q}) := \mathbf{C}  ~ \mathbf{z}^{(i)} - \mathcal{N}(\mathbf{D}  \mathbf{z}^{(i)},\mathbf{P},\mathbf{Q}), 
    \label{eq:minFun1}
\end{equation}
\item satisfying the IE, which according to Eq.~\eqref{eq:InvTrans} implies the minimization of the loss function:
\begin{equation}
\mathcal{L}_2 (\mathbf{z}^{(i)},\mathbf{C},\mathbf{D},\mathbf{P};\mathbf{Q}):= 
 \left( \mathbf{C} - \nabla_{\mathbf{y}} \mathcal{N}(\mathbf{D}  \mathbf{z}^{(i)},\mathbf{P},\mathbf{Q})  ~\mathbf{D} \right)  \mathbf{F}(\mathbf{z}^{(i)}),
    \label{eq:minFun2}  
\end{equation}
\item satisfying the pinning conditions of the transformation in Eq.~\eqref{eq:PinTrans} implemented as  minimization of the loss functions:
\begin{equation*}
\mathcal{L}_3^C (\mathbf{C},\mathbf{D}):=     \sum_{m=1}^M C^j_m + \sum_{d=1}^{N-M} D^j_d -1, \quad 
\mathcal{L}_3^R (\mathbf{C},\mathbf{D}):=
\begin{cases}
    \sum_{j=1}^N C^j_m-1 \\[2pt] 
    \sum_{j=1}^N D^j_d-1
\end{cases}, 
\end{equation*}
\begin{equation}
\mathcal{L}_3^{SC} (\mathbf{C},\mathbf{D}):= D_d^j \sum^M_{m=1} C_m^j, \quad 
\mathcal{L}_3^{FC} (\mathbf{C},\mathbf{D}):= C_m^j \sum^{N-M}_{d=1} D_d^j. 
\label{eq:minFunPC}  
\end{equation}
\end{enumerate}
The ANN function $\mathcal{N}(\cdot):=\mathcal{N}(\mathbf{y},\mathbf{P},\mathbf{Q}):\mathbb{R}^{N-M}\rightarrow\mathbb{R}^M$ in Eqs.~(\ref{eq:minFun1},\ref{eq:minFun2}) approximates the mapping $\mathbf{x} = \mathbf{h}(\mathbf{y})$ of the SIM in Eq.~\eqref{eq:SIM}, since it takes as inputs the transformed slow variables $\mathbf{y}=\mathbf{D} ~ \mathbf{z} \in \mathbb{R}^{N-M}$ and outputs an $M$-dim. vector, which, due to the minimization in Eq.~\eqref{eq:minFun1}, approximates the transformed fast variables $\mathbf{x}=\mathbf{C} ~ \mathbf{z} \in \mathbb{R}^M$.~$\mathcal{N}(\cdot)$ contains the parameters $\mathbf{P}$ of the PINN scheme (weights and biases of each layer) and the hyper-parameters $\mathbf{Q}$ such as the type and parameters of activation function, the learning rate, the number of epochs, etc.

In summary, the numerical approximation of the SIM via the proposed PINN methodology is obtained by the solution of the optimization problem:
\begin{equation}
    \min_{\mathbf{C}, \mathbf{D}, \mathbf{P},\mathbf{Q}} \mathcal{L} (\mathbf{C}, \mathbf{D},  \mathbf{P}; \mathbf{Q}) :=
     \sum_{i=1}^{n} \left( \big\lVert \mathcal{L}_1 (\mathbf{z}^{(i)},\mathbf{C},\mathbf{D},\mathbf{P};\mathbf{Q})  \big\rVert^2 +   \big\lVert \mathcal{L}_2 (\mathbf{z}^{(i)},\mathbf{C},\mathbf{D},\mathbf{P};\mathbf{Q})  \big\rVert^2\right) +  \omega \big\lVert \mathcal{L}_3 (\mathbf{C},\mathbf{D})  \big\rVert^2,
    \label{eq:minFunPI}
\end{equation}
where $\mathcal{L}_3 (\mathbf{C},\mathbf{D}) = \mathcal{L}^C_3 (\mathbf{C},\mathbf{D})+\mathcal{L}^R_3 (\mathbf{C},\mathbf{D})+\mathcal{L}^{SC}_3 (\mathbf{C},\mathbf{D})+\mathcal{L}^{FC}_3 (\mathbf{C},\mathbf{D})$ contains the $N^2+2N$ pinning conditions in Eq.~\eqref{eq:minFunPC} and $\omega\ge0$ is a weight for the contribution $\mathcal{L}_3$ w.r.t. $\mathcal{L}_1$ and $\mathcal{L}_2$ in the optimization task.~Note that the solution of the optimization problem in Eq.~\eqref{eq:minFunPI}, requires the derivation of the first and second order derivatives of $\mathcal{N}(\mathbf{y},\mathbf{P},\mathbf{Q})$ w.r.t. $\mathbf{y}$ and the parameters in $\mathbf{P}$, which can be obtained by numerical (e.g., using finite differences), symbolic or automatic differentiation \citep{baydin2018automatic,lu2021deepxde}.~Among the various structures that $\mathcal{N}(\cdot)$ can handle, here we consider single-layer feedforward ANNs with sigmoid activation functions, following \cite{patsatzis2023slow}.

\subsubsection{PINN implementation}
Here, we used single-layer feedforward neural networks (SLFNN).
The input and output dimensions of the SLFNN are $D=N-M$, and, $M$ respectively, as the dimensions of the slow and output fast variables in $\mathbf{y}\in\mathbb{R}^{N-M}$ and $\mathbf{x}\in\mathbb{R}^{M}$; the notation $\mathbf{x}/\mathbf{y}$ is kept for a more compact presentation.~For $L$ neurons in the hidden layer, the $m$-th output of the SLFNN ($m=1,\ldots,M$) can be written in the form:
\begin{equation}
    \mathcal{N}^{(m)}(\mathbf{y},\mathbf{P}^{(m)}) = \sum_{l=1}^L w^{o(m)}_{l} \mathcal{\phi}_l\left(\sum_{d=1}^{D} w^{(m)}_{ld} y_d + b^{(m)}_l\right)+ b^{o(m)} = \mathbf{w}^{o(m)\top} \phi \left( \mathbf{W}^{(m)} \mathbf{y}  + \mathbf{b}^{(m)} \right)  + b^{o(m)},
    \label{eq:SFLNNsumf}
\end{equation}
where $\mathbf{P}^{(m)} = [\mathbf{w}^{o(m)},b^{o(m)},\mathbf{W}^{(m)},\mathbf{b}^{(m)}]^\top\in\mathbb{R}^{L(D+2)+1}$ includes: (i) the output weights $\mathbf{w}^{o(m)} = [w^{o(m)}_1, \ldots, w^{o(m)}_L]^\top\in \mathbb{R}^{L}$ of the neurons between the hidden and the output layer, (ii) the bias $b^{o(m)} \in \mathbb{R}$ of the output layer, (iii) the internal weights $\mathbf{W}^{(m)}\in\mathbb{R}^{L\times D}$ between the input and the hidden layer, the columns $\mathbf{w}^{(m)}_l = [w_{l1}^{(m)}, \ldots, w_{lD}^{(m)}]^\top\in \mathbb{R}^{D}$ of which correspond to the weights between the input neurons and the $l$-th neuron in the hidden layer, and (iv) the internal biases  $\mathbf{b}^{(m)} = [b_1^{(m)}, \ldots, b_L^{(m)}]^\top\in\mathbb{R}^{L}$ of the neurons in the hidden layer.~Here, we chose as activation function $\phi_l(\cdot) = \phi(\mathbf{w}_l^{(m)\top} \mathbf{y}  + b_l^{(m)})$ the logistic sigmoid function which allows symbolic differentiation of the required, for the optimization problem, derivatives.

The solution of the PINN optimization problem in Eq.~\eqref{eq:minFunPI} requires the minimization of the three loss functions, $\mathcal{L}_1$, $\mathcal{L}_2$, $\mathcal{L}_3$ w.r.t. the parameters of the transformation and the SLFNN, which can be stacked into the column vector:
\begin{equation}
    \mathbf{R} = \big[\mathbf{C}_{1},\ldots,\mathbf{C}_{m},\ldots,\mathbf{C}_{M},\mathbf{D}_{1},\ldots,\mathbf{D}_{d},\ldots,\mathbf{D}_{N-M},\mathbf{P}^{(1)},\ldots,\mathbf{P}^{(m)},\ldots,\mathbf{P}^{(M)}\big]^\top \in\mathbb{R}^{N^2+(L(D+2)+1)M}.
    \label{eq:TotPar}
\end{equation}
To minimize the loss functions for every input point $i=1,\ldots,n$, the minimization of the following non-linear residuals is required:   
\begin{subequations}
\begin{align}
    \mathcal{F}_q^1(\mathbf{R}) & = \sum_{j=1}^N C_m^j z^{(i)}_j-\mathcal{N}^{(m)}( \mathbf{D}  \mathbf{z}^{(i)},\mathbf{P}^{(m)}),  \label{eq:ResL1} \\
    \mathcal{F}_q^2(\mathbf{R}) & = \sum_{j=1}^N \left(C_m^j -\sum_{d=1}^{N-M} \dfrac{\partial \mathcal{N}^{(m)}(\mathbf{D}  \mathbf{z}^{(i)},\mathbf{P}^{(m)})}{\partial y_d} D_d^j\right) F_j(\mathbf{z}^{(i)}) \label{eq:ResL2} \\
    \mathcal{F}^{3C}_j(\mathbf{R}) & = \omega\Bigg(\sum_{m=1}^M C^j_m + \sum_{d=1}^{N-M} D^j_d - 1\Bigg), \quad 
    \mathcal{F}^{3R}_{m}(\mathbf{R}) = \omega   \Bigg(\sum_{j=1}^N C^j_m - 1 \Bigg), \quad 
    \mathcal{F}^{3R}_{M+d}(\mathbf{R}) = \omega \Bigg(\sum_{j=1}^N D^j_d - 1 \Bigg)     \nonumber \\
    \mathcal{F}^{3SC}_{q1}(\mathbf{R}) & = \omega D_d^j \sum^M_{m=1} C_m^j, \qquad \qquad \qquad \quad~
    \mathcal{F}^{3FC}_{q2}(\mathbf{R}) = \omega C_m^j \sum^{N-M}_{d=1} D_d^j \label{eq:ResL3},
\end{align}
\label{eq:ResALL}
\end{subequations}
where $q=m+(i-1)M$, $q1=(d-1)N+j$ and $q2=(m-1)N+j$.~The term $F_j(\cdot)$ in Eq.~\eqref{eq:ResL2} corresponds to the $j$-th component of the analytically known from Eq.~\eqref{eq:gen} vector field $\mathbf{F}=[F_1,\ldots,F_j,\ldots,F_N]^\top$.~Clearly, the $M\times n$ residuals in Eqs.~(\ref{eq:ResL1}, \ref{eq:ResL2}) correspond to $\mathcal{L}_1$ and $\mathcal{L}_2$ and the $N^2+2N$ ones in Eq.~\eqref{eq:ResL3} correspond to the pinning conditions in $\mathcal{L}_3$, multiplied by $\omega$.

For computing the residuals in Eq.~\eqref{eq:ResALL}, it is required to calculate the $M\times (N-M)$ derivatives $\partial\mathcal{N}^{(m)}(\mathbf{D} \mathbf{z}^{(i)},\mathbf{P}^{(m)})/\partial y_d$; derivation with symbolic differentiation is enabled by using the derivative of the sigmoid function, $\phi'=\phi(1-\phi)$.~Thus, according to Eq.~\eqref{eq:SFLNNsumf}, the derivative of the $m$-th SLFNN output w.r.t. the $d$-th slow variable is:
\begin{equation}
    \dfrac{\partial \mathcal{N}^{(m)}(\mathbf{D} \mathbf{z}^{(i)},\mathbf{P}^{(m)})}{ \partial y_d} = \sum_{l=1}^L  w^{o(m)}_l w^{(m)}_{ld}\left[ \phi_l(\cdot) ( 1 - \phi_l (\cdot ) ) \right], \qquad \phi_l(\cdot) = \phi_l(\mathbf{w}_{l}^{(m)\top} \cdot (\mathbf{D} \mathbf{z}^{(i)})  + b^{(m)}_l),
    \label{eq:SFLNNder_sumf}
\end{equation}
for $m=1,\ldots,M$ and $d=1,\ldots,N-M$.~Given Eq.~\eqref{eq:SFLNNsumf} and \eqref{eq:SFLNNder_sumf}, all the components to calculate the non-linear residuals in Eq.~\eqref{eq:ResALL} are now available for every input point $i=1,\ldots,n$.

\subsubsection{Computation of the unknown PINN parameters}

As discussed, the solution of the PINN optimization problem in Eq.~\eqref{eq:minFunPI} is reduced to the minimization of the non-linear residuals in Eq.~\eqref{eq:ResALL} w.r.t. the unknown parameters $\mathbf{R}$ in Eq.~\eqref{eq:TotPar}.~For $n$ input points, the non-linear residuals can be stacked into the $\boldsymbol{\mathcal{F}}(\mathbf{R}) \in \mathbb{R}^{2Mn+N^2+N}$ column vector:
\begin{align}
\boldsymbol{\mathcal{F}}(\mathbf{R}) = [ & \mathcal{F}^1_1,\ldots,\mathcal{F}^1_q,\ldots,\mathcal{F}^1_{Mn},\mathcal{F}^2_1,\ldots,\mathcal{F}^2_q,\ldots,\mathcal{F}^2_{Mn},\mathcal{F}^{3C}_1,\ldots,\mathcal{F}^{3C}_j,\ldots,\mathcal{F}^{3C}_N, \mathcal{F}^{3R}_1,\ldots,\mathcal{F}^{3R}_m,\ldots,\mathcal{F}^{3R}_M,
\nonumber \\ 
& \mathcal{F}^{3R}_{M+1},\ldots,\mathcal{F}^{3R}_{M+d},\ldots,\mathcal{F}^{3R}_N,\mathcal{F}^{3SC}_1,\ldots,\mathcal{F}^{3SC}_{q1},\ldots,\mathcal{F}^{3SC}_{N(N-M)},\mathcal{F}^{3FC}_1,\ldots,\mathcal{F}^{3FC}_{q2},\ldots,\mathcal{F}^{3FC}_{NM} ]^\top.
\label{eq:ResVec}
\end{align}
In general, the minimization problem is overdetermined, requiring iterative methods for computing the unknown parameters.~Here, we utilized the Levenberg-Marquardt (LM) algorithm \citep{hagan1994training}.

For the implementation of the LM algorithm, the Jacobian matrix of the residuals w.r.t. $\mathbf{R}$ is required, which can be calculated through symbolic differentiation, as presented next.~Recall that $\boldsymbol{\mathcal{F}(\mathbf{R})}$ in Eq.~\eqref{eq:ResVec} is expressed using the indices $m$, $d$ and $j$.~Since, however, each parameter is not used for every residual, we introduce, only for differentiation, the indices $r=1,\ldots,M$, $h=1,\ldots,N-M$ and $k=1,\ldots,N$ (instead of  $m$, $d$ and $j$ respectively) to distinguish between the parameters of $\mathbf{R}$ over which the differentiation is performed.~Then, the first $Mn$ rows of the Jacobian $\nabla_{\mathbf{R}} \boldsymbol{\mathcal{F}} \in \mathbb{R}^{(2Mn+N^2+N) \times (N^2+M(L(N-M+2)+1))}$, corresponding to the derivatives of $\mathcal{F}^1_q$ in Eq.~\eqref{eq:ResL1}, are formed by the elements:
\begin{align}
\dfrac{\partial \mathcal{F}_q^1}{\partial C_r^k} & = 
 z_k^{(i)}, \quad \dfrac{\partial \mathcal{F}_q^1}{\partial D_h^k} = - \dfrac{\partial \mathcal{N}^{(m)}(\mathbf{D} \mathbf{z}^{(i)},\mathbf{P}^{(m)})}{\partial y_h} z_k^{(i)}
 \nonumber \\
\dfrac{\partial \mathcal{F}_q^1}{\partial w_l^{o(r)}} & = -\phi_l(\cdot), \quad 
\dfrac{\partial \mathcal{F}_q^1}{\partial b^{o(r)}} = -1, \quad 
\dfrac{\partial \mathcal{F}_q^1}{\partial w_{lh}^{(r)}} = -w_l^{o(r)} \left( \sum_{j=1}^{N} D_h^j z_j^{(i)} \right) \phi_l'(\cdot), \quad 
\dfrac{\partial \mathcal{F}_q^1}{\partial b_{l}^{(r)}} = - w_l^{o(r)} \phi_l'(\cdot),
\label{eq:Res1_ders}
\end{align}
if $r=m$ and zeros if $r\neq m$; recall that $q=m+(i-1)M$.~Note here that the derivatives w.r.t. $D_h^k$ are never zero, indicating that all the slow variables $\mathbf{y}$ are, in general, contributing to all $m$ components of the SIM functional.~In Eq.~\eqref{eq:Res1_ders}, $\phi_l'(\cdot)= \phi_l(\cdot) ( 1 - \phi_l (\cdot ) )$ is the derivative of the sigmoid function with $\phi_l(\cdot) = \phi_l(\mathbf{w}_{l}^{(r)\top} \mathbf{D} \mathbf{z}^{(i)}  + b^{(r)}_l)$.~Now, for the derivatives of $\mathcal{F}^2_q$ in Eq.~\eqref{eq:ResL2}, which are more complicated, we first present the elements of the Jacobian corresponding to the transformation (i.e., the first $N^2$ columns), reading: 
\begin{equation}
\dfrac{\partial \mathcal{F}_q^2}{\partial C_r^k} = F_k(\mathbf{z}^{(i)}), \quad \dfrac{\partial \mathcal{F}_q^2}{\partial D_h^k} = - \dfrac{\partial \mathcal{N}^{(m)}(\mathbf{D}\mathbf{z}^{(i)},\mathbf{P}^{(m)})}{\partial y_h} F_k(\mathbf{z}^{(i)})-\sum_{j=1}^N\left(\sum_{d=1}^{N-M} \dfrac{\partial^2 \mathcal{N}^{(m)}(\mathbf{D}\mathbf{z}^{(i)},\mathbf{P}^{(m)})}{\partial y_d \partial D_h^k} D_d^j \right)F_j(\mathbf{z}^{(i)}),
\label{eq:Res2_ders1}
\end{equation}
if $r=m$ and zeros if $r\neq m$.~As expected, all the elements of the second term are, in theory, non-zeros, since:
\begin{equation}
\dfrac{\partial^2 \mathcal{N}^{(m)}(\mathbf{D}\mathbf{z}^{(i)},\mathbf{P}^{(m)})}{\partial y_d \partial D_h^k} = \sum_{l=1}^L  w^{o(m)}_l w^{(m)}_{ld} w^{(m)}_{lh} \phi_l''(\cdot) z_k^{(i)},  \qquad \phi_l''(\cdot)= \phi_l(\cdot) ( 1 - \phi_l (\cdot ) )  ( 1 - 2 \phi_l (\cdot ) ).
\end{equation}
The remaining elements of the second $Mn$ rows of the Jacobian, corresponding to the SLFNN weights and biases, are expressed over a parameter $p\in \mathbf{P}$ as: 
\begin{equation}
\dfrac{\partial \mathcal{F}_q^2}{\partial p} = - \sum_{j=1}^N\left(\sum_{d=1}^{N-M} \dfrac{\partial^2 \mathcal{N}^{(m)}(\mathbf{D}\mathbf{z}^{(i)},\mathbf{P}^{(m)})}{\partial y_d \partial p} D_d^j \right) F_j(\mathbf{z}^{(i)}),
\label{eq:Res2_ders2}
\end{equation}
where, by differentiation of Eq.~\eqref{eq:SFLNNder_sumf} w.r.t $p$, the derivatives corresponding to (i) the output weights, for $p=w_l^{o(r)}$, and (ii) the output biases, for $p=b^{o(r)}$, read:
\begin{equation}
\dfrac{\partial^2 \mathcal{N}^{(m)}(\mathbf{D}\mathbf{z}^{(i)},\mathbf{P}^{(m)})}{\partial y_d \partial w_l^{o(r)}} = 
    \begin{cases}
        w_{ld}^{(r)}\phi_l'(\cdot), & \text{if } r = m\\
        0, & \text{if } r \neq m
    \end{cases} ,
\qquad
\dfrac{\partial^2 \mathcal{N}^{(m)}(\mathbf{D}\mathbf{z}^{(i)},\mathbf{P}^{(m)})}{\partial y_d \partial b^{o(r)} } = 0, 
\end{equation}
(iii) the input weights, for $p=w_{lh}^{(r)}$, read:
\begin{equation}
\dfrac{\partial^2 \mathcal{N}^{(m)}(\mathbf{D}\mathbf{z}^{(i)},\mathbf{P}^{(m)})}{\partial y_d \partial w_{lh}^{(r)} } = w_l^{o(r)} \cdot 
    \begin{cases}
         w_{ld}^{(r)}\left( \sum_{j=1}^{N} D_h^j z_j^{(i)} \right)\phi_l''(\cdot)+\phi_l'(\cdot), & \text{if } r=m\ , \ h=d  \\
         w_{ld}^{(r)}\left( \sum_{j=1}^{N} D_h^j z_j^{(i)} \right)\phi_l''(\cdot), & \text{if } r=m\ , \ h \neq d \\
         0, & \text{if } r \neq m
    \end{cases},
\end{equation}
and (iv) the input biases, for $p=b_l^{(r)}$, read:
\begin{equation} 
\dfrac{\partial^2 \mathcal{N}^{(m)}(\mathbf{y},\mathbf{P}^{(m)})}{\partial y_d \partial b_l^{(r)}} = \begin{cases}
        w_l^{o(r)} w_{ld}^{(r)} \phi_l''(\cdot), & \text{if } r=m \\
         0, & \text{if } r \neq m
    \end{cases}, 
\end{equation}
where $\phi_l'(\cdot)= \phi_l(\cdot) ( 1 - \phi_l (\cdot ) )$ and $\phi_l''(\cdot)= \phi_l(\cdot) ( 1 - \phi_l (\cdot ) )( 1 - 2\phi_l (\cdot ) )$ are the first and second derivatives of the sigmoid function with $\phi_l(\cdot) = \phi_l(\mathbf{w}_{l}^{(r)\top} \mathbf{D} \mathbf{z}^{(i)}  + b^{(r)}_l)$.~Finally, the last $N^2+N$ rows of the Jacobian matrix, corresponding to the derivatives of  $\mathcal{F}^{3C}_j$, $\mathcal{F}^{3R}_m$, $\mathcal{F}^{3R}_{M+d}$, $\mathcal{F}^{3SC}_{q1}$ and $\mathcal{F}^{3FC}_{q2}$ in Eq.~\eqref{eq:ResL3}, are formed by the elements:
\begin{align}
    \dfrac{\partial \mathcal{F}_j^{3C}}{\partial C_r^k} & =  \omega \begin{cases}
        1, & \text{if } k=j \text{ and } \forall r  \\
        0, & \text{if } k \neq j \text{ and } \forall r
    \end{cases}, 
    & \dfrac{\partial \mathcal{F}_j^{3C}}{\partial D_h^k} & = \omega \begin{cases}
        1, & \text{if } k=j \text{ and } \forall h  \\
        0, & \text{if } k \neq j \text{ and } \forall h
    \end{cases}, \nonumber \\
    \dfrac{\partial \mathcal{F}_m^{3R}}{\partial C_r^k} & = \omega \begin{cases}
        1, & \text{if } r=m \text{ and } \forall k \\
        0, & \text{if } r \neq m \text{ and } \forall k
    \end{cases},
    & \dfrac{\partial \mathcal{F}_{M+d}^{3R}}{\partial D_h^k} & = \omega \begin{cases}
        1, & \text{if } h=d \text{ and } \forall k \\
        0, & \text{if } h \neq d \text{ and } \forall k
    \end{cases}, \nonumber \\
    \dfrac{\partial \mathcal{F}^{3SC}_{q1}}{\partial C_r^k} & = \omega \begin{cases}
        D_d^k, & \text{if } k=j \text{ and } \forall r \\
        0, & \text{if } k \neq j \text{ and } \forall r
    \end{cases}, 
    & \dfrac{\partial \mathcal{F}^{3SC}_{q1}}{\partial D_h^k} & = \omega \begin{cases}
        \sum_{m=1}^M C_m^k, & \text{if } k=j \text{ and } h=d \\
        0, & \text{if } k \neq j \text{ or } h\neq d
    \end{cases}, \nonumber \\
    \dfrac{\partial \mathcal{F}^{3FC}_{q2}}{\partial C_r^k} & = \omega \begin{cases}
        \sum_{d=1}^{N-M} D_d^r, & \text{if } k=j \text{ and } r=m \\
        0, & \text{if } k \neq j \text{ or } r\neq m
    \end{cases}, 
    & \dfrac{\partial \mathcal{F}^{3FC}_{q2}}{\partial D_h^k} & = \omega \begin{cases}
        C_m^k, & \text{if } k=j \text{ and } \forall h \\
        0, & \text{if } k \neq j \text{ and } \forall h
    \end{cases},
    \label{eq:Res3_ders}
\end{align}
and all other derivatives w.r.t. the SLFNN weights and biases are zeros; recall that $q1=(d-1)N+j$ and $q2=(m-1)N+j$.~Using Eqs.~(\ref{eq:Res1_ders}, \ref{eq:Res2_ders1}, \ref{eq:Res2_ders2}, \ref{eq:Res3_ders}), all the components to form the Jacobian matrix $\nabla_{\mathbf{R}}\boldsymbol{\mathcal{F}}$ are now available.~Note that the Jacobian matrix can be alternatively obtained by numerical differentiation schemes (Finite Differences) or automatic differentiation.

For the computation of the unknown parameters $\mathbf{R}$ with the LM iterative algorithm, we start the optimization process with an initial guess of parameters, say $\mathbf{R}_0$, and an initial damping factor, say $\lambda_0$.~At the $\nu$-th iteration, we compute the residual vector $\boldsymbol{\mathcal{F}}(\mathbf{R}_\nu)$ in Eq.~\eqref{eq:ResVec} using Eqs.~\eqref{eq:ResALL}
and the Jacobian matrix $
\nabla_{\mathbf{R}_\nu}\boldsymbol{\mathcal{F}}$ using Eqs.~(\ref{eq:Res1_ders}, \ref{eq:Res2_ders1}, \ref{eq:Res2_ders2}, \ref{eq:Res3_ders}).~Estimating the Hessian matrix as $(\nabla_{\mathbf{R}_\nu}\boldsymbol{\mathcal{F}})^\top \nabla_{\mathbf{R}_\nu}\boldsymbol{\mathcal{F}}$, the LM algorithm determines the search direction (an $(N^2+(L(N-M + 2) + 1)M)$-dim. $\mathbf{d}_\nu$ vector)  by solving the linearized system:
\begin{equation}
\left( \left(\nabla_{\mathbf{R}_{\nu}} \boldsymbol{\mathcal{F}}\right)^\top \nabla_{\mathbf{R}_{\nu}} \boldsymbol{\mathcal{F}} + \lambda_\nu diag(\left(\nabla_{\mathbf{R}_{\nu}} \boldsymbol{\mathcal{F}}\right)^\top \nabla_{\mathbf{R}_{\nu}} \boldsymbol{\mathcal{F}})\right) \mathbf{d}_\nu = - \left( \nabla_{\mathbf{R}_{\nu}} \boldsymbol{\mathcal{F}}\right)^\top \boldsymbol{\mathcal{F}}(\mathbf{R}_{\nu})
    \label{eq:LMsd}
\end{equation}
where $diag(\cdot)$ denotes the diagonal matrix.~Subsequently, the unknown parameters $\mathbf{R}_{\nu+1}$ and the damping factor $\lambda_{\nu+1}$ at the next iteration are updated according to the following conditions:
\begin{itemize}
    \item successful step: if $\lVert \boldsymbol{\mathcal{F}} (\mathbf{R}_\nu+\mathbf{d}_\nu)\rVert_{l^2}<\lVert\boldsymbol{\mathcal{F}}(\mathbf{R}_\nu)\rVert_{l^2}$, then $\mathbf{R}_{\nu+1} = \mathbf{R}_\nu+\mathbf{d}_\nu$ and $\lambda_{\nu+1}=\lambda_{\nu}/10$, or
    \item unsuccessful step: if $\lVert \boldsymbol{\mathcal{F}}(\mathbf{R}_\nu+\mathbf{d}_\nu)\rVert_{l^2}\ge\lVert \boldsymbol{\mathcal{F}}(\mathbf{R}_\nu)\rVert_{l^2}$, then $\mathbf{R}_{\nu+1} = \mathbf{R}_\nu$ and $\lambda_{\nu+1}=10\lambda_{\nu}$.
\end{itemize}
The convergence of the algorithm is considered achieved when either of the relative following stopping criteria are met: (i) $\lVert \boldsymbol{\mathcal{F}}(\mathbf{R}_{\nu+1}) - \boldsymbol{\mathcal{F}}(\mathbf{R}_{\nu}) \rVert_{l^2}<tol_F (1+\lVert\boldsymbol{\mathcal{F}}(\mathbf{R}_{\nu})\rVert_{l^2})$ and (ii) $\lVert \mathbf{R}_{\nu+1} - \mathbf{R}_{\nu} \rVert_{l^2}<tol_R (1+\lVert \mathbf{R}_{\nu}\rVert_{l^2})$

\subsection{Implementation of the PINN approach and the QSSA/PEA/CSP methods for the SIM approximations}

In this section, we provide all practical information for the implementation of the proposed PINN methodology, as well as regarding the GSPT model reduction methodologies of QSSA, PEA and CSP, which were employed here for the derivation of SIM approximations.~We provide a pseudo-algorithm in \cref{alg:Outline} outlining the basic steps for employing the PINN methodology to learn SIM approximations and evaluate their numerical accuracy.

\begin{algorithm}[!h] 
\footnotesize
\caption{Training and evaluation of SIM approximations constructed via the PINN methodology} 
\label{alg:Outline}
\begin{algorithmic}[1]
\Require $d\mathbf{z}/dt=\mathbf{F}(\mathbf{z})$ and $\Omega\subset\mathbb{R}^N$ \Comment{General form in Eq.~\eqref{eq:gen} and domain to approximate the SIM.}
\State Sample collocation points $\mathbf{z}^{(i)}\in \Omega$ from trajectories lying in the same ($N$-$M$) dim. SIM. \Comment{Use Eq.~\eqref{eq:CSPcrit} to determine $M$}
\State Split to training/validation sets, select hyper-parameters $\mathbf{Q}$ and initialize parameters $\mathbf{C}$ and $\mathbf{D}$ of linear transformations and $\mathbf{P}$ of the SLFNNs.
\State For the training set, minimize the PINN optimization problem in Eq.~\eqref{eq:minFunPI}, initially for $\omega\neq0$, as follows: 
\Repeat
\State Form the residual vectors $\boldsymbol{\mathcal{F}}(\mathbf{R}_\nu)$ corresponding to all three loss functions $\mathcal{L}_1$, $\mathcal{L}_2$ and $\mathcal{L}_3$, using Eq.~\eqref{eq:ResALL}.
\State Form the Jacobian matrix $\nabla_{\mathbf{R}_\nu} \boldsymbol{\mathcal{F}}$.~For symbolic differentiation of $\mathcal{L}_1$, $\mathcal{L}_2$ and $\mathcal{L}_3$, use Eqs.~\eqref{eq:Res1_ders}, (\ref{eq:Res2_ders1}, \ref{eq:Res2_ders2}) and \eqref{eq:Res3_ders}.  
\State Update the parameters $\mathbf{R}_\nu$ with the Levenberg-Marquardt iteration; solve the linearized system in Eq.~\eqref{eq:LMsd}. 
\Until{convergence}
\State Having learned $\mathbf{C}$ and $\mathbf{D}$, set $\omega=0$ and continue training for $\mathbf{P}$ to get higher accuracy of the SLFNN-based SIM functional.
\State Sample test data sets $\mathbf{z}^{(i)}\in \Omega$ from trajectories; sample points with the same $M$ as in training. 
\State Evaluate the numerical accuracy of the PINN methodology by estimating the SIM approximation error $\lVert \mathbf{C}\mathbf{z}^{(i)}-\mathcal{N}(\mathbf{D} \mathbf{z}^{(i)})\rVert$.
\end{algorithmic}
\end{algorithm}

\subsubsection{Training process of the PINN}
\label{sbsb:trainPINN}

For learning the approximations of the SIMs with the proposed PINN-based methodology, data in the desired domain $\mathbf{z}^{(i)}\in \Omega \subset \mathbb{R}^N$ are required, serving as collocation points for $i=1,\ldots,n$.~However, there is no guarantee that for every point in the domain $\Omega$ there exists a transformation of $\mathbf{z}^{(i)}$, so that the SIM in Eq.~\eqref{eq:SIM} exists.~Thus, we collected the collocation points $\mathbf{z}^{(i)}$ from numerically derived trajectories.~In particular, for all problems considered, we selected $5^N$ initial conditions $\mathbf{z}_0$, randomly chosen outside $\Omega$ within specific - for each problem - bounds, and generated the trajectories by numerically solving the IVP in Eq.~\eqref{eq:gen}.~Then, we located the periods during which a constant (and same for each trajectory) $M$ was observed, using the criterion in Eq.~\eqref{eq:CSPcrit} with $\epsilon_{rel}=0.05$ and $\boldsymbol{\epsilon}_{abs}=10^{-10} \cdot \mathbf{1}$ and recorded within them 500 points per trajectory at equidistant - in time - points.~Finally, from the resulting datasets, we randomly chose $n$ points $\mathbf{z}^{(i)} \in \Omega$ to serve as collocation points for the PINN methodology.~To keep a similar ratio of number of residuals over number of unknown parameters for every problem under study, we set $n=500$ for the MM mechanism and $n=700$ for the TMDD and the fCSI mechanisms.     

The PINN schemes were trained using an 80\% uniformly random sample of the data sets described above, while the rest 20\% of the data set was used for validation.~For all problems, we used $L=20$ neurons in the hidden layer, with sigmoid activation functions, as already mentioned.~For determining the unknown parameters using the LM algorithm, we initialized the damping factor to $\lambda_0=0.01$ and the unknown parameters $\mathbf{R}_0$ in Eq.~\eqref{eq:TotPar} as follows: (i) the transformation related parameters were set to $C_m^j=D_d^j=1/N$ for $m=1,\ldots,M$, $d=1,\ldots,N-M$ and $j=1,\ldots,N$ and (ii) the SLFNN weights and biases were initialized using a uniform Xavier/Glorot initialization \cite{glorot2010understanding} for avoiding exploding or vanishing gradients.~Initially, we solved the optimization problem in Eq.~\eqref{eq:minFunPI} by setting $\omega=1$ (i.e., same contribution is provided by $\mathcal{L}_1$, $\mathcal{L}_2$ and $\mathcal{L}_3$) with the tolerances set to $tol_F=10^{-6}$ and $tol_R=10^{-4}$.~Then, using the parameters learned, we solved again the optimization problem with $\omega=100$ and the same tolerances, in order to learn more accurately the transformation; a maximum number of iterations $max_i=100$ was used.~Finally, from the resulting learned parameters, we rounded the elements of the transformation to the third decimal point and solved again the optimization problem, this time with $\omega=0$, $tol_F=10^{-10}$ and $tol_R=10^{-8}$, for optimizing the parameters of the SLFNN. 

To quantify the uncertainty in  the training performance of the PINNs, we performed 100 runs with different randomly sampled training and validation sets, and report the error of the residuals $\lVert \boldsymbol{\mathcal{F}} \rVert^2_2$, as well as the computational times needed for the training process.

\subsubsection{The SIM approximations provided by QSSA, PEA and CSP methodologies}
\label{sbsb:CSP_PEA_QSSA}

To assess the efficiency of the proposed scheme, we compared the SIM approximations provided by the proposed PINN scheme with the ones provided by the QSSA, PEA and CSP methods.~For deriving the latter, we followed the unifying CSP framework \cite{goussis2012quasi,goussis2015model,patsatzis2023algorithmic},  according to which the QSSA, PEA and CSP approximations are constructed by appropriate sets of basis vectors, spanning the fast and slow subspaces of the tangent space along which the dynamics of the original system in Eq.~\eqref{eq:gen} evolves.~The algorithmic procedure for the derivation of these sets is presented in \cref{app:CSP_SIM}.~We note here that the available in literature SIM approximations are the QSSA, PEA and CSP with one iteration for the first two problems \citep{patsatzis2023algorithmic,patsatzis2019new,patsatzis2016asymptotic} and the QSSA ones for the last problem \citep{rubinow1970time,pedersena2007total}, while all the rest SIM approximations were derived in this work; for the complete presentation of the SIM approximation of each problem considered, see \cref{app:MM_SIMs,app:TMDD_SIMs,app:Inh_SIMs}.

The derivation of SIM approximations in the context of QSSA, PEA and CSP methods requires specific assumptions.~In particular, the QSSA requires assumptions about the $M$ components of $\mathbf{z}$ that are considered fast (the ones being in quasi steady-state).~The same applies for the PEA with regard to the fast reactions (the ones being in partial equilibrium), so that assumptions on the components of the vector field $\mathbf{F}(\mathbf{z})$ are imposed; see \cref{app:sbPEAQSSA} for details.~Albeit erroneous assumptions of the above will not alter the algorithmic procedure, the resulting SIM expressions are deemed to be inaccurate \cite{goussis2012quasi}.~Secondly, even with the correct assumptions, the SIM approximations provided by QSSA and PEA might be inaccurate; accuracy is guaranteed from the  satisfaction of the conditions imposed by the CSP algorithmic criteria  \cite{goussis2012quasi,patsatzis2023algorithmic}.~Regarding the employment of CSP, no assumptions are required for deriving an accurate SIM approximation.~However, since CSP works on an iterative fashion, a good/bad initialization of the CSP algorithm will require less/more iterations for an accurate SIM approximation.~In particular, providing the correct assumption for the fast variables, one iteration of CSP guarantees the accuracy of the SIM approximation.~In the case where the assumption is erroneous, an additional iteration is required to guarantee accuracy \cite{goussis2012quasi}.~To avoid dependence on the fast variables assumption, here we perform CSP with two iterations.

%

On a final note, we highlight here that there is no guarantee that any of the QSSA, PEA and CSP approximations can be written as explicit functionals of the form $\mathbf{x}=\mathbf{h}(\mathbf{y})$.~In fact, for almost all the cases under study, the SIM expression derived via CSP with two iterations (that is the most accurate) is in an implicit form (i.e., $\mathbf{h}(\mathbf{x},\mathbf{y})=\mathbf{0}$), as well as the QSSA, PEA and CSP implementations for the third problem.

\subsection{Numerical accuracy of the SIM approximations}
\label{sb:NumAcc}

For assessing the numerical accuracy of the SIM approximations provided by the proposed PINN methodology, as well as the those obtained by the QSSA, PEA and CSP expressions, we constructed test data sets in the domain of interest $\Omega\subset \mathbb{R}^N$, consisting of numerically derived trajectories lying exclusively on the underlying $M$-dim. SIM.~In particular, we considered $5^N$ random initial conditions $\mathbf{z}_0$ chosen outside $\Omega$ within specific - for each problem - bounds, and generated the trajectories by numerically solving the IVP in Eq.\eqref{eq:gen}.~Then, we located the periods when a SIM of constant $M$ emerges, using the criterion in Eq.~\eqref{eq:CSPcrit} with $\epsilon_{rel}=0.05$ and $\boldsymbol{\epsilon}_{abs}=10^{-10} \cdot \mathbf{1}$, and recorded $500$ points per trajectory at equidistant - in time - points.~From the resulting datasets, we randomly chose $n_t=100\cdot5^N$ points $\mathbf{z}^{(i)} \in \Omega$ for $i=1,\ldots,n_t$ to form the test data sets. 

For the numerical SIM approximation accuracy of the PINN schemes, we evaluated the error through the output of the SLFNN function $\mathcal{N}(\cdot)$.~In particular, for every $\mathbf{z}^{(i)}$ in the test set, we compute the errors $\lVert \mathcal{X}(\mathbf{z}^{(i)})-\mathcal{N}(\mathcal{Y}(\mathbf{z}^{(i)})) \rVert = \lVert \mathbf{C}\mathbf{z}^{(i)}-\mathcal{N}(\mathbf{D} \mathbf{z}^{(i)}) \rVert$, which encapsulate the errors introduced by both the SLFNN function and the transformations. 

For the QSSA, PEA and CSP constructed SIM approximations, recall that the fast and slow variables, $\mathbf{x}$ and $\mathbf{y}$, are assumed to be some of the components of the original variables in $\mathbf{z}$.~Hence, to assess the numerical SIM approximation accuracy of the GSPT expressions, we first decomposed the observations of the test set  $\mathbf{z}^{(i)}$ to $\mathbf{x}^{(i)}\in\mathbb{R}^M$  and $\mathbf{y}^{(i)}\in\mathbb{R}^{(N-M)}$.
~In the case where the QSSA, PEA and CSP approximations provided explicit functionals $\mathbf{x}=\mathbf{h}(\mathbf{y})$, we computed the  errors $\lVert \mathbf{x}^{(i)}-\mathbf{h}( \mathbf{y}^{(i)}) \rVert$.~In the case where the SIM approximations were derived in the implicit form $\mathbf{h}(\mathbf{x},\mathbf{y})$, in order to provide a straightforward comparison, we numerically obtained a numerical estimation of $\tilde{\mathbf{x}}^{(i)}$ by solving the implicit functionals w.r.t. the, assumed, fast variables for every point $i=1,\ldots,n_t$.~In particular, we employed Newton-Raphson iterations by providing the derivatives of the SIM expressions w.r.t. $\mathbf{x}$ (symbolic differentiation for the analytic SIM expressions and finite differences for the numerical SIM expressions) and an initial guess, which is not always trivial for guaranteeing convergence.~After numerically estimating the roots $\tilde{\mathbf{x}}^{(i)}$ of the implicit functional $\mathbf{h}(\tilde{\mathbf{x}}^{(i)},\mathbf{y}^{(i)})$, we computed the explicit error $\lVert \mathbf{x}^{(i)}-\tilde{\mathbf{x}}^{(i)} \rVert$, for $i=1,\ldots,n_t$. 

\section{The benchmark problems}
\label{sec:Prob}

The efficiency of the proposed PINN numerical scheme is demonstrated via three benchmark problems.~First, we consider the Michaelis-Menten (MM)  mechanism, which is arguably the most extensively studied system in the context of SPT and GSPT; e.g. \citep{bowen1963singular,heineken1967mathematical,segel1989quasi,patsatzis2019new,schnell2000enzyme,gear2005projecting,borghans1996extending} just to name a few.~The $N=2$-dim. MM mechanism is an appropriate example to demonstrate the adaptiveness of the proposed framework, since it exhibits SIMs with different fast variables in different parameter regimes.~For our second problem, we consider the $N=3$-dim. Target Mediated Drug Disposition (TMDD) mechanism \cite{levy1994pharmacologic,mager2001general} to demonstrate generalization of the PINN scheme and showcase its ability to compute SIMs that do not contain the stable equilibrium of the system \citep{patsatzis2016asymptotic}, as is the case for the MM mechanism.~Both the MM and TMDD mechanisms have been studied in the context of GSPT (CSP in particular), thus allowing a comparison with the existing results \cite{patsatzis2019new,schnell2000enzyme,segel1989quasi,kristiansen2019geometric,patsatzis2016asymptotic,peletier2012dynamics}.~Finally, we consider a $N=4$-dim. enzyme reaction mechanism, the fully competitive substrate-inhibitor (fCSI) mechanism, which has been analyzed in the context of SPT, using the original system \citep{segel1988validity,schnell2000time,rubinow1970time}, but also an intuition-based transformed one \citep{pedersena2007total,bersani2017tihonov}, that can provide higher QSSA approximation accuracy.~By implementing the proposed methodology in both systems, we demonstrate that the accuracy of the PINN-derived SIM approximations are not affected by the form of the system and always provide higher approximation accuracy than the GSPT-generated SIMs; the PEA and CSP approximations are derived for the first time, in our work here.~In the following, we describe the three problems and provide the GSPT approximations of the SIM, derived on the basis of QSSA, PEA and CSP with one and two iterations, which are later used for comparison with the PINN scheme.

\subsection{The Michaelis-Menten mechanism}
\label{sub:MMdes}

The Michaelis-Menten (MM) mechanism is a chemical reaction scheme describing the basic mechanism of enzyme action \citep{michaelis1913kinetik,henri1903lois}.~According to the MM mechanism, an enzyme $E$ reversibly binds to a substrate $S$ for the formation of a complex $C$.~The complex can then decay irreversibly to a product $P$ and the same enzyme $E$, which is now free to bind to another substrate molecule.~The resulting MM reaction scheme is:
\begin{equation*}
    \ce{S + E <=>[k_{1f}][k_{1b}] C ->[k_2] E + P}, 
\end{equation*}
where $k_{1f}$ and $k_{1b}$ are the formation and dissociation rate constants and $k_2$ is the catalysis rate constant.~Using the chemical kinetics law of mass action and the conservation laws for the enzyme and the substrate, the MM mechanism takes the form of Eq.~\eqref{eq:gen} for the concentrations of the substrate $s$ and the complex $c$ as:
\begin{equation}
   \dfrac{d}{dt} \begin{bmatrix} s \\ c \end{bmatrix} = \begin{bmatrix} -k_{1f} (e_0-c) s + k_{1b} c~~~~~~~~~ \\ ~~~~~k_{1f} (e_0-c) s - k_{1b} c - k_2 c \end{bmatrix}, \qquad \qquad \begin{matrix}
       s(0) = s_0 \\
       c(0) = c_0
   \end{matrix} ~~,
   \label{eq:MM}
\end{equation}
where $e_0$, $s_0$ and $c_0<e_0$ are the initial concentrations of the enzyme, substrate and complex, respectively.~The $N=2$-dim. IVP in Eq.~\eqref{eq:MM} exhibits multi-scale character almost everywhere in the phase and parameter spaces, resulting to the emergence of an $M=1$-dim. SIM \citep{patsatzis2019new}.~The MM mechanism has been extensively studied in the context of GSPT and various analytic QSSA, PEA and CSP approximations of the SIM have been derived \citep{bowen1963singular,heineken1967mathematical,segel1989quasi,patsatzis2019new,schnell2000enzyme,gear2005projecting,borghans1996extending}.~However, these SIM approximations are accurate within specific regions of the parameter and phase spaces, depending on the fast dynamics (fast variable, $s$ or $c$, and fast reaction, the binding or the product formation one) of the MM mechanism in these regions \citep{patsatzis2019new}.

For comparison with the proposed PINN scheme, we considered all possible assumptions on the $M=1$ fast variables/reactions, to derive two QSSA, one PEA and five CSP approximations of the SIM, as discussed in Section~\ref{sbsb:CSP_PEA_QSSA}.~The rQSSA, CSP$_s$(1) and CSP$_s$(2) - with one and two iterations - approximations are constructed by assuming $s$ as the fast variable, while the sQSSA, CSP$_c$(1) and CSP$_c$(2) ones, are constructed by assuming $c$ as the fast variable.~The PEA approximation requires the assumption of the reversible, first, reaction to be fast (no assumption on fast variable), while the CSP$_e$ approximation, constructed on the basis of the eigenvectors, requires no assumptions.~Table~\ref{tb:MM_SIMs} enlists the above assumptions for the SIM approximations of the MM mechanism, along with their analytic implicit expressions, as constructed in \cref{app:MM_SIMs}.~It is further indicated that some of the SIM approximations cannot be solved explicitly w.r.t. the system's variables $s$ and $c$; thus, requiring as explained, numerical estimation with fixed-point (Newton) iterations.
\begin{table}[!h]
    \centering
    \resizebox{\textwidth}{!}{
    \begin{tabular}{l|c c c c c c c c}
    \toprule
    SIM approx. & rQSSA & sQSSA & PEA & CSP$_s$(1) & CSP$_s$(2) & CSP$_c$(1) & CSP$_c$(2) & CSP$_e$ \\
    \midrule
    Assumptions & $s$ fast & $c$ fast & 1st reac. fast & $s$ fast & $s$ fast & $c$ fast & $c$ fast & $-$ \\
    Expression & \eqref{eq:MM_rQSSA_imp} & \eqref{eq:MM_sQSSA_imp} & \eqref{eq:MM_PEA_imp} & \eqref{eq:MM_CSPs11_imp} & \eqref{eq:MM_CSPs21_imp} & \eqref{eq:MM_CSPc11_imp} & \eqref{eq:MM_CSPc21_imp} & \eqref{eq:MM_CSPe_imp}\\
    Derived in & \citep{schnell2000enzyme} & \citep{bowen1963singular,heineken1967mathematical} & \citep{patsatzis2019new} & here & here & here & here & \citep{patsatzis2019new} \\
    \midrule
    Explicit fun. & $s$, $c$ & $s$, $c$ & $s$, $c$ & $s$, $c$ & $s$ & $s$, $c$ & $-$ & $s$ \\
    Expression & (\ref{eq:MM_rQSSA_expS}, \ref{eq:MM_rQSSA_expC}) & (\ref{eq:MM_sQSSA_expS}, \ref{eq:MM_sQSSA_expC}) & (\ref{eq:MM_PEA_expS}, \ref{eq:MM_PEA_expC}) & (\ref{eq:MM_CSPs11_expS}, \ref{eq:MM_CSPs11_expC}) & \eqref{eq:MM_CSPs21_expS} & (\ref{eq:MM_CSPc11_expS}, \ref{eq:MM_CSPc11_expC}) & $-$ & \eqref{eq:MM_CSPe_expS} \\
    Requires Newton & & & & & $c$ & & $s$, $c$& $c$ \\ 
    \bottomrule
    \end{tabular}}
    \caption{SIM approximations for the MM mechanism, constructed on the basis of QSSA (for $s$ and $c$), PEA (for the first revirsible reaction) and CSP (assuming $s$ and $c$ fast with one and two iterations and using the eigenvectors) methods.~The assumptions made for constructing each approximation, the implicit functional expression and the reference where it is derived, are enlisted.~In addition, the explicit functional form w.r.t. either $s$ or $c$ is provided; when the latter is not available analytically, Newton iterations are required to numerically solve the SIM approximation w.r.t. $s$ or $c$.}
    \label{tb:MM_SIMs}
\end{table}

\begin{figure}[!h]
    \centering
    \subfigure[MM1 case]{
    \includegraphics[width=0.32\textwidth]{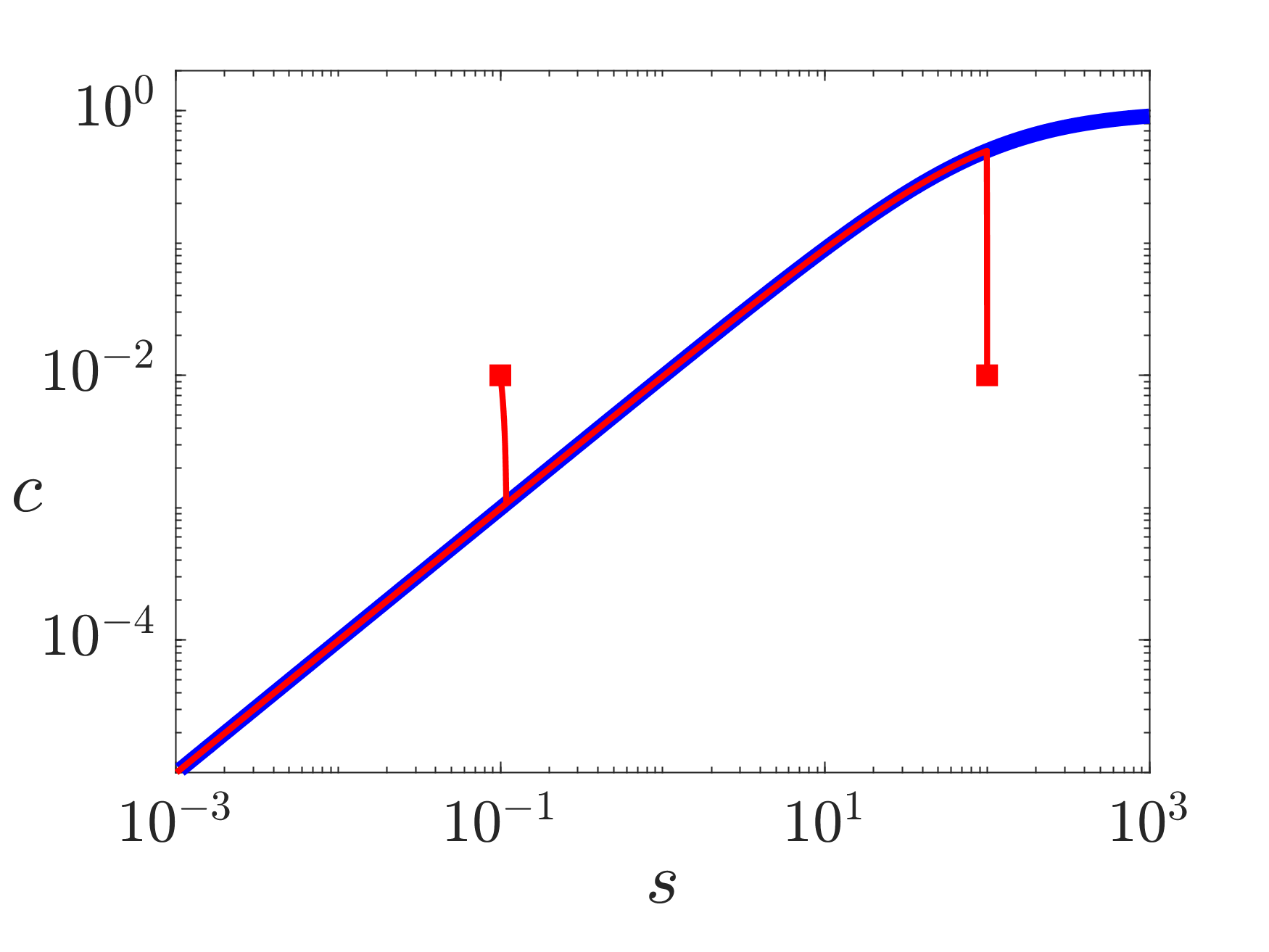}}
    \subfigure[MM2 case]{
    \includegraphics[width=0.32\textwidth]{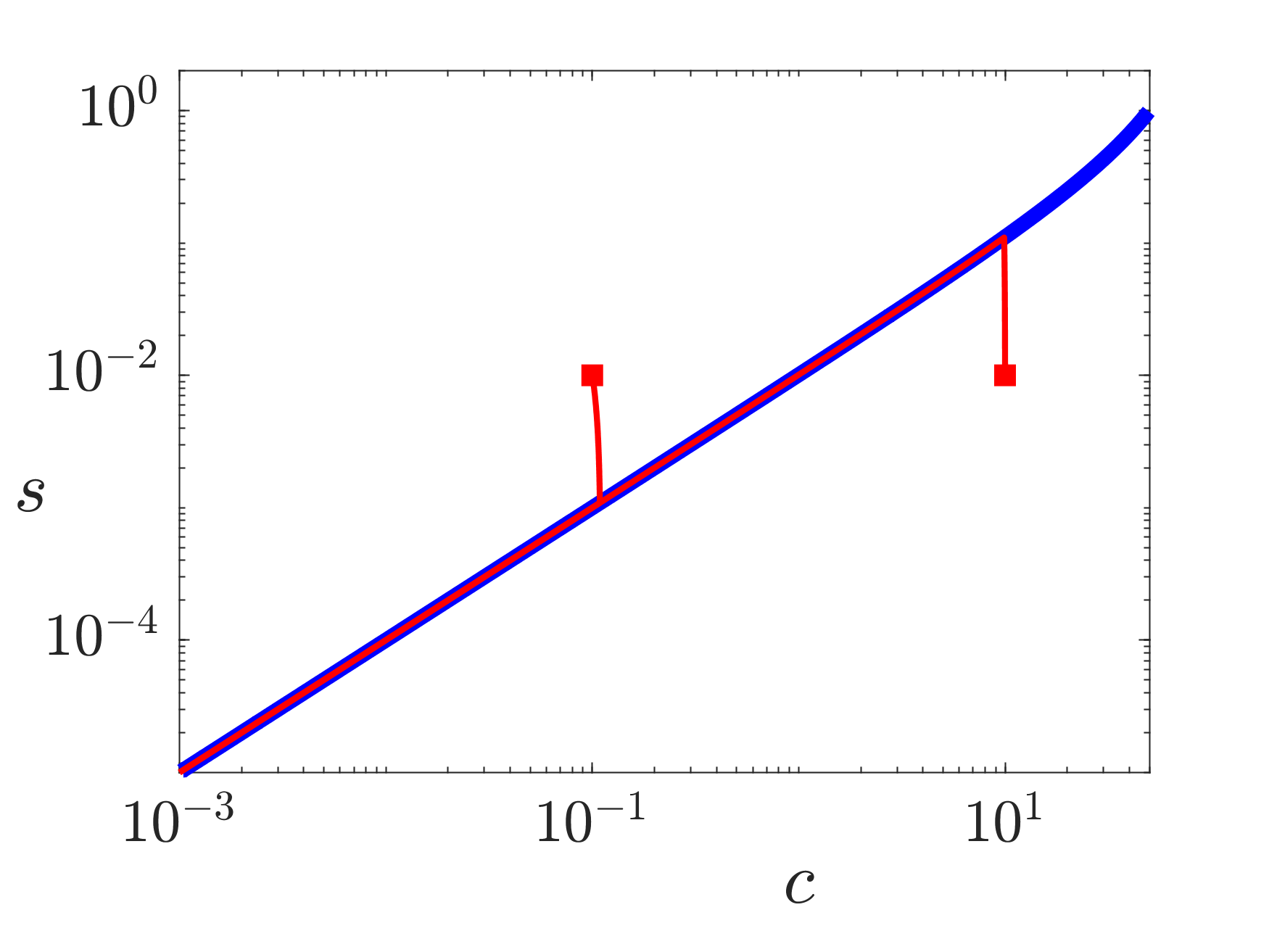}}
    \subfigure[MM3 case]{
    \includegraphics[width=0.32\textwidth]{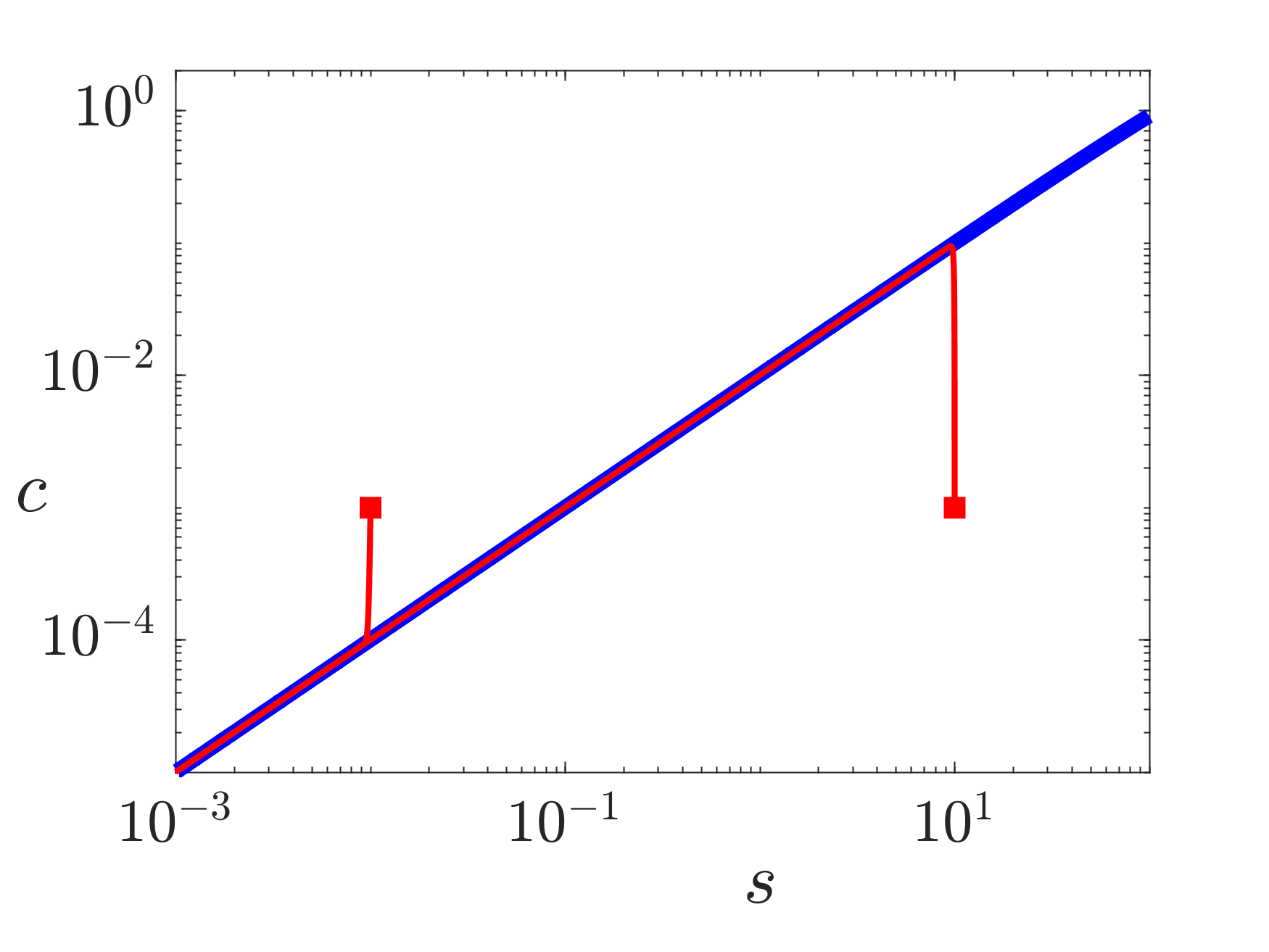}}
    \caption{The $M=1$-dim. SIM hypersurface arising in the phase space of the MM mechanism in Eq.~\eqref{eq:MM} for the three cases considered.~The red trajectories (starting at red squares) are attracted to the SIM along the direction of the fast variable (on y-axis) and then evolve on it towards reaching the fixed point $(0,0)$.}
    \label{fig:MM_SIMs}
\end{figure}

Here, to demonstrate the adaptability of the proposed PINN scheme to different regions of the phase and parameter spaces, we selected three indicative cases of, known from \citep{patsatzis2019new,patsatzis2023algorithmic}, different fast/slow dynamics.~For each case, different analytic expressions of the SIM in Table~\ref{tb:MM_SIMs} are expected to be accurate.~In particular, we consider:
\begin{enumerate}
    \item the MM1 case, where the parameter set is $k_{1f}=1$, $k_{1b}=10^2$, $k_2=1$, $e_0=1$ and the region of interest in the phase space is $\Omega = [10^{-3}, 10^3] \times [10^{-5}, 1]$.~In this case, the variable $c$ and the first binding reaction are fast.~Thus, we expect accurate SIM approximations from sQSSA, PEA, CSP$_c$(1), CSP$_c$(2), CSP$_s$(2) and CSP$_e$ expressions.
    \item the MM2 case, where the parameter set is $k_{1f}=1$, $k_{1b}=1$, $k_2=10^{-2}$, $e_0=10^2$ and the region of interest in the phase space is $\Omega = [10^{-5}, 1] \times [10^{-3}, 50]$.~In this case, the variable $s$ and the first binding reaction are fast; rQSSA, PEA, CSP$_s$(1), CSP$_s$(2), CSP$_c$(2) and CSP$_e$ expressions are expected to provide accurate SIM approximations.
    \item the MM3 case, where the parameter set is $k_{1f}=1$, $k_{1b}=1$, $k_2=10^3$, $e_0=10$ and the region of interest in the phase space is $\Omega = [10^{-3}, 10^2] \times [10^{-5}, 1]$.~Here, the variable $c$ is fast, but not the first binding reaction.~Thus, we expect accurate SIM approximations from sQSSA, CSP$_c$(1), CSP$_c$(2), CSP$_s$(2) and CSP$_e$ expressions.
\end{enumerate}
The representation of the SIM emerging in the phase space, in the three cases of the MM mechanism described above, is shown in Fig.~\ref{fig:MM_SIMs}.~In all three cases, the trajectories are attracted to the $M=1$-dim. SIM along the direction of the fast variable, and then evolve on it until reaching the stable stationary point $(0,0)$ of the MM mechanism.

\subsection{The Target Mediated Drug Disposition mechanism}
\label{sub:TMDDdes}
The Target Mediated Drug Disposition (TMDD) mechanism is a pharmacokinetic-pharmacodynamic reaction scheme  describing the binding process of a drug, with high affinity, to its pharmacological target and the compound outcome of this interaction \citep{levy1994pharmacologic,mager2001general}.~The simplest form of the TMDD mechanism is one-compartmental \citep{mager2001general}, which can be effectively described by the reaction scheme: 
\begin{equation*}
    \ce{L + R <=>[k_{on}][k_{off}] RL}, \qquad \ce{L ->[k_{el}] }, \qquad \ce{->[k_{syn}] R ->[k_{deg}] }, \qquad \ce{RL ->[k_{int}] }, 
\end{equation*}
where $L$, $R$ and $RL$ are the ligand (drug), its pharmacological target (receptor) and the ligand-target complex (active form of the drug in the body), respectively.~The one-compartmental TMDD reaction scheme includes the binding (formation and dissociation; rate constants $k_{on}$ and $k_{off}$), the ligand elimination (first-order; $k_{el}$), the receptor synthesis and degradation (zeroth and first-order; $k_{syn}$ and $k_{deg}$, respectively) and the complex internalization (first-order; $k_{int}$) processes.~According to the above, the $N=3$-dim. system of the TMDD model is written in the form of Eq.~\eqref{eq:gen}, as:
\begin{equation}
    \dfrac{d}{dt}\begin{bmatrix} L~~ \\ R~~ \\ RL
    \end{bmatrix} = \begin{bmatrix} -k_{on}L.R+k_{off}RL-k_{el} L \qquad \quad ~~~\\ 
    -k_{on}L.R+k_{off}RL+k_{syn}-k_{deg}R \\
    k_{on}L.R-k_{off}RL - k_{int}RL\quad ~~~~
    \end{bmatrix}, \qquad \quad 
    \begin{matrix}
        L(0)=L_0~~~ \\
        R(0)=R_0~~~ \\
        ~~RL(0)=RL_0
    \end{matrix} ~~,
    \label{eq:TMDD}
\end{equation}
where $L_0$, $R_0$ and $RL_0$ are the initial concentrations of the drug, its target and their complex, respectively, and the parameter values $k_{on}=0.091$, $k_{off}=0.001$, $k_{el}=0.0015$, $k_{syn}=0.11$, $k_{deg}=0.0089$ and $k_{int}=0.003$ are adopted, according to  \citep{peletier2015challenges,peletier2012dynamics}.

\begin{figure}[!h]
    \centering
    \subfigure[SIM $\mathcal{M}_1$]{
    \includegraphics[width=0.33\textwidth]{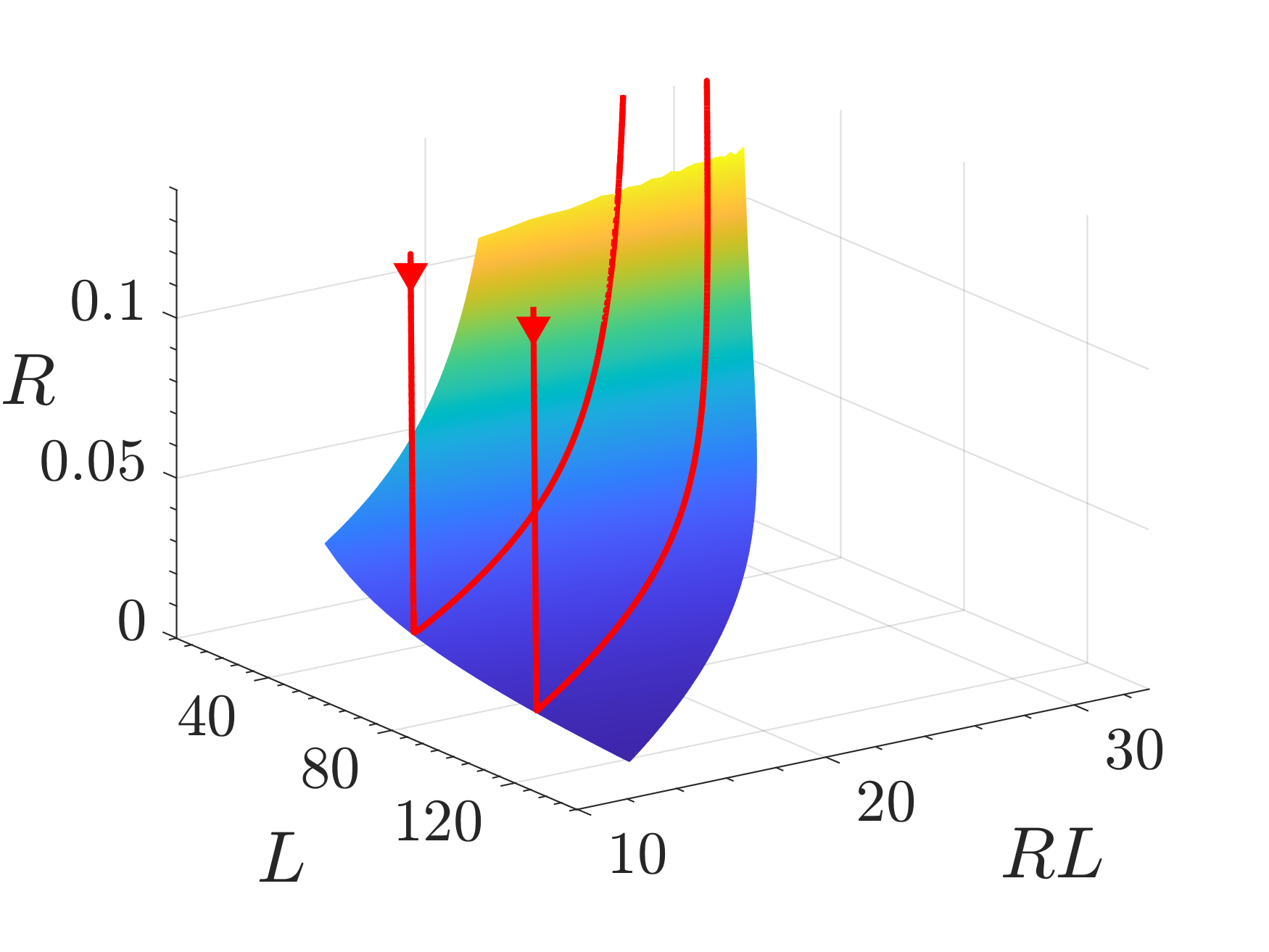}}
    \subfigure[SIM $\mathcal{M}_2$]{
    \includegraphics[width=0.33\textwidth]{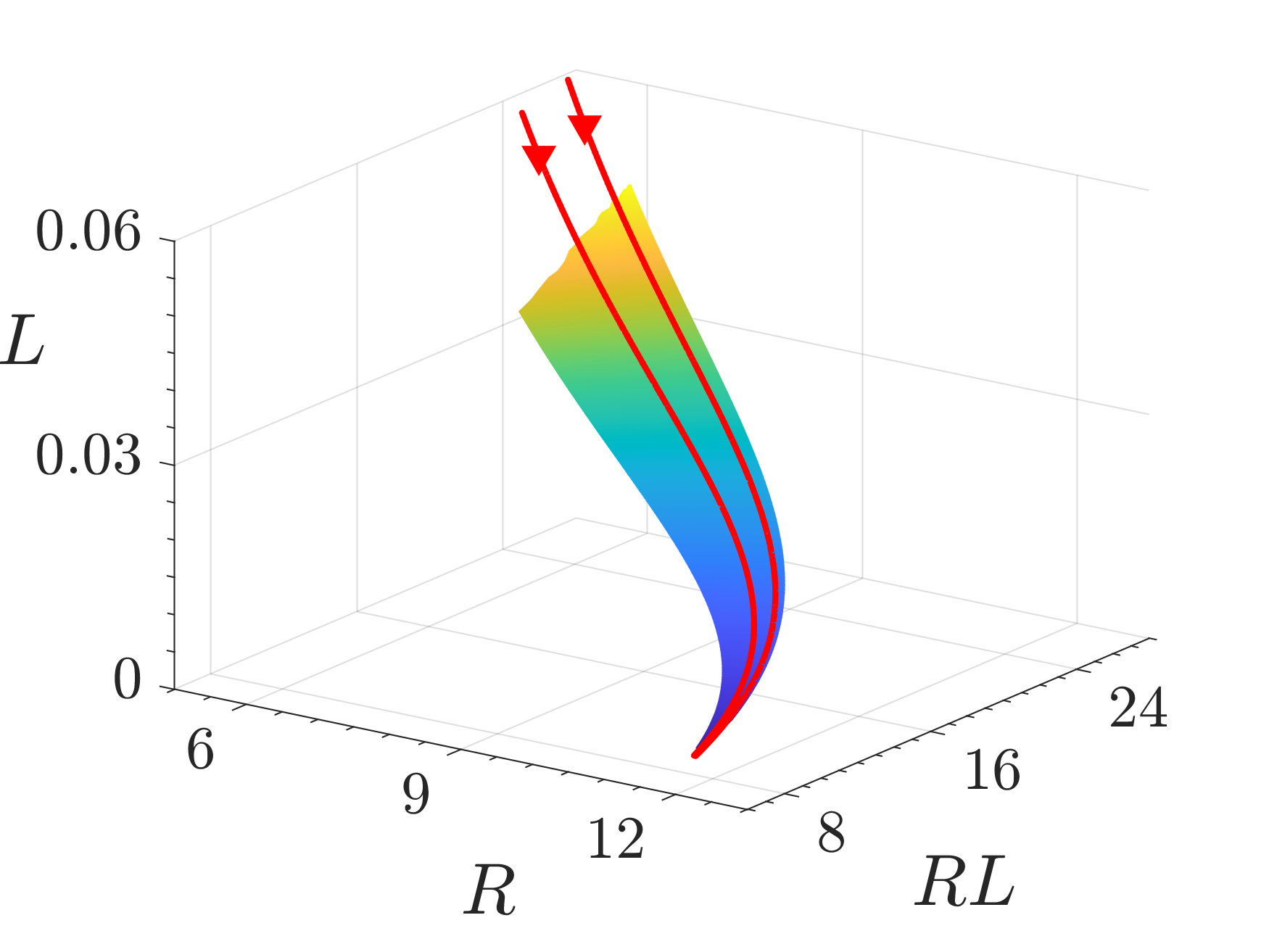}}
    \caption{The $M=1$-dim. SIM hypersurfaces of the TMDD mechanism in Eq.~\eqref{eq:TMDD}.~The red trajectories (initialized with $L_0=80/120$ at equilibrium $R_*=12$, $RL_*=0$) are attracted towards the SIM $\mathcal{M}_1$ and then, after evolving on it, they degenerate.~Then,  trajectories approach $\mathcal{M}_2$, along which they evolve until reaching the fixed point $(0,12,0)$.~The attraction towards $\mathcal{M}_1$ and $\mathcal{M}_2$ (red triangles) is mainly along the direction of the fast variable (z-axis).}
    \label{fig:TMDD_SIMs}
\end{figure}

The TMDD model in Eq.~\eqref{eq:TMDD} has been studied under dose-dependent initial conditions employed at equilibrium (i.e., for selected $L_0$ at $R_*=k_{syn}/k_{deg}=12$ and $RL_*=0$) \citep{peletier2015challenges,van2016topics} and it has been shown to exhibit multi-scale character in the phase space.~As a result, it has been extensively studied in the context of GSPT and various analytic QSSA and PEA approximations of the SIM have been proposed
\citep{kristiansen2019geometric,peletier2012dynamics,aston2011mathematical,gibiansky2008approximations,van2016topics,mager2001general,peletier2009dynamics,ma2012theoretical,mager2005quasi,patsatzis2016asymptotic}.~The systematic analysis in \citep{patsatzis2016asymptotic,kristiansen2019geometric} revealed that, under these initial conditions, the system initially undergoes an evolution along a $M=1$-dim. SIM, say $\mathcal{M}_1$, the fast dynamics of which relates to $R$ and the binding reaction.~However, this SIM does not include the equilibrium of the TMDD model $(L_*,R_*,RL_*) = (0,12,0)$ and, as a result, the system goes off  it, in order to subsequently approach another SIM, say $\mathcal{M}_2$, that leads the system to its equilibrium.~The fast dynamics on $\mathcal{M}_2$, which is initially $M=1$-dim., relates to $L$ and the binding reaction.~A representation of the two SIMs $\mathcal{M}_1$ and $\mathcal{M}_2$ arising in the phase space of the TMDD mechanism, as well as the dynamics of selected trajectories, are shown in Fig.~\ref{fig:TMDD_SIMs}.

For the construction of SIM approximations on the basis of GSPT, since both SIMs $\mathcal{M}_1$ and $\mathcal{M}_2$ are $M=1$-dim., we considered all possible assumptions on the $M=1$ fast variables/reactions, as discussed in Section~\ref{sbsb:CSP_PEA_QSSA}.~Assuming each of the three variables to be fast, we derived a SIM approximation on the basis of QSSA and two on the basis of CSP with one and two iterations, reading QSSA$_x$, CSP$_x$(1) and CSP$_x$(2) for $x=L$, $R$ or $RL$.~In addition, we derived the PEA approximation which requires only the assumption of the first, binding and reversible, reaction to be fast.~All the ten SIM approximations of the TMDD mechanism, along with their analytical implicit expressions, as constructed in \cref{app:TMDD_SIMs}, are listed in Table~\ref{tb:TMDD_SIMs}.~Since we know that either $L$ or $R$ are fast in any of $\mathcal{M}_1$ and $\mathcal{M}_2$, we did not solve the SIM expressions for $RL$, but only for $L$ and $R$.~We highlight here, that all CSP with two iterations approximations cannot be solved explicitly w.r.t. $L$ or $R$ and thus, require numerical estimation with fixed point (Newton) iterations; see last row of Table~\ref{tb:TMDD_SIMs}.
\begin{table}[!h]
    \centering
    \resizebox{\textwidth}{!}{
    \begin{tabular}{l|c c c c c c c c c c}
    \toprule
    SIM approx.	&	QSSA$_L$	&	QSSA$_R$	&	QSSA$_{RL}$	&	PEA	&	CSP$_L$(1)	&	CSP$_L$(2)	&	CSP$_R$(1)	&	CSP$_R$(2)	&	CSP$_{RL}$(1)	&	CSP$_{RL}$(2)	\\
    \midrule
    Assumptions	&	$L$ fast	&	$R$ fast	&	$RL$ fast	&	1st reac. fast	&	$L$ fast	&	$L$ fast	&	$R$ fast	&	$R$ fast	&	$RL$ fast	&	$RL$ fast	\\
    Expression	& \eqref{eq:TMDD_QSSAL_imp}	&	\eqref{eq:TMDD_QSSAR_imp}	&	\eqref{eq:TMDD_QSSARL_imp}	&	\eqref{eq:TMDD_PEA_imp}	  &	\eqref{eq:TMDD_CSPL11_imp}	&	\eqref{eq:TMDD_CSPL21_imp}	&	\eqref{eq:TMDD_CSPR11_imp}	&	\eqref{eq:TMDD_CSPR21_imp}	&	\eqref{eq:TMDD_CSPRL11_imp}	&	\eqref{eq:TMDD_CSPRL21_imp}	\\
    Derived in	&	\citep{peletier2012dynamics}	&	\citep{aston2011mathematical}	&	\citep{gibiansky2008approximations}	&	\citep{mager2001general}	&	\citep{patsatzis2016asymptotic}	&	here	&	\citep{patsatzis2016asymptotic}	&	here	&	here	&	here	\\
    \midrule
    Explicit fun.	&	$L$, $R$	&	$L$, $R$	&	$L$, $R$	&	$L$, $R$	&	$L$, $R$	&	$-$	&	$L$, $R$	&	$-$	&	$L$, $R$	&	$-$	\\
    Expression	&	(\ref{eq:TMDD_QSSAL_expL}, \ref{eq:TMDD_QSSAL_expR})	&	(\ref{eq:TMDD_QSSAR_expL}, \ref{eq:TMDD_QSSAR_expR})	&	(\ref{eq:TMDD_QSSARL_expL}, \ref{eq:TMDD_QSSARL_expR})	&	(\ref{eq:TMDD_PEA_expL}, \ref{eq:TMDD_PEA_expR})	   &	(\ref{eq:TMDD_CSPL11_expL}, \ref{eq:TMDD_CSPL11_expR})	&	$-$	&	(\ref{eq:TMDD_CSPR11_expL}, \ref{eq:TMDD_CSPR11_expR}) &	$-$ &	(\ref{eq:TMDD_CSPRL11_expL}, \ref{eq:TMDD_CSPRL11_expR}) & $-$	\\
    Requires Newton	&		&		&		&		&		&	$L$, $R$	&		&	$L$, $R$	&		&	$L$, $R$	\\ 
    \bottomrule
    \end{tabular}}
    \caption{SIM approximations for the TMDD mechanism, constructed on the basis of QSSA (for $L$, $R$ and $RL$), PEA (for the first binding reaction) and CSP (assuming $L$, $R$ or $RL$ fast with one and two iterations) methods.~The assumptions made for constructing each approximation, the implicit functional expression and the reference where it is derived, are enlisted.~In addition, the explicit functional form w.r.t. either $L$ or $R$ is provided; when the latter is not available analytically, Newton iterations are required to numerically solve the SIM approximation w.r.t. $L$ or $R$.}
    \label{tb:TMDD_SIMs}
\end{table}

Here, we approximated both $\mathcal{M}_1$ and $\mathcal{M}_2$ SIMs, and located the regions of the phase space where they arise, using the criterion in in Eq.~\eqref{eq:CSPcrit}; i.e., we detect, along the trajectory, the periods of $M=1$ separated by periods of $M=0$.~Following the above, we consider the domains $\Omega=[10,140]\times[0,0.12]\times[10,33]$ for $\mathcal{M}_1$ and $\Omega=[0,0.05]\times[6,13]\times[6,27]$ for $\mathcal{M}_2$.~Finally, regarding the GSPT expressions, we expect accurate SIM approximations provided by QSSA$_R$, PEA, CSP$_R$(1) and CSP$_R$(2) expressions in $\mathcal{M}_1$ (since $R$ is the fast variable there) and by QSSA$_L$, PEA, CSP$_L$(1) and CSP$_L$(2) expressions in $\mathcal{M}_2$ (since $L$ is the fast variable there) \citep{patsatzis2016asymptotic,kristiansen2019geometric}.

\subsection{The fully Competitive Substrate-Inhibitor mechanism}
\label{sub:CompInhdes}

Competitive substrate inhibition is the phenomenon of two or more substrates competing for binding with the same enzyme molecule \citep{cornish2013fundamentals}.~As a natural extension of the MM mechanism in Section~\ref{sub:MMdes} for complex intracellular networks, we consider the fully Competitive Substrate-Inhibitor (fCSI) mechanism where two substrates $S_1$ and $S_2$ can bind to the same site of an enzyme molecule $E$, but not in the same time \citep{cornish2013fundamentals}; hence $S_1$ inhibits $S_2$ and vice versa.~Both substrates follow the MM kinetics, thus resulting to the reaction scheme:
\begin{equation*}
    \ce{S_1 + E <=>[k_{1f}][k_{1b}] C_1 ->[k_2] E + P_1}, \qquad
    \ce{S_2 + E <=>[k_{3f}][k_{3b}] C_2 ->[k_4] E + P_2},
\end{equation*}
where $C_1$ and $C_2$ are the enzyme-substrate complexes and $P_1$ and $P_2$ the related products.~Using the chemical kinetics law of mass action and the conservation laws for the enzyme and the two substrates, the fCSI mechanism is formulated for the concentrations of the two substrates $s_1$, $s_2$ and the two complexes $c_1$, $c_2$ in the form of Eq.~\eqref{eq:gen} as:
\begin{equation}
    \dfrac{d}{dt} \begin{bmatrix} s_1 \\ c_1 \\ s_2 \\ c_2   \end{bmatrix} = 
    \begin{bmatrix}
        -k_{1f}(e_0-c_1-c_2) s_1 + k_{1b}c_1 \qquad \quad ~~~  \\ k_{1f}(e_0-c_1-c_2) s_1 - k_{1b}c_1-k_2c_1 \\
        -k_{3f}(e_0-c_1-c_2) s_2 + k_{3b}c_2 \qquad \quad ~~~ \\ k_{3f}(e_0-c_1-c_2) s_2 - k_{3b}c_2-k_4c_2
    \end{bmatrix}, \qquad \quad 
    \begin{matrix}
        s_1(0) = s_{10} \\
        c_1(0) = c_{10} \\
        s_2(0) = s_{20} \\
        c_2(0) = c_{20} 
    \end{matrix} ~~, 
    \label{eq:Inh_st}
\end{equation}
where $s_{10}$, $c_{10}$, $s_{20}$ and $c_{20}$ are the initial concentrations, with $c_{10}+c_{20}<e_0$.~Correspondingly to the MM mechanism, the parameters $k_{1f}$ ($k_{3f}$), $k_{1b}$ ($k_{3b}$) and $k_2$ ($k_4$) denote the formation, dissociation and catalysis rate constants of the first (second) substrate, respectively.

Due to its similarity to the MM mechanism, the system in Eq.~\eqref{eq:Inh_st} has been studied in the context of SPT, and QSSA was implemented with the assumption that $c_1$ and $c_2$ are the fast variables of the emergent $M=2$-dim. SIM \citep{segel1988validity,schnell2000time,rubinow1970time}; corresponding to sQSSA in the MM mechanism.~However, it has been shown that the implementation of the QSSA in Eq.~\eqref{eq:Inh_st} is not accurate in high enzyme concentrations or in high enzyme affinities \citep{pedersena2007total}.~In such cases, the tQSSA approach  \citep{borghans1996extending,tzafriri2003michaelis} is followed, according to which the original system is first linearly transformed, using $\bar{s}_1=s_1+c_1$ and $\bar{s}_2=s_2+c_2$, which results in the system:
\begin{equation}
    \dfrac{d}{d} \begin{bmatrix} \bar{s}_1 \\ c_1 \\ \bar{s}_2 \\ c_2   \end{bmatrix} = 
    \begin{bmatrix}
        -k_2 c_1 \\ k_{1f}(e_0-c_1-c_2) (\bar{s}_1-c_1) - k_{1b}c_1-k_2c_1 \\
        -k_4 c_2 \\ k_{3f}(e_0-c_1-c_2) (\bar{s}_2-c_2) - k_{3b}c_2-k_4c_2
    \end{bmatrix}, \qquad \quad \begin{matrix}
        \bar{s}_1(0) = s_{10}+ c_{10} \\
        c_1(0) = c_{10}~~~~~~~~~~ \\
        \bar{s}_2(0) = s_{20} + c_{20} \\
        c_2(0) = c_{20}~~~~~~~~~~
    \end{matrix} ~~. 
    \label{eq:Inh_tr}
\end{equation}
The transformed system in Eq.~\eqref{eq:Inh_tr} has been studied in the context of GSPT \citep{pedersena2007total,bersani2017tihonov}, and the QSSA was shown more accurate in comparison to the one employed for the original system in Eq.~\eqref{eq:Inh_st}.

For comparison purposes to the proposed PINN scheme, considering all possible combinations of $M=2$ fast variables/reactions in the context of GSTP, that was performed in the previous problems, would result to 19 SIM approximations for each of the original and transformed systems.~Thus, we only considered the $M=2$-dim. SIM approximations that can be constructed with the assumption of $c_1$ and $c_2$ being the fast variables; this allows for the consideration of the QSSA$_{c1c2}$ approximations reported in literature \citep{segel1988validity,pedersena2007total}.~On the basis of the above assumption, we further constructed (i) the CSP approximations with one and two iterations, namely CSP$_{c1c2}$(1) and CSP$_{c1c2}$(2), respectively, and (ii) the PEA$_{13}$ approximations, by further assuming that the reversible, 1st and 3rd, binding reactions are also fast.~Table~\ref{tb:Inh_SIMs} enlists the above SIM approximations as constructed in \cref{app:Inh_SIMs}, along with assumptions required for their construction and their analytic implicit expressions.
\begin{table}[!h]
    \centering
    \resizebox{\textwidth}{!}{
    \begin{tabular}{l|c c c c | c c c c}
    \toprule
    & \multicolumn{4}{c|}{original system in Eq.~\eqref{eq:Inh_st}} & \multicolumn{4}{c}{tranformed system in Eq.~\eqref{eq:Inh_tr}} \\
    \midrule
    SIM approx.	&	QSSA$_{c1c2}$	&	PEA$_{13}$	&	CSP$_{c1c2}$(1)	&	CSP$_{c1c2}$(2)	&	QSSA$_{c1c2}$	&	PEA$_{13}$	&	CSP$_{c1c2}$(1)	&	CSP$_{c1c2}$(2)	\\
    \midrule
    Assumptions	&	$c_1$, $c_2$ fast	&	1st, 3rd reac. fast	&	$c_1$, $c_2$ fast	&	$c_1$, $c_2$ fast	&	$c_1$, $c_2$ fast	&	1st, 3rd reac. fast	&	$c_1$, $c_2$ fast	&	$c_1$, $c_2$ fast	\\
    Expression	&	(\ref{eq:InhSt_QSSAc1c2_imp})	&	(\ref{eq:InhSt_PEA13_imp1}, \ref{eq:InhSt_PEA13_imp2})	&	(\ref{eq:InhSt_CSP11_imp1}, \ref{eq:InhSt_CSP11_imp2})	&	numerical	&	(\ref{eq:InhTr_QSSAc1c2_imp1}, \ref{eq:InhTr_QSSAc1c2_imp2})	&		(\ref{eq:InhTr_PEA13_imp1}, \ref{eq:InhTr_PEA13_imp2})	&	(\ref{eq:InhTr_CSP11_imp1}, \ref{eq:InhTr_CSP11_imp2})	&	numerical	\\
    Derived in	&	\citep{rubinow1970time}	&	here	&	here	&	here	&	\citep{pedersena2007total}	&	here	&	here	&	here	\\
    \midrule
    Explicit fun.	&	$c_1$ and $c_2$	&	$-$	&	$-$	&	$-$	&	$-$	&	$-$	&	$-$	&	$-$	\\
    Expression	&	(\ref{eq:InhSt_QSSAc1c2_exp1}, \ref{eq:InhSt_QSSAc1c2_exp2})	&	$-$	&	$-$	&	$-$	&	$-$	&	$-$	&	$-$	&	$-$	\\
    Requires Newton	&		&	$c_1$ and $c_2$	&	$c_1$ and $c_2$	&	$c_1$ and $c_2$	&	$c_1$ and $c_2$	&	$c_1$ and $c_2$	&	$c_1$ and $c_2$	&	$c_1$ and $c_2$	\\
    \bottomrule
    \end{tabular}}
    \caption{SIM approximations for the original and transformed systems in Eq.~(\ref{eq:Inh_st}, \ref{eq:Inh_tr}) of the fCSI mechanism, constructed on the basis of QSSA (for $c_1$ and $c_2$ fast), PEA (for the first and third binding reactions) and CSP (for $c_1$ and $c_2$ fast with one and two iterations) methods.~The assumptions made for constructing each approximation, the system of implicit functional expressions and the reference where it is derived, are enlisted.~Only QSSA$_{c1c2}$ for the original system provides a system of explicit functional forms w.r.t. $c_1$ and $c_2$, while all the remaining expressions, require Newton iterations to numerically solve the SIM approximation w.r.t. $c_1$ and $c_2$.~Note that the CSP expressions for two iterations are numerically derived.}
    \label{tb:Inh_SIMs} 
\end{table}
We further indicate that all SIM approximations, except from the QSSA$_{c1c2}$ of the original system, require Newton iterations (for a $M=2$-dim. system, in this case) for the numerical estimation of the SIM explicit functional w.r.t. the assumed fast variables $c_1$ and $c_2$.~We highlight here that the CSP with two iterations approximations are derived numerically for each point in the phase space, and, since they are implicit, the employment of Newton iterations requires the Jacobian, which is computed with finite differences; thus, requiring high computational cost.

\begin{figure}[!h]
    \centering
    \subfigure[Original system, projection to $c_1$]{
    \includegraphics[width=0.33\textwidth]{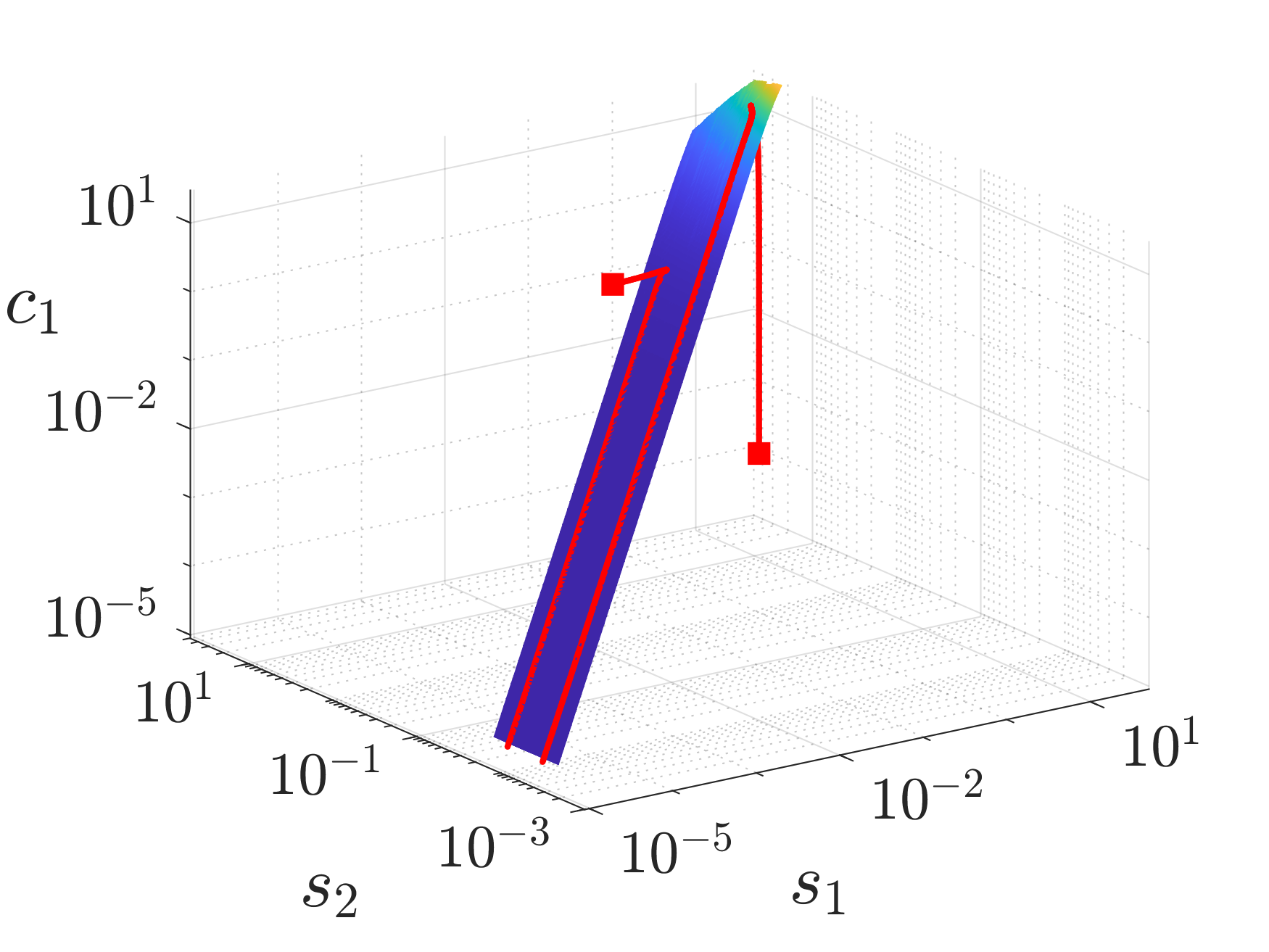}}
    \subfigure[Original system, projection to $c_2$]{
    \includegraphics[width=0.33\textwidth]{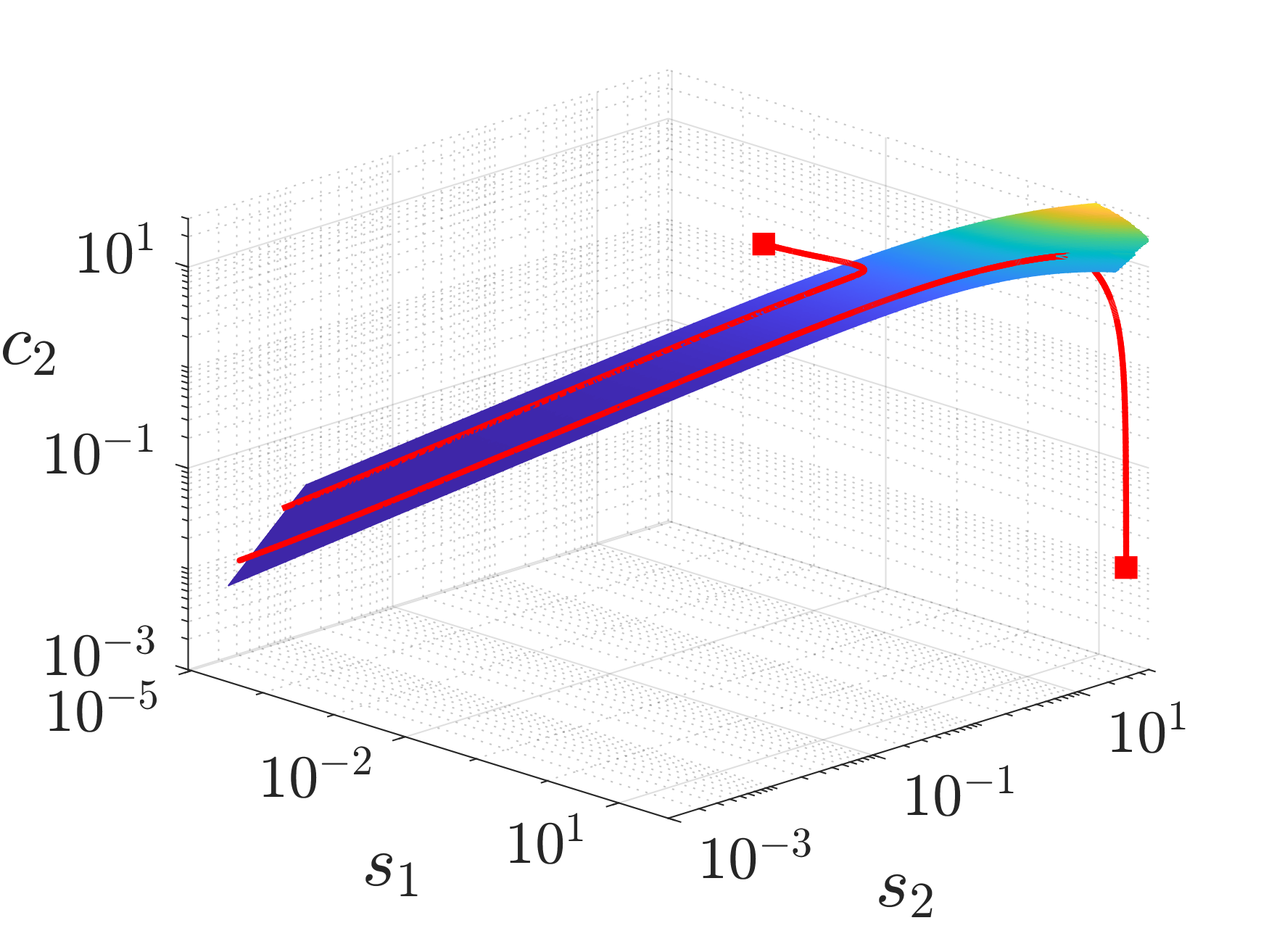}} \\
    \subfigure[Transformed system, projection to $c_1$]{
    \includegraphics[width=0.33\textwidth]{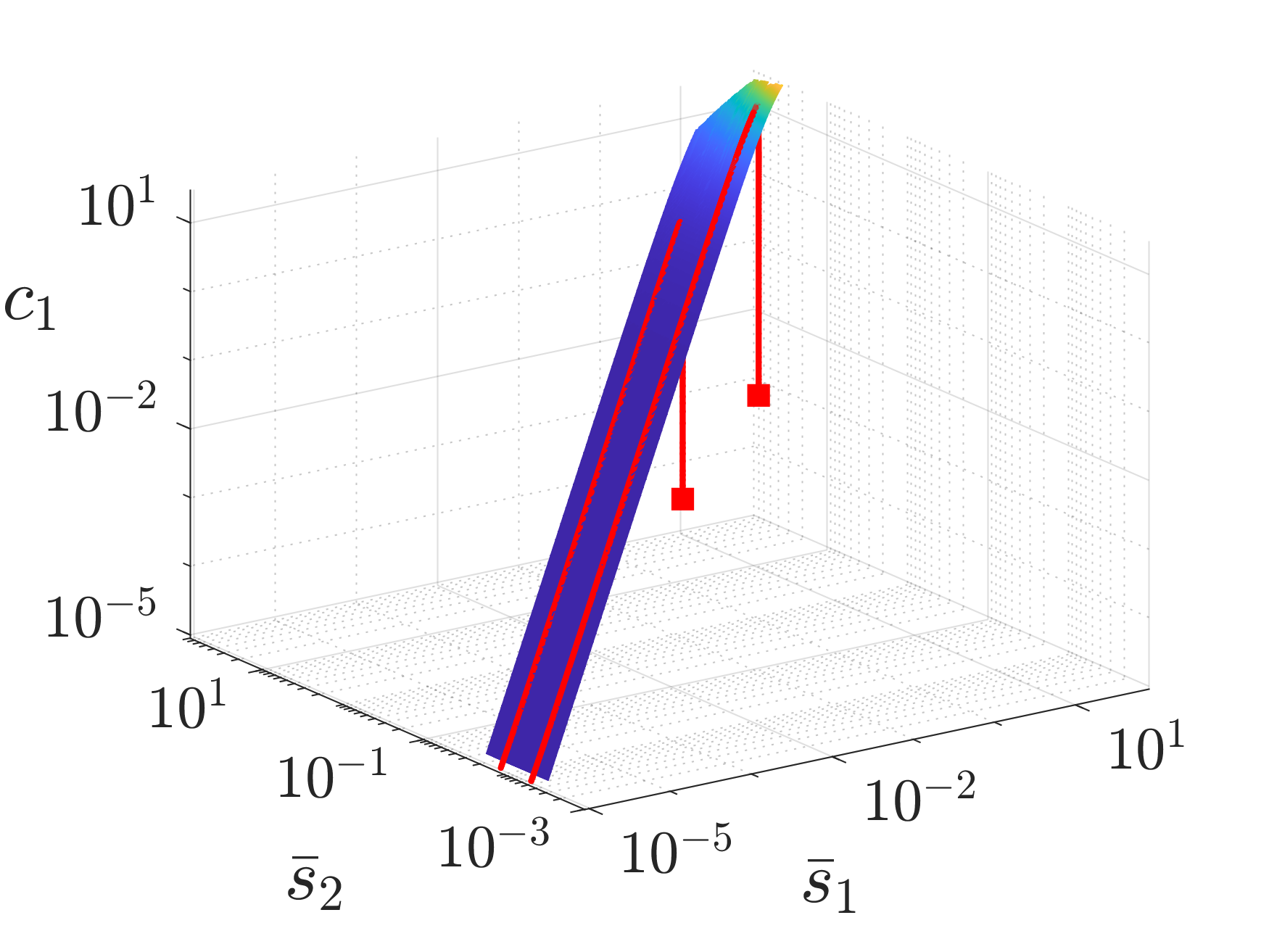}}
    \subfigure[Transformed system, projection to $c_2$]{
    \includegraphics[width=0.33\textwidth]{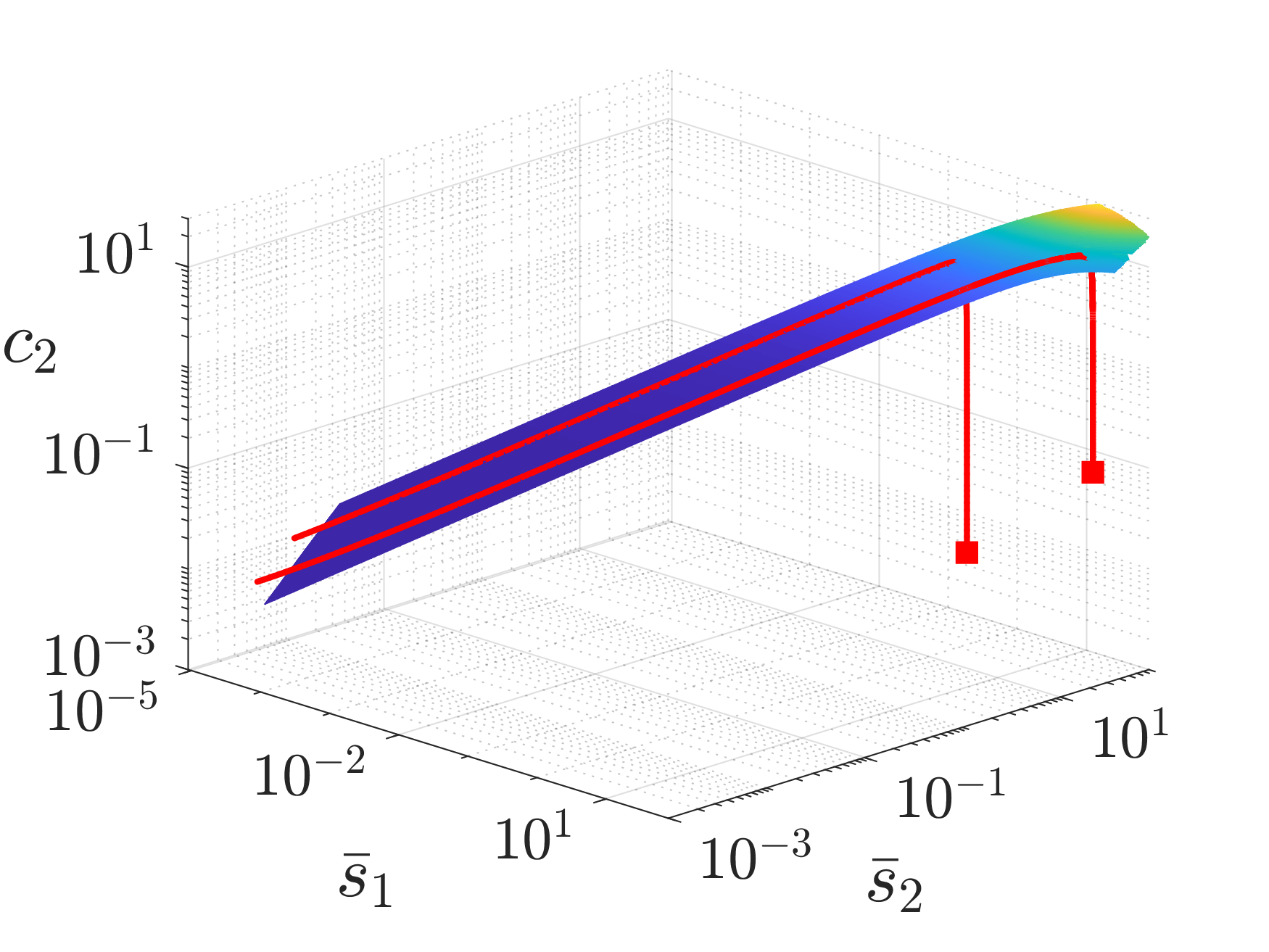}}
    \caption{The $M=2$-dim. SIM hypersurface arising in the phase space of the original (panels a and b) and the transformed (panels c and d) fCSI mechanism in Eq.~(\ref{eq:Inh_st}, \ref{eq:Inh_tr}), respectively.~Since the system is $N=4$-dim., projections to the fast variables $c_1$ and $c_2$ are shown in panels (a, c) and (b, d), respectively.~The red trajectories (starting at red squares) are attracted to the SIM and then evolve on it towards reaching the fixed point $(0,0,0,0)$.}
    \label{fig:Inh_SIMs}
\end{figure}

Here, to demonstrate that the PINN scheme is not affected by the transformation of the fCSI system, we studied both the original and transformed systems in Eqs.~(\ref{eq:Inh_st}, \ref{eq:Inh_tr}), under the parameter set $k_{1f} = 0.1522$, $k_{1b}=2.8$, $k_2=0.7$, $k_{3f}=0.0833$, $k_{3b}=1.667$, $k_4=0.4167$ and $e_0=50$ used in  \citep{pedersena2007total}.~In particular, we computed the $M=2$-dim. SIM in the domain $\Omega=[10^{-5},50] \times [10^{-5}, 30] \times [10^{-3}, 50] \times [10^{-3}, 30]$ for the original system, and in the domain $\Omega=[10^{-5},80] \times [10^{-5}, 30] \times [10^{-3}, 80] \times [10^{-3}, 30]$ for the transformed system; the conditions for the accuracy of the QSSA$_{c1c2}$ are satisfied in those domains, according to \citep{segel1988validity}.~A representation of the SIM emerging in the phase space of the original and transformed fCSI system is shown in Fig.~\ref{fig:Inh_SIMs}, where projections of the SIM onto $c_1$ or $c_2$ fast variables are depicted.~As shown, the trajectories are attracted to the SIM, and then evolve on it until reaching the stable stationary point $(0,0,0,0)$; similar behavior to that of the MM mechanism in Section~\ref{sub:MMdes}.

\section{Numerical Results}
\label{sec:NR}

We assess the efficiency of the proposed PINN-based methodology on learning SIM approximations for the three problems described in Section~\ref{sec:Prob}; both the numerical approximation accuracy and the computational cost are examined.~Furthermore, we compare the approximations provided by the PINN scheme with the GSPT-derived approximations of the SIM, as obtained on the basis of QSSA, PEA and CSP with one and two iterations.\par
For training the PINN scheme, we followed the procedure described in Section~\ref{sbsb:trainPINN} and solved the optimization problem as described in Section~\ref{sec:Meth}.~Next, we assessed the numerical SIM approximation accuracy of the PINN scheme and the GPST expressions, as described in Section~\ref{sb:NumAcc}.~All simulations were carried out with a CPU Intel(R) Xeon(R) CPU E5-2630 v4 @ 2.20GHz (2 processors), RAM 64.0 GB using MATLAB R2022b.

\subsection{The Michaelis-Menten mechanism}

We computed the SIM approximation of the MM mechanism for the three cases considered in Section~\ref{sub:MMdes} in the domain  (i) $\Omega=[10^{-3},10^3]\times[10^{-5},1]$ for the MM1 case, (ii) $\Omega=[10^{-5},1]\times [10^{-3},50]$ for the MM2 case, and (iii)  $\Omega=[10^{-3},10^2]\times[10^{-5},1]$ for the MM3 case.~To construct the training and validation sets for each case, we collected $n=500$ points $\mathbf{z}^{(i)}\in \Omega$, $i=1,\ldots,n$ lying on the $M=1$-dim. SIM from numerically derived trajectories of the IVP in Eq.~\eqref{eq:MM}, as described in Section~\ref{sbsb:trainPINN}; the initial conditions were chosen from the uniform random distributions $[1000,1100] \times [0,1]$ for the MM1 case, $[0,1] \times [50,51]$ for the MM2 case and $[100,110] \times [9,10]$ for the MM3 case.~For a visualization of the training and validation sets, see Supplementary Material Fig.~SF1.

The training results of the PINN scheme for all three cases of the MM model are shown in Table~\ref{tb:MM_train}.
\begin{table}[!h]
    \centering
    \begin{tabular}{l| c c | c c c}
    \toprule
    &	\multicolumn{2}{c}{\bf Loss Function $\lVert \boldsymbol{\mathcal{F}} \rVert^2_2$}	&		\multicolumn{3}{|c}{\bf Computational times (s)} \\
    case	&	Training	&	Validation	&	mean	&	min	&	max	\\
    \midrule
    MM1	&	8.99E$-$11	&	2.61E$-$11	&	1.42E$+$01	&	1.23E$+$01	&	2.27E$+$01	\\
    MM2	&	5.82E$-$06	&	2.67E$-$06	&	1.79E$+$01	&	7.66E$+$00	&	1.98E$+$01	\\
    MM3	&	2.03E$-$06	&	4.85E$-$07	&	2.19E$+$01	&	1.04E$+$01	&	2.82E$+$01	\\
    \bottomrule
    \end{tabular}
    \caption{MM system in \cref{eq:MM} for the MM1, MM2, and MM3 cases.~Loss function $\lVert \boldsymbol{\mathcal{F}} \rVert^2_2$ of the PINN scheme for the training and validation sets and computational times required (in seconds).~The results are obtained by averaging over 100 runs.}
    \label{tb:MM_train}
\end{table}
In particular, we report the average residuals $\lVert \boldsymbol{\mathcal{F}} \rVert^2_2$ of the loss function, as well as the mean, minimum and maximum computational time required over 100 runs with different randomly sampled training and validation sets.~In all 100 runs, the PINN scheme ``discovers'' the following transformations:
\begin{align*}
    \mathbf{C} & = \begin{bmatrix} 0 & 1 \end{bmatrix}, & \mathbf{C} & = \begin{bmatrix} 1 & 0 \end{bmatrix}, & \mathbf{C} & = \begin{bmatrix} 0 & 1 \end{bmatrix}, \\ 
    \mathbf{D} & = \begin{bmatrix} 1 & 0 \end{bmatrix}, & \mathbf{D} & = \begin{bmatrix} 0 & 1 \end{bmatrix}, & \mathbf{D} & = \begin{bmatrix} 1 & 0 \end{bmatrix},
\end{align*}
for the MM1, MM2, and MM3 cases, respectively.~Hence, the proposed PINN scheme correctly identifies that the fast (slow) variable is $x=c$ ($y=s$) in the MM1 and MM3 cases, while it is $x=s$ ($y=c$) in the MM2 case.

Next, to assess the numerical accuracy of the SIM approximations, we constructed the test sets for each case, as described in Section~\ref{sb:NumAcc}.~In particular, using the same range of initial conditions as with the training sets, we generated $25$ trajectories of the IVP in Eq.~\eqref{eq:MM} and collected $100$ equidistant -in time- points per trajectory lying exclusively on the $M=1$-dim. SIM in the respective domain $\Omega$.~A visualization of the test sets is provided for the MM1, MM2 and MM3 cases in Supplementary Material Fig.~SF1.

Table~\ref{tb:MM_test} reports the $l^2$, $l^{\infty}$ and MSE approximation errors on the test set for the three MM cases.~Given the data points $\mathbf{z}^{(i)}=[x^{(i)},y^{(i)}]^\top$ for $i=1,\ldots,n_t$, the PINN  approximation error is obtained as $\lVert \mathbf{C} \mathbf{z}^{(i)} - \mathcal{N}(\mathbf{D} \mathbf{z}^{(i)}) \rVert$, while the GSPT one as $\lVert x^{(i)} - h(y^{(i)}) \rVert$, thus requiring knowledge about the fast variables ($c$ in MM1 and MM3 cases and $s$ in MM2 cases).~For the cases where an implicit functional is provided by the GSPT expressions (see Table~\ref{tb:MM_SIMs}), we numerically solved the corresponding SIM expressions using Newton iterations w.r.t. the fast variable $\tilde{x}^{(i)}$ and obtained the SIM approximation error as $\lVert x^{(i)} - \tilde{x}^{(i)} \rVert$.~For the latter cases, we also report the computational time required for Newton iterations (CTN) in Table~\ref{tb:MM_test}. 
\begin{table}[!h]
    \centering
    \resizebox{\textwidth}{!}{
    \begin{tabular}{l l | c c c c c c c c c}
    \toprule
    &	\textbf{Error} 	&	\textbf{PINN}	&	\textbf{rQSSA}	&	\textbf{sQSSA}	&	\textbf{PEA}	&	\textbf{CSP$_s$(1)}	&	\textbf{CSP$_s$(2)}	&	\textbf{CSP$_c$(1)}	&	\textbf{CSP$_c$(2)}	&	\textbf{CSP$_e$}	\\
    \midrule
    \multirow{4}{*}{\rotatebox[origin=c]{90}{MM1}} &	$l^2$	&	4.80E$-$06	&	1.44E$-$01	&	4.03E$-$04	&	2.92E$-$06	&	2.92E$-$06	&	1.95E$-$09	&	1.29E$-$08	&	1.95E$-$09	&	9.81E$-$09	\\
    &	$l^\infty$	&	1.68E$-$07	&	2.48E$-$03	&	1.03E$-$05	&	7.90E$-$08	&	7.90E$-$08	&	3.18E$-$10	&	4.67E$-$10	&	3.18E$-$10	&	3.54E$-$10	\\
    &	MSE	&	1.90E$-$15	&	1.71E$-$06	&	1.34E$-$11	&	7.01E$-$16	&	7.01E$-$16	&	3.14E$-$22	&	1.36E$-$20	&	3.12E$-$22	&	7.93E$-$21	\\
    \cmidrule{2-11}
    &	CTN &	N/A	&	N/A	&	N/A	&	N/A	&	N/A	&	2.67E$-$02	&	N/A	&	2.06E$-$02	&	1.89E$-$02	\\
    \midrule
    \multirow{4}{*}{\rotatebox[origin=c]{90}{MM2}} &	$l^2$	&	2.18E$-$04	&	5.23E$-$03	&	1.85E$-$01	&	6.57E$-$06	&	6.57E$-$06	&	6.98E$-$06	&	2.99E$-$05	&	6.98E$-$06	&	6.61E$-$06	\\
    &	$l^\infty$	&	3.97E$-$05	&	3.82E$-$04	&	9.61E$-$03	&	1.82E$-$06	&	1.82E$-$06	&	2.03E$-$06	&	3.66E$-$06	&	2.03E$-$06	&	1.89E$-$06	\\
    &	MSE	&	3.82E$-$12	&	2.20E$-$09	&	2.75E$-$06	&	3.46E$-$15	&	3.46E$-$15	&	3.91E$-$15	&	7.16E$-$14	&	3.91E$-$15	&	3.51E$-$15	\\
    \cmidrule{2-11}
    &	CTN	&	N/A	&	N/A	&	N/A	&	N/A	&	N/A	&	N/A	&	N/A	&	2.62E$-$02	&	N/A	\\
    \midrule
    \multirow{4}{*}{\rotatebox[origin=c]{90}{MM3}} &	$l^2$	&	3.30E$-$04	&	6.15E$+$02	&	1.84E$-$01	&	5.64E$+$00	&	5.64E$+$00	&	5.22E$-$02	&	1.45E$-$03	&	3.43E$-$04	&	3.50E$-$04	\\
    &	$l^\infty$	&	7.45E$-$05	&	9.38E$+$00	&	6.89E$-$03	&	8.64E$-$02	&	8.64E$-$02	&	8.19E$-$04	&	1.17E$-$04	&	7.63E$-$05	&	6.60E$-$05	\\
    &	MSE	&	8.72E$-$12	&	3.03E$+$01	&	2.71E$-$06	&	2.55E$-$03	&	2.55E$-$03	&	2.18E$-$07	&	1.69E$-$10	&	9.42E$-$12	&	9.79E$-$12	\\
    \cmidrule{2-11}
    &	CTN&	N/A	&	N/A	&	N/A	&	N/A	&	N/A	&	3.53E$-$02	&	N/A	&	2.34E$-$02	&	2.53E$-$02	\\
    \bottomrule
    \end{tabular}}
    \caption{MM system in Eq.~\eqref{eq:MM} for the MM1, MM2 and MM3 cases.~SIM approximation accuracy over all points of the test set, in terms of $l^2$, $l^{\infty}$ and MSE approximation errors, resulting from the PINN scheme and the GSPT approximations enlisted in Table~\ref{tb:MM_SIMs}.~The numerical accuracy of each SIM approximation is compared with the numerical solution $\mathbf{z}^{(i)}$ of the MM system.~In the cases of implicit GSPT expressions, we report the computational time (in $s$) required for Newton iterations (CTN) over all points to obtain an estimation of the fast variable on the SIM.}
    \label{tb:MM_test}
\end{table}

Examining the accuracy of the GSPT expressions, Table~\ref{tb:MM_test} shows that in all MM cases, the SIM approximations constructed with the correct assumption on the fast variable/reaction (see Table~\ref{tb:MM_SIMs}) are accurate, as expected \citep{patsatzis2023algorithmic}.~In particular, the accurate SIM approximations are provided by (i) sQSSA, PEA, CSP$_c$(1) and CSP$_c$(2) in the MM1 case, (i) rQSSA, PEA, CSP$_s$(1) and CSP$_s$(2) in the MM2 case, and (iii) sQSSA, CSP$_c$(1) and CSP$_c$(2) in the MM3 case.~In addition, it is shown that the CSP approximations generated with the wrong assumption on the fast variables (CSP$_s$ in MM1 and MM3 cases and CSP$_c$ in MM2 case), required at least two iterations for achieving similarly high (or even lesser; e.g. in MM3 case) approximation accuracy with those generated with the correct assumption.~Last but not least, the CSP$_e$ approximation is in every case one of the most accurate SIM approximations.\par
Regarding the PINN scheme, Table~\ref{tb:MM_test} shows that a high SIM approximation accuracy is attained in all MM cases considered (the overall $l^\infty$ error is less than 1E$-$04).~In addition, the proposed PINN scheme is more accurate than the QSSA approximations for all cases.~For the MM3 case particularly, the PINN scheme provides a similarly high approximation accuracy as the CSP$_c(2)$ and CSP$_e$ approximations, without additionally requiring Newton iterations.\par
Finally, we provide a visualization of the SIM approximation accuracy for all the points of the test set in Supplementary Material Figs.~SF5, SF6 and SF7 for the MM1, MM2 and MM3 cases, respectively.~It is therein shown, that higher approximation accuracy is provided by the GSPT expressions for points that are closer to the fixed point $(0,0)$ of the MM system.~The proposed PINN scheme, maintains almost the same level of accuracy allover the test set.

\subsection{The Target Mediated Drug Disposition mechanism}
\label{sub:TMDDres}

For the TMDD mechanism, we computed the SIM approximations for both $\mathcal{M}_1$ and $\mathcal{M}_2$, in the domains $\Omega=[10,140]\times[0,0.12]\times[10,33]$ and $\Omega=[0,0.05]\times[6,13]\times[6,27]$, respectively, as described in Section~\ref{sub:TMDDdes}.~To construct the training and validation sets in the above domains, we produced a number of numerical trajectories of the TMDD model in Eq.~\eqref{eq:TMDD}; 125 initial conditions were realized from the uniform random distributions $[50,150] \times [10,15] \times [0,1]$.~Then, from all the resulting trajectories, we randomly sampled $n=700$ points $\mathbf{z}^{(i)}\in \Omega$, $i=1,\ldots,n$, lying on the respective SIM $\mathcal{M}_1$ or $\mathcal{M}_2$; for detection we used the criterion in Eq.~\eqref{eq:CSPcrit}.~A visualization of the training and validation sets is provided in Supplementary Material Fig.~SF2.

Table~\ref{tb:TMDD_train} summarizes the average squared $l^2$ residuals $\lVert \boldsymbol{\mathcal{F}} \rVert^2_2$ of the loss functions on the training and validation sets of the PINN scheme, corresponding to both $\mathcal{M}_1$ and $\mathcal{M}_2$ SIMs.
\begin{table}[!h]
    \centering
    \begin{tabular}{l| c c | c c c}
    \toprule
    &	\multicolumn{2}{c}{\bf Loss Function $\lVert \boldsymbol{\mathcal{F}} \rVert^2_2$}	&		\multicolumn{3}{|c}{\bf Computational times (s)} \\
    SIM	&	Training	&	Validation	&	mean	&	min	&	max	\\
    \midrule
    $\mathcal{M}_1$	&	3.52E$-$10	&	1.11E$-$10	&	1.19E$+$01	&	8.94E$+$00	&	1.48E$+$01	\\
    $\mathcal{M}_2$	&	7.67E$-$10	&	2.53E$-$10	&	1.39E$+$01	&	1.16E$+$01	&	1.49E$+$01	\\
    \bottomrule
    \end{tabular}
    \caption{TMDD system in \cref{eq:TMDD} for the SIMs $\mathcal{M}_1$ and $\mathcal{M}_2$.~Loss function $\lVert \boldsymbol{\mathcal{F}} \rVert^2_2$ of the PINN scheme for the training and validation sets and computational times required (in seconds).~The results are obtained by averaging over 100 runs.}
    \label{tb:TMDD_train}
\end{table}
~The mean, minimum and maximum computational time required over 100 runs with different randomly sampled training and validation sets are also included in Table~\ref{tb:TMDD_train}.~The transformations discovered from the PINN scheme in each of the 100 random runs, are of the form:
\begin{align*}
    \mathbf{C} & = \begin{bmatrix} 0 & 1 & 0 \end{bmatrix}, & \mathbf{C} & = \begin{bmatrix} 1 & 0 & 0 \end{bmatrix},  \\ 
    \mathbf{D} & = \begin{bmatrix} \alpha_1 & 0 & 1-\alpha_1 \\ 1-\alpha_1 & 0 & \alpha_1\end{bmatrix}, & \mathbf{D} & = \begin{bmatrix} 0 & \alpha_2  & 1-\alpha_2 \\ 0 & 1-\alpha_2 & \alpha_2 \end{bmatrix},
\end{align*}
for $\mathcal{M}_1$ and $\mathcal{M}_2$, respectively, where $\alpha_1\in[0.23,0.88]$ and $\alpha_2\in[0.17,0.22]\cup [0.77,0.88]$.~Thus, the proposed PINN scheme correctly discovers that the fast variable in $\mathcal{M}_1$ ($\mathcal{M}_2$) is $x=R$ ($x=L$), while the slow variables $\mathbf{y}$ are linear combinations of $[L,RL]$ ($[R,RL]$).

For evaluating the numerical SIM approximation accuracy, we collected test sets for $\mathcal{M}_1$ and $\mathcal{M}_2$, as described in Section~\ref{sb:NumAcc}.~In particular, using random initial conditions in the same range with the ones of the train sets, we generated $125$ trajectories of the TMDD model in Eq.~\eqref{eq:TMDD}.~From the resulting trajectories, we collected $100$ equidistant -in time- points, lying exclusively on the $M=1$-dim. SIM in the respective domain $\Omega$.~A visualization of the test sets for the $\mathcal{M}_1$ and $\mathcal{M}_2$ SIMs is provided in Supplementary Material Fig.~SF2.

Table~\ref{tb:TMDD_test} reports the $l^2$, $l^{\infty}$ and MSE approximation errors overall the test set for $\mathcal{M}_1$ and $\mathcal{M}_2$ SIMs of the TMDD mechanism.~For the PINN scheme, the approximation error is obtained as $\lVert \mathbf{C} \mathbf{z}^{(i)} - \mathcal{N}(\mathbf{D} \mathbf{z}^{(i)}) \rVert$, while for the GSPT expressions as $\lVert x^{(i)} - h(\mathbf{y}^{(i)}) \rVert$; the latter requiring knowledge about the fast variables ($R$ for $\mathcal{M}_1$ and $L$ for $\mathcal{M}_2$).~For the CSP with two iterations expressions, all three of which provide implicit functionals (see Table~\ref{tb:TMDD_SIMs}); in this case, we used Newton iterations w.r.t. the fast variable $\tilde{x}^{(i)}$.~In Table~\ref{tb:TMDD_test} we report the respective SIM approximation errors $\lVert x^{(i)} - \tilde{x}^{(i)} \rVert$, along with the computational time required for Newton iterations (CTN). 
\begin{table}[!h]
    \centering
    \resizebox{\textwidth}{!}{
    \begin{tabular}{l l | c c c c c c c c c c c}
    \toprule
    &	\textbf{Error} 	&	\textbf{PINN}	&	\textbf{QSSA$_L$}	&	\textbf{QSSA$_R$}	&	\textbf{QSSA$_{RL}$}	&	\textbf{PEA}	&	\textbf{CSP$_L$(1)}	&	\textbf{CSP$_L$(2)}	&	\textbf{CSP$_R$(1)}	&	\textbf{CSP$_R$(2)}	&	\textbf{CSP$_{RL}$(1)}	&	\textbf{CSP$_{RL}$(2)}	\\
    \midrule
    \multirow{4}{*}{\rotatebox[origin=c]{90}{$\mathcal{M}_1$}} &	$l^2$	&	2.97E$-$05	&	1.38E$+$01	&	5.55E$-$02	&	3.07E$+$00	&	9.30E$-$04	&	1.11E$-$02	&	4.86E$-$05	&	5.85E$-$04	&	1.49E$-$05	&	5.06E$-$03	&	8.87E$-$06	\\
    &	$l^\infty$	&	1.52E$-$06	&	1.07E$-$01	&	9.51E$-$04	&	3.42E$-$02	&	2.06E$-$05	&	1.53E$-$04	&	1.50E$-$06	&	1.35E$-$05	&	1.47E$-$06	&	1.01E$-$04	&	1.48E$-$06	\\
    &	MSE	&	1.41E$-$14	&	3.03E$-$03	&	4.93E$-$08	&	1.50E$-$04	&	1.38E$-$11	&	1.98E$-$09	&	3.78E$-$14	&	5.48E$-$12	&	3.53E$-$15	&	4.09E$-$10	&	1.26E$-$15	\\
    \cmidrule{2-13}
    &	CTN	&	N/A	&	N/A	&	N/A	&	N/A	&	N/A	&	N/A	&	8.93E$-$01	&	N/A	&	4.14E$-$01	&	N/A	&	1.73E$-$01	\\
    \midrule
    \multirow{4}{*}{\rotatebox[origin=c]{90}{$\mathcal{M}_2$}} &	$l^2$	&	6.73E$-$05	&	5.59E$-$02	&	7.45E$+$00	&	1.42E$+$01	&	1.86E$-$03	&	1.98E$-$03	&	1.50E$-$04	&	9.53E$-$02	&	8.05E$-$04	&	5.47E$-$02	&	1.55E$-$04	\\
    &	$l^\infty$	&	3.15E$-$06	&	1.10E$-$03	&	9.95E$-$02	&	1.43E$-$01	&	5.08E$-$05	&	5.37E$-$05	&	4.91E$-$06	&	1.54E$-$03	&	1.21E$-$05	&	7.23E$-$04	&	1.52E$-$06	\\
    &	MSE	&	7.24E$-$14	&	5.00E$-$08	&	8.89E$-$04	&	3.22E$-$03	&	5.53E$-$11	&	6.30E$-$11	&	3.61E$-$13	&	1.45E$-$07	&	1.04E$-$11	&	4.79E$-$08	&	3.84E$-$13	\\
    \cmidrule{2-13}
    &	CTN	&	N/A	&	N/A	&	N/A	&	N/A	&	N/A	&	N/A	&	5.32E$-$01	&	N/A	&	5.68E$-$01	&	N/A	&	1.62E$-$01	\\
    \bottomrule
    \end{tabular}}
    \caption{TMDD system in Eq.~\eqref{eq:TMDD} for $\mathcal{M}_1$ and $\mathcal{M}_2$ SIMs.~SIM approximation accuracy over all points of the test set, in terms of $l^2$, $l^{\infty}$ and MSE approximation errors, resulting from the PINN scheme and the GSPT approximations enlisted in Table~\ref{tb:TMDD_SIMs}.~The numerical accuracy of each SIM approximation is compared with the numerical solution $\mathbf{z}^{(i)}$ of the TMDD system.~In the cases of implicit GSPT expressions, we report the computational time (in $s$) required for Newton iterations (CTN) over all points to obtain an estimation of the fast variable on the SIM.}
    \label{tb:TMDD_test}
\end{table}

For the accuracy of the GSPT expressions, Table~\ref{tb:TMDD_test} shows that the QSSA and CSP expressions constructed with the correct assumption on the fast variable ($R$ for $\mathcal{M}_1$ and $L$  for $\mathcal{M}_2$), provide higher approximation accuracy than the rest of the same kind.~High approximation accuracy is also provided by the PEA expression; in both SIMs the assumption on the first reaction being fast is satisfied.~In addition, although the CSP expressions constructed with the wrong assumption on the fast variable, are accurate (CSP requires at least two iterations for achieving similarly high accuracy with those generated with the correct assumption).~For example, to achieve an order of $1E-03$ $l^2$ error, either the CSP$_R$(1) expression or the CSP$_L$(2) or CSP$_{RL}$(2) expressions for $\mathcal{M}_1$ should be adopted; the latter require Newton iterations.

The high numerical SIM approximation accuracy provided by the proposed PINN scheme is also shown Table~\ref{tb:TMDD_test}; the overall $l^\infty$ error is of the order $1E-06$.~In particular, the PINN provides similar or even higher accuracy than that provided by CSP with two iterations, the most accurate GSPT expressions considered, which also require additional computational time for the Newton iterations.\par
Finally, a visualization of the approximation accuracy, in both $\mathcal{M}_1$ and $\mathcal{M}_2$ SIMs, for all the points of these sets is provided in Supplementary Material Figs.~SF8 and SF9, respectively.~It is therein shown that the accurate GSPT expressions result to rather poor approximation accuracy close to the boundaries of the SIM, particularly as the trajectories exit $\mathcal{M}_1$ and as they enter $\mathcal{M}_2$.~This is not the case for the PINN scheme, as expected by \citep{patsatzis2023slow}, which maintains a similar level of accuracy all over the range.

\subsection{The fully Competitive Substrate-Inhibitor mechanism}

As described in Section~\ref{sub:CompInhdes}, we computed SIM approximations of the original fCSI system in Eq.~\eqref{eq:Inh_st} in the domain $\Omega=[10^{-5},50] \times [10^{-5}, 30] \times [10^{-3}, 50] \times [10^{-3}, 30]$, and those of the transformed fCSI system in Eq.~\eqref{eq:Inh_tr} in the domain $\Omega=[10^{-5},80] \times [10^{-5}, 30] \times [10^{-3}, 80] \times [10^{-3}, 30]$.~In both systems, we constructed the training and validation points in the above domains from numerically derived trajectories of the fCSI system; $625$ initial conditions are randomly selected from uniform distributions $s_{10}, s_{20}\in[50,150]$ and $c_{10}, c_{20}\in[0,1]$.~From the resulting trajectories, we detected, with the criterion in Eq.~\eqref{eq:CSPcrit}, the points where $M=2$, and randomly sampled $n=700$ of them, so that $\mathbf{z}^{(i)}\in \Omega$, $i=1,\ldots,n$; for a visualization, see Supplementary Material Figs.~SF3 and SF4. 

The loss functions $\lVert \boldsymbol{\mathcal{F}} \rVert^2_2$ obtained with the PINN scheme in the training data set are shown in Table~\ref{tb:Inh_train} for both the original and transformed systems.
\begin{table}[!h]
    \centering
    \begin{tabular}{l| c c | c c c}
    \toprule
    &	\multicolumn{2}{c}{\bf Loss Function $\lVert \boldsymbol{\mathcal{F}} \rVert^2_2$}	&		\multicolumn{3}{|c}{\bf Computational times (s)} \\
    system	&	Training	& 	Validation	& 	mean	& 	min	& 	max	\\
     \midrule
    original	&	1.91E$-$04	& 	7.83E$-$05	& 	6.62E$+$01	& 	3.24E$+$01	& 	1.33E$+$02	\\
    transformed	&	5.63E$-$05	& 	1.99E$-$05	& 	1.04E$+$02	& 	4.98E$+$01	& 	1.34E$+$02	\\
    \bottomrule
    \end{tabular}
    \caption{Original and transformed systems of the fCSI mechanism in \cref{eq:Inh_st,eq:Inh_tr}.~Loss function $\lVert \boldsymbol{\mathcal{F}} \rVert^2_2$ of the PINN scheme for the training and validation sets and computational times required (in seconds).~The results are obtained by averaging over 100 runs.}
    \label{tb:Inh_train}
\end{table}
The PINN optimization scheme converges (since the average loss function $\lVert \boldsymbol{\mathcal{F}} \rVert^2_2$ of training set is slightly higher than that of the validation set) to higher errors in comparison to the MM and TMDD problems (see Tables~\ref{tb:MM_train}, \ref{tb:TMDD_train}).~This is because here, the number of residuals is two-fold higher than the other cases, since $M=2$ in this case.~Table~\ref{tb:Inh_train} also reports the computational time required for training; mean, minimum and maximum over 100 runs with different randomly sampled training and validation sets.
The transformations discovered by the PINN scheme, are of the form:
\begin{align*}
    \mathbf{C} & = \begin{bmatrix} 0 & \alpha_1 & 0 & 1-\alpha_1 \\
    0 & 1-\alpha_1 & 0 & \alpha_1 \end{bmatrix}, & \mathbf{C} & = \begin{bmatrix} 0 & \alpha_1' & 0 & 1-\alpha_1' \\
    0 & 1-\alpha_1' & 0 & \alpha_1' \end{bmatrix},  \\ 
    \mathbf{D} & = \begin{bmatrix} \alpha_2 & 0 & 1-\alpha_2 & 0\\ 1-\alpha_2 & 0 & \alpha_2 & 0\end{bmatrix}, & \mathbf{D} & = \begin{bmatrix} \alpha_2' & 0 & 1-\alpha_2' & 0 \\ 1-\alpha_2' & 0 & \alpha_2' & 0 \end{bmatrix},
\end{align*}
where $\alpha_1=0.5$ and $\alpha_2\in[0.32,0.35]\cup [0.66,0.77]$ for the original system, while $\alpha_1'\in[0.49,0.51]$ and $\alpha_2'\in[0.44,0.58]$ for the transformed one.~In both systems, the proposed PINN sheme discovers that the fast variables $x_1$ and $x_2$ are a linear combination of the $c_1$ and $c_2$; the slow variables are a linear combination of $s_1$/$\bar{s}_1$ and $s_2$/$\bar{s}_2$ for the original/transformed system.~\textcolor{blue}
~The above finding is in agreement with the assumptions of fast variables in \citep{segel1988validity,schnell2000time,rubinow1970time}.

Next, for assessing the numerical SIM approximation accuracy, we formed the test sets as described in Section~\ref{sb:NumAcc}.~Using the same range of initial conditions with that of the training sets, we generated $625$ trajectories of the original and transformed systems of the fCSI systems in \cref{eq:Inh_st,eq:Inh_tr}.~From the resulting trajectories, we sampled $100$ equidistant -in time- points, lying exclusively on the $M=2$-dim. SIM in the respective domain $\Omega$.~A visualization of the test sets is provided in Supplementary Material Figs.~SF3 and SF4. 

Based on the test sets, we report in Table~\ref{tb:Inh_test} the $l^2$, $l^{\infty}$ and MSE approximation errors, overall the test set, for both the original and transformed systems of the fCSI mechanism.~Note that since, in this case, the SIM is $M=2$-dim., we report separately the errors for the  fast variables $x_1$ and $x_2$.~For the PINN scheme, the approximation error is obtained as $\lVert \mathbf{C} \mathbf{z}^{(i)} - \mathcal{N}(\mathbf{D} \mathbf{z}^{(i)}) \rVert$, while for the GSPT expressions it would be calculated as $\lVert \mathbf{x}^{(i)} - \mathbf{h}(\mathbf{y}^{(i)}) \rVert$ using the knowledge of the fast variables ($x_1=c_1$ and $x_2=c_2$).~However, this only applies for the QSSA$_{c1c2}$ expression of the original system, since all other GSPT expressions provide implicit functionals (see Table~\ref{tb:Inh_SIMs}).~For these expressions, we numerically solved for the corresponding $M=2$-dim. system of the implicit SIM expressions using Newton iterations w.r.t. the fast variables, the estimation of which is denoted as $\tilde{\mathbf{x}}=[\tilde{c}_1, \tilde{c}_2]$.~Then, we calculated the SIM approximation errors as $\lVert \mathbf{x}^{(i)} - \tilde{\mathbf{x}}^{(i)} \rVert$ and report them in Table~\ref{tb:Inh_test}, along with the computational time required for the Newton iterations (CTN). 
\begin{table}[!h]
    \centering
    \resizebox{\textwidth}{!}{
    \begin{tabular}{l l | c c c c c c c c c c}
    \toprule
    &	\textbf{Error} 	&	\multicolumn{2}{c}{\textbf{PINN}}	&	\multicolumn{2}{c}{\textbf{QSSA$_{c1c2}$}}	&	\multicolumn{2}{c}{\textbf{PEA$_{13}$}}	&	\multicolumn{2}{c}{\textbf{CSP$_{c1c2}$(1)}}	&	\multicolumn{2}{c}{\textbf{CSP$_{c1c2}$(2)}}	\\
    \midrule
	&		&	$x_1$	&	$x_2$	&	$x_1$	&	$x_2$	&	$x_1$	&	$x_2$	&	$x_1$	&	$x_2$	&	$x_1$	&	$x_2$	\\
	\midrule	
    \multirow{4}{*}{\rotatebox[origin=c]{90}{original}} &	$l^2$	&	5.73E$-$02	&	5.69E$-$02	&	4.85E$+$01	&	9.00E$+$01	&	3.75E$+$00	&	6.42E$+$00	&	2.06E$+$00	&	4.06E$+$00	&	5.80E$-$01	&	7.43E$-$01	\\
	&	$l^\infty$	&	9.39E$-$03	&	9.38E$-$03	&	6.99E$-$01	&	9.97E$-$01	&	6.45E$-$02	&	7.43E$-$02	&	8.95E$-$02	&	1.13E$-$01	&	7.94E$-$02	&	9.58E$-$02	\\
	&	MSE	&	5.25E$-$08	&	5.18E$-$08	&	3.76E$-$02	&	1.30E$-$01	&	2.24E$-$04	&	6.60E$-$04	&	6.79E$-$05	&	2.64E$-$04	&	5.37E$-$06	&	8.83E$-$06	\\
    \cmidrule{2-12}
	&	CTN	&		\multicolumn{2}{c}{N/A}&		\multicolumn{2}{c}{N/A}	&	\multicolumn{2}{c}{4.05E$+$00}	&	\multicolumn{2}{c}{4.11E$+$00}	&	\multicolumn{2}{c}{5.63E$+$01}	\\
	\midrule					
	\multirow{4}{*}{\rotatebox[origin=c]{90}{transformed}}	&	$l^2$	&	7.02E$-$02	&	6.90E$-$02	&	2.83E$+$01	&	4.29E$+$01	&	2.56E$+$00	&	3.61E$+$00	&	9.45E$-$01	&	2.08E$+$00	&	4.25E$-$01	&	5.49E$-$01	\\
	&	$l^\infty$	&	8.55E$-$03	&	8.55E$-$03	&	4.10E$-$01	&	5.31E$-$01	&	4.71E$-$02	&	4.88E$-$02	&	5.73E$-$02	&	7.80E$-$02	&	5.32E$-$02	&	6.60E$-$02	\\
	&	MSE	&	7.87E$-$08	&	7.62E$-$08	&	1.28E$-$02	&	2.95E$-$02	&	1.05E$-$04	&	2.09E$-$04	&	1.43E$-$05	&	6.95E$-$05	&	2.89E$-$06	&	4.82E$-$06	\\
     \cmidrule{2-12}
     &	CTN	&	\multicolumn{2}{c}{N/A}	&	\multicolumn{2}{c}{3.98E$+$00}	&	\multicolumn{2}{c}{3.87E$+$00}	&	\multicolumn{2}{c}{4.06E$+$00}	&	\multicolumn{2}{c}{5.36E$+$01}		\\
    \bottomrule
    \end{tabular}}
    \caption{Original and transformed systems of the fCSI mechanism in \cref{eq:Inh_st,eq:Inh_tr}.~SIM approximation accuracy over all points of the test set, in terms of $l^2$, $l^{\infty}$ and MSE approximation errors, resulting from the PINN scheme and the GSPT approximations enlisted in Table~\ref{tb:Inh_SIMs}.~The numerical accuracy of each SIM approximation is compared with the numerical solution $\mathbf{z}^{(i)}$ of the fCSI system for each fast variable $x_1$ and $x_2$.~In the cases of implicit GSPT expressions, we report the computational time (in $s$) required for Newton iterations (CTN) to obtain an estimation of the fast variables on the SIM, over all points of the test set.}
    \label{tb:Inh_test}
\end{table}

For the evaluation of the GSPT approximations, Table~\ref{tb:Inh_test} shows that in both systems, the lower numerical SIM approximation accuracy is provided by QSSA$_{c1c2}$.~Higher numerical accuracy is progressively provided by the PEA$_{13}$, CSP$_{c1c2}$(1) and CSP$_{c1c2}$(2) expressions, derived in this work with the same assumptions as QSSA$_{c1c2}$ about the fast variables.~Comparing the two systems, it is shown in Table~\ref{tb:Inh_test}, that a slightly higher approximation accuracy is provided by the GSPT expressions derived for the transformed system (in agreement with the results about QSSA reported in \citep{pedersena2007total,bersani2017tihonov}).~We re-iterate here that all, but the QSSA$_{c1c2}$ expressions for the original system, require Newton iterations to obtain an estimation of the fast variables on the SIM.\par
Table~\ref{tb:Inh_test} reports the high numerical SIM approximation accuracy provided by the proposed PINN scheme (similarly high for both original and transformed systems).~In particular, the accuracy of the PINN scheme is at least one order higher, in terms of every approximation error, than that of the CSP$_{c1c2}$(2) expression, which is the most accurate among the GSPT expressions considered.~Note here that the computational time required for training the PINN scheme (see Table~\ref{tb:Inh_train}) is comparable to the computational time needed for CSP$_{c1c2}$(2) to obtain an explicit functional.~This is because CSP with two iterations is derived numerically in this problem, and the employment of Newton iterations requires additional computational cost, due to the numerical differentiation required for the Jacobian matrix.\par
Finally, we provide a visualization of the SIM approximation accuracy over all points of the test sets, for both the original and transformed systems, in Supplementary Material Figs.~SF10 and SF11, respectively.~It is therein shown that higher approximation
accuracy is provided by the GSPT expressions for points that are closer to the fixed point $(0, 0 , 0, 0)$ of the fCSI system.~This is not the case for the proposed PINN scheme, which maintains almost homogeneous accuracy over all the points of the test set.

\section{Conclusions}

The discovery and approximation of slow invariant manifolds (SIMs) is an open and challenging problem of utmost importance across diverse disciplines ranging from biology and neuroscience to ecology and process engineering, just to name a few. SIMs not only facilitate, through the construction of ROMs, the efficient - in terms of accuracy and computational cost - completion of tasks such as the numerical bifurcation analysis for the identification of tipping points, the (epistemic) uncertainty quantification of model predictions, the optimization and control, the efficient. They also provide useful physical insight about the essential/parsimonious mechanisms that govern the observed dynamics, thus advancing our understanding of the emergent behavior of complex and multiscale systems. \par
In recent years, physics-informed machine learning has raised as a promising tool to construct a synergism between state-of-the-art numerical analysis and machine learning. Exploiting this synergism, we proposed a PINN approach to find explicit closed forms of SIMs and their dimension for the general class of stiff systems of ODEs. As the resulting SIMs are constructed in a closed form, they may in turn be used in a straightforward manner to build ROMs. Our approach generalizes previous studies, which either focused on the data-driven ``black-box'' construction of ROMs without firstly identifying the underlying SIM in a closed form, or others that require the a priori knowledge of the splitting of the fast and slow variables. \par 
The proposed approach can be seen as a hybrid one, bridging the CSP-inspired criterion for the  ``on the fly'' identification of the dimension of the SIM with PINNs for the simultaneous transformation of the state-space into fast and slow components, and the solution of the PDE arising from the invariance equation within the framework of GSPT. The transformation that we use here is a linear one, in analogy to the GSPT methods that have been primarily developed with a focus on the simplification of high-dimensional chemical kinetics. For such problems, but also for many problems in other disciplines (e.g., in biology), one is usually interested in discovering the essential physical variables per se or a linear combination of them, in order to get a physical insight to the essential mechanisms that govern the long-term evolution.\par

The comparison of the proposed ML approach with other traditional GSPT methods (QSSA, PEA and CSP) for the three multi-scale problems considered here, showed that the proposed method results in similarly high SIM approximation accuracy to the CSP with two iterations.~In addition, it was shown through the TMDD problem in Section~\ref{sub:TMDDres} that the PINN approach  results in significantly higher approximation accuracy close to the boundaries of the SIM.~We highlight here that in contrast to QSSA, PEA and CSP with one iteration, the PINN approach does not require any assumption on the fast variables (or reactions) in order to provide accurate SIM approximations.~On a final note, since the PINN provides - by construction - explicit functionals of the SIM, it does not require additional numerical root-finding algorithms to estimate the values of the fast variables on the SIM, which might be the case for the other GSPT methods; see \cref{tb:MM_SIMs,tb:TMDD_SIMs,tb:Inh_SIMs} for the problem specific SIM approximations.~Hence, as highlighted above, the resulting, by  the propose PINN approach, SIM approximations, can in turn be used in
a straightforward manner to build ROMs.  

The proposed PINN framework can be extended and serve as a basis for further developments in several directions. These include, indicatively, the construction of a non-linear map of the high-dimensional space to the fast and slow components using for example either autoencoders, or manifold learning coupled with a method of the solution of the pre-image problem (see for example \cite{patsatzis2023data,evangelou2023double}), and the epistemic uncertainty quantification of the computed SIM; the latter task can be attempted using for example Gaussian process regression models or ensembles of ANNs (for a comprehensive discussion see for example in \cite{psaros2023uncertainty}.

\section*{APPENDIX}
\appendix
\renewcommand{\theequation}{A.\arabic{equation}}
\renewcommand{\thefigure}{A.\arabic{figure}}
\setcounter{equation}{0}
\setcounter{figure}{0}
\section{The CSP framework to derive SIM approximations on the basis of QSSA, PEA and CSP with two iterations}
\label{app:CSP_SIM}

Computational Singular Perturbation (CSP) is an algorithmic methodology employed in the context of GSPT for the derivation of SIMs and reduced order models and the identification of their geometrical properties \citep{lam1989understanding,goussis1992study,hadjinicolaou1998asymptotic,valorani2005higher,goussis2012quasi,goussis2006efficient}.~Given a system in the form of Eq.~\eqref{eq:gen}, CSP follows an iterative procedure to locally approximate the \emph{fast} and \emph{slow} subspaces of the \emph{tangent} space, along which the solution of the IVP evolves.~In this way, the fast and slow dynamics of the system (i.e., the fast and slow timescales and modes) is decomposed to the respective subspaces.~Then, when the solution of the IVP evolves on the SIM, the projection of the vector field in Eq.~\eqref{eq:gen} onto the fast subspace is tangential to the SIM, thus providing an expression for the SIM approximation \citep{kuehn2015multiple,goussis2011model,kaper1999introduction,zagaris2004analysis}.~To allow this decompositions, the CSP algorithmic procedure is specifically designed to approximate the set of \emph{fast} and \emph{slow} basis vectors spanning the fast and slow subspaces.~At each iteration, the CSP fast and slow vectors approximate the ``true'' (and unknown) ones by $\mathcal{O}(\epsilon)$, where $\epsilon$ is the underlying timescale gap \cite{zagaris2004analysis,valorani2005higher,kaper2015geometry}.~In addition, it has been shown that the QSSA and PEA approximations introduce their own sets of vectors to approximate the fast and slow subspaces \cite{goussis2012quasi}, which however are only accurate when specific criteria are satisfied \cite{goussis2015model,goussis2012quasi,
patsatzis2023algorithmic}.~In what follows, we briefly present the basic concepts of CSP for the derivation of SIM approximations through fast/slow vectors, and then describe the iterative procedure for computing the CSP vectors, as well as the ones introduced by QSSA and PEA.

Let us consider the $N$-dim. autonomous system of ODEs in Eq.~\eqref{eq:gen} and assume (i) the existence of $M$ fast timescales, that are much faster than the remaining $N-M$ slow ones, and (ii) the dissipative nature of the $M$ fast timescales.~Under these assumptions, also assumed in Section~\ref{sec:Meth}, the system in Eq.~\eqref{eq:gen} exhibits an attracting SIM that attracts all the neighbouring trajectories (under the action of the fast $M$ timescales), which then evolve on the SIM with a slow flow (characterized by the remaining $N-M$ slow timescales) \citep{fenichel1979geometric,zagaris2004analysis,kaper1999introduction,hek2010geometric,lizarraga2020computational}.~According to the CSP framework, the vector field $\mathbf{F}(\mathbf{z})$ at every point of the trajectory can be decomposed in two components $\mathbf{F}_r(\mathbf{z})$ and $\mathbf{F}_s(\mathbf{z})$, such that:
\begin{equation}
    \dfrac{d\mathbf{z}}{dt} = \mathbf{F}(\mathbf{z}) = \mathbf{F}_r(\mathbf{z}) + \mathbf{F}_s(\mathbf{z}),
    \label{eq:VF1}
\end{equation}
where $\mathbf{F}_r,\mathbf{F}_s:\mathbb{R}^N\rightarrow\mathbb{R}^N$ are the local projections of the vector field to the fast and slow subspaces.~Now, dropping the dependency on $\mathbf{z}$, the $M$-dim. fast subspace is spanned by the $M$ column basis vectors $\mathbf{a}_i\in\mathbb{R}^N$ for $i=1,\ldots,M$ and the slow subspace is spanned by the $N-M$ column basis vectors $\mathbf{a}_j\in\mathbb{R}^N$ for $j=M+1,\ldots,N$.~Collecting all the basis vectors, the matrices:
\begin{equation}
    \mathbf{A}_r = \begin{bmatrix} \mathbf{a}_1 & \ldots & \mathbf{a}_M \end{bmatrix} \in \mathbb{R}^{N \times M}, \qquad \mathbf{A}_s = \begin{bmatrix} \mathbf{a}_{M+1} & \ldots & \mathbf{a}_N \end{bmatrix} \in \mathbb{R}^{N \times (N-M)},
    \label{eq:A_BVs}
\end{equation}
are formed, in which each column basis vector is associated to a timescale, ordered according to their magnitude; i.e., $\mathbf{a}_1$ corresponds to the fastest timescale, while $\mathbf{a}_N$ corresponds to the slowest one.~In addition, let a second set of $N$ row basis vectors, defined as the dual basis vectors $\mathbf{b}^i$ for $i=1,\ldots,M$ and $\mathbf{b}^j$ for $j=M+1,\ldots,N$, which is also collected in the matrices:
\begin{equation}
    \mathbf{B}^r = \begin{bmatrix} \mathbf{b}^1; & \ldots ; & \mathbf{b}^M \end{bmatrix}^\top \in \mathbb{R}^{M \times N}, \qquad \mathbf{B}^s = \begin{bmatrix} \mathbf{b}^{M+1}; & \ldots; & \mathbf{b}^N \end{bmatrix}^\top \in \mathbb{R}^{(N-M) \times N},
    \label{eq:B_BVs}
\end{equation}
where the symbol $;$ is used for indicating a new row.~Due to orthogonality of the basis vectors and their duals, the following matrix relations are implied:
\begin{equation}
\mathbf{B}^r \mathbf{A}_r = \mathbf{I}^M_M, \quad \mathbf{B}^r \mathbf{A}_s = \mathbf{0}^M_{N-M}, \quad
\mathbf{B}^s \mathbf{A}_r = \mathbf{0}^{N-M}_M, \quad  \mathbf{B}^s \mathbf{A}_s = \mathbf{I}^{N-M}_{N-M}, \quad
\mathbf{A}_r \mathbf{B}^r + \mathbf{A}_s \mathbf{B}^s = \mathbf{I}^N_N,
\label{eq:Orth}
\end{equation}
where $\mathbf{I}_q^q$ and $\mathbf{0}_p^q$ are the identity and zero matrices; superscripts and subscripts denoting row and column dimensions, respectively.

Using the above matrix definition of the fast and slow basis vectors, the vector field in Eq.~\eqref{eq:VF1} can be written as a projection onto the fast and slow subspaces as:
\begin{equation}
    \dfrac{d\mathbf{z}}{dt} = \mathbf{F}_r + \mathbf{F}_s = \mathbf{A}_r  \boldsymbol{f}^r + \mathbf{A}_s \boldsymbol{f}^s,
    \label{eq:VF2}
\end{equation}
where the $M$ and $(N-M)$-dim. vectors $\boldsymbol{f}^r=[f^1,\ldots,f^M]$ and $\boldsymbol{f}^s=[f^{M+1},\ldots,f^N]$ contain the amplitudes $f^i$, defined as the projections of the vector field to the dual basis vectors; i.e., $f^i = \mathbf{b}^i \cdot \mathbf{F}$ for $i=1,\ldots, N$.~When the system evolves on the SIM, the projections of the vector field to the fast subspace are negligible \citep{lam1989understanding,goussis1992study}, so that
\begin{equation}
    \boldsymbol{f}^r = \mathbf{B}^r \cdot \mathbf{F} \approx \mathbf{0}_1^M,
    \label{eq:VF3}
\end{equation}
which introduces $M$ algebraic equations for the derivation of the SIM approximation.~Note that both $\mathbf{B}^r$ and $\mathbf{F}$ are functions of $\mathbf{z}$.~Hence, in general, the $M$ algebraic equation of the SIM approximation are of implicit form, in contrast to the ones in Eq.~\eqref{eq:SIM} that are explicit.~In turn, when on the SIM, Eq.~\eqref{eq:VF3} implies that the governing equations of the slow flow, now contain projections only into the slow subspace, such that
\begin{equation}
\dfrac{d\mathbf{z}}{dt} \approx \mathbf{F}_s =  \mathbf{A}_s \boldsymbol{f}^s = (\mathbf{I}_N^N - \mathbf{A}_r \mathbf{B}^r) \mathbf{F}. 
\label{eq:VF4}
\end{equation}
Equations~(\ref{eq:VF3}, \ref{eq:VF4}) constitute the reduced system, which clearly relies only in the appropriate approximation of the fast basis vectors $\mathbf{A}_r$ and their duals $\mathbf{B}^r$.~In the following, we present the CSP, PEA and QSSA vectors for the approximation the ``true'' fast/slow basis vectors described above.

\subsection{The CSP basis vectors}
\label{app:sbCSP}

To approximate the fast/slow basis vectors in Eqs.~(\ref{eq:A_BVs}, \ref{eq:B_BVs}), CSP follows an iterative procedure that makes use of two complementary refinements; the $\mathbf{A}_r$- and the $\mathbf{B}^r$-refinements \citep{lam1989understanding,hadjinicolaou1998asymptotic,zagaris2004analysis,valorani2005higher}.~The $\mathbf{A}_r$-refinement updates the vectors in $\mathbf{A}_r$ and $\mathbf{B}^s$, leaving unaffected the vectors in $\mathbf{A}_s$ and $\mathbf{B}^r$; implementation of one such refinement is sufficient to guarantee the numerical stability of the reduced system in Eq.~\eqref{eq:VF4}.~The $\mathbf{B}^r$-refinement updates the vectors in $\mathbf{A}_s$ and $\mathbf{B}^r$, leaving unaffected the vectors in $\mathbf{A}_r$ and $\mathbf{B}^s$, and is associated to the numerical accuracy of the SIM approximation in Eq.~\eqref{eq:VF3} and thus, the numerical accuracy of the system in Eq.~\eqref{eq:VF4}.~Each CSP iteration consists of a $\mathbf{B}^r$-refinement followed by an $\mathbf{A}_r$-refinement; at each iteration the accuracy improvement of all the CSP vectors is of $\mathcal{O}(\epsilon)$ \citep{zagaris2004analysis,valorani2005higher}, where $\epsilon$ is the characteristic fast/slow timescale gap.~In this work, a maximum of two CSP iteration were employed, and since we are only interested in the SIM approximation, the  $\mathbf{A}_r$-refinement of the second iteration is omitted.~Thus, in the following, we present the iterative procedure followed at each time point of the trajectory, and use the notation $(k,m)$ to denote the output after the $k$-th $\mathbf{B}^r$-refinement and the $m$-th $\mathbf{A}_r$-refinement; $k=0,1,2$ and $m=0,1$.
\begin{enumerate}[font={\bfseries},label={Step \arabic*}]
    \item \emph{Initialization} \\ Guess initial set of vectors satisfying the orthogonality conditions in Eq.~\eqref{eq:Orth} for given $M$: 
    \begin{flalign*}
        \mathbf{A}_r(0,0), \qquad \mathbf{A}_s(0,0), \qquad \mathbf{B}^r(0,0), \qquad \mathbf{B}^s(0,0). & &
    \end{flalign*}
    \item \emph{$\mathbf{B}^r$-refinement} \\
    Using the Jacobian $\mathbf{J}$ of the system, compute:
    \begin{flalign*}
        \boldsymbol{\Lambda}(0,0) & = \mathbf{B}^r(0,0) \mathbf{J} \mathbf{A}_r(0,0), & & \\
        \mathbf{T}(0,0) & = \left( \boldsymbol{\Lambda}(0,0)\right)^{-1}.
    \end{flalign*}
    Update:
    \begin{flalign*}
        \mathbf{B}^r(1,0) & = \mathbf{T}(0,0) \mathbf{B}^r(0,0) \mathbf{J},  & & \\
        \mathbf{A}_r(1,0) & = \mathbf{A}_r(0,0), & & \\
        \mathbf{B}^s(1,0) & = \mathbf{B}^s(0,0), & & \\
        \mathbf{A}_s(1,0) & = \left( \mathbf{I}^N_N - \mathbf{A}_r(1,0) \mathbf{B}^r(1,0) \right) \mathbf{A}_s(0,0). & &
    \end{flalign*}
    \item \emph{$\mathbf{A}_r$-refinement} \\
    Using the Jacobian $\mathbf{J}$ of the system, compute:
    \begin{flalign*}
        \boldsymbol{\Lambda}(1,0) & = \mathbf{B}^r(1,0) \mathbf{J} \mathbf{A}_r(1,0), & & \\
        \mathbf{T}(1,0) & = \left( \boldsymbol{\Lambda}(1,0)\right)^{-1}.
    \end{flalign*}
    Update:
    \begin{flalign*}
        \mathbf{A}_r(1,1) & =\mathbf{J} \mathbf{A}^r(1,0) \mathbf{T}(1,0),  & & \\
        \mathbf{B}^r(1,1) & = \mathbf{B}^r(1,0), & & \\
        \mathbf{A}_s(1,1) & = \mathbf{A}_s(1,0), & & \\
        \mathbf{B}^s(1,1) & = \mathbf{B}^s(1,0) \left( \mathbf{I}^N_N - \mathbf{A}_r(1,1) \mathbf{B}^r(1,1) \right). & &
    \end{flalign*}
    \item \emph{$\mathbf{B}^r$-refinement} \\
    Using the Jacobian $\mathbf{J}$ of the system and the derivative of $\mathbf{B}^r(1,1)$ in time, compute:
    \begin{flalign*}
        \boldsymbol{\Lambda}(1,1) & = \left( \dfrac{d \mathbf{B}^r(1,1)}{dt}+\mathbf{B}^r(1,1) \mathbf{J} \right) \mathbf{A}_r(1,1), & & \\
        \mathbf{T}(1,1) & = \left( \boldsymbol{\Lambda}(1,1)\right)^{-1}.
    \end{flalign*}
    Update:
    \begin{flalign*}
        \mathbf{B}^r(2,1) & = \mathbf{T}(1,1) \left( \dfrac{d \mathbf{B}^r(1,1)}{dt}+\mathbf{B}^r(1,1) \mathbf{J} \right), & & \\
        \mathbf{A}_r(2,1) & = \mathbf{A}_r(1,1), & & \\
        \mathbf{B}^s(2,1) & = \mathbf{B}^s(1,1), & & \\
        \mathbf{A}_s(2,1) & = \left( \mathbf{I}^N_N - \mathbf{A}_r(2,1) \mathbf{B}^r(2,1) \right) \mathbf{A}_s(1,1). & &
    \end{flalign*}
\end{enumerate}
As evident by the updates in Steps 2-4, after every update the resulting set of basis vectors satisfies the orthogonality conditions in Eq.~\eqref{eq:Orth}.

For the initialization in Step 1, any arbitrary set of vectors $\mathbf{A}_r(0,0)$, $\mathbf{A}_s(0,0)$ and their duals $\mathbf{B}^r(0,0)$, $\mathbf{B}^s(0,0)$ satisfying the orthogonality conditions in Eq.~\eqref{eq:Orth} would suffice to approximate the ``true'' basis vectors, but it may require more than two CSP iterations \cite{goussis2012quasi,valorani2005higher}.~In fact, it has been shown that leading order accuracy and stability of the generated reduced model in Eq.~\eqref{eq:VF4} is guaranteed with only one CSP iteration if the vectors in $\mathbf{A}_r(0,0)$ and $\mathbf{A}_s(0,0)$ have negligible components along the slow and fast subspace, respectively \cite{lam1994csp,goussis2012quasi}.~To ensure accurate SIM approximations, we performed two CSP iterations for every problem considered, using as initial vectors:
\begin{equation}
\mathbf{A}_r(0,0) = \begin{bmatrix}
        \mathbf{I}^M_M \\[2pt]
        \mathbf{0}^{N-M}_M
    \end{bmatrix}, ~~     
    \mathbf{A}_s(0,0) = \begin{bmatrix}
		\mathbf{0}^M_{N-M} \\[2pt]
        \mathbf{I}^{N-M}_{N-M}
    \end{bmatrix}, ~~
    \mathbf{B}^r(0,0) = \begin{bmatrix}
         \mathbf{I}^M_M &  \mathbf{0}_{N-M}^M
    \end{bmatrix}, ~~
    \mathbf{B}^s(0,0) = \begin{bmatrix}
         \mathbf{0}^{N-M}_M &  \mathbf{I}_{N-M}^{N-M}
    \end{bmatrix}.
    \label{eq:ICbv}
\end{equation}
This initialization is natural when the first $M$ components of $\mathbf{z}$ are considered fast and the following $N-M$ are slow, thus guaranteeing, in these cases, leading order SIM approximations with only one CSP iteration; for all other cases, a re-ordering of the initialization is required to comply with the assumed $M$ fast components.

For the construction of the CSP basis vectors in Steps 2 and 3, the calculation of the Jacobian matrix $\mathbf{J}=\partial \mathbf{F}(\mathbf{z})/\partial \mathbf{z}$ of the system at each time point of the trajectory is required.~After the initialization in Eq.~\eqref{eq:ICbv}, the CSP vectors $\mathbf{A}_r(1,1)$ and $\mathbf{B}^r(1,1)$ approximate to leading order the ``true'' fast/slow basis vectors.~For the employment of the second $\mathbf{B}^r$-refinement in Step 4, except from the Jacobian matrix, the derivative of $\mathbf{B}^r(1,1)$ w.r.t. time is required, which can be calculated as:
\begin{equation}
    \dfrac{d\mathbf{B}^r(1,1)}{dt} = \mathbf{T}(0,0) \mathbf{B}^r(0,0) \dfrac{d \mathbf{J}}{dt} \left( \mathbf{I}^N_N - \mathbf{A}_r(0,0) \mathbf{B}^r(1,0)\right),
\end{equation}
where the derivative of the Jacobian $d \mathbf{J}/dt$  is used.~As a result, the second $\mathbf{B}^r$-refinement accounts for the curvature of the SIM, which was not considered in the first CSP iteration (Steps 2 and 3).~Hence, the resulting CSP vectors $\mathbf{A}_r(2,1)$ and $\mathbf{B}^r(2,1)$ provide higher order corrections to the ``true'' fast/slow basis vectors \citep{zagaris2004analysis,valorani2005higher}.

In this work, we performed two CSP iterations (two $\mathbf{B}^r$-refinement and one $\mathbf{A}_r$-refinement) for deriving the SIM approximations.~Following Eq.~\eqref{eq:VF3}, the resulting SIM approximations after one and two CSP iterations are:
\begin{equation}
    \boldsymbol{f}^r(1) = \mathbf{B}^r(1,1) \cdot \mathbf{F} = \mathbf{0}^M_1, \qquad \boldsymbol{f}^r(2) = \mathbf{B}^r(2,1) \cdot \mathbf{F}= \mathbf{0}^M_1,
\end{equation}
where the number $i=1,2$ in the amplitudes $\boldsymbol{f}^r(i)$ denotes the number of CSP iterations.~Note that the above expressions are, in general, formulated in an implicit form.

Finally, we note that, as shown in \cite{valorani2005higher}, carrying out multiple times the  first iteration in Steps 2 and 3, results in CSP vectors $\mathbf{a}_j$ and their duals $\mathbf{b}^j$ that are equivalent to the right, say $\mathbf{v}_j$, and left, say $\mathbf{u}^j$, eigenvectors of the Jacobian matrix (sorted by decreasing eigenvalues $\lambda_j$), for $j=1,\ldots,N$.~Note that (i) SIM approximations can be derived on the basis of the eigenvectors instead of the CSP vectors, which is the case of the CSP$_e$ SIM approximation for the MM mechanism \citep{patsatzis2019new}, considered in \cref{app:MM_SIMs}, and (ii) the criterion for the identification of $M$ in Eq.~\eqref{eq:CSPcrit} makes use of the first $M$ eigenvectors.

\subsection{The PEA and QSSA basis vectors}
\label{app:sbPEAQSSA}

For introducing the QSSA and PEA vectors, we first cast the system in Eq.~\eqref{eq:gen} in the form:
\begin{equation}
\dfrac{d\mathbf{z}}{dt} = \mathbf{F}(\mathbf{z}) = \mathbf{S}_1 R^1(\mathbf{z}) + \ldots \mathbf{S}_k R^k(\mathbf{z}) + \ldots + \mathbf{S}_K R^K(\mathbf{z}) = \mathbf{S}~ \mathbf{R}(\mathbf{z}), 
\label{eq:genSRform}
\end{equation} 
with the vector field $\mathbf{F}(\mathbf{z})$ written as a sum of $K$ terms $\mathbf{S}_k R^k(\mathbf{z})$, where $\mathbf{S}_k$ denotes the $N$-dim. stoichiometry vector of the $k$-th reaction and $R^k(\mathbf{z})$ denotes the respective reaction rate.~The last compact form of Eq.~\eqref{eq:genSRform} is formed by collecting all stoichiometric column vectors in the matrix $\mathbf{S}=[\mathbf{S}_1,\ldots,\mathbf{S}_K]\in \mathbb{R}^{N\times K}$ and all reaction rate terms in the column vector $\mathbf{R}(\mathbf{z})=[R^1(\mathbf{z}), \ldots,R^K(\mathbf{z})]^\top\in \mathbb{R}^{K}$; the dependency on $\mathbf{z}$ is dropped from now on.

Now, as already discussed in Section \ref{sbsb:CSP_PEA_QSSA}, to derive the QSSA and PEA approximations we need to make assumptions about the fast variables and the fast reactions; i.e., the variables in $\mathbf{z}$ and reactions in $
\mathbf{R}$ that are associated the most to the fast timescales $\tau_i$ for $i=1,\ldots,M$.~Here, we consider only the case were $M$ variables and/or $M$ reactions are related to the $M$ fast timescales; see \cite{patsatzis2023algorithmic} for more than $M$ fast variables/reactions.~We assume that the $M$ fast variables are included in the first $\mathbf{z}^M$ elements of $\mathbf{z}$ and that the $M$ fast reactions are included in the first $\mathbf{R}^M$ elements of $\mathbf{R}$; all other cases reduce to this one with a proper re-ordering of the system.~Then, the system in Eq.~\eqref{eq:genSRform} can be decomposed to the form:
\begin{equation}
\dfrac{d}{dt} \begin{bmatrix} \mathbf{z}^M ~~~~~ \\ \mathbf{z}^{N-M} \end{bmatrix}  = \mathbf{S}_M \mathbf{R}^M + \mathbf{S}_{K-M} \mathbf{R}^{K-M} = \begin{bmatrix} \mathbf{S}^M_M ~~~~~ \\ \mathbf{S}^{N-M}_M \end{bmatrix} \mathbf{R}^M + \begin{bmatrix} \mathbf{S}^M_{K-M} \\ \mathbf{S}^{N-M}_{K-M} \end{bmatrix} \mathbf{R}^{K-M}, 
\label{eq:genSR_dec}
\end{equation} 
where $\mathbf{z}^M=[z_1,\ldots,z_M]^\top$ contains the $M$ fast variables and $\mathbf{z}^{N-M}=[z_{M+1},\ldots,z_N]^\top$ the remaining slow ones, $\mathbf{R}^M=[R^1,\ldots,R^M]^\top$ contains the $M$ fast reactions and $\mathbf{R}^{N-M}=[R^{M+1},\ldots,R^N]^\top$ the remaining slow ones.~According to above decomposition, the stoichiometry matrix $\mathbf{S}$ is decomposed to its $N\times M$ left block and $N\times (K-M)$ right block matrices:
\begin{equation*}
\mathbf{S}_M=[\mathbf{S}_1,\ldots,\mathbf{S}_M]=\begin{bmatrix} \mathbf{S}^M_M ~~~~~ \\ \mathbf{S}^{N-M}_M \end{bmatrix}, \qquad 
\mathbf{S}_{K-M}=[\mathbf{S}_{M+1},\ldots,\mathbf{S}_K]=\begin{bmatrix} \mathbf{S}^M_{K-M} \\ \mathbf{S}^{N-M}_{K-M} \end{bmatrix}. 
\end{equation*}

Given the decomposed form of the system in Eq.~\eqref{eq:genSR_dec}, the PEA vectors for $M$ fast reactions are \cite{goussis2012quasi,goussis2015model}:
\begin{equation}
\mathbf{A}_r = \begin{bmatrix} \mathbf{I}^M_M \\ \mathbf{a}^{N-M}_M \end{bmatrix}, ~~ 
\mathbf{A}_s  = \begin{bmatrix} -\mathbf{V}^M_{N-M} \\ \mathbf{I}_{N-M}^{N-M} \end{bmatrix} \mathbf{Y}^{N-M}_{N-M}, ~~
\mathbf{B}^r = \mathbf{Y}^M_M \begin{bmatrix} \mathbf{I}^M_M & \mathbf{V}_{N-M}^M \end{bmatrix}, ~~
\mathbf{B}^s = \begin{bmatrix} -\mathbf{a}_M^{N-M} & \mathbf{I}_{N-M}^{N-M} \end{bmatrix}, 
\label{eq:PEA_BVs}
\end{equation}
where
\begin{align*}
\mathbf{a}_M^{N-M} & = \mathbf{S}^{N-M}_M \left( \mathbf{S}^M_M \right)^{-1} & \mathbf{V}_{N-M}^M & = \left( \dfrac{\partial \mathbf{R}^M}{\partial \mathbf{z}^M}\right)^{-1} \dfrac{\partial \mathbf{R}^M}{\partial \mathbf{z}^{N-M}} \\
\mathbf{Y}^{N-M}_{N-M} & = \left( \mathbf{I}^{N-M}_{N-M} + \mathbf{a}_M^{N-M} \mathbf{V}_{N-M}^M\right)^{-1} & \mathbf{Y}^{M}_{M} & = \left( \mathbf{I}^M_M + \mathbf{V}^M_{N-M} \mathbf{a}^{N-M}_M\right)^{-1}. 
\end{align*}
Substituting to Eq.~\eqref{eq:VF3} and using Eq.~\eqref{eq:genSR_dec}, the SIM approximation on the basis of PEA for the $M$ fast reactions reads:
\begin{equation}
\mathbf{S}_M^M \mathbf{R}^M + \mathbf{Y}^M_M\left( \mathbf{S}^M_{K-M}+\mathbf{V}^M_{N-M} \mathbf{S}^{N-M}_{K-M}\right) \mathbf{R}^{K-M} \approx \mathbf{0}^M_1.
\label{eq:PEA_SIMgen}
\end{equation}

In the case of QSSA, the QSSA vectors for $M$ fast variables are \cite{goussis2012quasi}:
\begin{equation}
\mathbf{A}_r = \begin{bmatrix} \mathbf{I}^M_M \\ \mathbf{a}^{N-M}_M \end{bmatrix}, \quad
\mathbf{A}_s  = \begin{bmatrix} \mathbf{0}^M_{N-M} \\ \mathbf{I}_{N-M}^{N-M} \end{bmatrix}, \quad
\mathbf{B}^r =  \begin{bmatrix} \mathbf{I}^M_M & \mathbf{0}_{N-M}^M \end{bmatrix}, \quad
\mathbf{B}^s = \begin{bmatrix} -\mathbf{a}_M^{N-M} & \mathbf{I}_{N-M}^{N-M} \end{bmatrix}. 
\label{eq:QSSA_BVs}
\end{equation} 
Substituting to Eq.~\eqref{eq:VF3} and using Eq.~\eqref{eq:genSR_dec}, the SIM approximation on the basis of QSSA for the $M$ fast variables reads:
\begin{equation}
\mathbf{S}_M^M \mathbf{R}^M +  \mathbf{S}^M_{K-M} \mathbf{R}^{K-M} \approx \mathbf{0}^M_1.
\label{eq:QSSA_SIMgen}
\end{equation}

\renewcommand{\theequation}{B.\arabic{equation}}
\renewcommand{\thefigure}{B.\arabic{figure}}
\setcounter{equation}{0}
\setcounter{figure}{0}

\section{The SIM approximations of the MM mechanism}
\label{app:MM_SIMs}

The SIM approximations underlying the slow evolution of the MM mechanism in Eq.~\eqref{eq:MM} are presented herein, as derived by the employment of the CSP, PEA and QSSA methods.~A summary of the eight, in total, SIM approximations for the MM mechanism (five constructed with CSP, one with PEA, and two with QSSA) is presented in Table~\ref{tb:MM_convSIMs}. 

\begin{table}[!h]
\centering
\resizebox{\textwidth}{!}{
\begin{tabular}{l c c c c c c c c c c c c c c c c}
\toprule
SIM approximation & \multicolumn{2}{c}{rQSSA} & \multicolumn{2}{c}{sQSSA} & \multicolumn{2}{c}{PEA} & \multicolumn{2}{c}{CSP$_s$(1)} & \multicolumn{2}{c}{CSP$_s$(2)}  & \multicolumn{2}{c}{CSP$_c$(1)} & \multicolumn{2}{c}{CSP$_c$(2)} & \multicolumn{2}{c}{CSP$_e$} \\
\midrule
Assumptions & \multicolumn{2}{c}{$s$ fast}  &  \multicolumn{2}{c}{$c$ fast}  &  \multicolumn{2}{c}{1st reaction fast}  & \multicolumn{2}{c}{$s$ fast} & \multicolumn{2}{c}{$s$ fast} & \multicolumn{2}{c}{$c$ fast} & \multicolumn{2}{c}{$c$ fast} & \\
\midrule
Solved for	& $s$	& $c$ & $s$	& $c$  & $s$	& $c$ 	& $s$	& $c$ & $s$	& $c$ & $s$	& $c$ & $s$	& $c$ & $s$	& $c$ \\
Explicit/Implicit  & E & E & E & E & E & E & E & E & E & I & E & E & I & I & E & I \\
\midrule
Equation & \eqref{eq:MM_rQSSA_expS} & \eqref{eq:MM_rQSSA_expC}  & \eqref{eq:MM_sQSSA_expS} & \eqref{eq:MM_sQSSA_expC}   & \eqref{eq:MM_PEA_expS} & \eqref{eq:MM_PEA_expC}   &  \eqref{eq:MM_CSPs11_expS} & \eqref{eq:MM_CSPs11_expC} & \eqref{eq:MM_CSPs21_expS} & \eqref{eq:MM_CSPs21_imp} & \eqref{eq:MM_CSPc11_expS} & \eqref{eq:MM_CSPc11_expC} & \eqref{eq:MM_CSPc21_imp} & \eqref{eq:MM_CSPc21_imp} & \eqref{eq:MM_CSPe_expS} & \eqref{eq:MM_CSPe_imp} \\
Requires Newton	& N & N & N & N & N & N & N & N & N & Y & N & N & Y & Y & N & Y\\
\bottomrule
\end{tabular}}
\caption{SIM approximations for the MM mechanism, constructed on the basis of QSSA, PEA and CSP methods.~For each approximation, we enlist (i) the assumptions made for constructing it, (ii) the functional form (E/I for explicit/implicit) when solved for the system's variables $s$ or $c$, and (iii) the respective equation and whether Newton iterations are required (Y/N for yes/no) to numerically solve the SIM approximation for each of the system's variables.}
\label{tb:MM_convSIMs} 
\end{table} 

\subsection{The CSP expressions with one and two iterations}
Here, we present the SIM approximations constructed on the basis of the CSP methodology described in Section~\ref{app:sbCSP}.~Given that $M=1$, we account for the two distinct cases where $s$ or $c$ is the fast variable.~In both cases, we construct the SIM approximations using one and two CSP iterations.~We re-iterate that when the assumption on the fast variable is correct, one CSP iteration is sufficient to guarantee the accuracy of the SIM approximation.~However, when this assumption is not correct, a second CSP iteration is required \cite{goussis2012quasi,valorani2005higher}.

We begin by the case where the fast variable is assumed to be $s$; i.e., the first $M=1$ variable in the state variable vector $[s,c]^\top$ of the system in Eq.~\eqref{eq:MM}.~According to Eq.~\eqref{eq:ICbv}, we initialize the CSP iterative procedure by setting the CSP vectors as:  
\begin{equation}
\mathbf{A}_r(0,0) = \begin{bmatrix}
        1 \\[2pt]
        0
    \end{bmatrix}, \quad     
    \mathbf{A}_s(0,0) = \begin{bmatrix}
		0 \\[2pt]
        1
    \end{bmatrix}, \quad
    \mathbf{B}^r(0,0) = \begin{bmatrix}
         1 &  0
    \end{bmatrix}, \quad
    \mathbf{B}^s(0,0) = \begin{bmatrix}
         0 & 1
    \end{bmatrix}.
    \label{eq:ICbv_MM1}
\end{equation}
Then, after performing the calculations in Steps 2 and 3 of the algorithm in Section~\ref{app:sbCSP}, the SIM approximation after one CSP iteration reads:
\begin{equation}
f^1(1)= k_{1b}c - k_{1f}(e_0-c)s + \dfrac{(k_{1b}+k_{1f} s)(c (k_{1b}+k_2+k_{1f} s) -e_0 k_{1f} s )}{k_{1f}(e_0-c) }=0, 
\label{eq:MM_CSPs11_imp}
\end{equation}
The above implicit expression can be solved analytically for both variables resulting to the CSP$_s$(1) SIM approximations; solving for $s$ the explicit expression is:
\begin{align}
\hat{s} = \dfrac{1}{2 k_{1f} (c-e_0)} \left(  \vphantom{\left(k_2^2\right)^{1/2}} c^2 k_{1f} - c (2 k_{1b}+2 e_0 k_{1f}+k_2)+e_0 (k_{1b}+e_0 k_{1f}) -\left((e_0 k_{1b}+(c-e_0)^2 k_{1f})^2- \right. \right. \nonumber \\
\left. \left. -2 c (-e_0 k_{1b}+(c-e_0)^2 k_{1f}) k_2 + c^2 k_2^2 \right)^{1/2}\right).
\label{eq:MM_CSPs11_expS}
\end{align}
and solving for $c$ the explicit expression is:
\begin{align}
\hat{c} = \dfrac{1}{2 k_{1f} (k_{1b}+k_{1f} s)} \left( \vphantom{\left(k_2^2\right)^{1/2}} k_{1b}^2+k_{1f} s (2 e_0 k_{1f}+k_2+k_{1f} s)+k_{1b} (e_0 k_{1f}+k_2+2 k_{1f} s) - \left( \vphantom{k_2^2} (-4 e_0 k_{1f}^2 s (k_{1b}+\right. \right. \nonumber \\
\left. \left. +k_{1f} s) (k_{1b}+k_{1f} (e_0+s))+(k_{1b}^2+k_{1f} s (2 e_0 k_{1f}+k_2+k_{1f} s)+k_{1b} (e_0 k_{1f}+k_2+2 k_{1f} s))^2\right)^{1/2} \right).
\label{eq:MM_CSPs11_expC}
\end{align}

For performing the second CSP iteration, we implemented Step 4 of the algorithm in Section~\ref{app:sbCSP}.~The resulting CSP$_s$(2) SIM approximation after two CSP iterations, reads:
\begin{align}
%
f^1(2)=\dfrac{(k_{1f} (-c+e_0+s)+k_{1b}) \left( \vphantom{3 c^2 k_{1f}^2} (e_0 k_{1f} s-c (k_{1b}+k_{1f} s+k_2)) \left( \vphantom{c^2_0}2 c^2 k_{1f} (k_{1b}+k_{1f} s)+\right. \right.}{k_{1f} (c-e_0) \left(e_0 \left(3 c^2 k_{1f}^2+k_{1f} s (-6 c k_{1f}+k_{1b}+k_2)- \right. \right.} \nonumber \\
\dfrac{\left. \left. +e_0 k_{1f} s (-4 c k_{1f}+k_{1b}+k_2)+e_0 k_{1b} (-3 c k_{1f}+k_{1b}+k_2)+e_0^2 k_{1f} (k_{1b}+2 k_{1f} s)\right) \right.}{\left. \left.-5 c k_{1b} k_{1f}+k_{1b}^2+k_{1b} k_2 \vphantom{3 c^2 k_{1f}^2} \right)+c^2 k_{1f} (-c k_{1f}+3 k_{1b}+3 k_{1f} s)+ \right.} \nonumber \\
\dfrac{\left.-k_{1f} (c-e_0)^2 (k_{1f} (-c+e_0+s)+k_{1b}) (c (k_{1b}+k_{1f} s)-e_0 k_{1f} s)\vphantom{3 c^2 k_{1f}^2} \right)}{ \left. + e_0^2 k_{1f} (3 k_{1f} (s-c)+2 k_{1b})+e_0^3 k_{1f}^2 \right)}=0.
\label{eq:MM_CSPs21_imp}
\end{align}
The above implicit expression can be solved analytically for both variables.~However, the resulting expressions are very extensive and thus, instead of the analytic solution for $c$, we numerically solved the above implicit expression using Newton's iterations, to get explicit functionals of $c$.~The analytic solution for $s$ is more manageable, reading:
\begin{align}
\hat{s}=\dfrac{1}{2 k_{1f}^2 (c-e_0) \left(3 c^2 k_{1f}+e_0 (-6 c k_{1f}+3 e_0 k_{1f}+k_{1b}+k_2)\right)} \left( \vphantom{\left(\left( k^2_2\right)^2\right)^{1/2}} e_0 k_{1f} (k_{1b}+k_2) (e_0 k_{1b}-c (2 k_{1b}+k_2)) - \right. \nonumber \\
\left. -2 k_{1f}^2 (c-e_0)^2 (c (3 k_{1b}+k_2)-e_0 k_{1b})+k_{1f}^3 (c-e_0)^4-\left(\vphantom{\left( k^2_2\right)^2} k_{1f}^2 \left(4 c k_{1b} (c-e_0) \left(3 c^2 k_{1f}-e_0 (6 c k_{1f} -  \right. \right.  \right. \right. \nonumber \\
\left. \left. \left. \left.- 3 e_0 k_{1f} - k_{1b}-k_2)\right) \left(-e_0 \left(3 c^2 k_{1f}^2-5 c k_{1b} k_{1f}-3 c k_{1f} k_2+k_{1b}^2+2 k_{1b} k_2+k_2^2\right)+c^2 k_{1f} (c k_{1f} - \right. \right.  \right. \right. \nonumber \\
\left. \left. \left. \left. -3 k_{1b} -2 k_2) -  e_0^2 k_{1f} (-3 c k_{1f} +2 k_{1b}+k_2)-e_0^3 -k_{1f}^2\right)+\left(c^4 k_{1f}^2-2 c^3 k_{1f} (2 e_0 k_{1f}+3 k_{1b}+k_2)+ \right. \right.  \right. \right. \nonumber \\ 
\left. \left. \left. \left. + 2 c^2 e_0 k_{1f} (3 e_0 k_{1f}  + 7 k_{1b} +2 k_2) -c e_0 \left(4 e_0^2 k_{1f}^2+10 e_0 k_{1b} k_{1f}+2 e_0 k_{1f} k_2+2 k_{1b}^2+3 k_{1b} k_2+k_2^2\right)+ \right. \right.  \right. \right. \nonumber \\
\left. \left. \left. \left. + e_0^2 \left((e_0 k_{1f}+k_{1b})^2+k_{1b} k_2\right)\right)^2\right)\right)^{1/2}\right).
\label{eq:MM_CSPs21_expS}
\end{align}
which is the CSP$_s(2)$ explicit expression for $s$.

In the case where the fast variable is assumed to be $c$, one should reorder the system in Eq.~\eqref{eq:MM} - so that the first $M=1$ variable in the new state vector $[c,s]^\top$ is the fast one - and then realize the same choice of initial CSP vectors as in Eq.~\eqref{eq:ICbv_MM1}.~Equivalently, we considered the system in Eq.~\eqref{eq:MM} and selected the initial CSP vectors as:
\begin{equation}
\mathbf{A}_r(0,0) = \begin{bmatrix}
        0 \\[2pt]
        1
    \end{bmatrix}, \quad     
    \mathbf{A}_s(0,0) = \begin{bmatrix}
		1 \\[2pt]
        0
    \end{bmatrix}, \quad
    \mathbf{B}^r(0,0) = \begin{bmatrix}
         0 &  1
    \end{bmatrix}, \quad
    \mathbf{B}^s(0,0) = \begin{bmatrix}
         1 & 0
    \end{bmatrix},
    \label{eq:ICbv_MM2}
\end{equation}
which agree naturally with the assumption of the second variable $c$ being fast.~Then, after performing the calculations in Steps 2 and 3 in Section~\ref{app:sbCSP}, the SIM approximation after one CSP iteration reads:
\begin{equation}
f^1(1)= k_{1b} c +k_2c-(e0-c) k_{1f} s+\dfrac{(e_0-c) k_{1f} (c k_{1b}-(e_0-c) k_{1f} s)}{k_{1b}+k_2+k_{1f} s} = 0, 
\label{eq:MM_CSPc11_imp}
\end{equation}
which can be solved analytically for both variables resulting to the CSP$_c$(1) SIM approximation; solving for $s$ the explicit expression is:
\begin{align}
\hat{s}=\dfrac{1}{2 k_{1f} (c-e_0)} \left( 
\vphantom{\left( \left( k_2^2\right)^{1/2}\right)^2} 
c^2 k_{1f} (e_0-2 c) (e_0 k_{1f}+k_{1b}+k_2) - \left( \vphantom{\left( k_2^2\right)^{1/2}} 
2 k_2 \left(k_{1f} (e_0-2 c) (c-e_0)^2 + \right. \right. \right. \nonumber \\
   \left. \left. \left. +e_0^2 k_{1b}\right)+\left(k_{1f} (c-e_0)^2+e_0 k_{1b}\right)^2+e_0^2 k_2^2\right)^{1/2}\right),
\label{eq:MM_CSPc11_expS}
\end{align}
and solving for $c$ the explicit expression is:
\begin{align}
\hat{c}=\dfrac{1}{2 k_{1f} (k_{1b}+k_{1f} s)} \left( \vphantom{\left( \left( k_2^2\right)^{1/2}\right)^2} e_0  k_{1f} (k_{1b} +2 k_{1f} s ) +  (k_{1b} + k_{1f}s +k_2)^2 - \left( \vphantom{\left( k_2^2\right)^{1/2}} \left(e_0 k_{1b} k_{1f}+2 e_0 k_{1f}^2 s+k_{1b}^2+  \right. \right. \right. \nonumber \\
\left. \left. \left. + 2 k_{1b} (k_{1f} s+k_2)+(k_{1f} s+k_2)^2\right)^2-4 e_0 k_{1f}^2 s (k_{1b}+k_{1f} s) (k_{1f} (e_0+s)+k_{1b}+k_2) \right)^{1/2} \right).
\label{eq:MM_CSPc11_expC}
\end{align}

For performing the second CSP iteration, we implemented Step 4 of the algorithm in Section~\ref{app:sbCSP}.~The resulting CSP$_s$(2) SIM approximation after two CSP iterations reads:
\begin{align}
f^1(2) = \dfrac{(k_{1b}+k_{1f} s+k_2) \left[k_{1b} \left(c^2 k_{1f}^2 k_2- \right. \right.}{-\left(k_{1b} k_{1f} (-c+e_0+2 s)+k_{1f}^2 s (e_0-c)+k_{1b}^2+2 k_{1b} k_2+(k_{1f} s+k_2)^2\right) \left[c^3 k_{1f}^2 k_2 (k_{1b}+ \right. } \nonumber \\
\dfrac{\left. \left. -c k_{1f} \bigl(e_0 k_{1f} (k_{1f} s+2 k_2) +3 (k_{1f} s+k_2) (3 k_{1f} s+k_2) \bigr)+ e_0^2 k_{1f}^2 (k_{1f} s+k_2)+ \right. \right.}{\left. +k_{1f} s)-c^2 k_{1f} \left( \vphantom{\left((k_2^2)^2\right)^2} k_{1b}^2 (e_0 k_{1f}+9 k_{1f} s+6 k_2)+k_{1b} \left( \vphantom{(k_2^2)} e_0 k_{1f} ( k_{1f} s+2 k_2)+3 (k_{1f} s+k_2) (3 k_{1f} s+ \right. \right. \right. } \nonumber \\
\dfrac{\left. \left. + e_0 k_{1f} \bigl(8 k_{1f}^2 s^2 +9 k_{1f} k_2 s+2 k_2^2\bigr)+4 (k_{1f} s+k_2)^3\right)+ \right.}{\left. \left. \left. +k_2)\vphantom{(k_2^2)} \right)+3 k_{1f} s \left(e_0 k_{1f} k_2+(k_{1f} s+k_2)^2\right)+3 k_{1b}^3 \vphantom{\left((k_2^2)^2\right)^2} \right)+c \left(k_{1f} s \bigl( \vphantom{(k_{1f}^2)} k_2 (3 e_0^2 k_{1f}^2+4 e_0 k_{1f} k_2+4 k_2^2)+ \right. \right.  }  \nonumber \\
\dfrac{\left. + k_{1f}^3 k_2 s \left(c^2-2 c (e_0+3 s) +(e_0+s) (e_0+4 s)\right)+ \right.}{\left. \left.  +2 k_{1b}^2 (5 e_0 k_{1f}+6 k_2)+2 k_{1b} (e_0 k_{1f}+k_2) (e_0 k_{1f}+6 k_2)+4 k_{1b}^3\bigr)+2 k_{1f}^2 s^2 \left(7 e_0 k_{1b} k_{1f}+ \right. \right. \right. } \nonumber \\
\dfrac{\left. + k_{1b}^2 \left(-c k_{1f} (e_0 k_{1f}+9 k_{1f} s+6 k_2)+e_0^2 k_{1f}^2+e_0 k_{1f} (7 k_{1f} s+4 k_2)+6 (k_{1f} s+k_2)^2\right)+ \right.}{\left. \left. \left. +5 e_0 k_{1f} k_2+3 k_{1b}^2+6 k_{1b} k_2+3 k_2^2\right)+ 2 k_{1f}^3 s^3 (3 e_0 k_{1f}+2 k_{1b}+2 k_2)+(k_{1b}+k_2) \bigl( k_{1b} k_2 (2 e_0 k_{1f}+\right. \right.  } \nonumber \\
\dfrac{\left. + k_{1b}^3 \bigl(-3 c k_{1f}+2 e_0 k_{1f} +4 (k_{1f} s+k_2)\bigr)+k_{1f}^4 s^3 (-3 c+3 e_0+s)+ \right. }{\left. \left.  +3 k_{1b})+k_{1b} (e_0 k_{1f}+k_{1b})^2+3 k_{1b} k_2^2+k_2^3\bigr)+k_{1f}^4 s^4\right)+e_0 k_{1f} s \Bigl(-k_{1f} s \left(5 e_0 k_{1b} k_{1f}+4 e_0 k_{1f} k_2+ \right. \right. } \nonumber \\
\dfrac{\left. +k_{1f}^2 k_2^2 s (-3 c+2 e_0+6 s)+k_{1b}^4+4 k_{1f} k_2^3 s+ \right.}{\left. \left. +3 k_{1b}^2+6 k_{1b} k_2+3 k_2^2\right)-3 k_{1f}^2 s^2 (e_0 k_{1f}+k_{1b}+k_2)-(k_{1b}+k_2) \bigl(k_2 (e_0 k_{1f}+2 k_{1b})+ \right. } \nonumber \\
\dfrac{\left. +k_2^4\right]}{\left. +(e_0 k_{1f}+k_{1b})^2+ k_2^2\bigr)-k_{1f}^3 s^3\Bigr) \right]}.
\label{eq:MM_CSPc21_imp}
\end{align}
As evident, the above implicit expression is very extensive.~While it can be solved explicitly for both $s$ and $c$ variables, the resulting expressions are even more complicated.~Thus, we use the CSP$_c$(2) implicit expression in Eq.~\eqref{eq:MM_CSPc21_imp} and employ Newton's iterations to numerically acquire the explicit functional in terms of $s$ or $c$. 

Finally, in the case of the MM mechanism, we also considered the SIM approximation constructed in \cite{patsatzis2019new} via CSP by using the eigenvectors of the Jacobian, instead of the CSP basis vectors.~The resulting CSP$_e$ approximation of the SIM reads:
\begin{align}
f^1_e = \dfrac{1}{2} \Bigr(-c^2 k_{1f}+c (e_0 k_{1f}-k_{1b}-k_{1f} s-k_2)+2 e_0 k_{1f} s-c \bigl( k_{1f} (e_0-c+s)+ \nonumber \\
+ k_{1b}+k_2)^2+4 k_{1f} k_2 (c-e_0)\bigr)^{1/2} \Bigl).
\label{eq:MM_CSPe_imp}
\end{align}
The above implicit expression can be solved analytically for both variables.~The resulting explicit expression for $c$ is very extensive and thus, we solved the implicit expression Eq.~\eqref{eq:MM_CSPe_imp} numerically over $c$, using Newton's iterations.~On the other hand, solving for $s$ results to a simpler expression, yielding the CSP$_e$ explicit SIM approximation for $s$:
\begin{align}
\hat{s} = \dfrac{c}{2 e_0 k_{1f} (c-e_0)} \left[(c^2 k_{1f}-e_0 (2 c k_{1f}+k_{1b}+k_2)+e_0^2 k_{1f}-\Bigl(\bigl(k_{1f} (e_0-c)^2- \right. \nonumber \\
\left. - e_0 (k_{1b}+k_2)\bigr)^2+4 e_0 k_{1b} k_{1f} (c-e_0)^2\Bigr)^{1/2} \right].
\label{eq:MM_CSPe_expS}
\end{align}

\subsection{The PEA and QSSA expressions}

Here, we report the three existing PEA and QSSA approximations for the SIM of the MM mechanism; the procedure for their derivation, as described in Section~\ref{app:sbPEAQSSA}, was followed in  \cite{patsatzis2023algorithmic} in Sections 6.2-6.6.~For the employment of the PEA, we assume that the first binding (and reversible) reaction is fast, so that its partial equilibrium can be realized as suggested by \citep{michaelis1913kinetik}.~Following  Eq.~\eqref{eq:PEA_SIMgen}, the resulting SIM approximation \citep{patsatzis2019new} reads:
\begin{equation}
    (e_0-c)s - k_{1b}c/k_{1f} - \dfrac{k_{1b}/k_{1f}+s}{k_{1b}c/k_{1f}+s + e_0-c} k_2c/k_{1f} \approx 0, 
    \label{eq:MM_PEA_imp}
\end{equation}
which is an implicit expression that can be solved for either $s$ or $c$ for deriving explicit functionals of the PEA SIM approximations; solving for $s$ yields:
\begin{align}
\hat{s} = \dfrac{1}{2 k_{1f} (c-e_0)} \Bigl( c^2 k_{1f}-c (2 e_0 k_{1f}+2 k_{1b}+k_2)+e_0 (e_0 k_{1f}+k_{1b}) - \bigl(c^2 k_2^2-  \nonumber \\
 - 2 c k_2 \left(k_{1f} (c-e_0)^2-e_0 k_{1b}\right)+\left(k_{1f} (c-e_0)^2+e_0 k_{1b}\right)^2 \bigr)^{1/2}\Bigr),  \label{eq:MM_PEA_expS}
\end{align}
and solving for $c$ yields:
\begin{align}
\hat{c} = \dfrac{1}{2 k_{1f} (k_{1b}+k_{1f} s)} \Bigl( k_{1b} (e_0 k_{1f}+2 k_{1f} s+k_2)+k_{1f} s (2 e_0 k_{1f}+k_{1f} s+k_2)+k_{1b}^2 - \bigl( \left(k_{1b} (e_0 k_{1f}+ \right. \nonumber \\
\left. + 2 k_{1f} s+k_2)+k_{1f} s (2 e_0 k_{1f}+k_{1f} s+k_2)+k_{1b}^2\right)^2-4 e_0 k_{1f}^2 s (k_{1b}+k_{1f} s) (k_{1f} (e_0+s)+k_{1b})\bigr)^{1/2} \Bigr). \label{eq:MM_PEA_expC}
\end{align}

For the employment of QSSA, we first assume that the fast variable is $s$, so that quasi steady-state of $s$ can be realized.~According to Eq.~\eqref{eq:QSSA_SIMgen}, the resulting SIM approximation reads:
\begin{equation}
    k_{1f}(e_0-c)s - k_{1b} c \approx 0, 
    \label{eq:MM_rQSSA_imp}
\end{equation}
which is the well-known \emph{reverse} QSSA (rQSSA) SIM approximation \citep{schnell2000enzyme}.~To derive explicit functionals of the rQSSA SIM approximation, we solve the implicit expression in Eq.~\eqref{eq:MM_rQSSA_imp} for $s$, yielding:
\begin{align}
\hat{s} = \dfrac{k_{1b} c}{k_{1f}(e_0-c)},
\label{eq:MM_rQSSA_expS}
\end{align} 
and for $c$, implying:
\begin{align}
\hat{c} = \dfrac{k_{1f}e_0s}{k_{1b}+k_{1f}s}. \label{eq:MM_rQSSA_expC}
\end{align} 
 
Finally, we also assume that the fast variable is $c$ to enable the QSSA for $c$.~According to Eq.~\eqref{eq:QSSA_SIMgen}, the resulting SIM approximation reads:
\begin{equation}
    k_{1f}(e_0-c)s - k_{1b} c - k_2 c \approx 0, 
    \label{eq:MM_sQSSA_imp}
\end{equation}
which is the well-known \emph{standard} QSSA (sQSSA) SIM approximation \citep{bowen1963singular,heineken1967mathematical,segel1989quasi}.~For the explicit functionals of the sQSSA SIM approximation, we solve the implicit expression in Eq.~\eqref{eq:MM_sQSSA_imp} for $s$, resulting to:
\begin{align}
\hat{s} = \dfrac{(k_{1b}+k_2)c}{k_{1f}(e_0-c)},
\label{eq:MM_sQSSA_expS}
\end{align} 
and for $c$, yielding:
\begin{align}
\hat{c} = \dfrac{k_{1f}e_0s}{k_2+k_{1b}+k_{1f}s}. 
\label{eq:MM_sQSSA_expC}
\end{align}

\renewcommand{\theequation}{C.\arabic{equation}}
\renewcommand{\thefigure}{C.\arabic{figure}}
\setcounter{equation}{0}
\setcounter{figure}{0}
\section{The SIM approximations of the TMDD mechanism}
\label{app:TMDD_SIMs}

The SIM approximations underlying the slow evolution of the TMDD mechanism in Eq.~\eqref{eq:TMDD} are presented herein, as derived by the employment of the CSP, PEA and QSSA methods.~We consider all possible cases for the construction of an $M=1$-dim. SIM approximation; some of these expressions have been derived for $M=2$ as well \cite{patsatzis2016asymptotic}.~Here, we solve the resulting SIM expressions only for $L$ and $R$ variables, while we do not solve for $RL$ to avoid presenting lengthy expressions that are not used in the main text.~A summary of the ten, in total, SIM approximations for the TMDD mechanism (six constructed by CSP, one by PEA, and three by QSSA) is presented in Table~\ref{tb:TMDD_convSIMs}. 

\begin{table}[!h]
\centering
\footnotesize
\begin{tabular}{l c c c c c c c c}
\toprule
SIM approximation & \multicolumn{2}{c}{QSSA$_L$} & \multicolumn{2}{c}{QSSA$_R$} & \multicolumn{2}{c}{QSSA$_{RL}$} & \multicolumn{2}{c}{PEA} \\
\midrule
Assumptions & \multicolumn{2}{c}{$L$ fast}  &  \multicolumn{2}{c}{$R$ fast} &  \multicolumn{2}{c}{$RL$ fast}  &  \multicolumn{2}{c}{1st reaction fast}  \\
\midrule
Solved for	& $L$	& $R$ & $L$	& $R$ & $L$	& $R$ & $L$	& $R$  \\
Explicit/Implicit  & E & E & E & E & E & E & E & E  \\
\midrule
Equation & \eqref{eq:TMDD_QSSAL_expL} & \eqref{eq:TMDD_QSSAL_expR}  & \eqref{eq:TMDD_QSSAR_expL} & \eqref{eq:TMDD_QSSAR_expR} & \eqref{eq:TMDD_QSSARL_expL} & \eqref{eq:TMDD_QSSARL_expR} &  \eqref{eq:TMDD_PEA_expL} & \eqref{eq:TMDD_PEA_expL} \\
Requires Newton	& N & N & N & N & N & N & N & N \\
\bottomrule
\end{tabular} \\[2pt]
\resizebox{\textwidth}{!}{\begin{tabular}{l c c c c c c c c c c c c}
\toprule
SIM approximation &  \multicolumn{2}{c}{CSP$_L$(1)} & \multicolumn{2}{c}{CSP$_L$(2)}  & \multicolumn{2}{c}{CSP$_R$(1)} & \multicolumn{2}{c}{CSP$_R$(2)} & \multicolumn{2}{c}{CSP$_{RL}$(1)} & \multicolumn{2}{c}{CSP$_{RL}$(2)} \\
\midrule
Assumptions & \multicolumn{2}{c}{$L$ fast} & \multicolumn{2}{c}{$L$ fast} & \multicolumn{2}{c}{$R$ fast} & \multicolumn{2}{c}{$R$ fast} & \multicolumn{2}{c}{$RL$ fast} & \multicolumn{2}{c}{$RL$ fast} \\
\midrule
Solved for	& $L$	& $R$ & $L$	& $R$ & $L$	& $R$ & $L$	& $R$ & $L$	& $R$ & $L$	& $R$ \\
Explicit/Implicit  & E & E & I & I & E & E & I & I & E & E & I & I \\
\midrule
Equation & \eqref{eq:TMDD_CSPL11_expL} & \eqref{eq:TMDD_CSPL11_expR} & \eqref{eq:TMDD_CSPL21_imp} & \eqref{eq:TMDD_CSPL21_imp} & \eqref{eq:TMDD_CSPR11_expL} & \eqref{eq:TMDD_CSPR11_expR} & \eqref{eq:TMDD_CSPR21_imp} & \eqref{eq:TMDD_CSPR21_imp} & \eqref{eq:TMDD_CSPRL11_expL} & \eqref{eq:TMDD_CSPRL11_expR} & \eqref{eq:TMDD_CSPRL21_imp} & \eqref{eq:TMDD_CSPRL21_imp} \\
Requires Newton	& N & N & Y & Y & N & N & Y & Y & N & N & Y & Y\\
\bottomrule
\end{tabular}}
\caption{SIM approximations for the TMDD mechanism, constructed on the basis of QSSA, PEA and CSP methods.~For each approximation, we enlist (i) the assumptions made for constructing it, (ii) the functional form (E/I for explicit/implicit) when solved for the system's variables $L$ or $R$, and (iii) the respective equation and whether Newton iterations are required  (Y/N for yes/no) to numerically solve the SIM approximation for each of the system's variables.}
\label{tb:TMDD_convSIMs} 
\end{table} 

\subsection{The CSP expressions with one and two iterations}

First, we construct SIM approximations on the basis of the CSP methodology, as described in Section~\ref{app:sbCSP}.~Considering that $M=1$, we account for the three distinct cases where one of $L$, $R$ and $RL$ might the fast variable and construct the CSP approximations accordingly, using one and two CSP iterations.

We first consider the case where the fast variable is assumed to be $L$, the first variable in Eq.~\eqref{eq:TMDD}.~According to Eq.~\eqref{eq:ICbv}, the initial set of CSP vectors is set to:  
\begin{equation}
\mathbf{A}_r(0,0) = \begin{bmatrix}
        1 \\[2pt]
        0 \\[2pt]
        0
    \end{bmatrix}, \quad     
    \mathbf{A}_s(0,0) = \begin{bmatrix}
		0 & 0 \\[2pt]
        1 & 0 \\[2pt]
        0 & 1 
    \end{bmatrix}, \quad
    \mathbf{B}^r(0,0) = \begin{bmatrix}
         1 &  0 & 0
    \end{bmatrix}, \quad
    \mathbf{B}^s(0,0) = \begin{bmatrix}
         0 & 1 & 0 \\[2pt]
         0 & 0 & 1
    \end{bmatrix}.
    \label{eq:ICbv_TMDD1}
\end{equation}
With the above set, we performed Steps 2 and 3 of the Algorithm in Section~\ref{app:sbCSP} to get the CSP vectors after one iteration.~The resulting SIM approximation, according to Eq.~\eqref{eq:VF3}, reads:
\begin{multline}
f^1(1) = -k_{on} L.R+k_{off} RL-k_{el} L+ \\
+\dfrac{k_{on} L (k_{syn}-R (k_{deg}+k_{off}+k_{on} L))+k_{off} RL (k_{int}+k_{off}+k_{on} L)}{k_{el}+k_{on} R} = 0,
\label{eq:TMDD_CSPL11_imp}
\end{multline}
which is in full agreement with the expression provided in Eq.~(84) in \citep{patsatzis2016asymptotic}.~The above implicit expression of the SIM can be solved analytically for both $L$ and $R$, resulting to the CSP$_L$(1) SIM approximations; solving for $L$, the explicit expression is:
\begin{multline}
\hat{L} = \dfrac{1}{2 k_{on}^2 R} \bigg(k_{on} (-R (k_{deg}+k_{off}+k_{on} R)+k_{off} RL+k_{syn})-k_{el}^2-2 k_{el} k_{on} R + \Big(\big(k_{on} (R (k_{deg} +k_{off}+ \\
+k_{on} R)-k_{off} RL-k_{syn})+k_{el}^2+2 k_{el} k_{on} R\big)^2+4 k_{off} k_{on}^2 R.RL (k_{el}+k_{int}+k_{off}+k_{on} R)\Big)^{1/2} \bigg),
\label{eq:TMDD_CSPL11_expL}
\end{multline}
and solving for $R$ the explicit expression is:
\begin{multline}
\hat{R} = \dfrac{1}{2 k_{on} L} \Big(- L (k_{deg}+2 k_{el}+k_{off}+k_{on} L)+k_{off} RL + \big(L^2 ((k_{deg}+ k_{off}+k_{on} L) (k_{deg}+ \\
+4 k_{el}+k_{off}+k_{on} L)+4 k_{on} k_{syn})+2 k_{off} L. RL (-k_{deg}+2 k_{int}+k_{off}+k_{on} L)+k_{off}^2 RL^2\big)^{1/2} \Big).
\label{eq:TMDD_CSPL11_expR}
\end{multline}
For performing the second CSP iteration, we implemented Step 4 of the algorithm in Section~\ref{app:sbCSP}.~The resulting CSP$_L$(2) SIM approximation after two CSP iteration, reads:
\begin{multline}
f^1(2) = \dfrac{\left(k_{el}^2+2 k_{el} k_{on} R+k_{on} R (k_{off}+k_{on} (L+R))\right) \Big(k_{off} (RL (k_{int}+k_{off})-k_{on} L. R) \big(k_{on} (-k_{deg} R+ }{(k_{el}+k_{on} R) \Big(k_{el} k_{on} R \big(k_{on} (L (k_{deg}+k_{on} L) + }\\ 
\dfrac{+R (k_{int}+k_{off}+k_{on} R)+k_{off} RL+k_{syn})+k_{el}^2+k_{el} (k_{int}+k_{off}+k_{on} L+2 k_{on} R)\big)+ }{ +6 k_{on} L.R+4 k_{on} R^2)+k_{int} k_{off}+k_{off}^2+k_{off} k_{on} (2 L+4 R-RL)\big)+} \\
\dfrac{+k_{on} (-R (k_{deg}+k_{on} L)+k_{off} RL+k_{syn}) \big(L \left(k_{el} (k_{deg}+2 k_{el}+k_{off}+k_{on} L)+k_{on} R (4 k_{el}+k_{off})+ \right. }{+k_{on}^2 R \big(k_{off} (-k_{deg} R+R (k_{int}+k_{on} (L+2 R))+k_{on} RL (L-R)+k_{syn}) } \\
\dfrac{\left.  +2 k_{on}^2 R^2+k_{on} k_{syn}\right)-k_{off} RL (k_{el}+k_{on} (R-L))\big)-\big((k_{el}+k_{on} R) \left(k_{el}^2+2 k_{el} k_{on} R+ \right.}{ +k_{off}^2 (R+RL)+k_{on} \left(k_{on} R^2 (3 L+R)+k_{syn} L\right)\big)+k_{el}^4+4 k_{el}^3 k_{on} R+} \\
\dfrac{\left. +k_{on} R (k_{off}+k_{on} (L+R))\right) (k_{el} L-k_{off} RL+k_{on} L. R)\big)\Big)}{+k_{el}^2 k_{on} R (2 k_{off}+3 k_{on} (L+2 R))\Big)} = 0,
\label{eq:TMDD_CSPL21_imp}
\end{multline}
which is a very complicated implicit expression.~Although an explicit analytic solution for $L$ can be obtained (but not for $R$), the resulting expression is very complicated as well.~Hence, we solve the CSP$_L$(2) implicit expression in Eq.~\eqref{eq:TMDD_CSPL21_imp} numerically using Newton's iterations to acquire an explicit functional in terms of $L$ or $R$.

Next, we consider the case where the fast variable is assumed to be $R$.~Since this is the second variable of the state vector in Eq.~\eqref{eq:TMDD}, instead of reordering the system, we select the set of initial CSP vectors as:
\begin{equation}
\mathbf{A}_r(0,0) = \begin{bmatrix}
        0 \\[2pt]
        1 \\[2pt]
        0
    \end{bmatrix}, \quad     
    \mathbf{A}_s(0,0) = \begin{bmatrix}
		1 & 0 \\[2pt]
        0 & 0 \\[2pt]
        0 & 1 
    \end{bmatrix}, \quad
    \mathbf{B}^r(0,0) = \begin{bmatrix}
         0 &  1 & 0
    \end{bmatrix}, \quad
    \mathbf{B}^s(0,0) = \begin{bmatrix}
         1 & 0 & 0 \\[2pt]
         0 & 0 & 1
    \end{bmatrix},
    \label{eq:ICbv_TMDD2}
\end{equation}
that is consistent with the assumption of $R$ being the fast variable.~For performing the first CSP iteration, we follow Steps 2 and 3 of the Algorithm in Section~\ref{app:sbCSP}.~According to Eq.~\eqref{eq:VF3}, the new CSP vectors, result to the SIM approximation:
\begin{multline}
f^1(1) = -k_{deg} R+k_{off} RL-k_{on} L R+k_{syn} + \\
+ \dfrac{k_{off} RL (k_{int}+k_{off}+k_{on} R)-k_{on} L. R (k_{el}+k_{off}+k_{on} R)}{k_{deg}+k_{on} L} = 0,
\label{eq:TMDD_CSPR11_imp}
\end{multline}
which is in full agreement with the expression provided in Eq.~(79) in \citep{patsatzis2016asymptotic}.~The above implicit SIM approximation provided by CSP with one iteration, can be solved analytically for both $L$ and $R$.~The resulting analytic CSP$_R$(1) approximation, when solving for $L$, yields: 
\begin{multline}
\hat{L} = \dfrac{1}{2 k_{on} R} \Big(-2 k_{deg} R-k_{el} R-k_{off} R+k_{off} RL-k_{on} R^2+k_{syn}+\big(-R (2 k_{deg}+k_{el}+k_{off}+ \\
+k_{on} R)+k_{off} RL+k_{syn})^2+4 R (k_{off} RL (k_{deg}+k_{int}+k_{off}+k_{on} R)+k_{deg} (k_{syn}-k_{deg} R))\big)^{1/2}\Big),
\label{eq:TMDD_CSPR11_expL}
\end{multline}
and when solving for $R$, yields:
\begin{multline}
\hat{R} = \dfrac{1}{2 k_{on}^2 L} \bigg(-(k_{deg}+k_{on} L)^2-k_{on} L (k_{el}+k_{off})+k_{off} k_{on} RL + \Big(\big((k_{deg}+k_{on} L)^2 + k_{on} L (k_{el}+ \\
 + k_{off}) - k_{off} k_{on} RL\big)^2+4 k_{on}^2 L (k_{off} RL (k_{deg}+k_{int}+k_{off}+k_{on} L)+k_{syn} (k_{deg}+k_{on} L))\Big)^{1/2} \bigg).
 \label{eq:TMDD_CSPR11_expR}
\end{multline}
For the second CSP iteration, we implemented Step 4 of  the Algorithm in Section~\ref{app:sbCSP}, to obtain the resulting CSP$_R$(2) SIM approximation, as:
\begin{multline}
f^1(2) = \dfrac{\big( (k_{deg}+k_{on}L)^2+k_{on} L (k_{off}+k_{on} R)\big) \Big(-k_{off} (k_{on} L.R-RL (k_{int}+k_{off})) \big((k_{deg}+k_{on} L)^2+}{(k_{deg}+k_{on} L) \Big(k_{deg}^4+4 k_{deg}^3 k_{on} L+k_{deg}^2 k_{on} L (2 k_{off}+ } \\
\dfrac{+ k_{deg} (k_{int}+ k_{off}+k_{on} R)+k_{on} L (-k_{el}+k_{int}+k_{off})+k_{off} k_{on} RL\big)-k_{on} (k_{el} L-k_{off} RL+}{+3 k_{on} (2 L+R))+k_{deg} k_{on} L \big(k_{on} \left(k_{el} R+k_{on} (4 L^2+6 L.R+R^2)-k_{syn}\right)+ } \\
\dfrac{ +k_{on} L. R) \big(2 k_{deg}^2 R+k_{deg} (R (k_{el}+k_{off}+4 k_{on} L+k_{on} R)-k_{off} RL-k_{syn})+k_{on} (L.R (k_{off}+2 k_{on} L)+ }{ + k_{int} k_{off}+k_{off}^2+ k_{off} k_{on} (4 L+2 R-RL)\big) + k_{on}^2 L \big(L \left(-k_{el} k_{off}+k_{int} k_{off}+ \right.} \\
\dfrac{+ k_{off} RL (R-L)-k_{syn} L)\big)+(k_{deg}+k_{on} L) \big((k_{deg}+k_{on} L)^2+k_{on} L (k_{off} + }{ \left.  +k_{on} R (k_{off}+3 k_{on} L) +(k_{off}+k_{on} L)^2-k_{on} k_{syn}\right)+ } \\
\dfrac{+k_{on} R)\big) (-R (k_{deg}+k_{on} L)+k_{off} RL+k_{syn})\Big)}{+k_{off} RL (k_{off}+k_{on} (R-L))\big)\Big)} = 0,
\label{eq:TMDD_CSPR21_imp}
\end{multline} 
that is a very complicated implicit expression.~Although we can derive an analytic explicit solution for $R$ (but not for $L$), it is even more complicated than the above.~Hence, we solve the CSP$_R$(2) implicit expression in Eq.~\eqref{eq:TMDD_CSPR21_imp} numerically, using Newton's iterations, to get explicit solutions for $L$ or $R$.

Finally, we consider the case where the fast variable is assumed to be $RL$.~Again, instead of reordering the system in Eq.~\eqref{eq:TMDD}, we equivalently select the set of initial CSP vectors as:
\begin{equation}
\mathbf{A}_r(0,0) = \begin{bmatrix}
        0 \\[2pt]
        0 \\[2pt]
        1
    \end{bmatrix}, \quad     
    \mathbf{A}_s(0,0) = \begin{bmatrix}
		1 & 0 \\[2pt]
        0 & 1 \\[2pt]
        0 & 0 
    \end{bmatrix}, \quad
    \mathbf{B}^r(0,0) = \begin{bmatrix}
         0 &  0 & 1
    \end{bmatrix}, \quad
    \mathbf{B}^s(0,0) = \begin{bmatrix}
         1 & 0 & 0 \\[2pt]
         0 & 1 & 0
    \end{bmatrix}.
    \label{eq:ICbv_TMDD3}
\end{equation}
which is in agreement to the assumption of the 3rd variable in the state vector of Eq.~\eqref{eq:TMDD} being fast.~For the first CSP iteration, we employ Steps 2 and 3 of the Algorithm in Section~\ref{app:sbCSP}.~Using the new CSP vectors, the resulting from Eq.~\eqref{eq:VF3} SIM approximation reads:
\begin{multline}
f^1(1) = k_{on} L.R-k_{off} RL-k_{int} RL + \\ + \dfrac{k_{on} L (R (k_{deg}-k_{el}+k_{on} L-k_{on} R)-k_{syn})+k_{off} k_{on} RL (R-L)}{k_{int}+k_{off}} = 0,
\label{eq:TMDD_CSPRL11_imp}
\end{multline}
which is written in an implicit functional form.~To express the above CSP with one iteration expression in explicit form, we solve analytically for both $L$ and $R$.~The resulting analytic CSP$_{RL}$(1) approximation when solving for $L$ is:
\begin{multline}
\hat{L} = -\dfrac{1}{2 k_{on} R} \bigg( (k_{deg}+k_{el}+k_{int}+k_{off}) R-k_{off} RL+k_{on} R^2-k_{syn} -\Big(-R (k_{deg}+k_{el}+ \\ + k_{int}+k_{off}+k_{on} R)+k_{off} RL+k_{syn})^2+4 R.RL \left((k_{int}+k_{off})^2+k_{off} k_{on} R\right)\Big)^{1/2}\bigg),
\label{eq:TMDD_CSPRL11_expL}
\end{multline}
and when solving for $R$ is:
\begin{multline}
\hat{R} = -\dfrac{1}{2 k_{on} L} \bigg( (k_{deg}+k_{el} +k_{int} +k_{off}+k_{on} L) L-k_{off} RL -\Big(L (k_{deg}+k_{el}+k_{int}+k_{off}+k_{on} L)- \\
-k_{off} RL)^2+4 L \left(RL (k_{int}+k_{off})^2+k_{on} L (k_{off} RL+k_{syn})\right) \Big)^{1/2} \bigg).
\label{eq:TMDD_CSPRL11_expR}
\end{multline}
For the second CSP iteration, we obtained the CSP$_{RL}$(2) SIM approximation, by employing Step 4 of  the Algorithm in Section~\ref{app:sbCSP}, as:
\begin{multline}
f^1(2,1) = \dfrac{\big((k_{int}+k_{off})^2+k_{off} +k_{on} (L+R)\big) \Big(k_{on} (R (k_{deg}+k_{on} L)-k_{off} RL-k_{syn}) \big(L (k_{deg}+ k_{el}+ }{(k_{int}+k_{off}) \Big(k_{off} \big(k_{on} (R (k_{deg}+k_{el}+4 k_{on} L)+L (k_{deg}+k_{el})+ } \\
\dfrac{ +k_{int}+k_{off}+k_{on} (L+2 R) )-k_{off} RL\big)-k_{on} (k_{off} RL-L (k_{el}+k_{on} R)) \big(R (k_{deg}+k_{el}+k_{int}+k_{off} -k_{syn} + }{+k_{on} R^2-k_{syn}) + (k_{off}+k_{on} L)^2+2 k_{off} k_{on} (R-RL)\big) + } \\
\dfrac{+ k_{on} (R+2L))-k_{off} RL \big)+\left((k_{int}+ k_{off})^2+k_{off} k_{on} (L+R)\right) (k_{on} L. R-RL (k_{int}+k_{off}))\Big)}{+k_{int}^3+3 k_{int}^2 k_{off}+k_{int} k_{off} (3 k_{off}+2 k_{on} (L+R))\Big)}= 0.
\label{eq:TMDD_CSPRL21_imp}
\end{multline}
Although the above very complicated implicit expression can be analytically solved for both $L$ and $R$, the resulting explicit expressions are even more complicated.~Hence, we solve the CSP$_{RL}$(2) implicit expression in Eq.~\eqref{eq:TMDD_CSPRL21_imp} numerically for getting an explicit functional in terms of $L$ or $R$, via the utilization of Newton's iterations.

\subsection{The PEA and QSSA expressions}
Here, we employ the framework presented in Section~\ref{app:sbPEAQSSA} to derive all possible PEA and QSSA expressions for the approximation of the $M=1$-dim. SIM of the TMDD mechanism.~For the construction of the PEA expression, we assume that the first, binding reaction of the TMDD mechanism in Eq.~\eqref{eq:TMDD} is fast, since it is the only reversible reaction; thus, it can accommodate for \emph{partial equilibrium}.~For the construction of the QSSA expressions, we consider all three cases of $L$, $R$ or $RL$ variables being fast.~Note here that while the choice of fast reaction may change the QSSA reduced model in Eq.~\eqref{eq:VF4}, it is not changing the SIM approximation.~Hence, we will not consider any other fast reaction other than the first, for the derivation of SIM approximations by QSSA.

First, we consider the case where the first, binding reaction and $L$ variable are assumed fast.~Then, the TMDD system is written in the decomposed form of Eq.~\eqref{eq:genSR_dec}, as:
\begin{equation}
\dfrac{d}{dt} \begin{bmatrix} L~~ \\ R~~ \\ RL
    \end{bmatrix} = \begin{bmatrix}
    -1 \\ -1 \\ ~~~1
    \end{bmatrix}
    (k_{on}L.R-k_{off}RL)
    +
    \begin{bmatrix}
    -1 & ~~~0 & ~~~0 & ~~~0 \\
    ~~~0 & ~~~1 & -1 & ~~~0 \\
    ~~~0 & ~~~0 & ~~~0 & -1
    \end{bmatrix}
    \begin{bmatrix}
    k_{el} L~\quad \\
    k_{syn}~\quad \\
    k_{deg}R~~ \\
    k_{int}RL
    \end{bmatrix},
    \label{eq:TMDD_dec_1L}
\end{equation}
where: 
\begin{align}
\mathbf{z}^M & =L, & \mathbf{z}^{N-M} & = \begin{bmatrix}  R~~ \\ RL \end{bmatrix}, & 
\mathbf{R}^M & = k_{on}L.R-k_{off}RL, & \mathbf{R}^{K-M} & = \begin{bmatrix}   k_{el} L &    k_{syn} &    k_{deg}R &    k_{int}RL \end{bmatrix}^\top, \nonumber \\
    \mathbf{S}^M_M & = -1, & \mathbf{S}^{N-M}_M & = \begin{bmatrix} -1 \\ ~~~1 \end{bmatrix}, & \mathbf{S}^M_{K-M} & = \begin{bmatrix}
    -1 & 0 &  0 &  0 \end{bmatrix}, & \mathbf{S}^{N-M}_{K-M} &= \begin{bmatrix}
    0 & 1 & -1 & ~~~0 \\
    0 & 0 & ~~~0 & -1
    \end{bmatrix}.
    \label{eq:TMDD_dec_1Lm}
\end{align}
Given the above, we calculate the matrices required for the construction of the PEA vectors in Section~\ref{app:sbPEAQSSA}.~Then, following Eq.~\eqref{eq:PEA_SIMgen}, the resulting SIM approximation on the basis of PEA is:
\begin{equation}
k_{on} L.R-k_{off} RL-\dfrac{k_{on} L (k_{syn}-R (k_{deg}+k_{el}))+k_{int} k_{off} RL}{k_{off}+k_{on} (L+R)}=0,
\label{eq:TMDD_PEA_imp}
\end{equation}
which is in full agreement with Eq.~(85) in \citep{patsatzis2016asymptotic}.~The above implicit expression of the SIM can be analytically solved for both $L$ and $R$, resulting to the PEA explicit SIM approximations; solving for $L$, the expression reads:
\begin{multline}
\hat{L}=\dfrac{1}{2 k_{on} R} \Big(-R(k_{deg}+k_{el} +k_{off} +k_{on} R)+k_{off} RL+k_{syn}+ \big((-R (k_{deg}+k_{el}+k_{off}+k_{on} R)+\\
+k_{off} RL+k_{syn})^2+4 k_{off} R RL (k_{int}+k_{off}+k_{on} R)\big)^{1/2} \Big),
\label{eq:TMDD_PEA_expL}
\end{multline}
and solving for $R$, the explicit expression reads:
\begin{multline}
\hat{R} = -\dfrac{1}{2 k_{on} L} \big(L(k_{deg}+k_{el}+k_{off} +k_{on} L)-k_{off} RL-\big((L (k_{deg}+k_{el}+k_{off}+k_{on} L) -\\
-k_{off} RL)^2+4 L (k_{off} RL (k_{int}+k_{off}+k_{on} L)+k_{on} k_{syn} L)\big)^{1/2} \Big).
\label{eq:TMDD_PEA_expR}
\end{multline}
Further, considering the TMDD system in the form of Eq.~\eqref{eq:TMDD_dec_1L} accompanied by the expressions in Eq.~\eqref{eq:TMDD_dec_1Lm}, we derive the SIM approximation on the basis of QSSA when $L$ is the fast variable.~Following Eq.~\eqref{eq:QSSA_SIMgen}, the resulting QSSA$_L$ SIM approximation is:
\begin{equation}
k_{off} RL-L (k_{el}+k_{on} R) = 0,
\label{eq:TMDD_QSSAL_imp}
\end{equation}
which is in full agreement with Eq.~(86) in \citep{patsatzis2016asymptotic}.~The above implicit expression can be analytically solved for both $L$ and $R$.~Solving for $L$, the QSSA$_L$ explicit SIM approximation is:
\begin{equation}
\hat{L} = \dfrac{k_{off} RL}{k_{el}+k_{on} R},
\label{eq:TMDD_QSSAL_expL}
\end{equation}
and solving for $R$:
\begin{equation}
\hat{R} = \dfrac{k_{off} RL-k_{el} L}{k_{on} L}.
\label{eq:TMDD_QSSAL_expR}
\end{equation}

Next, we consider the case where the first, binding reaction and $R$ variable are assumed fast.~In this case, we reorder the system, so that the fast variable is first, to obtain  the decomposed form of the TMDD system in Eq.~\eqref{eq:genSR_dec}, as:
\begin{equation}
\dfrac{d}{dt} \begin{bmatrix} R~~ \\ L~~ \\ RL
    \end{bmatrix} = \begin{bmatrix}
    -1 \\ -1 \\ ~~~1
    \end{bmatrix}
    (k_{on}L.R-k_{off}RL)
    +
    \begin{bmatrix}
    ~~~0 & ~~~1 & -1 & ~~~0   \\
    -1 & ~~~0 & ~~~0 & ~~~0 \\
    ~~~0 & ~~~0 & ~~~0 & -1
    \end{bmatrix}
    \begin{bmatrix}
    k_{el} L~\quad \\
    k_{syn}~\quad \\
    k_{deg}R~~ \\
    k_{int}RL
    \end{bmatrix},
    \label{eq:TMDD_dec_1R}
\end{equation}
where, in this case:
\begin{align}
\mathbf{z}^M & =R, & \mathbf{z}^{N-M} & = \begin{bmatrix}  L~~ \\ RL \end{bmatrix}, & 
\mathbf{R}^M & = k_{on}L.R-k_{off}RL, & \mathbf{R}^{K-M} & = \begin{bmatrix}   k_{el} L &    k_{syn} &    k_{deg}R &    k_{int}RL \end{bmatrix}^\top, \nonumber \\
    \mathbf{S}^M_M & = -1, & \mathbf{S}^{N-M}_M & = \begin{bmatrix} -1 \\ ~~~1 \end{bmatrix}, & \mathbf{S}^M_{K-M} & = \begin{bmatrix}
    0 & 1 & -1 & 0  \end{bmatrix}, & \mathbf{S}^{N-M}_{K-M} &= \begin{bmatrix}
    -1 & 0 &  0 &  ~~~0 \\
    ~~~0 & 0 & 0 & -1
    \end{bmatrix}.
    \label{eq:TMDD_dec_1Rm}
\end{align}
Given the above decomposition, we constructed the SIM approximation on the basis of PEA, following Eq.~\eqref{eq:PEA_SIMgen}.~The resulting implicit expression (as well as the explicit functionals) is  equivalent with the PEA expression in Eq.~\eqref{eq:TMDD_PEA_imp}, where the variable $L$ was assumed fast.~This is due to the fact that for the construction of PEA approximations, the assumption on the fast variable is unnecessary \citep{goussis2012quasi,goussis2015model}.
However, this is not the case for the QSSA, which requires an assumption on the fast variables.~Thus, considering the TMDD system in the form of Eq.~\eqref{eq:TMDD_dec_1R} accompanied by the expressions in Eq.~\eqref{eq:TMDD_dec_1Rm}, we derive the SIM approximation on the basis of QSSA when $R$ is the fast variable.~Following Eq.~\eqref{eq:QSSA_SIMgen}, the resulting QSSA$_R$ SIM approximation is:
\begin{equation}
-R (k_{deg}+k_{on} L)+k_{off} RL+k_{syn}=0,
\label{eq:TMDD_QSSAR_imp}
\end{equation}
which is in full agreement with Eq.~(81) in \citep{patsatzis2016asymptotic}.~The above implicit expression can be analytically solved for both $L$ and $R$.~Solving for $L$, the QSSA$_R$ explicit SIM approximation is:
\begin{equation}
\hat{L} = \dfrac{-k_{deg} R+k_{off} RL+k_{syn}}{k_{on} R},
\label{eq:TMDD_QSSAR_expL}
\end{equation}
and solving for $R$:
\begin{equation}
\hat{R} = \dfrac{k_{off} RL+k_{syn}}{k_{deg}+k_{on} L}.
\label{eq:TMDD_QSSAR_expR}
\end{equation}

Finally, we consider the case where the first, binding reaction and $RL$ variable are assumed fast.~As in the previous case, we reorder the system, so that the fast variable is first, to obtain  the decomposed form of the TMDD system in Eq.~\eqref{eq:genSR_dec}, as:
\begin{equation}
\dfrac{d}{dt} \begin{bmatrix} RL \\ L~~ \\ R~~
    \end{bmatrix} = \begin{bmatrix}
    ~~~1 \\ -1 \\ -1
    \end{bmatrix}
    (k_{on}L.R-k_{off}RL)
    +
    \begin{bmatrix}
     ~~~0 & ~~~0 & ~~~0 & -1 \\
    -1 & ~~~0 & ~~~0 & ~~~0 \\
    ~~~0 & ~~~1 & -1 & ~~~0
    \end{bmatrix}
    \begin{bmatrix}
    k_{el} L~\quad \\
    k_{syn}~\quad \\
    k_{deg}R~~ \\
    k_{int}RL
    \end{bmatrix},
    \label{eq:TMDD_dec_1RL}
\end{equation}
where, in this case:
\begin{align}
\mathbf{z}^M & =RL, & \mathbf{z}^{N-M} & = \begin{bmatrix}  L \\ R \end{bmatrix}, & 
\mathbf{R}^M & = k_{on}L.R-k_{off}RL, & \mathbf{R}^{K-M} & = \begin{bmatrix}   k_{el} L &    k_{syn} &    k_{deg}R &    k_{int}RL \end{bmatrix}^\top, \nonumber \\
    \mathbf{S}^M_M & = 1, & \mathbf{S}^{N-M}_M & = \begin{bmatrix} -1 \\ -1 \end{bmatrix}, & \mathbf{S}^M_{K-M} & = \begin{bmatrix}
      0 & 0 & 0 & -1\end{bmatrix}, & \mathbf{S}^{N-M}_{K-M} &= \begin{bmatrix}
    -1 & 0 &  ~~~0 & 0 \\
    ~~~0 & 1 & -1 & 0
    \end{bmatrix}.
    \label{eq:TMDD_dec_1RLm}
\end{align}
Given the above decomposition, the SIM approximation on the basis of PEA is, as expected by  \citep{goussis2012quasi,goussis2015model}, resulting to an implicit expression that is  equivalent with the PEA expression in Eq.~\eqref{eq:TMDD_PEA_imp}; the same holds for the explicit functionals in Eqs.~\eqref{eq:TMDD_PEA_expL} and \eqref{eq:TMDD_PEA_expR}.~For the derivation of the SIM approximations on the basis of QSSA when $RL$ is the fast variable, we consider the decomposed form of the TMDD system in Eq.~\eqref{eq:TMDD_dec_1RL} accompanied by the expressions in Eq.~\eqref{eq:TMDD_dec_1RLm}.~Then, following Eq.~\eqref{eq:QSSA_SIMgen}, the resulting QSSA$_{RL}$ SIM approximation is:
\begin{equation}
k_{on} L R-RL (k_{int}+k_{off})=0.
\label{eq:TMDD_QSSARL_imp}
\end{equation}
The above implicit expression can be analytically solved for both $L$ and $R$, resulting to the explicit QSSA$_{RL}$ SIM approximations; solving for $L$, yields:
\begin{equation}
\hat{L} = \dfrac{RL (k_{int}+k_{off})}{k_{on} R},
\label{eq:TMDD_QSSARL_expL}
\end{equation}
and solving for $R$ yields:
\begin{equation}
\hat{R} = \dfrac{RL (k_{int}+k_{off})}{k_{on} L}.
\label{eq:TMDD_QSSARL_expR}
\end{equation}

\renewcommand{\theequation}{D.\arabic{equation}}
\renewcommand{\thefigure}{D.\arabic{figure}}
\setcounter{equation}{0}
\setcounter{figure}{0}
\section{The SIM approximations of the fCSI mechanism}
\label{app:Inh_SIMs}
Here, we present the SIM approximations underlying the slow evolution of the fCSI mechanism in Section~\ref{sub:CompInhdes}.~We realize both the original system in Eq.~\eqref{eq:Inh_st} and the intuition-based transformed system in Eq.~\eqref{eq:Inh_tr}, and derive the approximations on the basis of CSP, PEA and QSSA methods.~Considering all the possible cases for the construction of an $M=2$-dim. SIM approximation would require, for each system, the derivation of 12 CSP, 1 PEA and 6 QSSA expressions; 6 different pairs of fast variables and 1 pair of fast reactions can be considered.~Thus, here we consider only the case where the fast variables are the complexes $c_1$ and $c_2$ (which is in agreement with   \citep{segel1988validity,schnell2000time,rubinow1970time} and the findings of the proposed PINN scheme) and solve the resulting SIM expressions for $c_1$ and $c_2$ variables.

\subsection{The original system}
\label{app:InhSt_SIMs}
First, we consider the original system in Eq.~\eqref{eq:Inh_st} for which we construct SIM approximations on the basis of CSP with one and two iterations, PEA and QSSA.~A summary of the four, in total, SIM approximations is presented in Table~\ref{tb:InhSt_convSIMs}.

\begin{table}[!h]
\centering
\footnotesize
\begin{tabular}{l c c c c}
\toprule
SIM approximation & QSSA$_{c1c2}$ & PEA$_{13}$ & CSP$_{c1c2}$(1) & CSP$_{c1c2}$(2)  \\
\midrule
Assumptions & $c_1$, $c_2$ fast & 1st, 3rd reactions fast  & $c_1$, $c_2$ fast & $c_1$, $c_2$ fast \\
\midrule
Solved for	& $c_1$ and $c_2$  & $c_1$ and $c_2$ & $c_1$ and $c_2$ & $c_1$ and $c_2$ \\
Explicit/Implicit  & E & I & I & I  \\
\midrule
Equation & (\ref{eq:InhSt_QSSAc1c2_exp1}, \ref{eq:InhSt_QSSAc1c2_exp2})  & (\ref{eq:InhSt_PEA13_imp1}, \ref{eq:InhSt_PEA13_imp2}) & (\ref{eq:InhSt_CSP11_imp1}, \ref{eq:InhSt_CSP11_imp2}) & numerical$^\dagger$ \\
Requires Newton	& N & Y & Y & Y \\
\bottomrule
\end{tabular}
\caption{SIM approximations for the original system of the fCSI mechanism in Eq.~\eqref{eq:Inh_st}, constructed on the basis of QSSA, PEA and CSP methods.~For each approximation, we enlist (i) the assumptions made for constructing it, (ii) the functional form (E/I for explicit/implicit) whem solved for $c_1$ and $c_2$ variables, and (iii) the respective equations and whether Newton iterations are required  (Y/N for yes/no) to numerically solve the SIM approximation for $c_1$ and $c_2$.\\
$^\dagger$ the CSP$_{c1c2}$(2) approximation is computed numerically for every given point of the $[s_1,c_1,s_2,c_2]$ space.}
\label{tb:InhSt_convSIMs}
\end{table} 

\subsubsection{The CSP expressions with one and two iterations}
For the employment of the CSP methodology, we follow the procedure described in Section~\ref{app:sbCSP} by setting  $M=2$ and assuming that the fast variables are $c_1$ and $c_2$.~When the latter assumption is not correct, two iterations of CSP are required for the accuracy of the SIM approximation \cite{goussis2012quasi,valorani2005higher}.~To initialize the CSP iterative procedure, instead of reordering the system and taking an initial set of CSP vectors as in Eq.~\eqref{eq:ICbv}, we just set the initial set as:  
\begin{equation}
    \mathbf{A}_r(0,0) = \begin{bmatrix}
        0 & 0 \\
        1 & 0 \\
        0 & 0 \\
        0 & 1 
    \end{bmatrix}, \quad     
    \mathbf{A}_s(0,0) = \begin{bmatrix}
        1 & 0 \\
        0 & 0 \\
        0 & 1 \\
        0 & 0 
    \end{bmatrix}, \quad
    \mathbf{B}^r(0,0) = \begin{bmatrix}
        0 & 1 & 0 & 0 \\ 0 & 0 & 0 & 1
    \end{bmatrix}, \quad
    \mathbf{B}^s(0,0) = \begin{bmatrix}
        1 & 0 & 0 & 0 \\ 0 & 0 & 1 & 0
    \end{bmatrix},
    \label{eq:ICbv_InhSt}
\end{equation}
which is in agreement with the assumption of the 2nd and 4th variables, $c_1$ and $c_2$ respectively, being fast.~With the above set, we performed Steps 2 and 3 of the Algorithm in Section~\ref{app:sbCSP} to get the CSP vectors after one iteration.~Then, according to Eq.~\eqref{eq:VF3}, the resulting CSP$_{c1c2}$(1) SIM approximation is $M=2$-dimensional $\mathbf{f}^r(1)=[f^1(1), f^2(1)]^\top$, consisting of the components given by:
\begin{multline}
f^1(1) = k_{1f} s_1 (e_0-c_2)-c_1 (k_{1b}+k_{1f} s_1+k_2)+ k_{1f} (c_1+c_2-e_0) \cdot \\
\cdot \dfrac{(k_{3b}+k_{3f} s_2+k_4) (k_{1f} s_1 (c_1+c_2-e_0)+c_1 k_{1b})-k_{3f} s_1 (k_{3f} s_2 (c_1+c_2-e_0)+c_2 k_{3b})}{(k_{3b}+k_4) (k_{1b}+k_{1f} s_1+k_2)+k_{3f} s_2 (k_{1b}+k_2)}=0,
\label{eq:InhSt_CSP11_imp1}
\end{multline}
and
\begin{multline}
f^2(1) = k_{3f} s_2 (e_0-c_1)-c_2 (k_{3b}+k_{3f} s_2+k_4)+ k_{3f} (c_1+c_2-e_0) \cdot \\
\cdot \dfrac{(k_{1b}+k_{1f} s_1+k_2) (k_{3f} s_2 (c_1+c_2-e_0)+c_2 k_{3b})-k_{1f} s_2 (k_{1f} s_1 (c_1+c_2-e_0)+c_1 k_{1b})}{(k_{3b}+k_4) (k_{1b}+k_{1f} s_1+k_2)+k_{3f} s_2 (k_{1b}+k_2)}=0.
\label{eq:InhSt_CSP11_imp2}
\end{multline}
The above system of implicit expressions cannot be solved analytically for $c_1$ and $c_2$ variables.~Thus, in order to get explicit functionals of the CSP$_{c1c2}$(1) implicit SIM approximation, we numerically solve Eqs.~\eqref{eq:InhSt_CSP11_imp1} and \eqref{eq:InhSt_CSP11_imp2} using Newton's iterations (requiring the derivatives w.r.t. $c_1$ and $c_2$).

For performing the second CSP iteration, the implementation of Step 4 of the Algorithm in Section~\ref{app:sbCSP} is required.~However, the analytic computations lead to enormous expressions for $\mathbf{f}^r(2)=[f^1(2),f^2(2)]^\top$ of the SIM approximation.~Hence, we compute the CSP$_{c1c2}$(2) SIM approximation $\mathbf{f}^r(2)$ numerically for every given point of the $[s_1,c_1,s_2,c_2]$ space.~Naturally, the numerical SIM approximation is an implicit functional $\mathbf{f}^r(2)(s_1,c_1,s_2,c_2)=\mathbf{0}$ and thus, we solve it numerically using Newton iterations to get explicit estimations for $c_1$ and $c_2$ on the SIM.~Note that in this case, since $\mathbf{f}^r(2)$ is calculated numerically, its derivatives w.r.t. $c_1$ and $c_2$ are  also computed numerically via finite differences.

\subsubsection{The PEA and QSSA expressions}
For the employment of the PEA and QSSA methodologies, we follow the framework presented in Section~\ref{app:sbPEAQSSA}.~For the construction of the PEA expression, we assume that the first and third binding reactions of the fCSI mechanism in Eq.~\eqref{eq:Inh_st} are fast; these are the only reversible reactions that can accommodate for \emph{partial equilibrium}.~For the construction of the QSSA expression, we only consider the case of $c_1$ and $c_2$ being fast, as discussed.~With the above assumptions on the fast reactions and variables, we reorder the system in Eq.~\eqref{eq:Inh_st}, so that the fast reactions and variables are first, thus casting the system to the decomposed form of Eq.~\eqref{eq:genSR_dec}, as:
\begin{equation}
\dfrac{d}{dt} \begin{bmatrix} c_1 \\ c_2 \\ s_1 \\ s_2
    \end{bmatrix} = \begin{bmatrix}
    ~~~1 & ~~~0 \\ ~~~0 & ~~~1 \\ -1 & ~~~0 \\ ~~~0 & -1
    \end{bmatrix}
    \begin{bmatrix}
    k_{1f} (e_0-c_1-c_2) s_1 -k_{1b} c_1\\
    k_{3f} (e_0-c_1-c_2) s_2 -k_{3b} c_2\\
    \end{bmatrix}
    +
    \begin{bmatrix}
    -1 & ~~~0  \\
    ~~~0 & -1 \\
    ~~~0 & ~~~0 \\
    ~~~0 & ~~~0
    \end{bmatrix}
    \begin{bmatrix}
    k_2 c_1 \\
    k_4 c_2 
    \end{bmatrix},
    \label{eq:CInhOr_dec_13c1c2}
\end{equation}
where: 
\begin{align}
\mathbf{z}^M, & = \begin{bmatrix} c_1 \\ c_2 \end{bmatrix} & \mathbf{z}^{N-M} & = \begin{bmatrix}  s_1 \\ s_2 \end{bmatrix}, & 
\mathbf{R}^M & = \begin{bmatrix}
    k_{1f} (e_0-c_1-c_2) s_1 -k_{1b} c_1\\
    k_{3f} (e_0-c_1-c_2) s_2 -k_{3b} c_2\\
    \end{bmatrix}, & \mathbf{R}^{K-M} & = \begin{bmatrix}
    k_2 c_1 \\
    k_4 c_2 
    \end{bmatrix}, \nonumber \\
    \mathbf{S}^M_M & = \begin{bmatrix} 1 & 0 \\ 0 & 1 \end{bmatrix}, & \mathbf{S}^{N-M}_M & = \begin{bmatrix} -1 & ~~~0 \\ ~~~0 & -1 \end{bmatrix}, & \mathbf{S}^M_{K-M} & =  \begin{bmatrix} -1 & ~~~0  \\  ~~~0 & -1 \end{bmatrix}, & \mathbf{S}^{N-M}_{K-M} &= \begin{bmatrix} 0 & 0 \\ 0 & 0 \end{bmatrix}.
    \label{eq:CInhOr_dec_13c1c2m}
\end{align}
Given the above, we calculate the matrices required for the construction of the PEA vectors in Section~\ref{app:sbPEAQSSA}.~Then, following Eq.~\eqref{eq:PEA_SIMgen}, the resulting PEA$_{13}$ SIM approximation is composed by the following two components:
\begin{multline}
k_{1f} s_1 (e_0-c_2)-c_1 (k_{1b}+k_{1f} s_1+k_2) + \\
+ \dfrac{k_{1f} (c_1+c_2-e_0) \left(c_1^2 k_2 k_{3f}-c_1 k_2 (k_{3f} (-c_2+e_0+s_2)+k_{3b})+c_2 k_{3f} k_4 s_1\right)}{(k_{3b}-k_{3f} (c_1+c_2-e_0)) (k_{1f} (-c_1-c_2+e_0+s_1)+k_{1b})+k_{3f} s_2 (k_{1b}-k_{1f} (c_1+c_2-e_0))} = 0,
\label{eq:InhSt_PEA13_imp1}
\end{multline}
\begin{multline}
k_{3f} s_2 (e_0-c_1)-c_2 (k_{3b}+k_{3f} s_2+k_4) + \\
+ \dfrac{k_{3f} (c_1+c_2-e_0) (c_2 k_4 (k_{1f} (c_1+c_2-e_0-s_1)-k_{1b})+c_1 k_{1f} k_2 s_2)}{(k_{3b}-k_{3f} (c_1+c_2-e_0)) (k_{1f} (-c_1-c_2+e_0+s_1)+k_{1b})+k_{3f} s_2 (k_{1b}-k_{1f} (c_1+c_2-e_0))} = 0.
\label{eq:InhSt_PEA13_imp2}
\end{multline}
The above system of implicit expressions cannot be solved analytically for $c_1$ and $c_2$ variables.~Thus, to get the explicit functionals of the PEA$_{13}$ SIM approximation w.r.t. $c_1$ and $c_2$, we solve Eqs.~\eqref{eq:InhSt_PEA13_imp1} and \eqref{eq:InhSt_PEA13_imp2} numerically using Newton's iterations (requiring the derivatives w.r.t. $c_1$ and $c_2$).

Finally, considering the fCSI original system in the decomposed form of Eq.~\eqref{eq:CInhOr_dec_13c1c2} accompanied by the expressions in Eq.~\eqref{eq:CInhOr_dec_13c1c2m}, we derive the SIM approximation on the basis of QSSA when $c_1$ and $c_2$ are the fast variables.~Following Eq.~\eqref{eq:QSSA_SIMgen}, the resulting SIM approximation consists of the following two components:
\begin{equation}
k_{1f} s_1 (e_0-c_2)-c_1 (k_{1b}+k_{1f} s_1+k_2) = 0, \qquad
k_{3f} s_2 (e_0-c_1)-c_2 (k_{3b}+k_{3f} s_2+k_4) = 0.
\label{eq:InhSt_QSSAc1c2_imp}
\end{equation}
In this case, the above system of implicit expressions can be solved analytically for $c_1$ and $c_2$, thus resulting to the explicit QSSA$_{c1c2}$ SIM approximation:
\begin{align}
\hat{c}_1 & =\dfrac{e_0 k_{1f} s_1 (k_{3b}+k_4)}{(k_{3b}+k_4) (k_{1b}+k_{1f} s_1+k_2)+k_{3f} s_2 (k_{1b}+k_2)}, 
\label{eq:InhSt_QSSAc1c2_exp1} \\
\hat{c}_2 & =\dfrac{e_0 k_{3f} s_2 (k_{1b}+k_2)}{(k_{3b}+k_4) (k_{1b}+k_{1f} s_1+k_2)+k_{3f} s_2 (k_{1b}+k_2)}.
\label{eq:InhSt_QSSAc1c2_exp2}
\end{align}

\subsection{Intuition-based linearly transformed system}
\label{app:InhTr_SIMs}
Next, consider the transformed system of the fCSI mechanism in Eq.~\eqref{eq:Inh_tr} for which we construct SIM approximations on the basis of CSP with one and two iterations, PEA and QSSA.~A summary of the four, in total, SIM approximation is provided in Table~\ref{tb:InhTr_convSIMs}.

\begin{table}[!h]
\centering
\footnotesize
\begin{tabular}{l c c c c}
\toprule
SIM approximation & QSSA$_{c1c2}$ & PEA$_{13}$ & CSP$_{c1c2}$(1) & CSP$_{c1c2}$(2)  \\
\midrule
Assumptions & $c_1$, $c_2$ fast & 1st, 3rd reactions fast  & $c_1$, $c_2$ fast & $c_1$, $c_2$ fast \\
\midrule
Solved for	& $c_1$ and $c_2$  & $c_1$ and $c_2$ & $c_1$ and $c_2$ & $c_1$ and $c_2$ \\
Explicit/Implicit  & E & I & I & I  \\
\midrule
Equation & (\ref{eq:InhTr_QSSAc1c2_imp1}, \ref{eq:InhTr_QSSAc1c2_imp2})  & (\ref{eq:InhTr_PEA13_imp1}, \ref{eq:InhTr_PEA13_imp2}) & (\ref{eq:InhTr_CSP11_imp1}, \ref{eq:InhTr_CSP11_imp2}) & numerical$^\dagger$ \\
Requires Newton	& Y & Y & Y & Y \\
\bottomrule
\end{tabular}
\caption{SIM approximations for the original system of the fCSI mechanism in Eq.~\eqref{eq:Inh_tr}, constructed on the basis of QSSA, PEA and CSP methods.~For each approximation, we enlist (i) the assumptions made for constructing it, (ii) the functional form (E/I for explicit/implicit) when solved for $c_1$ and $c_2$ variables, and (iii) the respective equations and whether Newton iterations are required  (Y/N for yes/no) to numerically solve the SIM approximation for $c_1$ and $c_2$.\\
$^\dagger$ the CSP$_{c1c2}$(2) approximation is computed numerically for every given point of the $[\bar{s}_1,c_1,\bar{s}_2,c_2]$ space.}
\label{tb:InhTr_convSIMs}
\end{table} 

\subsubsection{The CSP expressions with one and two iterations}
First, we employ the CSP methodology described in  Section~\ref{app:sbCSP} to construct the CSP, with one and two iterations, SIM approximations.~Since $M=2$ and the fast variables of the transformed system remain $c_1$ and $c_2$, we again initialize the CSP iterative procedure with the set of CSP vectors in Eq.~\eqref{eq:ICbv_InhSt}.~Following this initialization, we obtain the CSP vectors after one iteration, following the Steps 2 and 3 of the Algorithm in Section~\ref{app:sbCSP}.~Then, according to Eq.~\eqref{eq:VF3}, the resulting CSP$_{c1c2}$(1) SIM approximation is $M=2$-dim. $\mathbf{f}^r(1)=[f^1(1), f^2(1)]^\top$, consists of the components given by:
\begin{multline}
f^1(1) = k_{1f} (c_1-\bar{s}_1) (c_1+c_2-e_0)-c_1 k_{1b}-c_1 k_2 + \dfrac{k_{1f} (e_0-c_1-c_2)\big(c_1^2 k_2 k_{3f}-}{k_{1f} k_{3f} (c_1-\bar{s}_1) (c_2-\bar{s}_2)-}\\
\dfrac{ -c_1 (k_2 (k_{3f} (e_0+\bar{s}_2-2 c_2)+k_{3b}+k_4)+c_2 k_{3f} k_4)+c_2 k_{3f} k_4 \bar{s}_1\big)}{-(k_{1f} (e_0+\bar{s}_1-2 c_1-c_2)+k_{1b}+k_2) (k_{3f} (e_0+\bar{s}_2-c_1-2 c_2)+k_{3b}+k_4)} = 0,
\label{eq:InhTr_CSP11_imp1}
\end{multline}
and
\begin{multline}
f^2(1) = k_{3f} (c_2-\bar{s}_2) (c_1+c_2-e_0)-c_2 k_{3b}-c_2 k_4 + \dfrac{k_{3f} (e_0-c_1-c_2) (c_2 k_4 (2 c_1 k_{1f}+}{k_{1f} k_{3f} (c_1-\bar{s}_1) (c_2-\bar{s}_2)-}\\
\dfrac{+c_2 k_{1f}-k_{1f} (e_0+\bar{s}_1)-k_{1b}-k_2)+c_1 k_{1f} k_2 (\bar{s}_2-c_2))}{-(k_{1f} (e_0+\bar{s}_1-2 c_1-c_2)+k_{1b}+k_2) (k_{3f} (e_0+\bar{s}_2-c_1-2 c_2)+k_{3b}+k_4)} = 0.
\label{eq:InhTr_CSP11_imp2}
\end{multline}
The above system of implicit expressions cannot be solved analytically for $c_1$ and $c_2$.~To get explicit functionals of the CSP$_{c1c2}$(1) SIM approximation, we numerically solve Eqs.~\eqref{eq:InhTr_CSP11_imp1} and \eqref{eq:InhTr_CSP11_imp2} w.r.t. $c_1$ and $c_2$, using Newton's iterations.

For the second CSP iteration, the implementation of Step 4 of the Algorithm in Section~\ref{app:sbCSP} results to to enormous analytic expressions for $\mathbf{f}^r(2)=[f^1(2),f^2(2)]^\top$ of the SIM approximation.~Hence, we compute the CSP$_{c1c2}$(2) SIM approximation $\mathbf{f}^r(2)$ numerically for every given point of the $[\bar{s}_1,c_1,\bar{s}_2,c_2]$ space.~However, this numerical SIM approximation is also an implicit functional $\mathbf{f}^r(2) (\bar{s}_1,c_1,\bar{s}_2,c_2)=\mathbf{0}$.~Therefore, to get explicit functionals of $c_1$ and $c_2$ we perform Newton iterations, which in this case require numerical computation (finite differences) of the derivatives of $\mathbf{f}^r(2)$ w.r.t. $c_1$ and $c_2$.

\subsubsection{The PEA and QSSA expressions}
For the derivation of the SIM approximations on the basis of the PEA and QSSA methodologies, we follow the framework presented in Section~\ref{app:sbPEAQSSA}.~For the PEA expression, we assume that the first and third binding reactions of the system in Eq.~\eqref{eq:Inh_tr} are fast, since only reversible reactions allow for realizing \emph{partial equilibrium}.~For the QSSA expression, we similarly consider only the case of $c_1$ and $c_2$ being fast, as discussed.~As in the case of the original system in Section~\ref{app:InhSt_SIMs}, we reorder the the system in Eq.~\eqref{eq:Inh_tr} so that the assumed fast reactions and variables are presented first.~Then, the system is cast to the decomposed form of Eq.~\eqref{eq:genSR_dec}, as:
\begin{equation}
\dfrac{d}{dt} \begin{bmatrix} c_1 \\ c_2 \\ \bar{s}_1 \\ \bar{s}_2  \end{bmatrix} = \begin{bmatrix}
    1 & 0 \\ 0 & 1 \\ 0 & 0 \\ 0 & 0
    \end{bmatrix}
    \begin{bmatrix}
    k_{1f} (e_0-c_1-c_2) (\bar{s}_1-c_1) -k_{1b} c_1\\
    k_{3f} (e_0-c_1-c_2) (\bar{s}_2-c_2) -k_{3b} c_2
    \end{bmatrix}
    +
    \begin{bmatrix}
    -1 & ~~~0  \\
    ~~~0 & -1 \\
    -1 & ~~~0  \\
    ~~~0 & -1
    \end{bmatrix}
    \begin{bmatrix}
    k_2 c_1 \\
    k_4 c_2 
    \end{bmatrix},
    \label{eq:CInhTr_dec_13c1c2}
\end{equation}
where: 
\begin{align}
\mathbf{z}^M & = \begin{bmatrix} c_1 \\ c_2 \end{bmatrix}, & \mathbf{z}^{N-M} & = \begin{bmatrix}  \bar{s}_1 \\ \bar{s}_2 \end{bmatrix}, & 
\mathbf{R}^M & = \begin{bmatrix}
    k_{1f} (e_0-c_1-c_2) (\bar{s}_1-c_1) -k_{1b} c_1\\
    k_{3f} (e_0-c_1-c_2) (\bar{s}_2-c_2) -k_{3b} c_2 
    \end{bmatrix}, & \mathbf{R}^{K-M} & = \begin{bmatrix}
    k_2 c_1 \\
    k_4 c_2 
    \end{bmatrix}, \nonumber \\
    \mathbf{S}^M_M & = \begin{bmatrix} 1 & 0 \\ 0 & 1 \end{bmatrix}, & \mathbf{S}^{N-M}_M & = \begin{bmatrix} 0 & 0 \\ 0 & 0 \end{bmatrix}, & \mathbf{S}^M_{K-M} & =  \begin{bmatrix} -1 & ~~~0  \\  ~~~0 & -1 \end{bmatrix}, & \mathbf{S}^{N-M}_{K-M} &=  \begin{bmatrix} -1 & ~~~0  \\  ~~~0 & -1 \end{bmatrix}.
    \label{eq:CInhTr_dec_13c1c2m}
\end{align}
Given the above, we calculate the matrices in Section~\ref{app:sbPEAQSSA} in order to construct the PEA vectors.~Then, following Eq.~\eqref{eq:PEA_SIMgen}, the resulting PEA$_{13}$ SIM approximation consists of the following two components:
\begin{multline}
k_{1f} (c_1-\bar{s}_1) (c_1+c_2-e_0)-c_1 k_{1b}-c_1 k_2 + \dfrac{k_{1f} (c_1+c_2-e_0) \big(c_1^2 k_2 k_{3f}-}{k_{1f} \big(2 c_1^2 k_{3f}-c_1 (k_{3f} (\bar{s}_1+\bar{s}_2-4 c_2+3 e_0)+2 k_{3b})+}\\ \dfrac{-c_1 (k_2 (k_{3f} (-2 c_2+e_0+\bar{s}_2)+k_{3b})+c_2 k_{3f} k_4)+}{+2 c_2^2 k_{3f}-c_2 (k_{3f} (3 e_0+\bar{s}_1+\bar{s}_2)+k_{3b})+(e_0+\bar{s}_1) (e_0 k_{3f}+k_{3b})+e_0 k_{3f} \bar{s}_2\big)+} \\ \dfrac{+c_2 k_{3f} k_4 \bar{s}_1\big)}{+k_{1b} (k_{3f} (e_0+\bar{s}_2-c_1-2 c_2)+k_{3b})} =0,
\label{eq:InhTr_PEA13_imp1}
\end{multline}
\begin{multline}
k_{3f} (c_2-\bar{s}_2) (c_1+c_2-e_0)-c_2 k_{3b}-c_2 k_4 - \dfrac{k_{3f} (c_1+c_2-e_0) (c_1 k_{1f} (c_2 (k_2-}{k_{1f} \big(2 c_1^2 k_{3f}-c_1 (k_{3f} (\bar{s}_1+\bar{s}_2-4 c_2+3 e_0)+2 k_{3b})+} \\ \dfrac{-2 k_4)-k_2 \bar{s}_2)+c_2 k_4 (k_{1f} (e_0+\bar{s}_1-c_2)}{+2 c_2^2 k_{3f}-c_2 (k_{3f} (3 e_0+\bar{s}_1+\bar{s}_2)+k_{3b})+(e_0+\bar{s}_1) (e_0 k_{3f}+k_{3b})+e_0 k_{3f} \bar{s}_2\big)+)} \\ \dfrac{+k_{1b}))}{+k_{1b} (k_{3f} (e_0+\bar{s}_2-c_1-2 c_2)+k_{3b}} =0.
\label{eq:InhTr_PEA13_imp2}
\end{multline}
The system of the above two implicit expressions cannot be solved analytically for $c_1$ and $c_2$ variables.~Thus, to get the explicit functionals of the PEA$_{13}$ SIM approximation w.r.t. $c_1$ and $c_2$, we solve Eqs.~\eqref{eq:InhTr_PEA13_imp1} and \eqref{eq:InhTr_PEA13_imp2} numerically using Newton's iterations (requiring the derivatives w.r.t. $c_1$ and $c_2$).

Finally, considering the transformed system of the fCSI mechanism in the decomposed form of Eq.~\eqref{eq:CInhTr_dec_13c1c2} accompanied by the expressions in Eq.~\eqref{eq:CInhTr_dec_13c1c2m}, we derive the SIM approximation on the basis of QSSA when $c_1$ and $c_2$ are the fast variables.~Following Eq.~\eqref{eq:QSSA_SIMgen}, the resulting QSSA$_{c1c2}$ SIM approximation consists of the following two components:
\begin{align}
k_{1f} (c_1-\bar{s}_1) (c_1+c_2-e_0)-c_1 k_{1b}-c_1 k_2 & = 0, 
\label{eq:InhTr_QSSAc1c2_imp1} \\
k_{3f} (c_2-\bar{s}_2) (c_1+c_2-e_0)-c_2 k_{3b}-c_2 k_4 & = 0.
\label{eq:InhTr_QSSAc1c2_imp2}
\end{align}
In contrast to the original system, the  QSSA$_{c1c2}$ SIM approximation in Eq.~\eqref{eq:InhTr_QSSAc1c2_imp1} and \eqref{eq:InhTr_QSSAc1c2_imp2} for the transformed system, cannot be solved for analytically for $c_1$ and $c_2$.~Therefore, to obtain the explicit functionals, we once again employ Newton iterations.


%
%
%
%

\end{document}


\maketitle

\section{Training and test data sets for the MM, TMDD and fCSI systems}

\begin{figure}[!h]
    \centering
    \subfigure[MM1 case]{
    \includegraphics[width=0.32\textwidth]{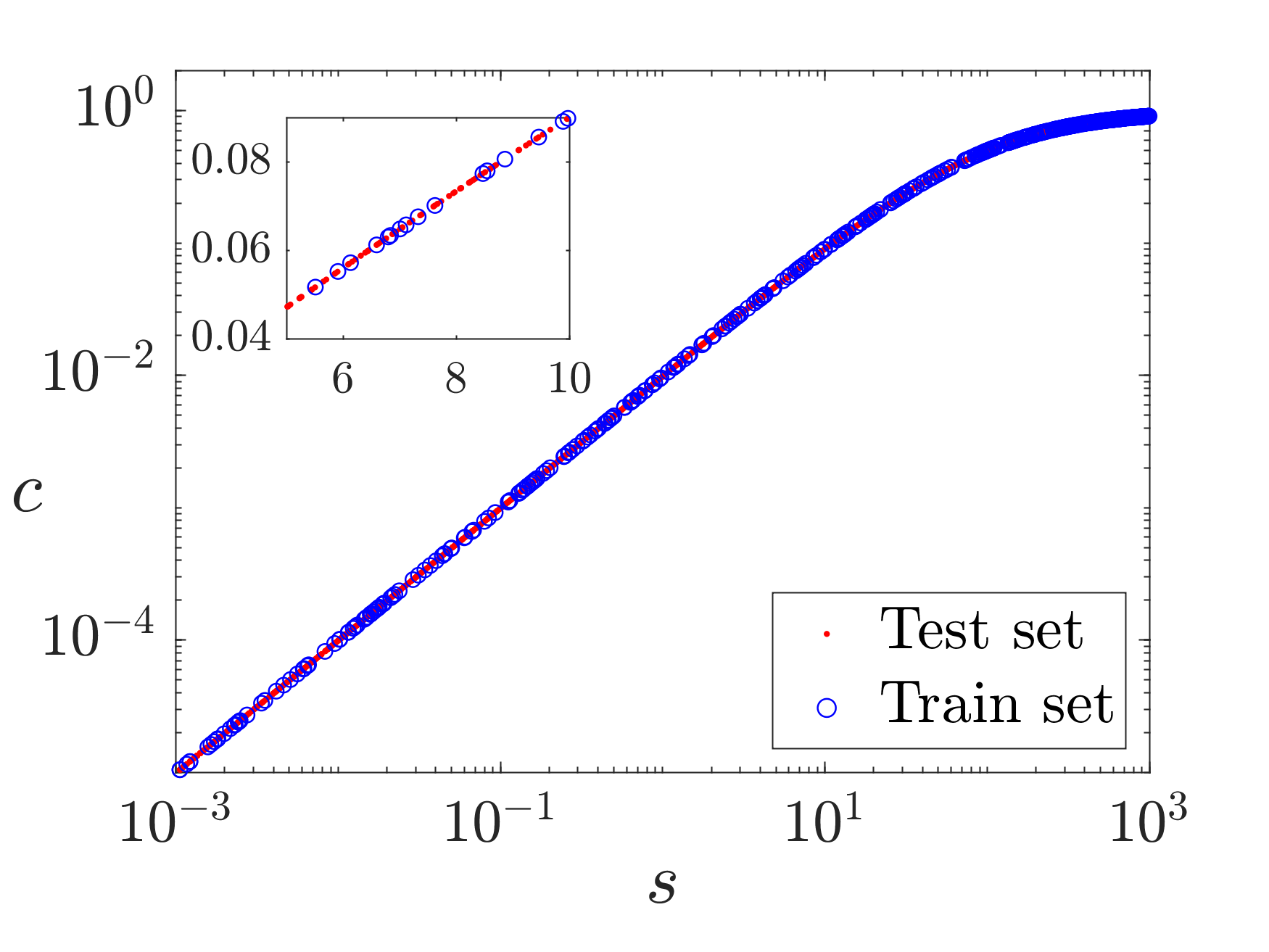}}
    \subfigure[MM2 case]{
    \includegraphics[width=0.32\textwidth]{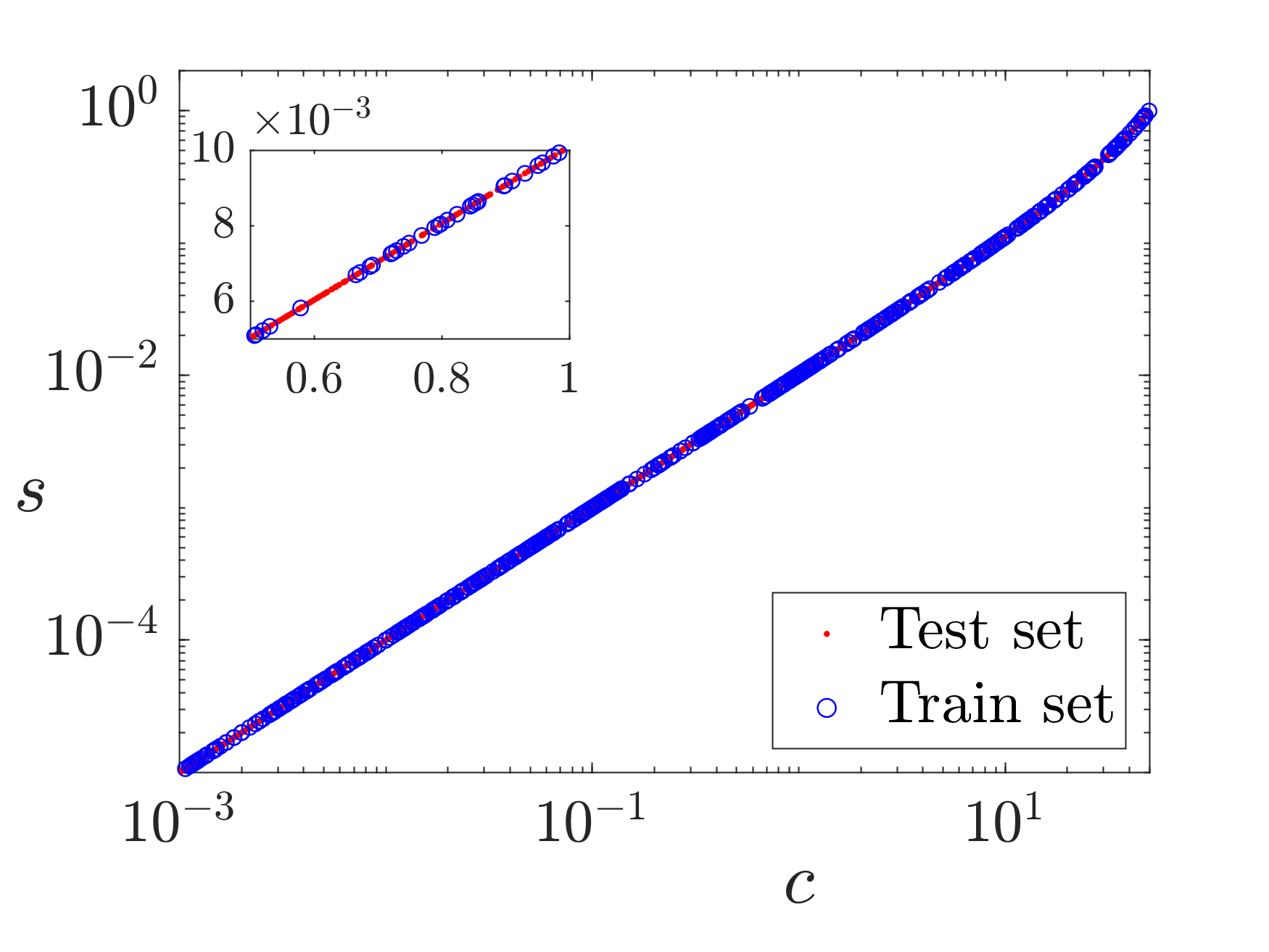}}
    \subfigure[MM3 case]{
    \includegraphics[width=0.32\textwidth]{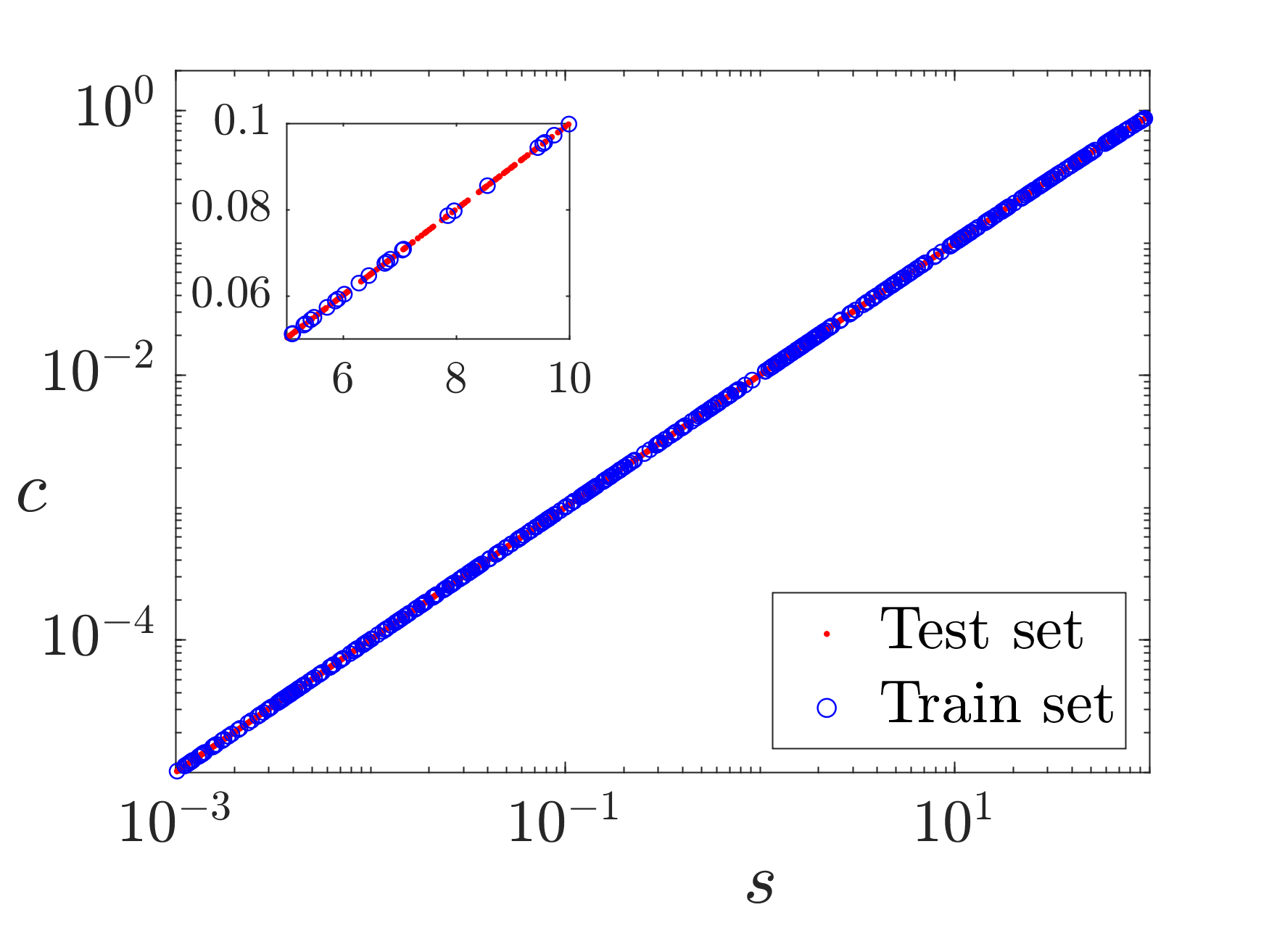}}
    \caption{MM system in Eq.~(32) for the three MM1, MM2 and MM3 cases.~Data sets used for training and validation of the PINN scheme (blue circles) for the three MM cases.~The test data sets (red dots) were used to assess the numerical approximation accuracy of the SIMs.}
    \label{SF:MM_Data}
\end{figure}

\begin{figure}[!h]
    \centering
    \subfigure[SIM $\mathcal{M}_1$]{
    \includegraphics[width=0.32\textwidth]{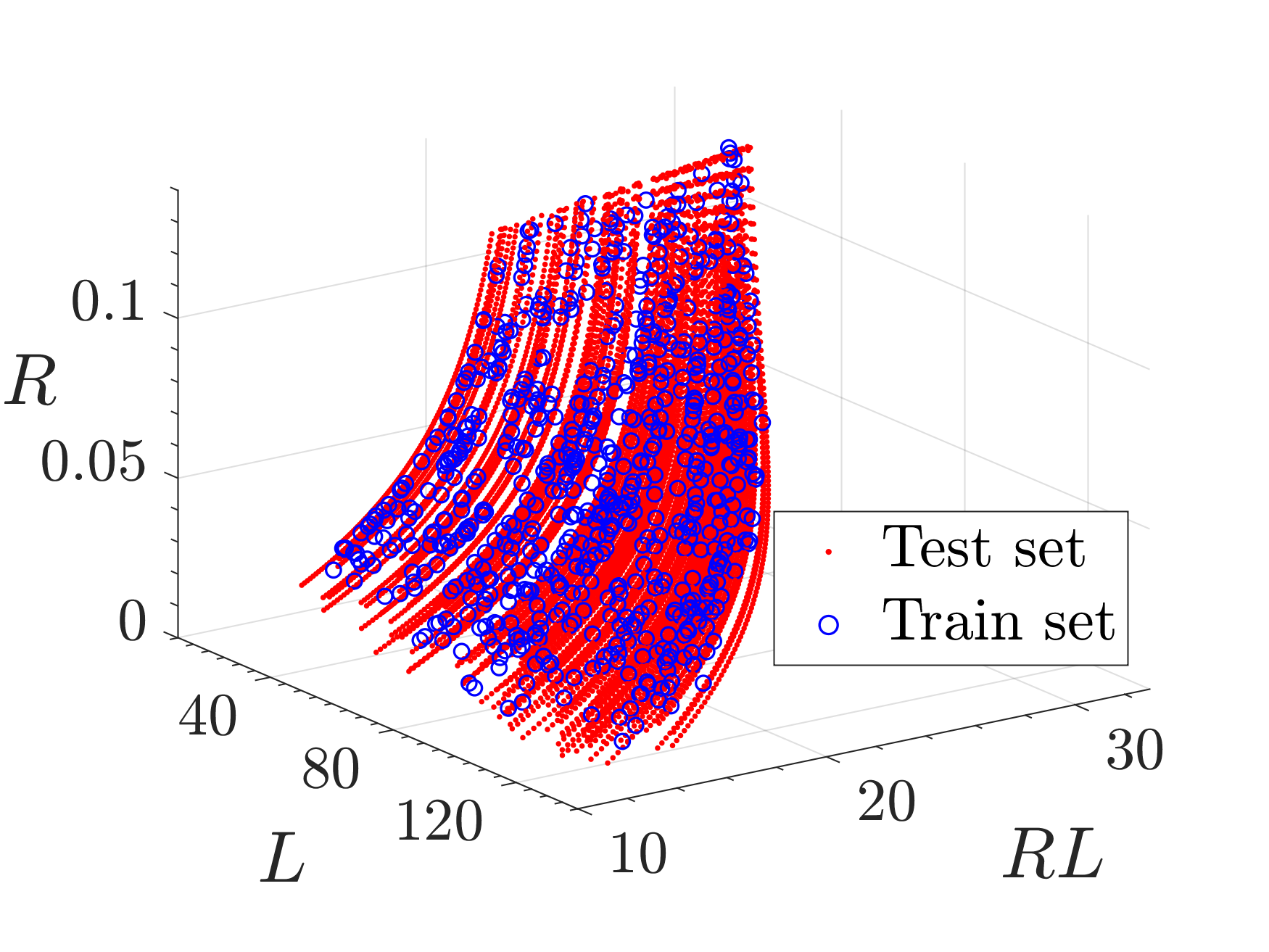}}
    \subfigure[SIM $\mathcal{M}_2$]{
    \includegraphics[width=0.32\textwidth]{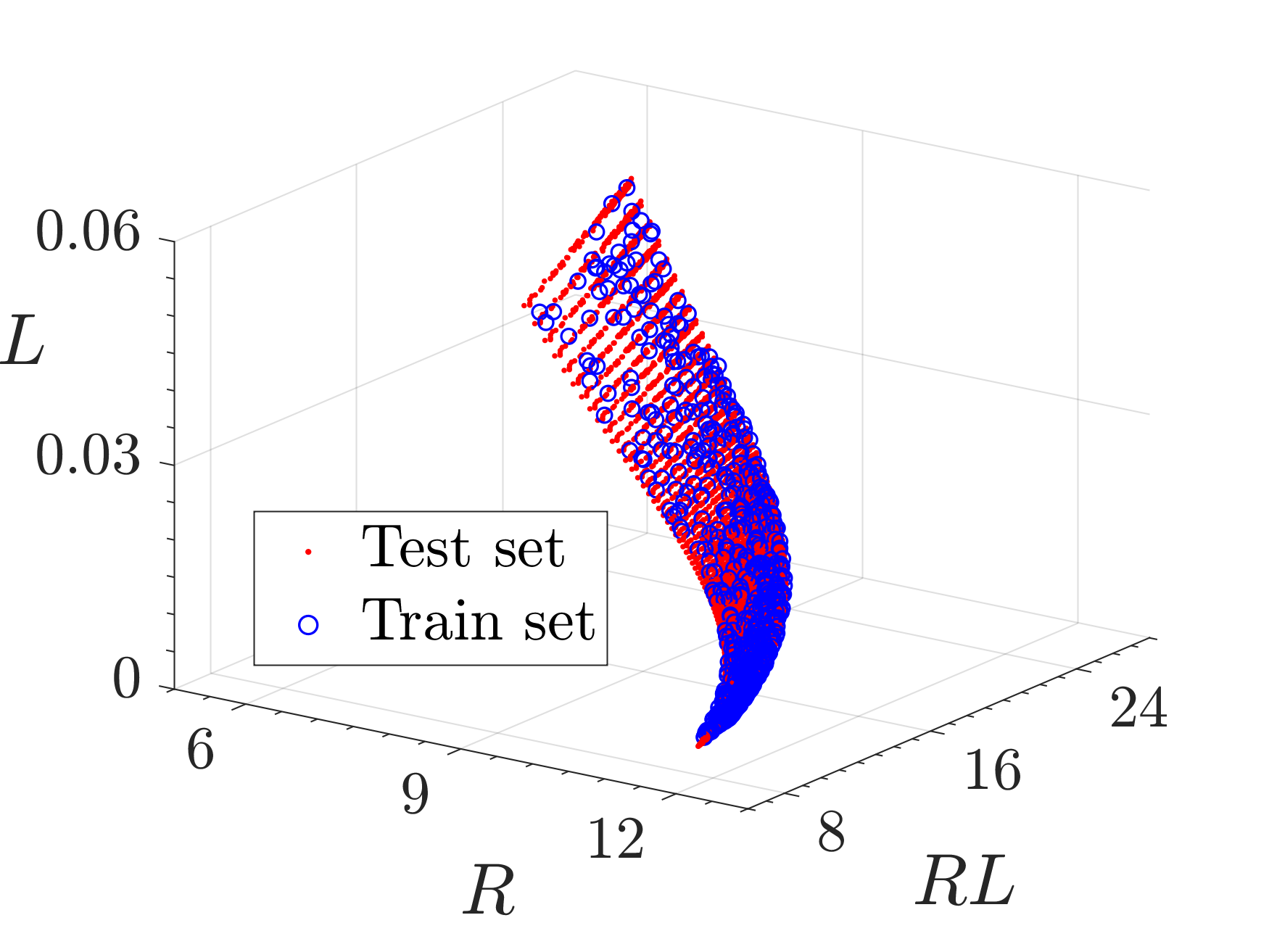}}
    \caption{TMDD system in Eq.~(33) for the SIMs $\mathcal{M}_1$ and $\mathcal{M}_2$.~Data sets used for training and validation of the PINN scheme (blue circles).~The test data sets (red dots) were used to assess the numerical approximation accuracy of the SIMs.}
    \label{SF:TMDD_Data}
\end{figure}

\begin{figure}[!h]
    \centering
    \subfigure[Projection to $c_1$]{
    \includegraphics[width=0.32\textwidth]{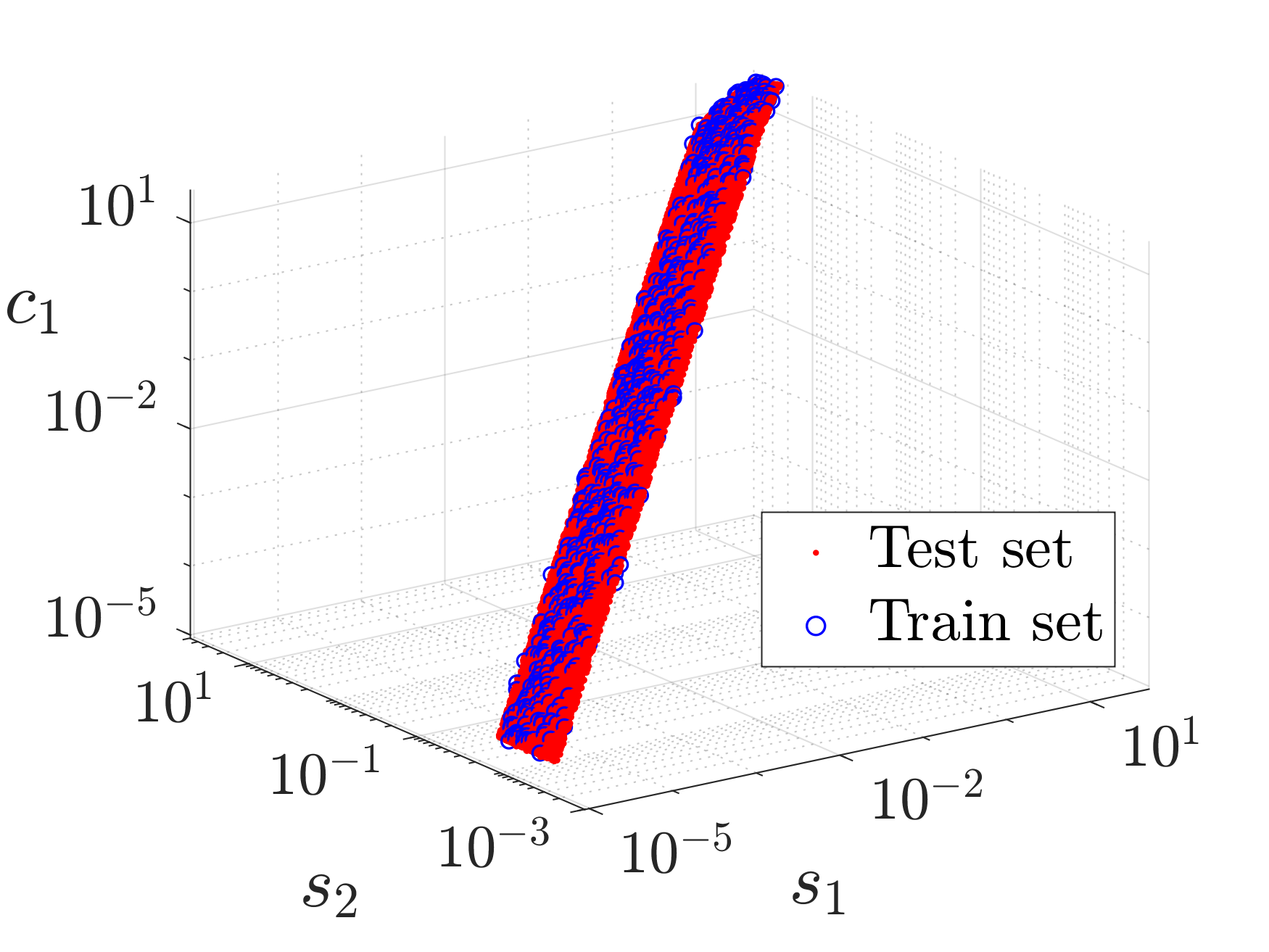}}
    \subfigure[Projection to $c_2$]{
    \includegraphics[width=0.32\textwidth]{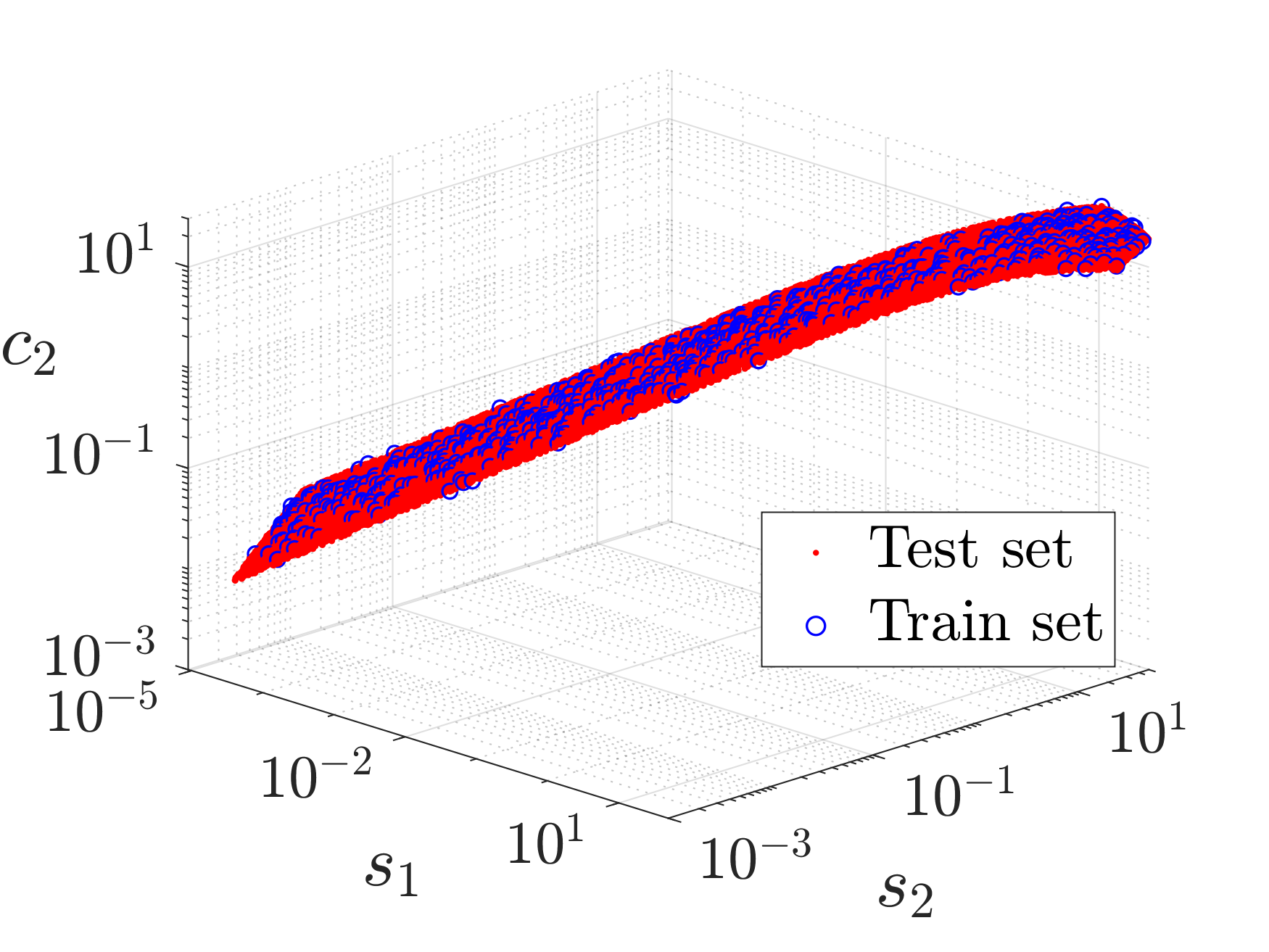}}
    \caption{Original fCSI system in Eq.~(34); projections to $c_1$ and $c_2$ fast variables.~Data sets used for training and validation of the PINN scheme (blue circles).~The test data sets (red dots) were used to assess the numerical approximation accuracy of the SIMs.}
    \label{SF:InhSt_Data}
\end{figure}

\begin{figure}[!h]
    \centering
    \subfigure[Projection to $c_1$]{
    \includegraphics[width=0.32\textwidth]{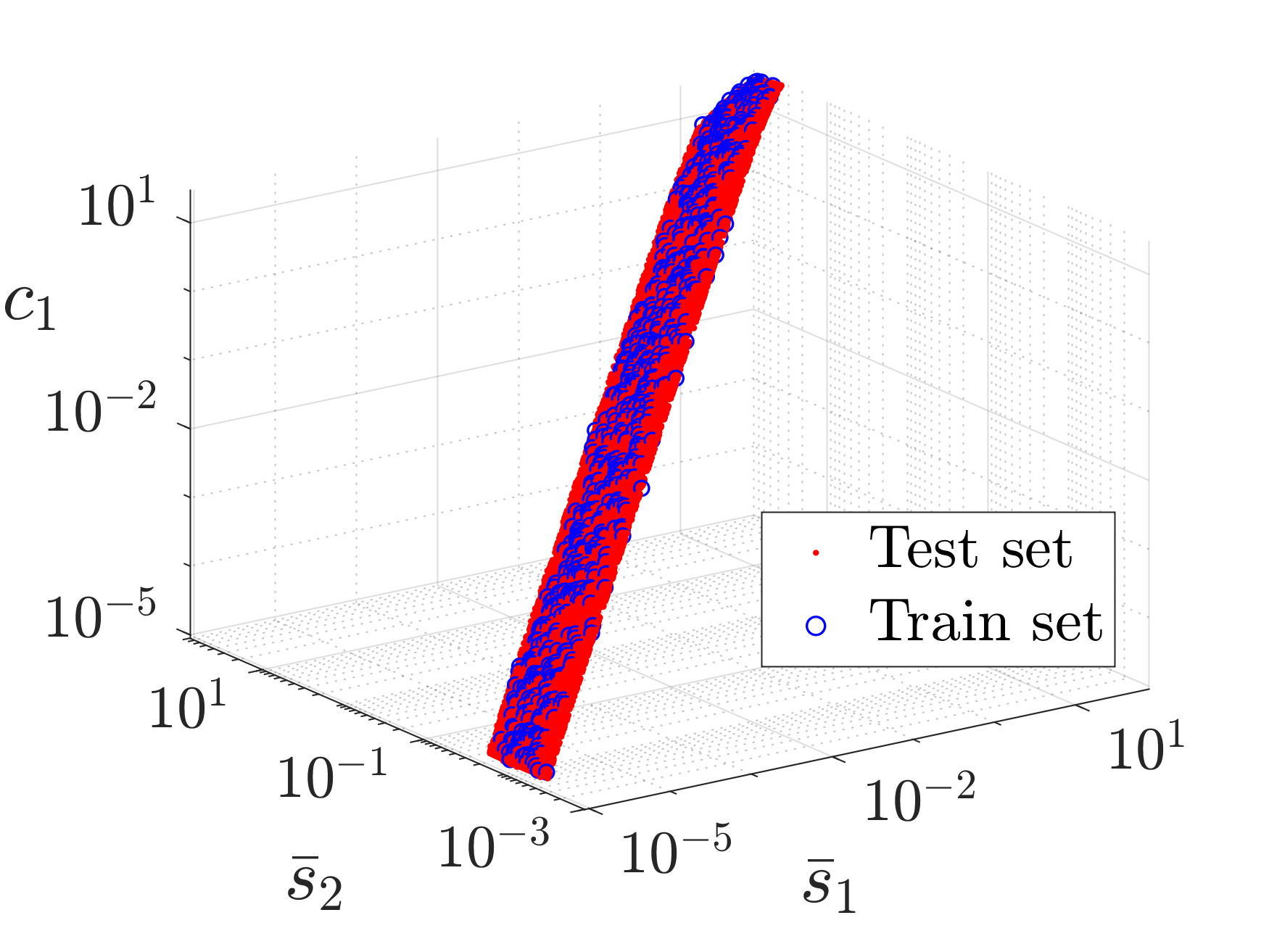}}
    \subfigure[Projection to $c_2$]{
    \includegraphics[width=0.32\textwidth]{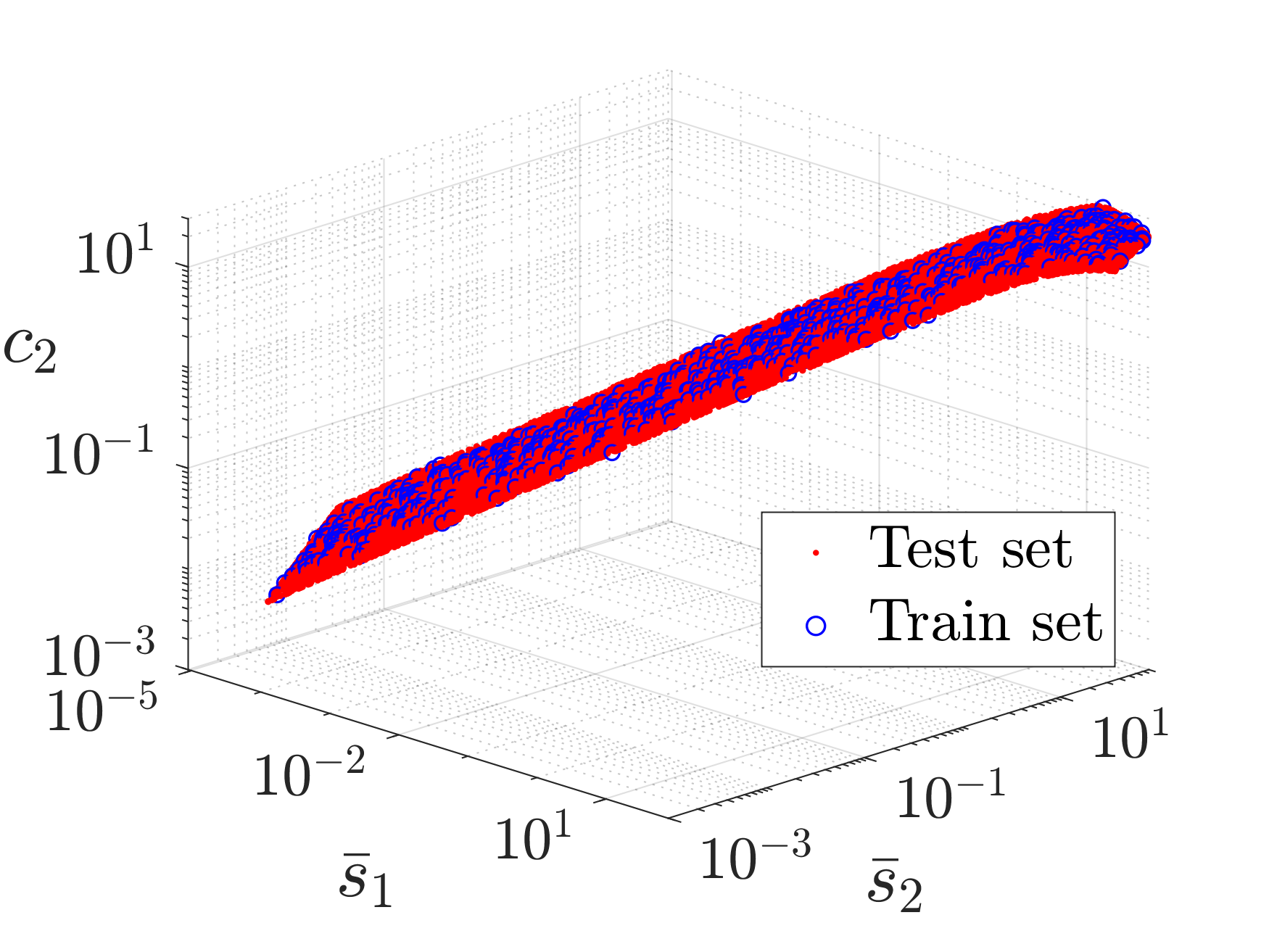}}
    \caption{Transformed fCSI system in Eq.~(35); projections to $c_1$ and $c_2$ fast variables.~Data sets used for training and validation of the PINN scheme (blue circles).~The test data sets (red dots) were used to assess the numerical approximation accuracy of the SIMs.}
    \label{SF:InhTr_Data}
\end{figure}

\newpage
\clearpage
\section{Numerical SIM approximation errors}

\begin{figure}[!h]
    \centering
    \subfigure[PIML]{
    \includegraphics[width=0.32\textwidth]{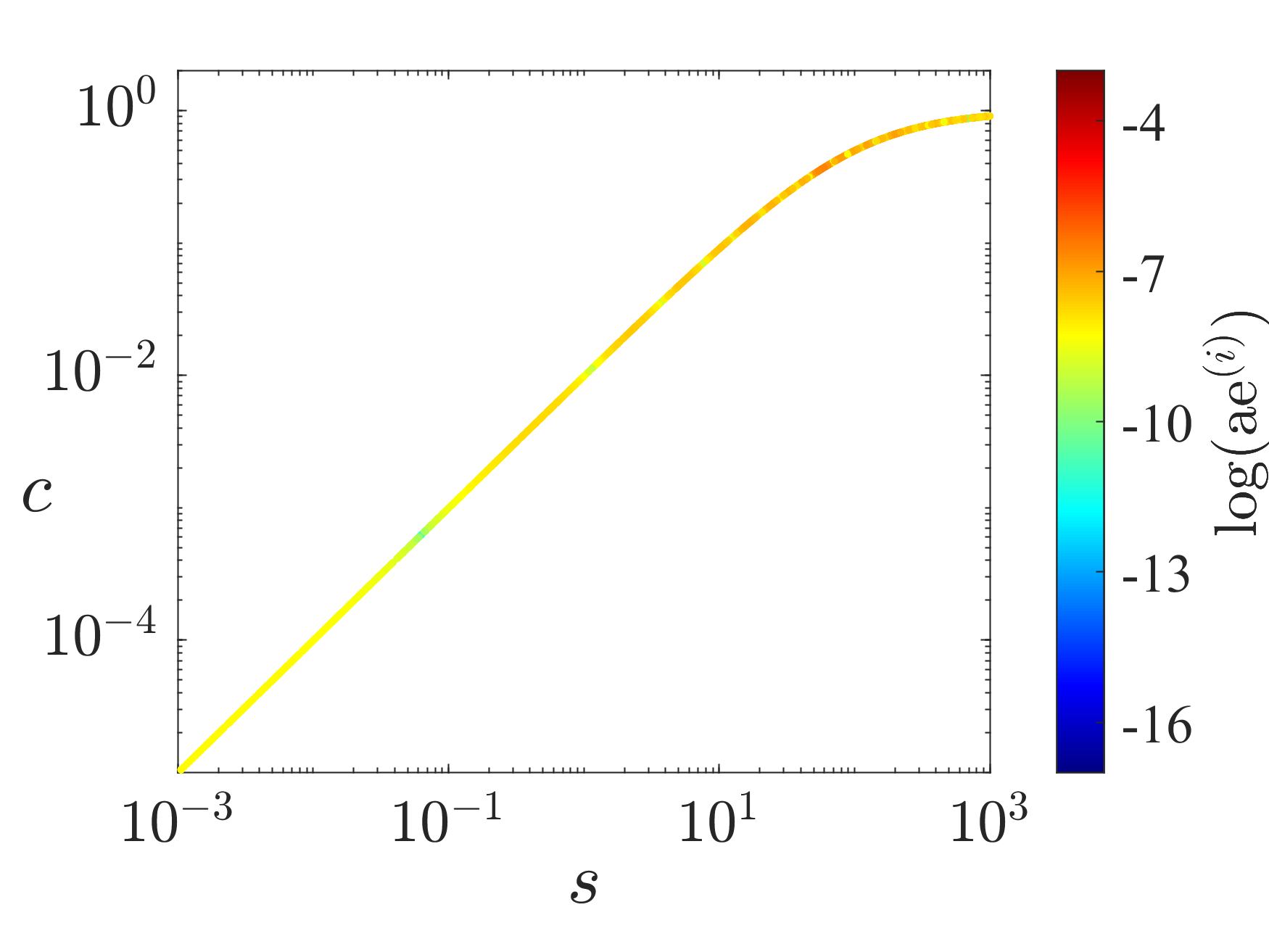}}
    \subfigure[PEA]{
    \includegraphics[width=0.32\textwidth]{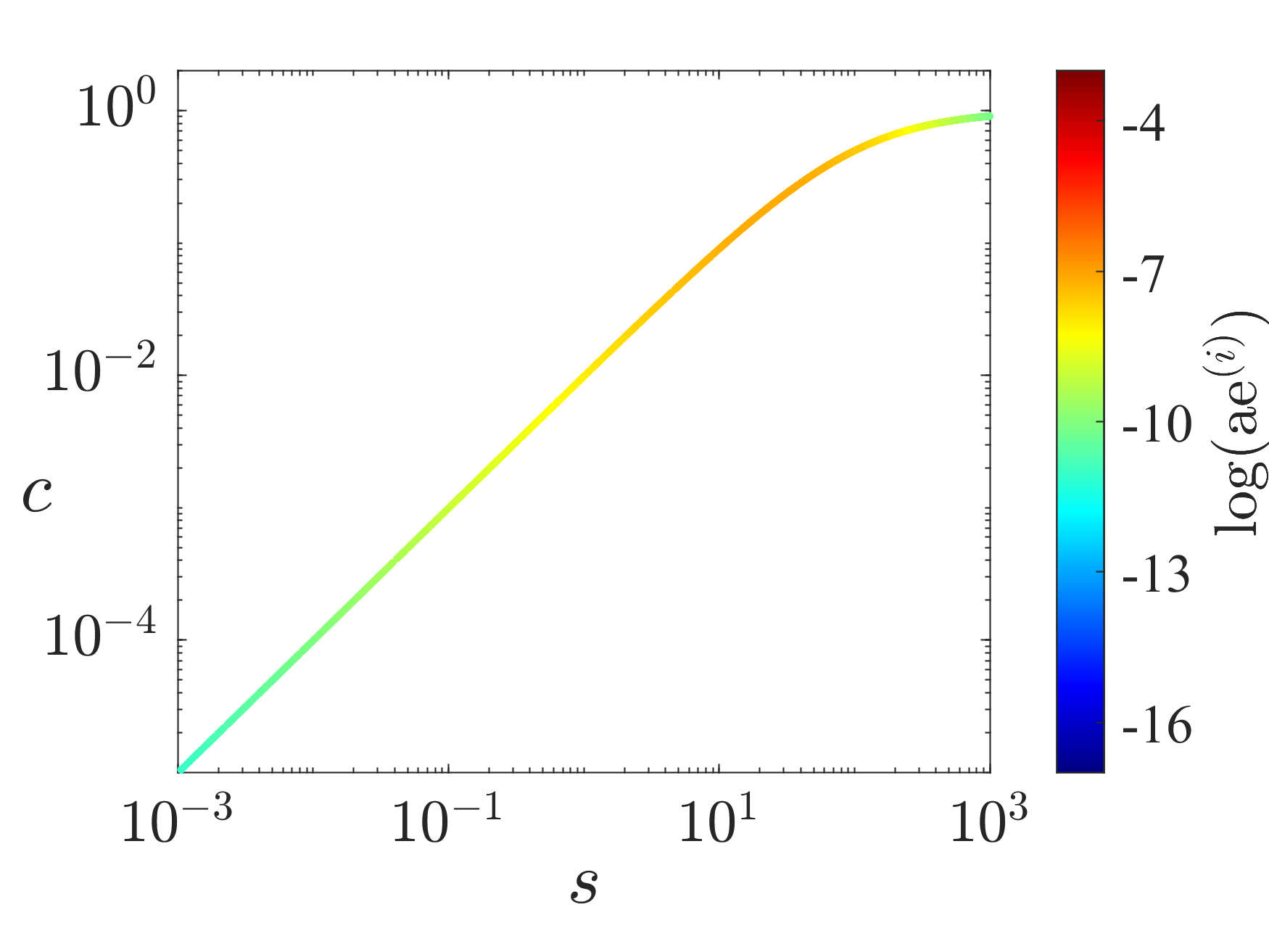}}
    \subfigure[CSP$_e$]{
    \includegraphics[width=0.32\textwidth]{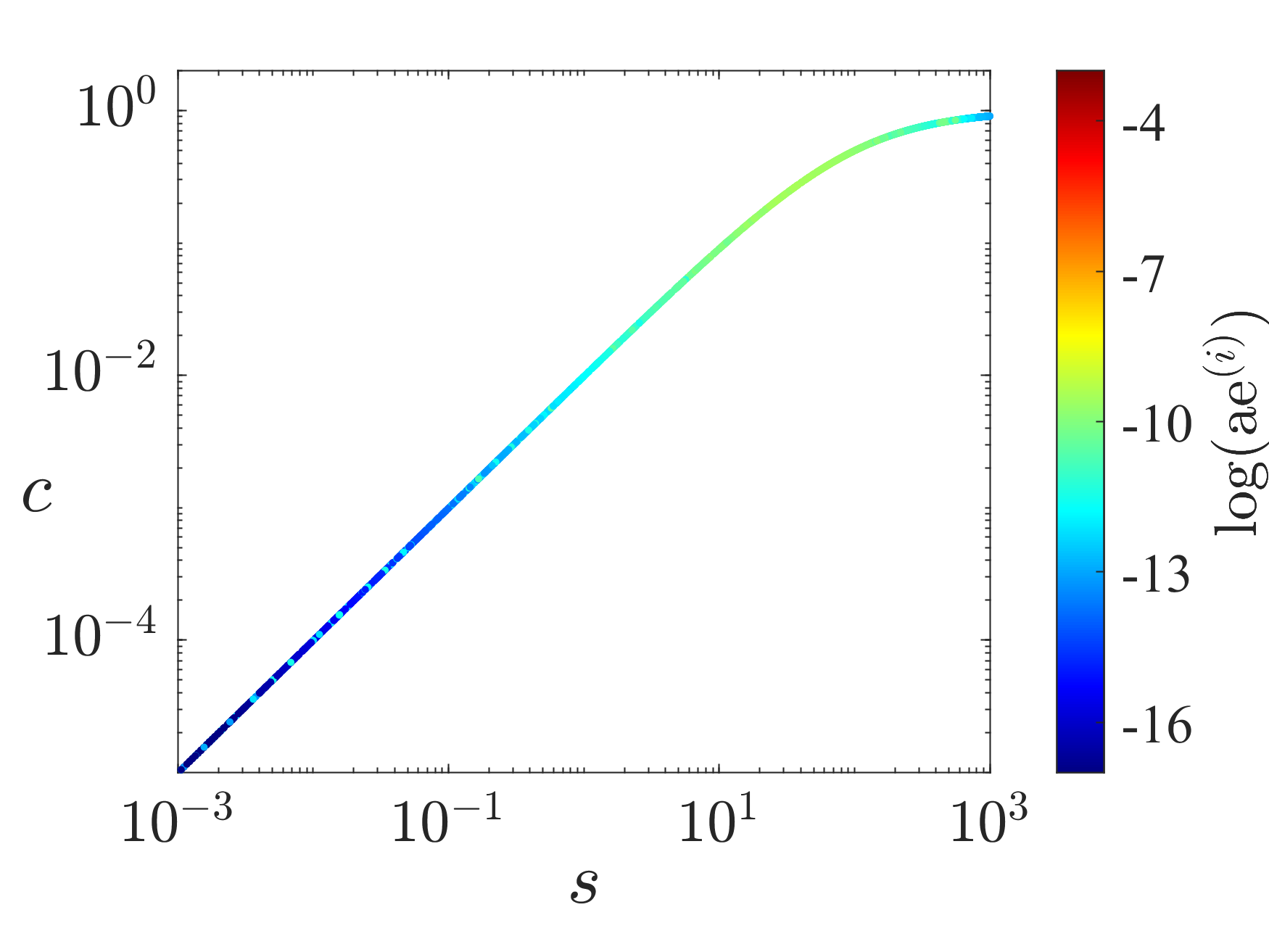}} \\ 
    \subfigure[rQSSA]{
    \includegraphics[width=0.32\textwidth]{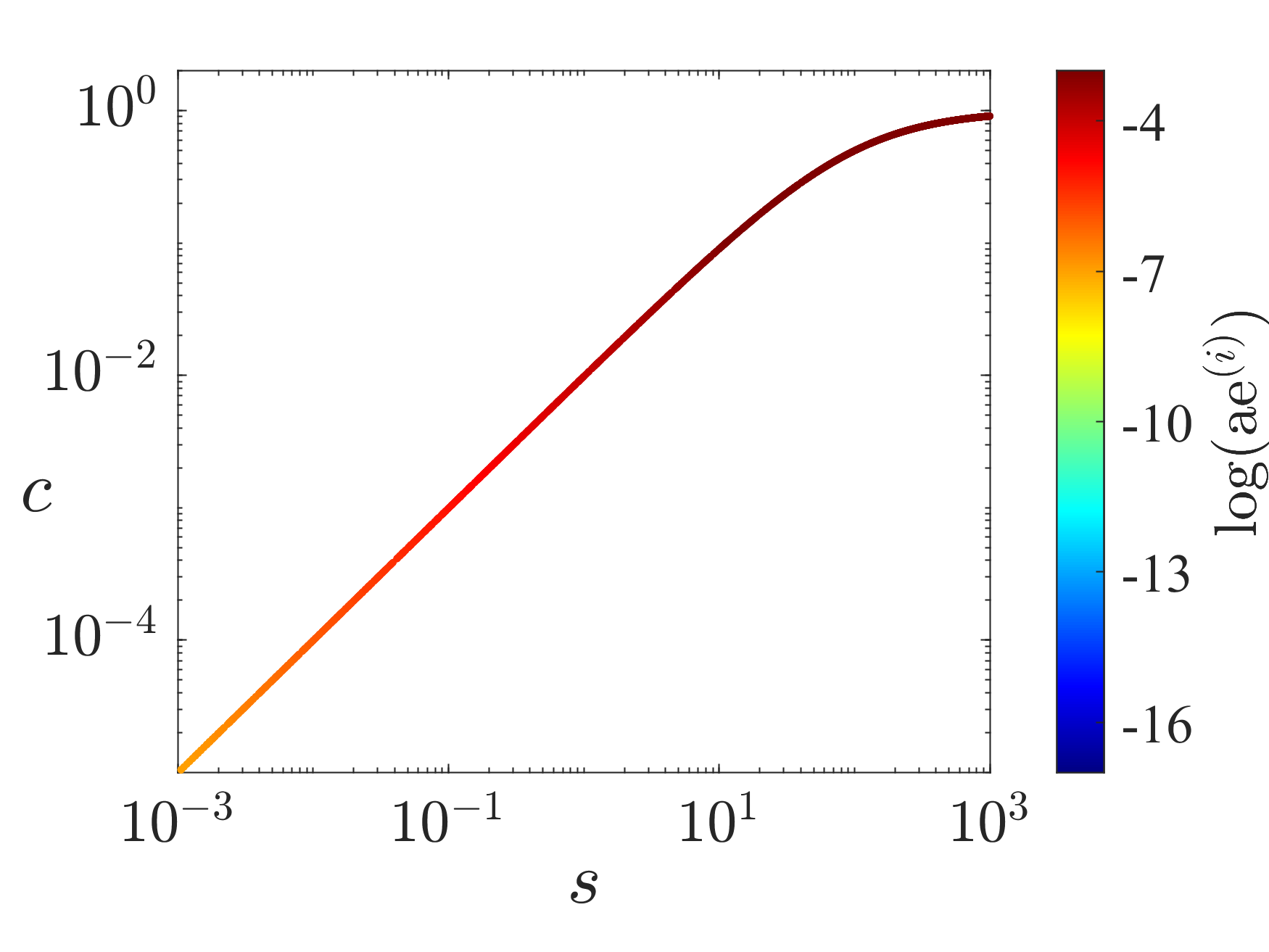}}
    \subfigure[CSP$_s$(1)]{
    \includegraphics[width=0.32\textwidth]{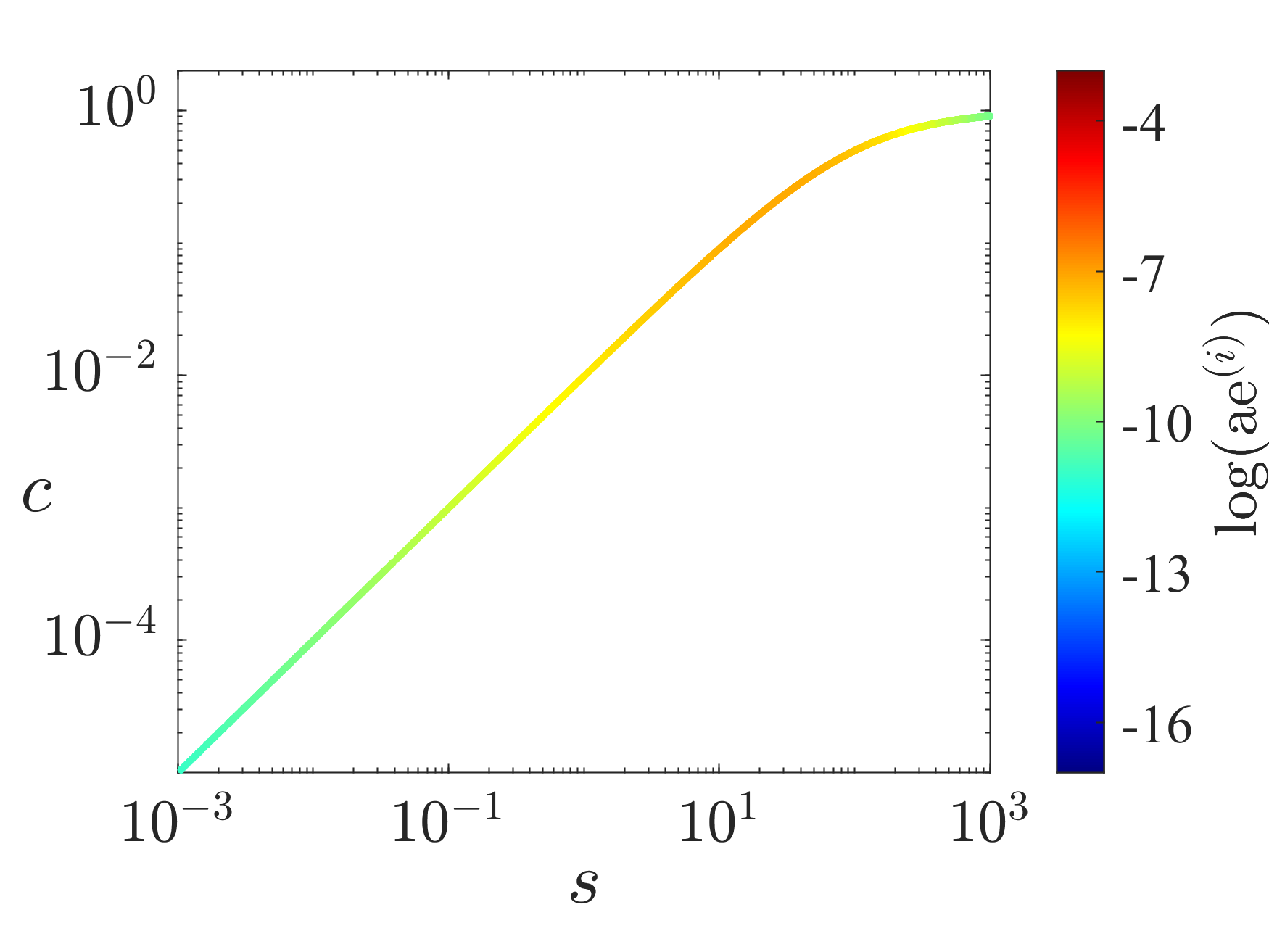}}
    \subfigure[CSP$_s$(2)]{
    \includegraphics[width=0.32\textwidth]{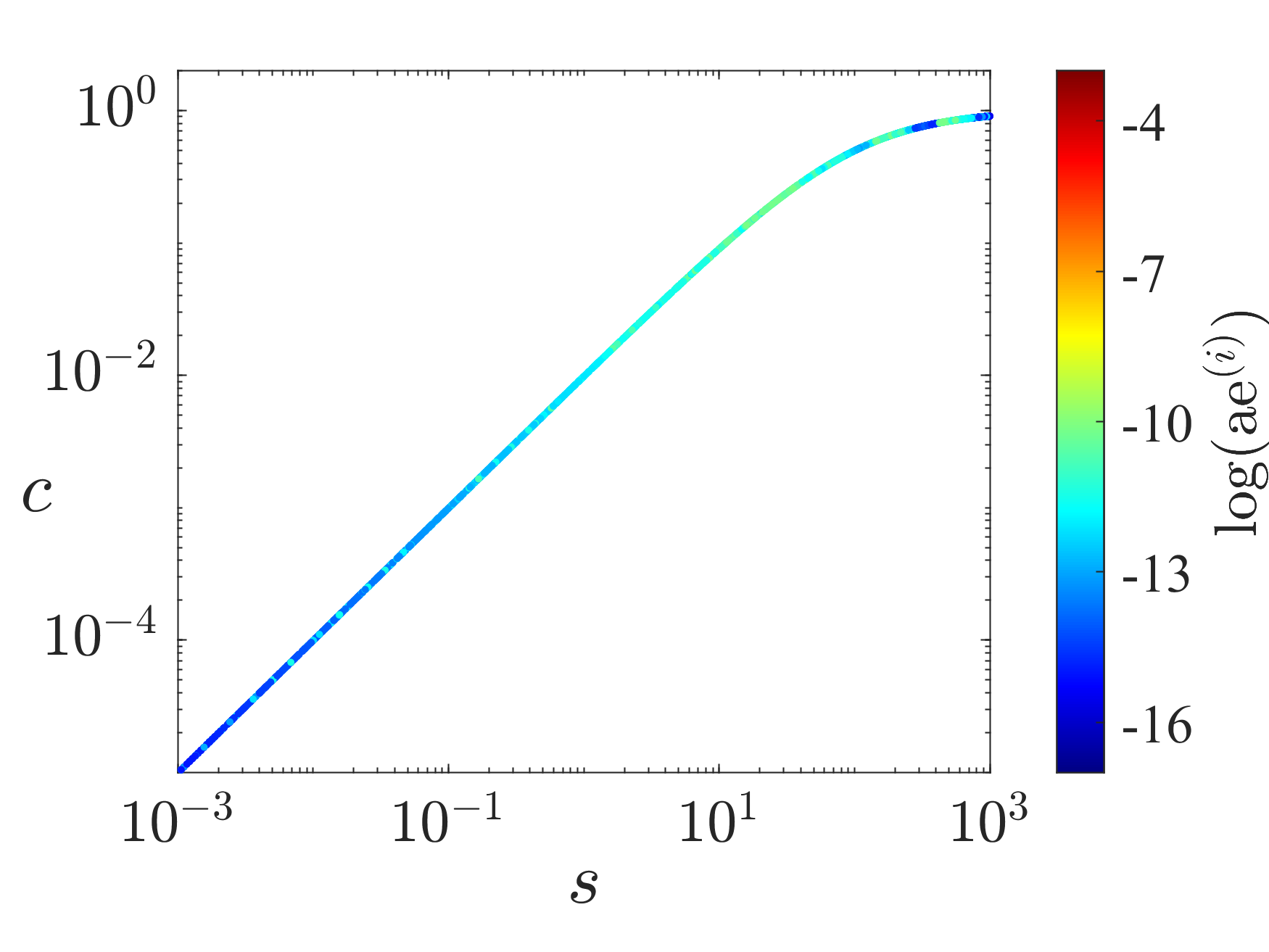}} \\
    \subfigure[sQSSA]{
    \includegraphics[width=0.32\textwidth]{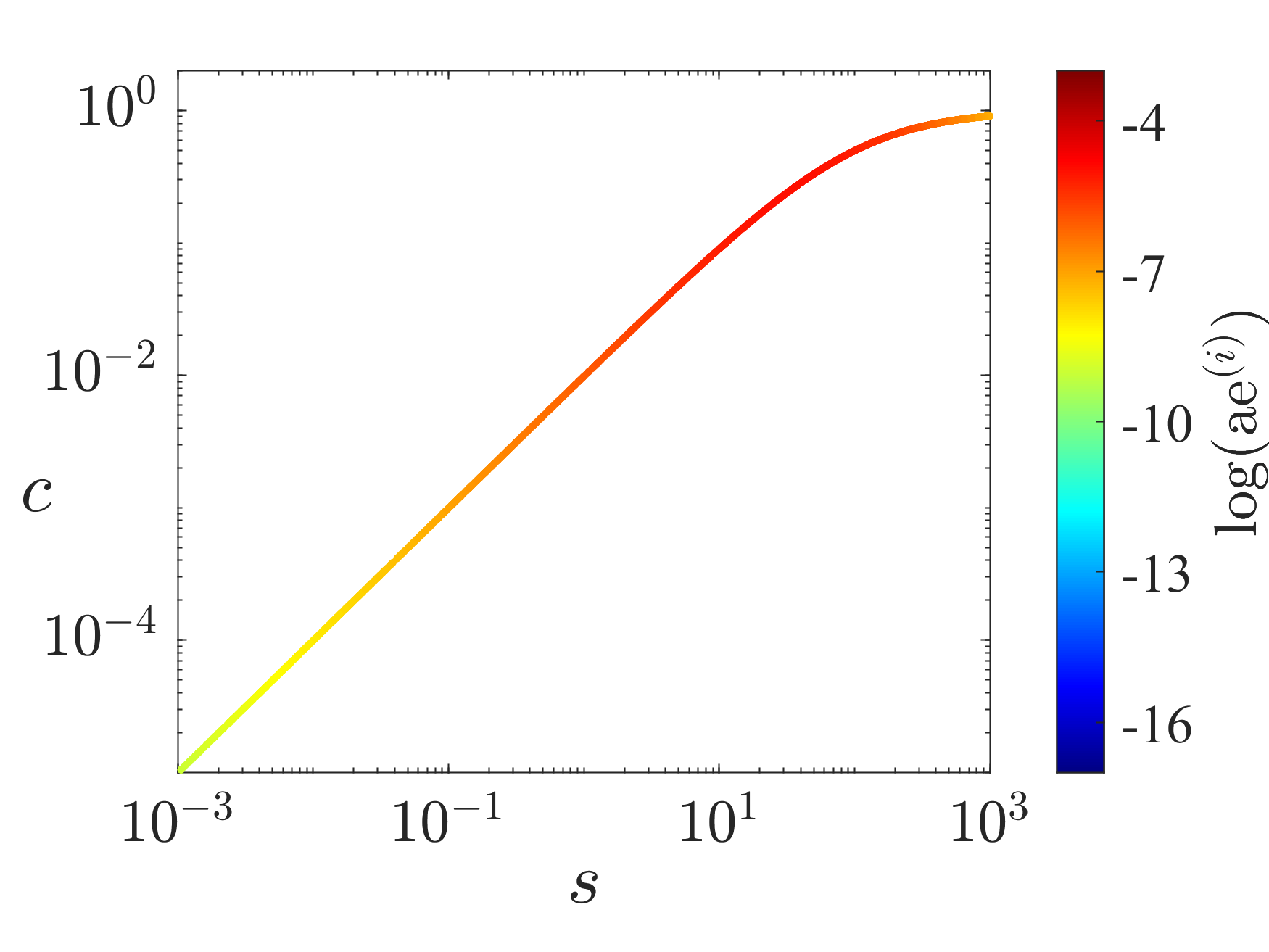}}
    \subfigure[CSP$_c$(1)]{
    \includegraphics[width=0.32\textwidth]{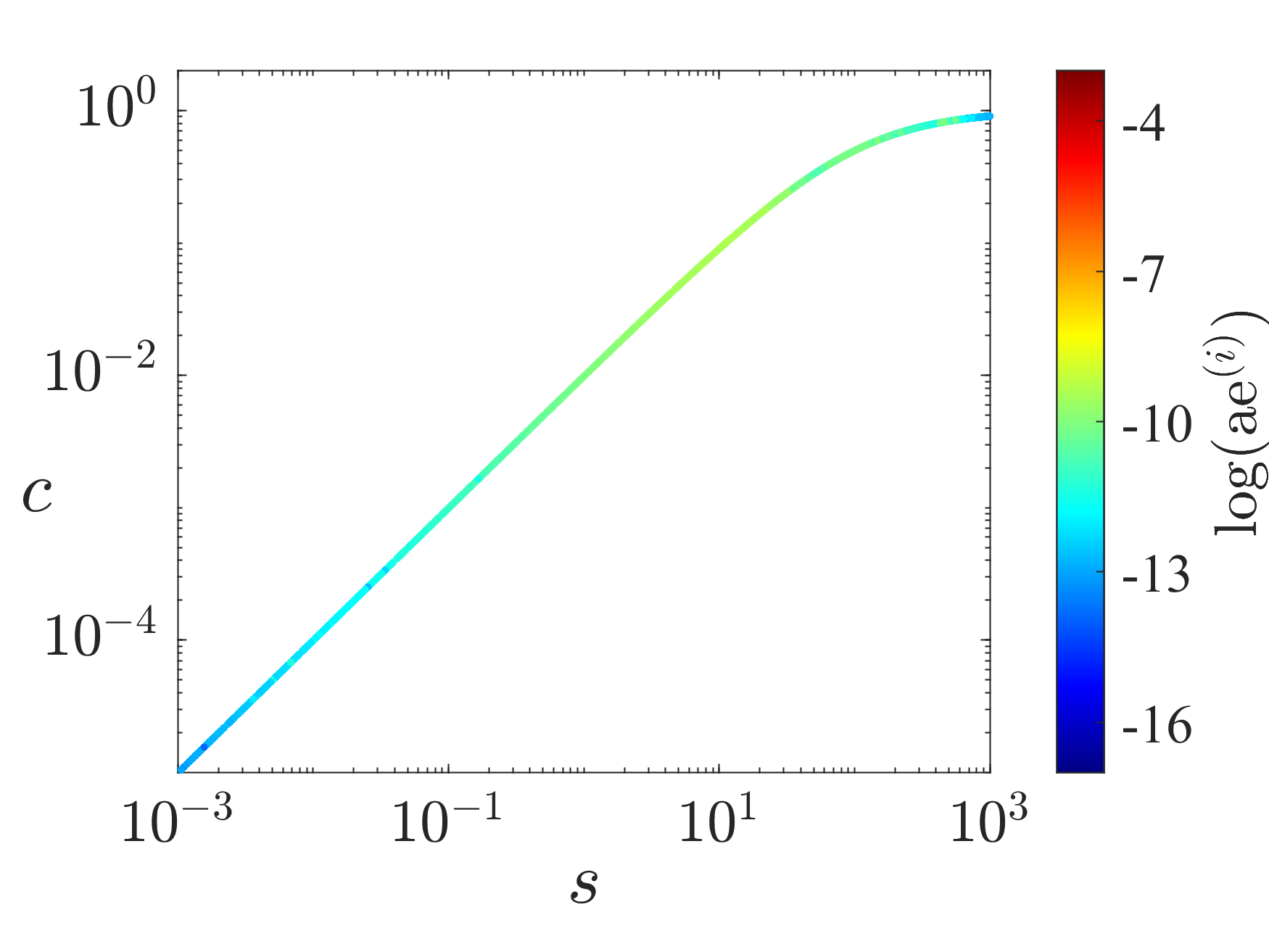}}
    \subfigure[CSP$_c$(2)]{
    \includegraphics[width=0.32\textwidth]{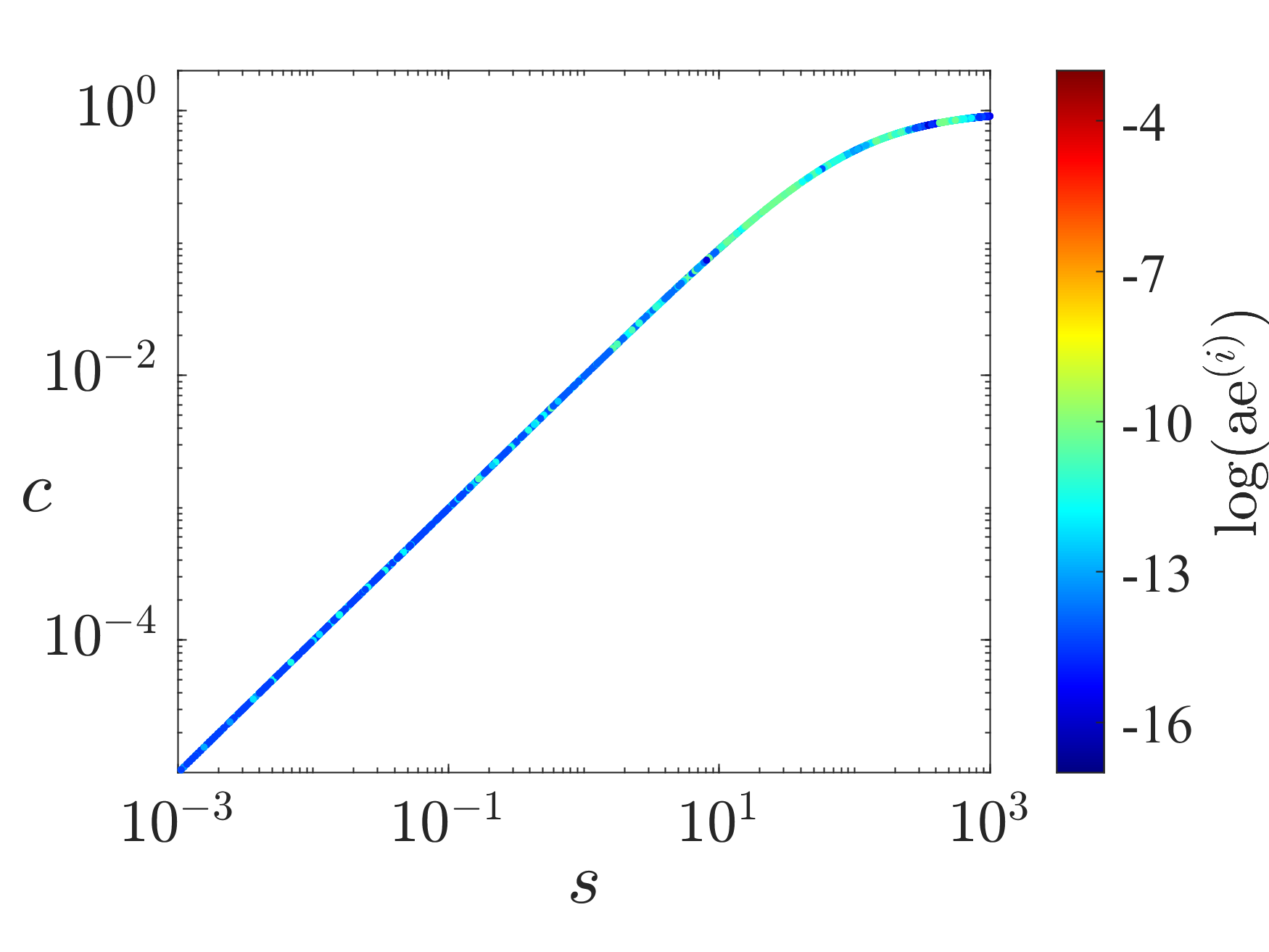}}
    \caption{MM system in Eq.~(32) for the MM1 case, where $x=c$ is the fast variable.~Absolute errors ($ae^{(i)}$) of the SIM approximations over all points of the test set, in comparison to the numerical solution $\mathbf{z}^{(i)}=[x^{(i)},y^{(i)}]^\top$ for $i=1,\ldots,n_t$.~Panel (a) depicts the $\lvert \mathbf{C} \mathbf{z}^{(i)} - \mathcal{N}(\mathbf{D} \mathbf{z}^{(i)}) \rvert$ of the PINN scheme, panels (b), (d), (e), (g) and (h) depict $\lvert x^{(i)} - h(y^{(i)}) \rvert$ of the PEA, rQSSA,  CSP$_s$(1), sQSSA and CSP$_c$(1) explicit functionals w.r.t. $x$, and panels (c), (f) and (i) depict $\lvert x^{(i)} - \hat{x}^{(i)}) \rvert$ of the CSP$_e$, CSP$_s$(2) and CSP$_c$(2) implicit functionals, solved numerically with Newton for $x$.~Note that the approximations of the second/third row were constructed with $s$/$c$ assumed as fast variable.}
    \label{SF:MM1_AE}
\end{figure}

\begin{figure}[!h]
    \centering
    \subfigure[PIML]{
    \includegraphics[width=0.32\textwidth]{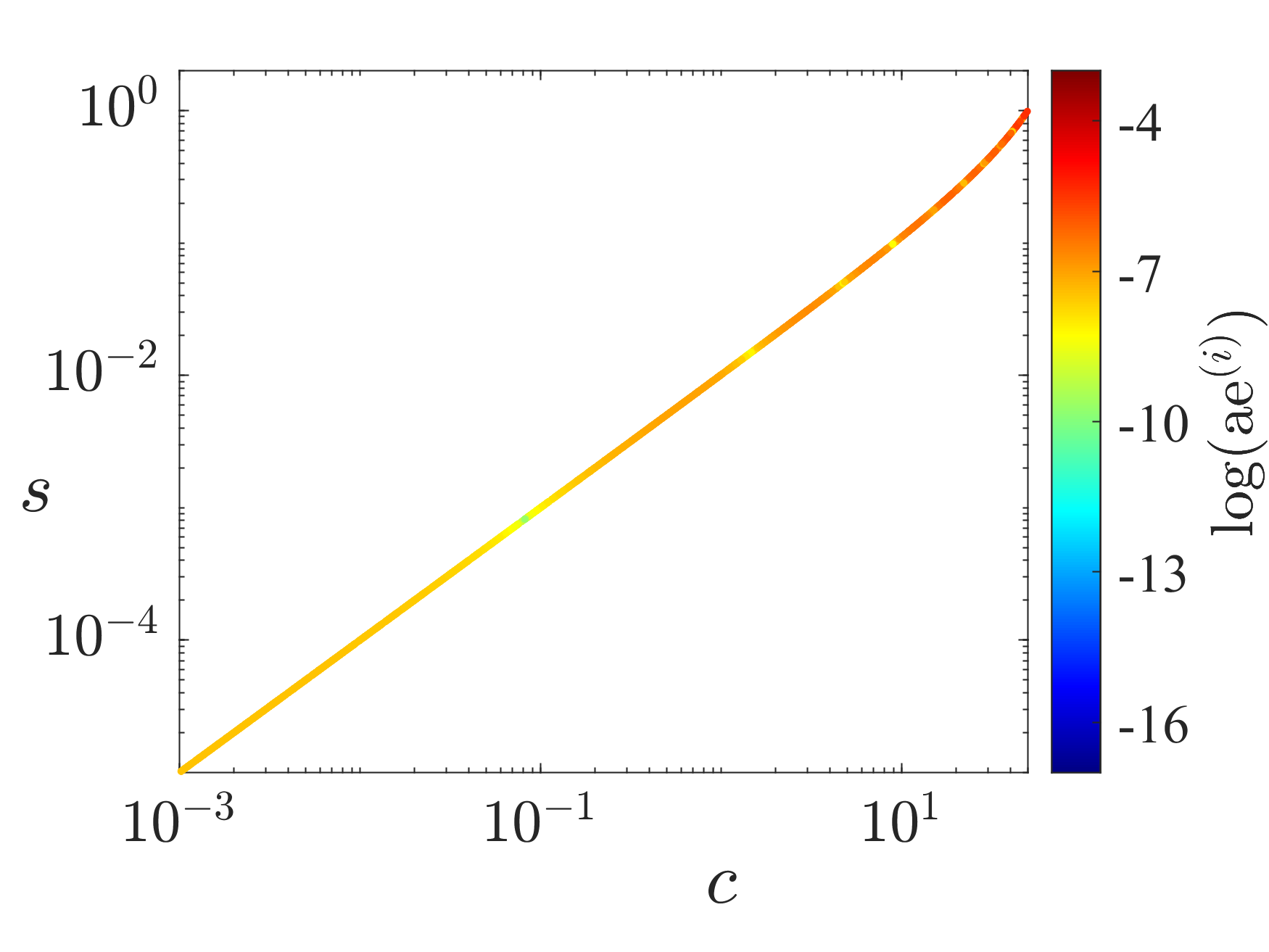}}
    \subfigure[PEA]{
    \includegraphics[width=0.32\textwidth]{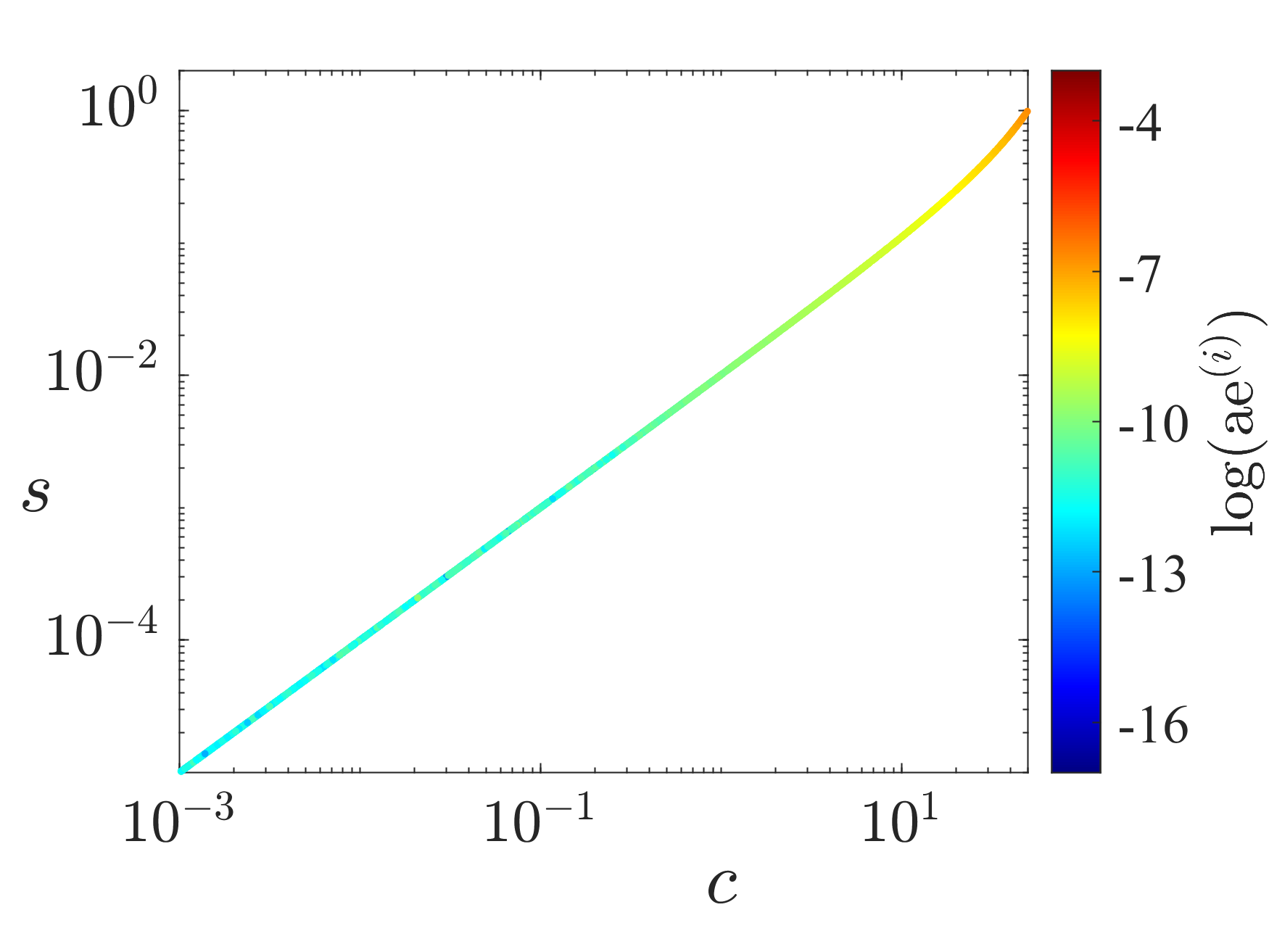}}
    \subfigure[CSP$_e$]{
    \includegraphics[width=0.32\textwidth]{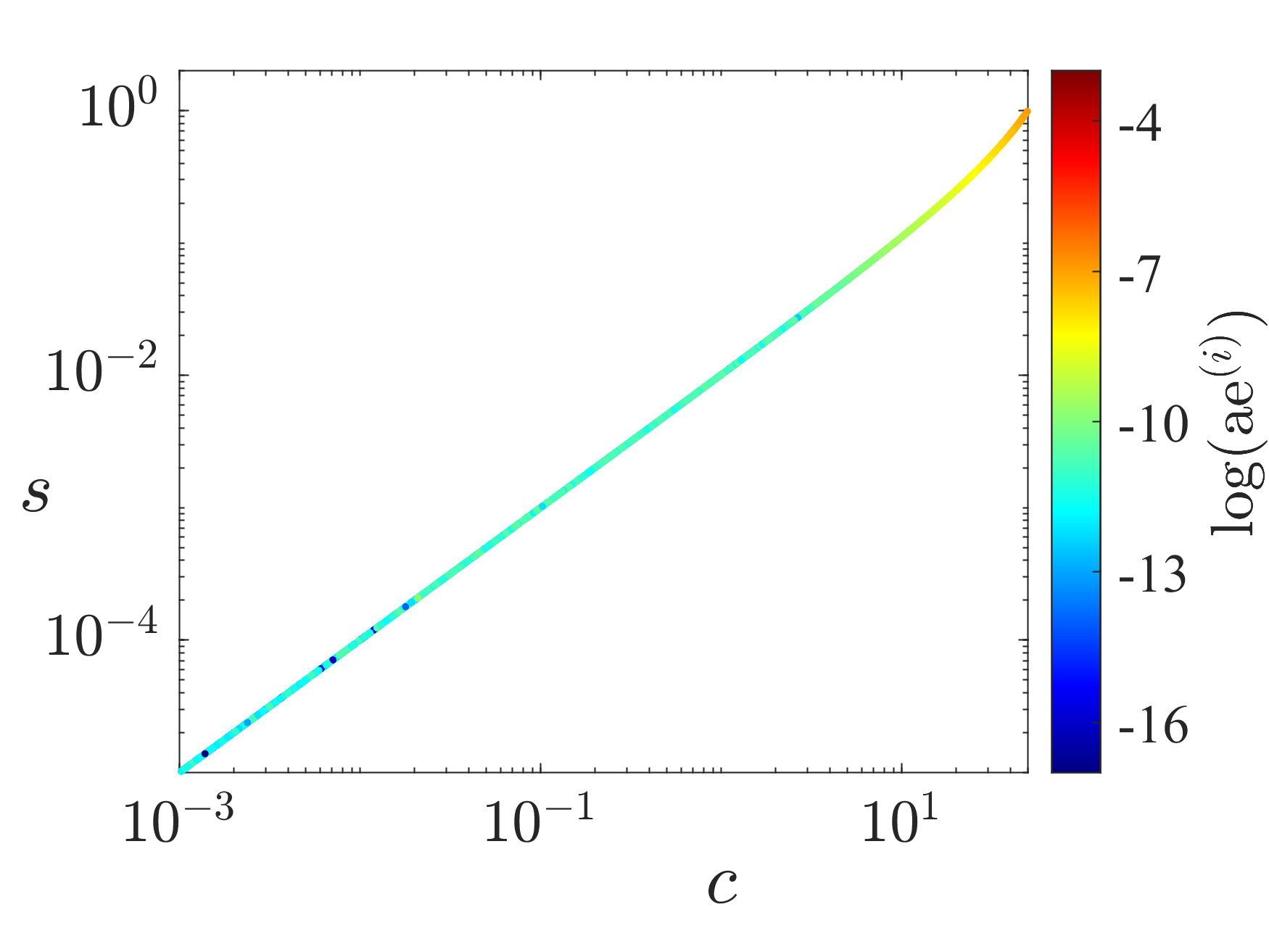}} \\ 
    \subfigure[rQSSA]{
    \includegraphics[width=0.32\textwidth]{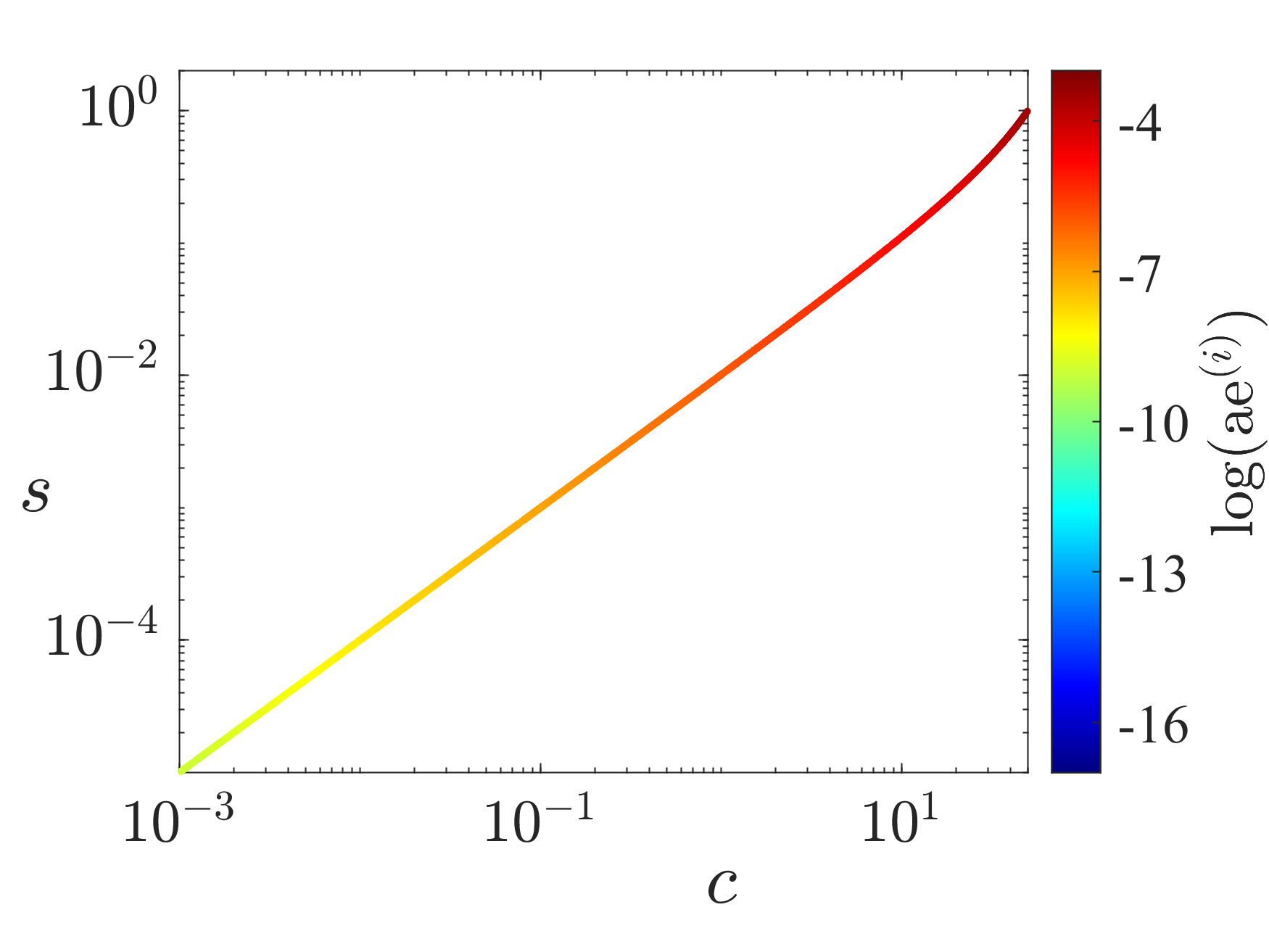}}
    \subfigure[CSP$_s$(1)]{
    \includegraphics[width=0.32\textwidth]{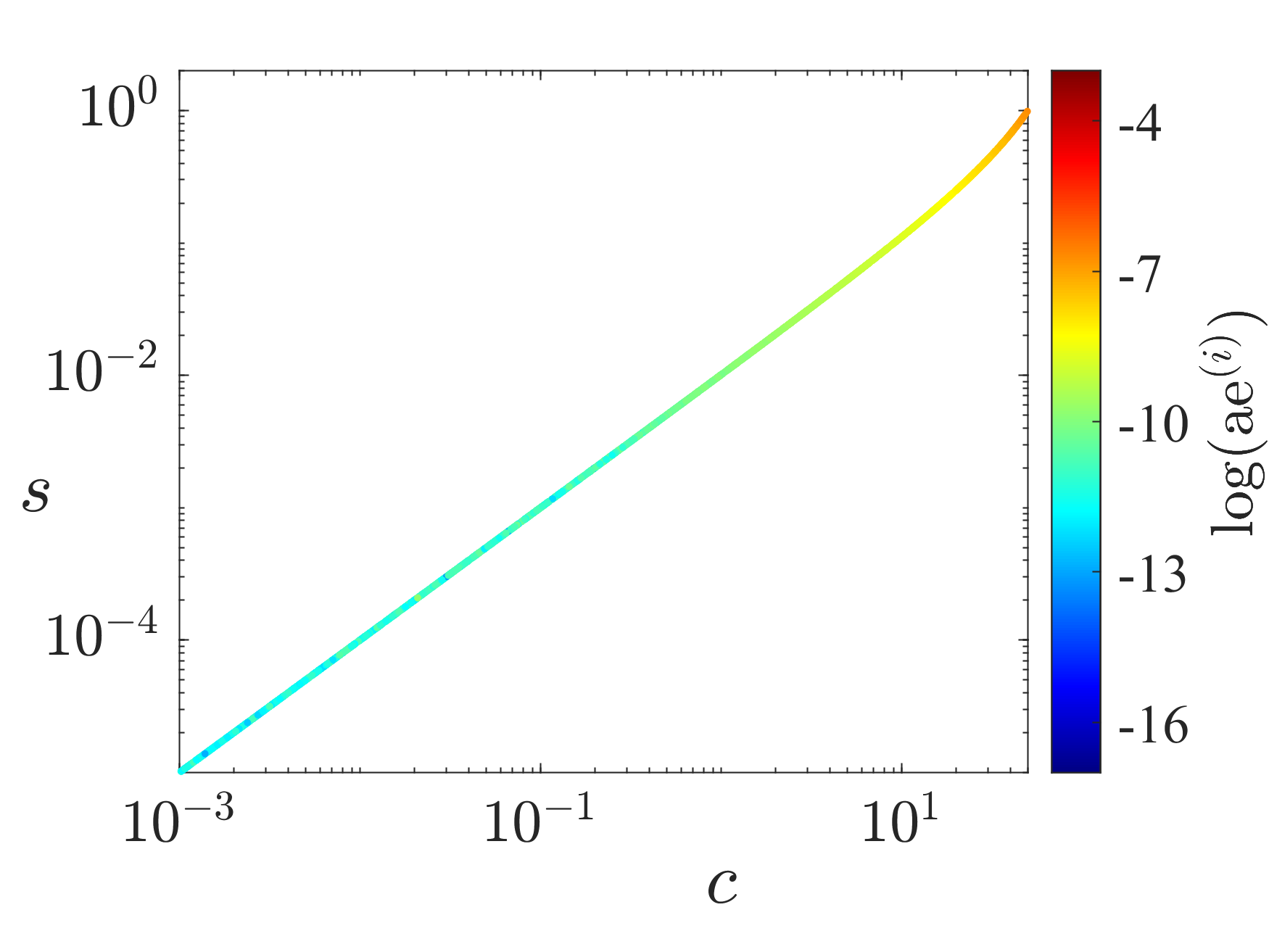}}
    \subfigure[CSP$_s$(2)]{
    \includegraphics[width=0.32\textwidth]{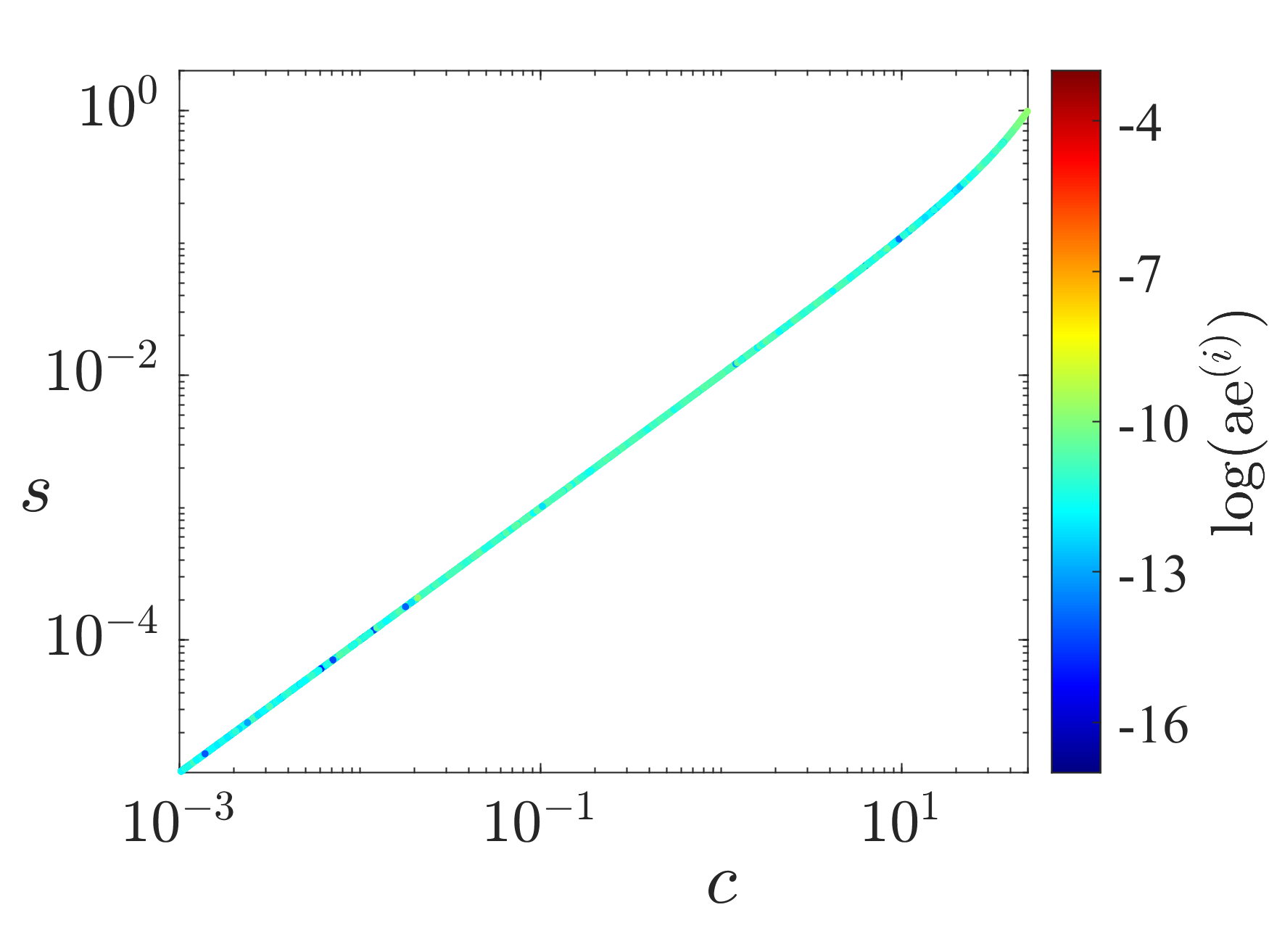}} \\
    \subfigure[sQSSA]{
    \includegraphics[width=0.32\textwidth]{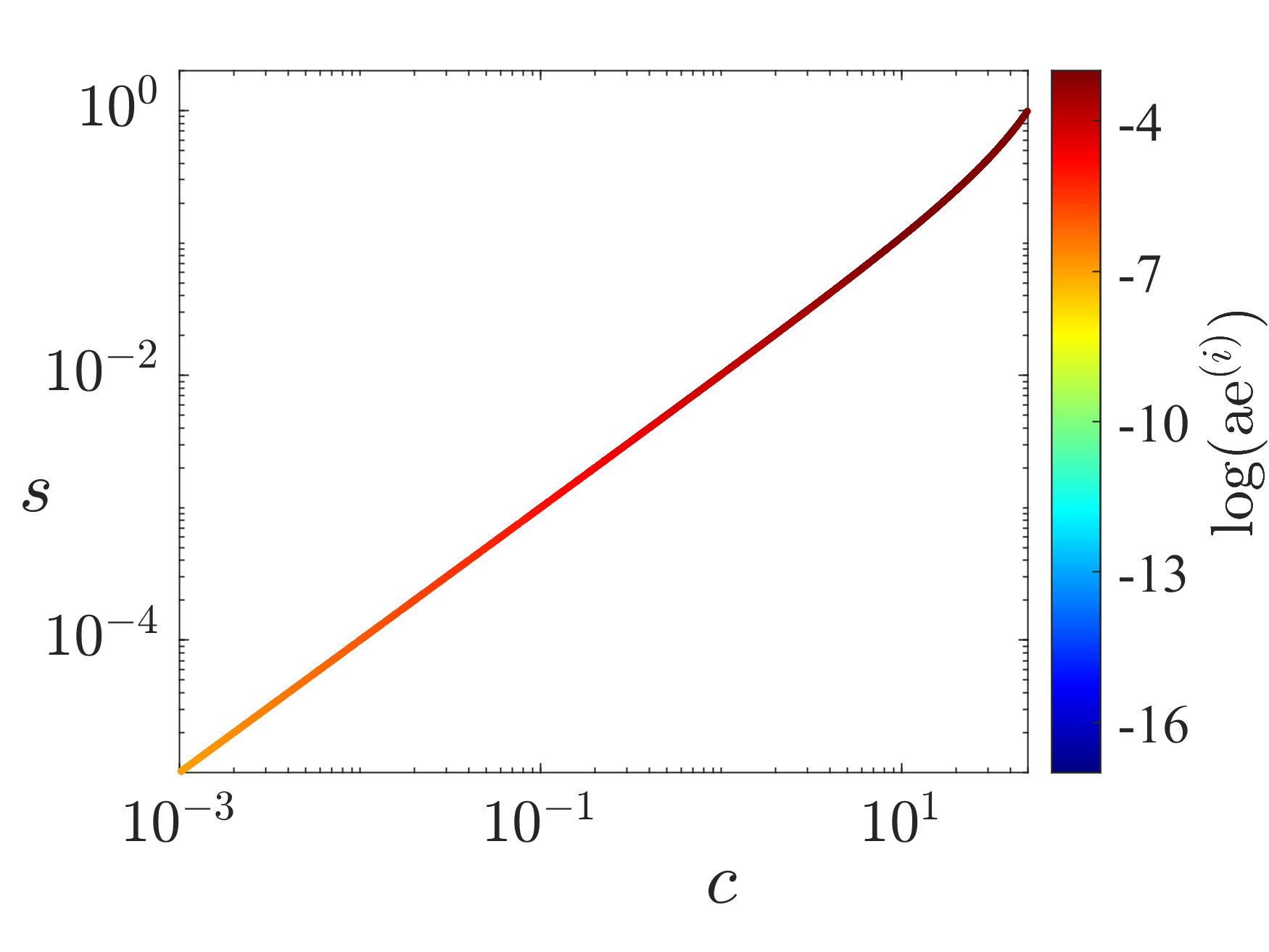}}
    \subfigure[CSP$_c$(1)]{
    \includegraphics[width=0.32\textwidth]{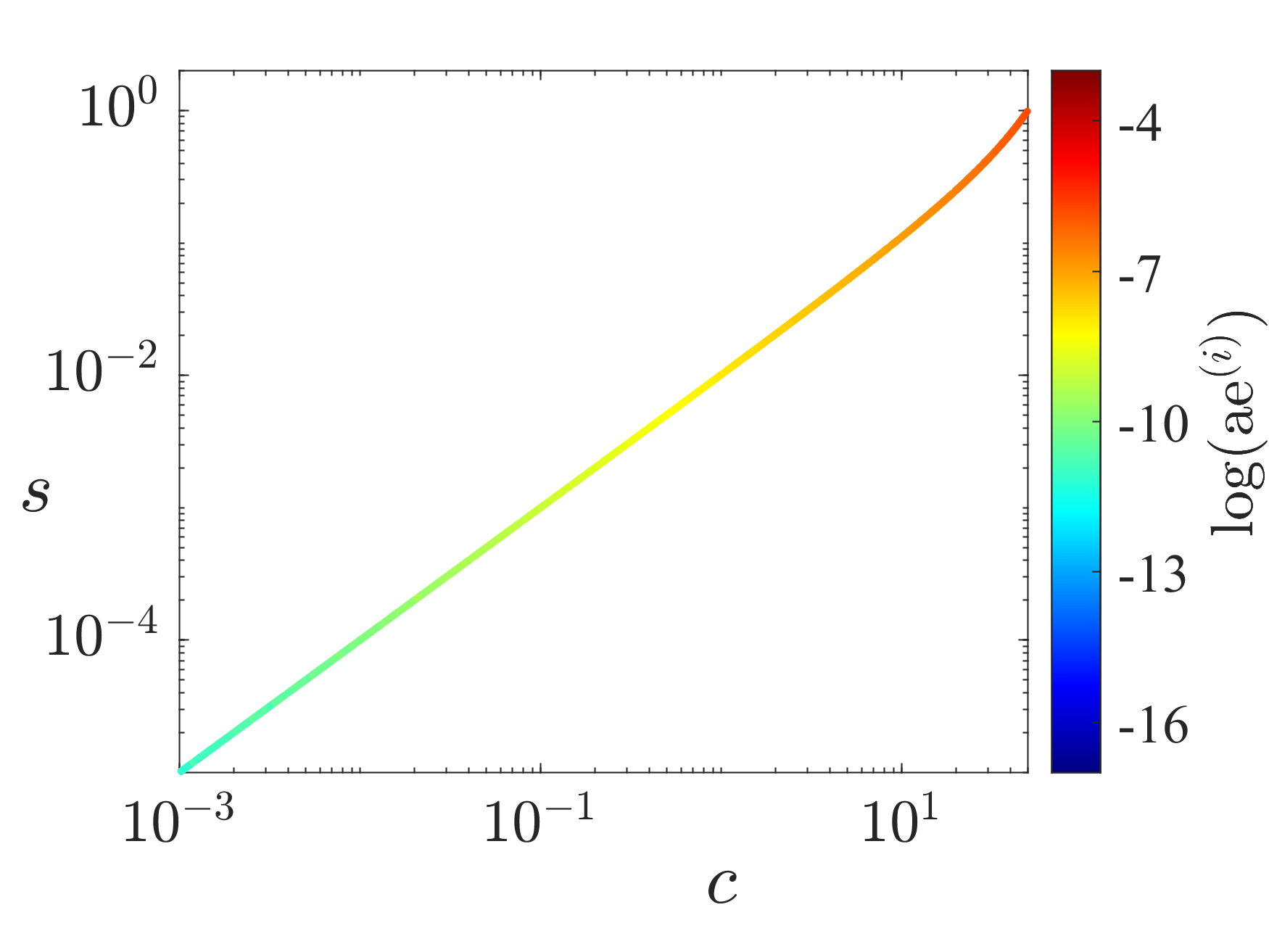}}
    \subfigure[CSP$_c$(2)]{
    \includegraphics[width=0.32\textwidth]{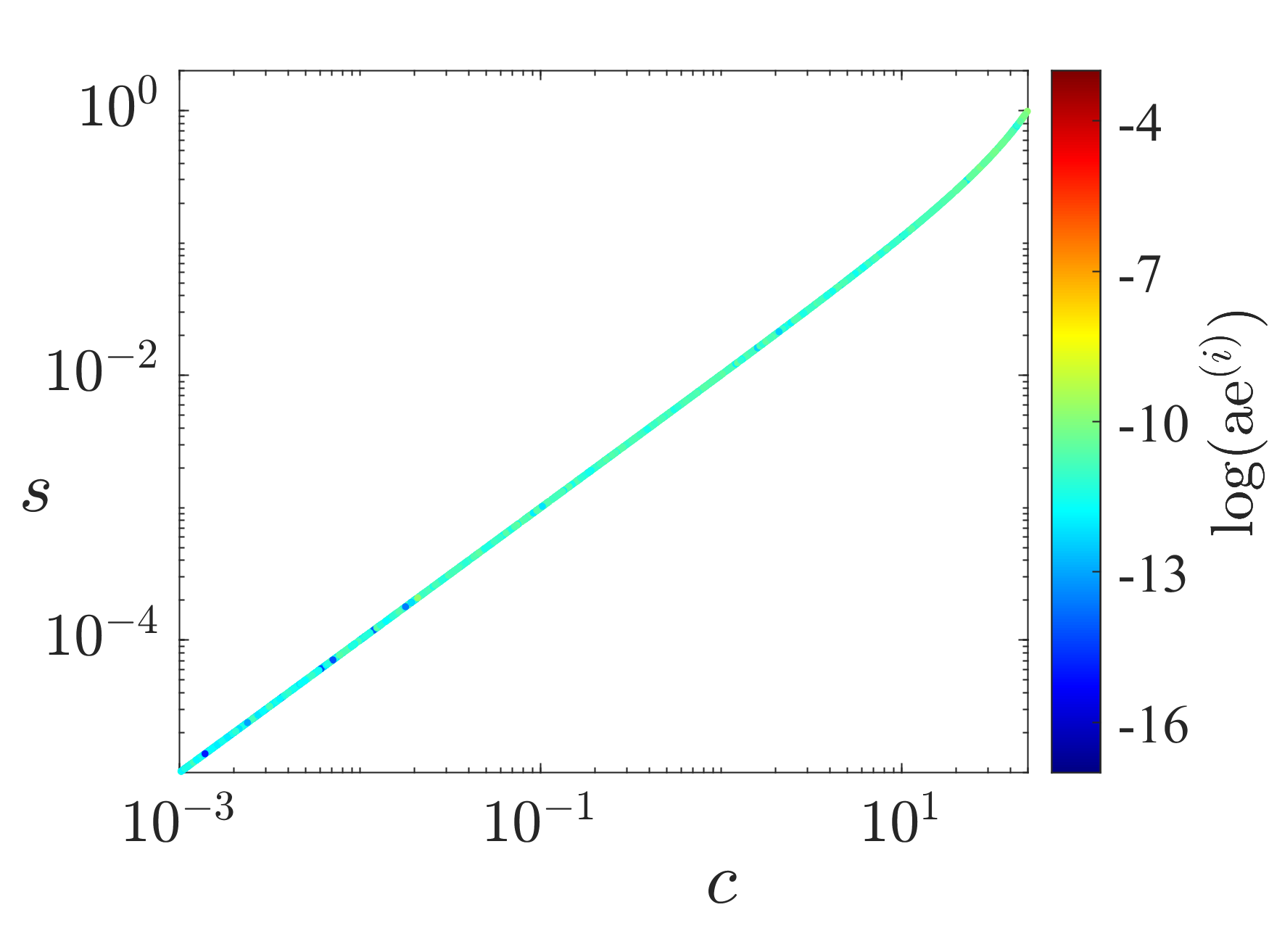}}
    \caption{MM system in Eq.~(32) for the MM2 case, where $x=s$ is the fast variable.~Absolute errors ($ae^{(i)}$) of the SIM approximations over all points of the test set, in comparison to the numerical solution $\mathbf{z}^{(i)}=[x^{(i)},y^{(i)}]^\top$ for $i=1,\ldots,n_t$.~Panel (a) depicts the $\lvert \mathbf{C} \mathbf{z}^{(i)} - \mathcal{N}(\mathbf{D} \mathbf{z}^{(i)}) \rvert$ of the PINN scheme, panels (b-h) depict $\lvert x^{(i)} - h(y^{(i)}) \rvert$ of the  PEA, CSP$_e$, rQSSA, CSP$_s$(1), CSP$_s$(2), sQSSA and CSP$_c$(1) explicit functionals w.r.t. $x$, and panel (i) depicts $\lvert x^{(i)} - \hat{x}^{(i)}) \rvert$ of the CSP$_c$(2) implicit functional, solved numerically with Newton for $x$.~Note that the approximations of the second/third row were constructed with $s$/$c$ assumed as fast variable.}
    \label{SF:MM2_AE}
\end{figure}

\begin{figure}[!h]
    \centering
    \subfigure[PIML]{
    \includegraphics[width=0.32\textwidth]{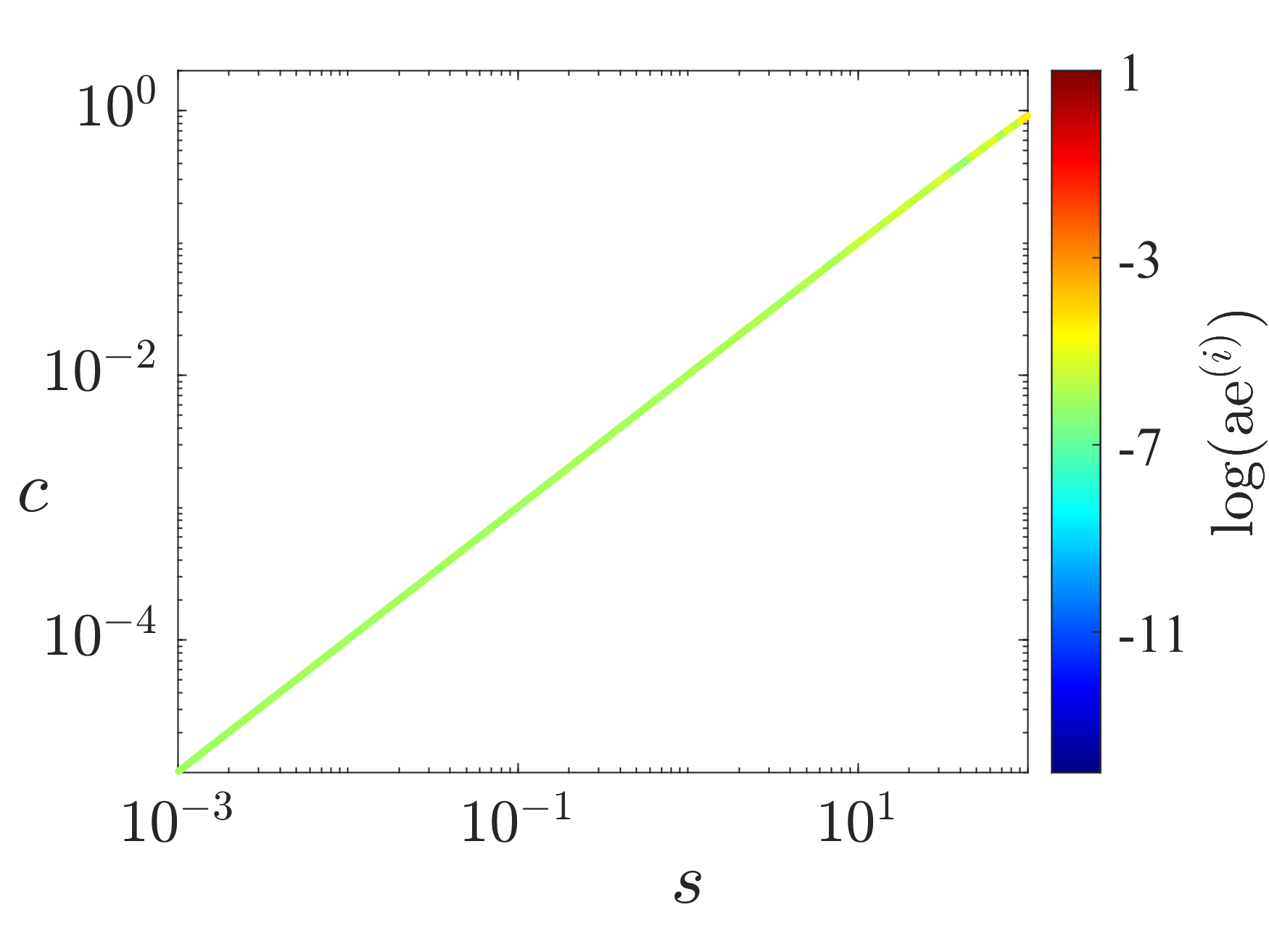}}
    \subfigure[PEA]{
    \includegraphics[width=0.32\textwidth]{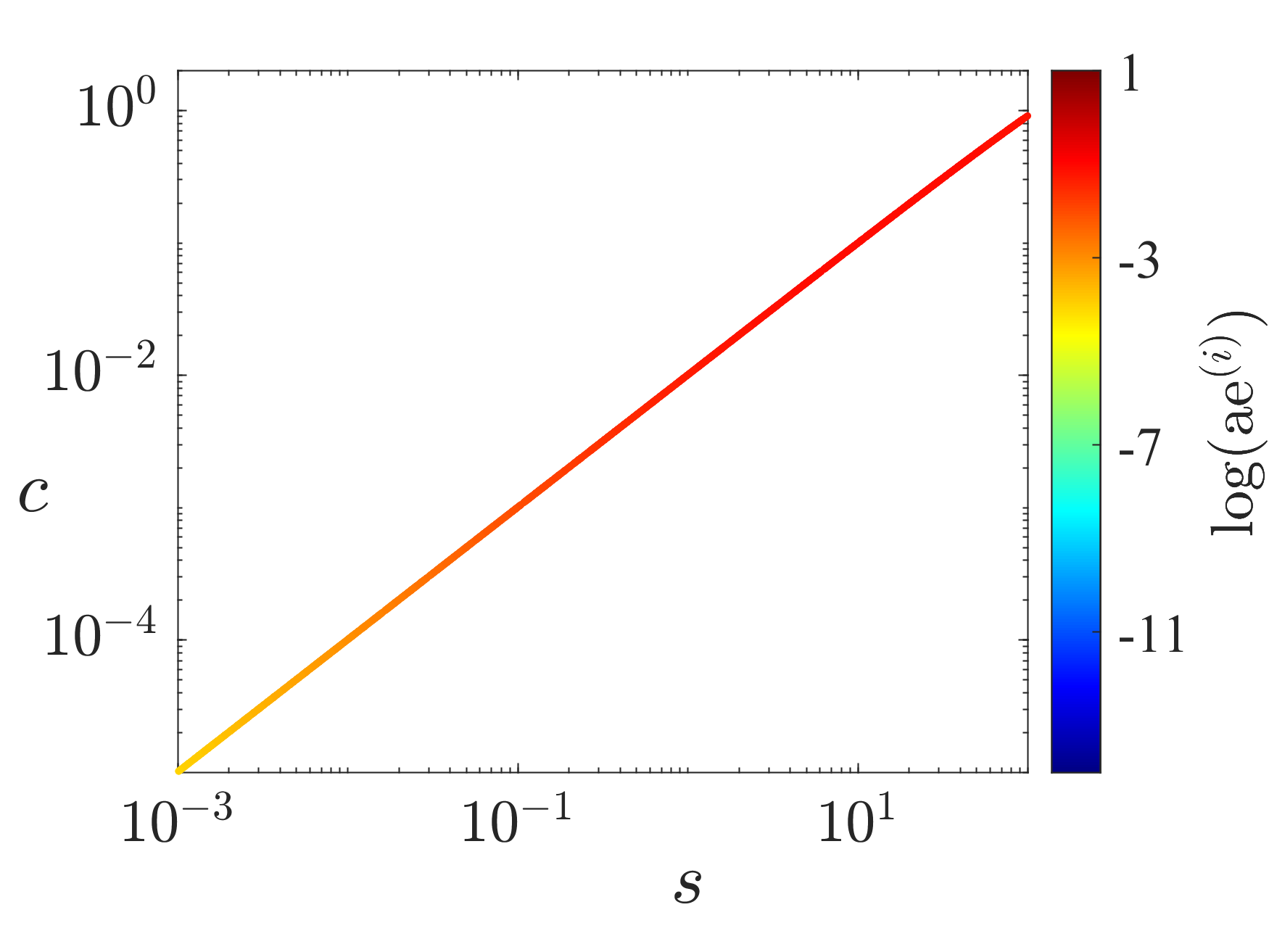}}
    \subfigure[CSP$_e$]{
    \includegraphics[width=0.32\textwidth]{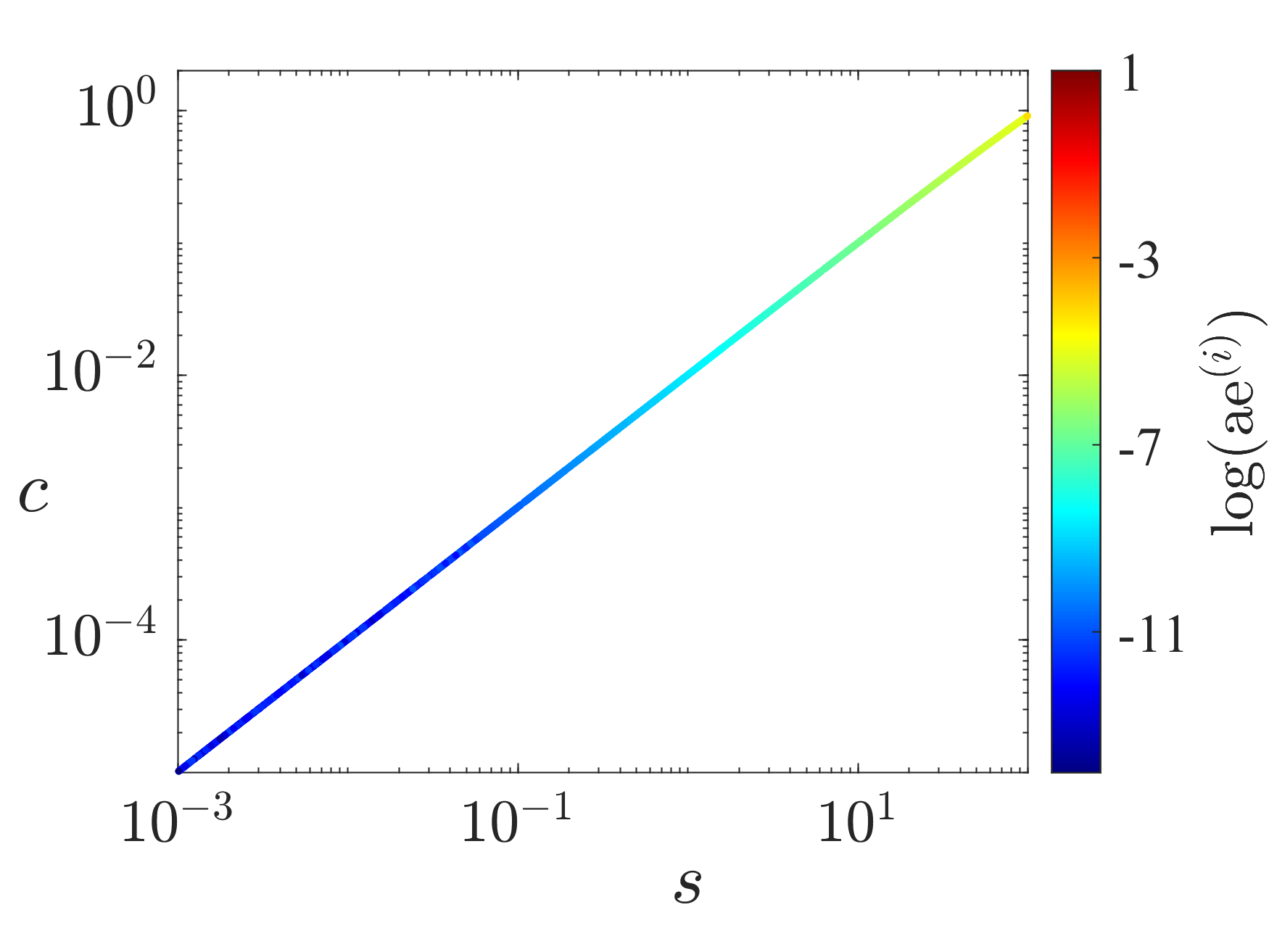}} \\ 
    \subfigure[rQSSA]{
    \includegraphics[width=0.32\textwidth]{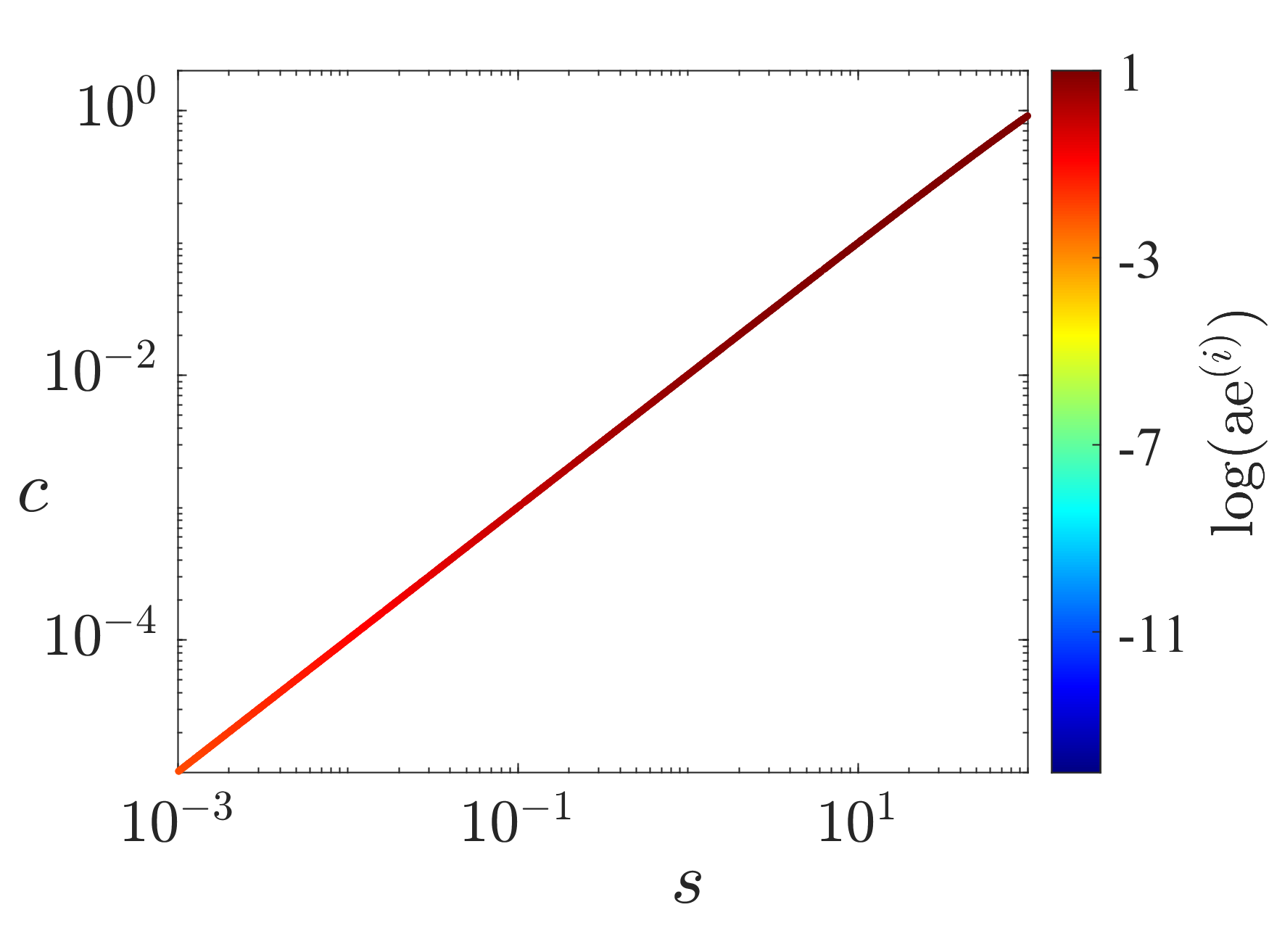}}
    \subfigure[CSP$_s$(1)]{
    \includegraphics[width=0.32\textwidth]{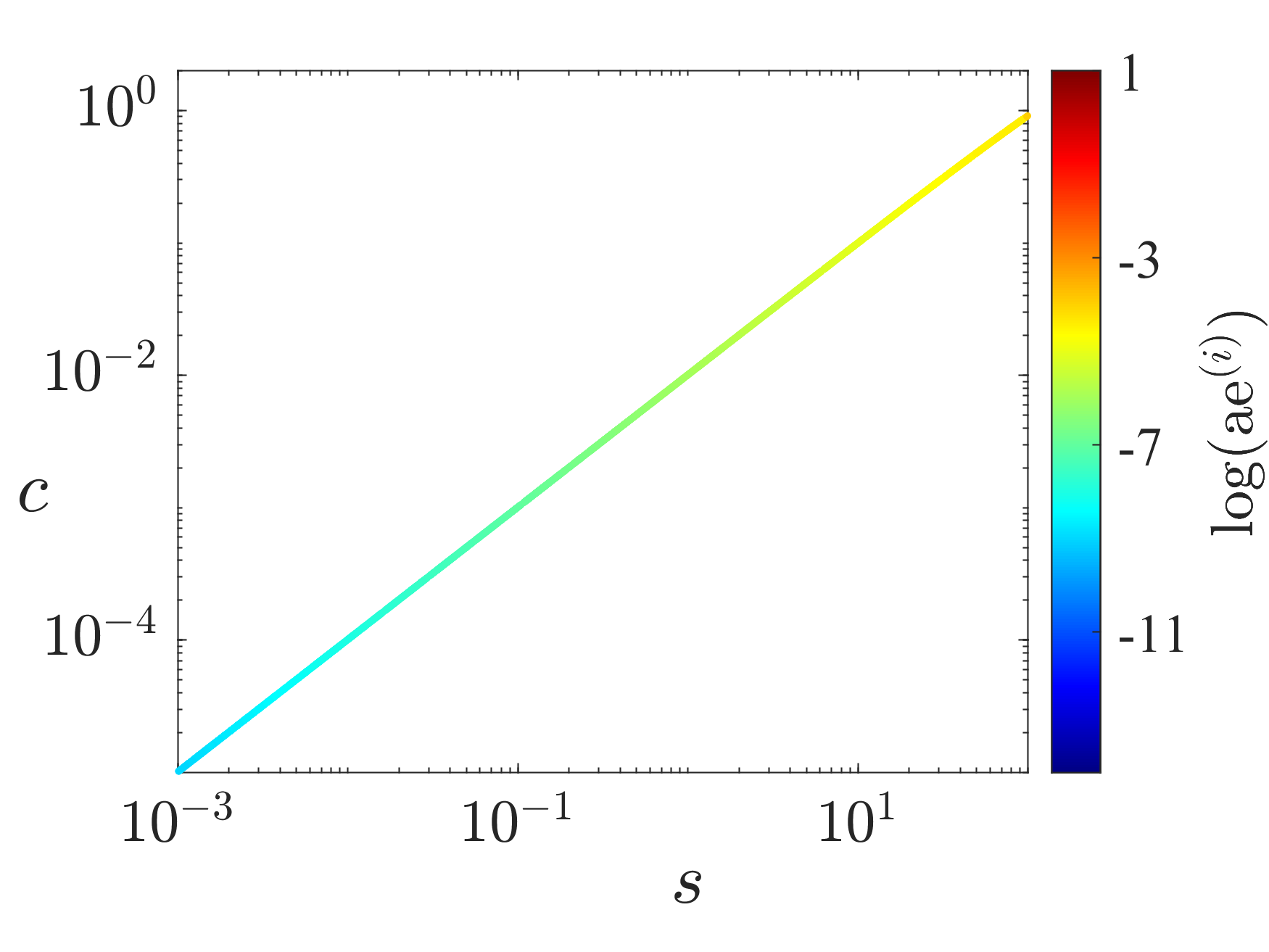}}
    \subfigure[CSP$_s$(2)]{
    \includegraphics[width=0.32\textwidth]{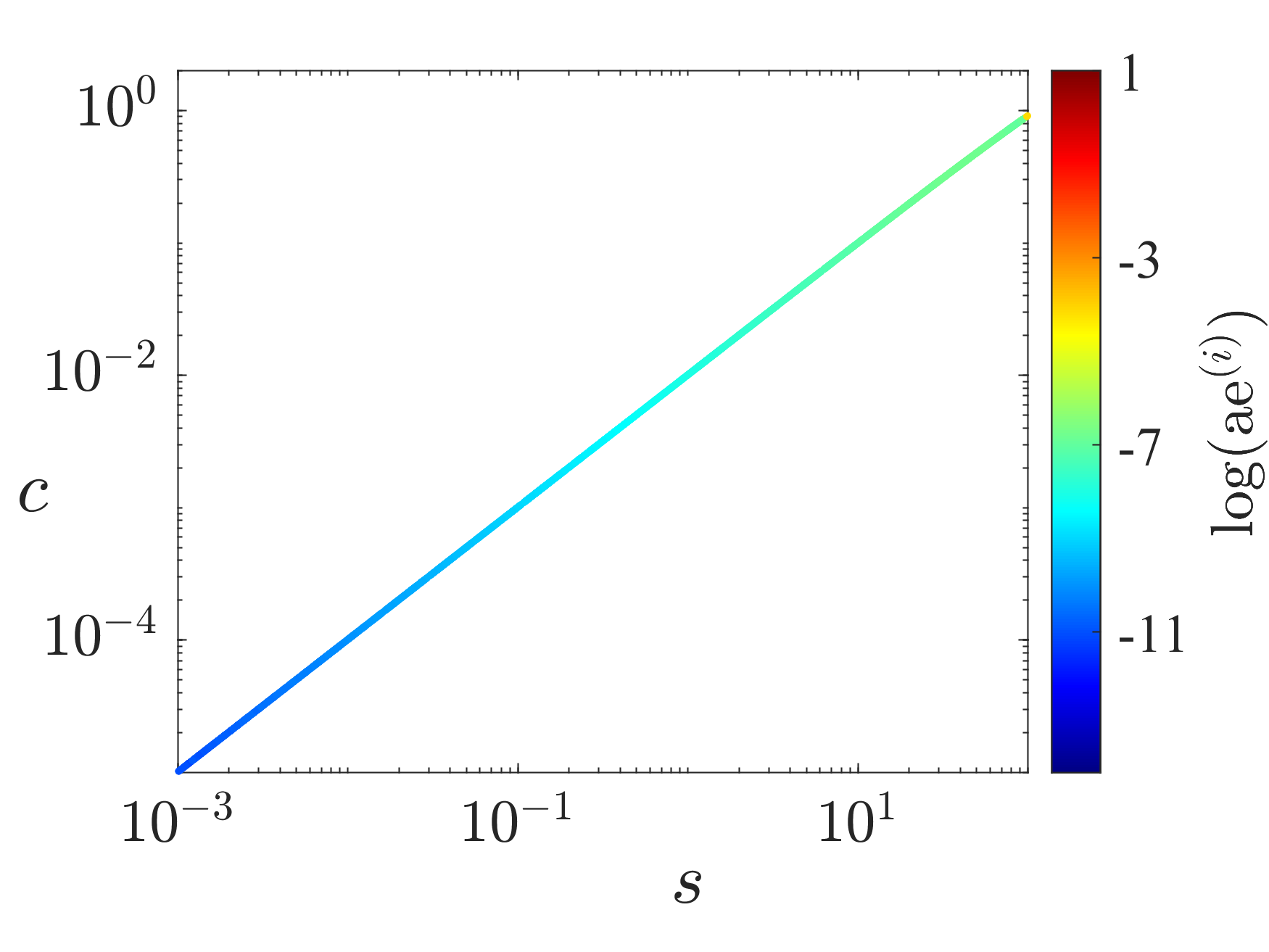}} \\
    \subfigure[sQSSA]{
    \includegraphics[width=0.32\textwidth]{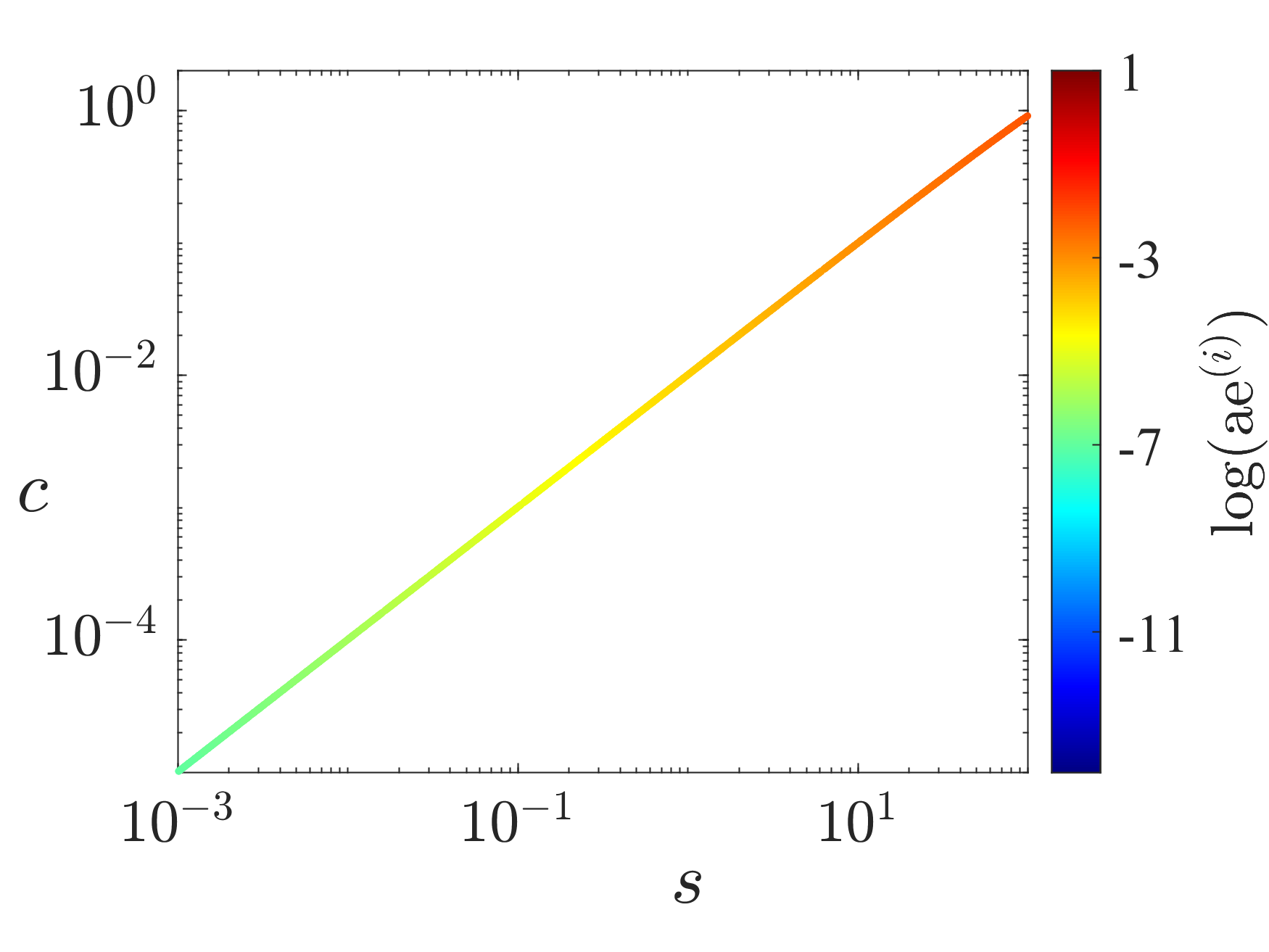}}
    \subfigure[CSP$_c$(1)]{
    \includegraphics[width=0.32\textwidth]{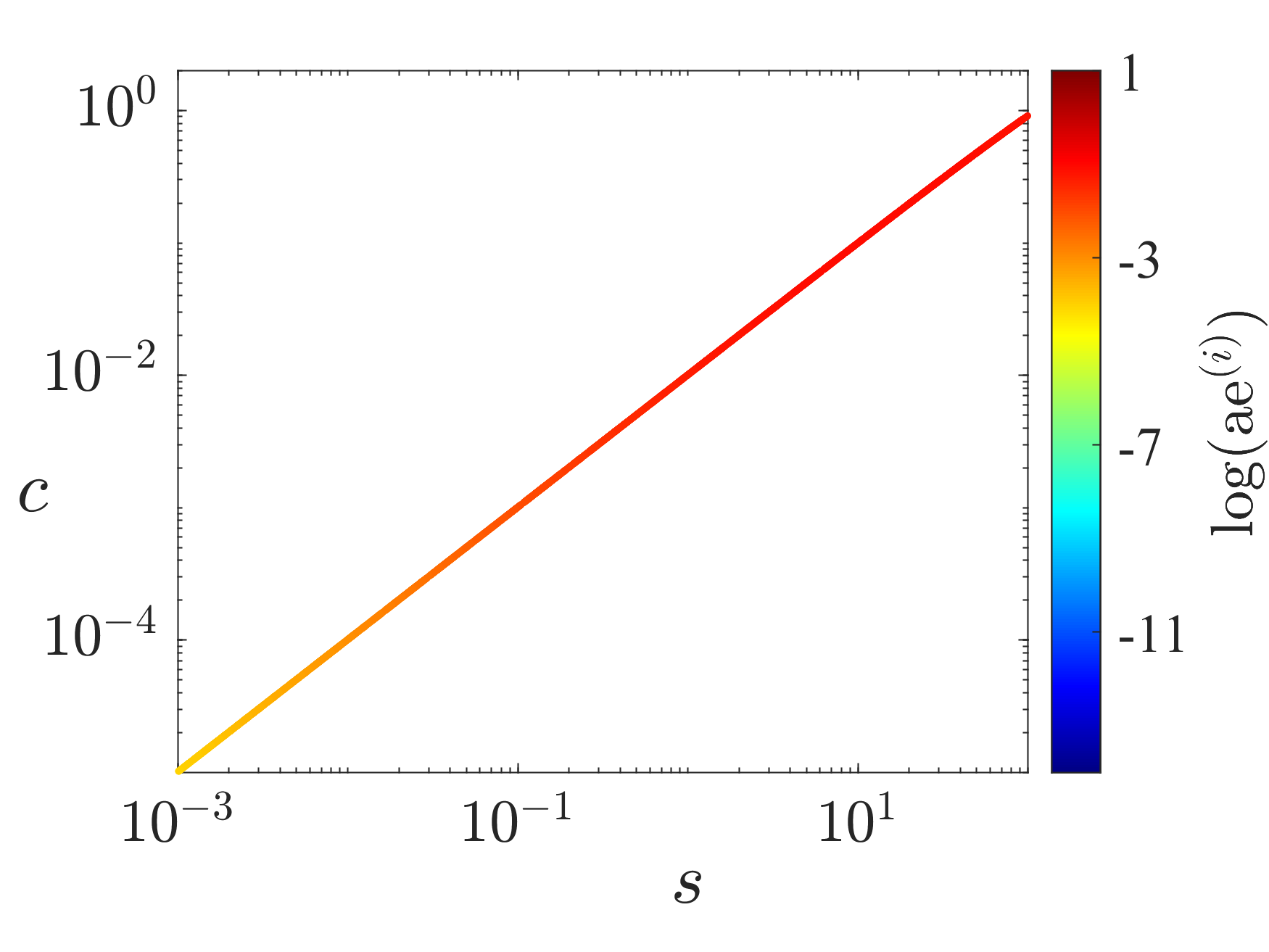}}
    \subfigure[CSP$_c$(2)]{
    \includegraphics[width=0.32\textwidth]{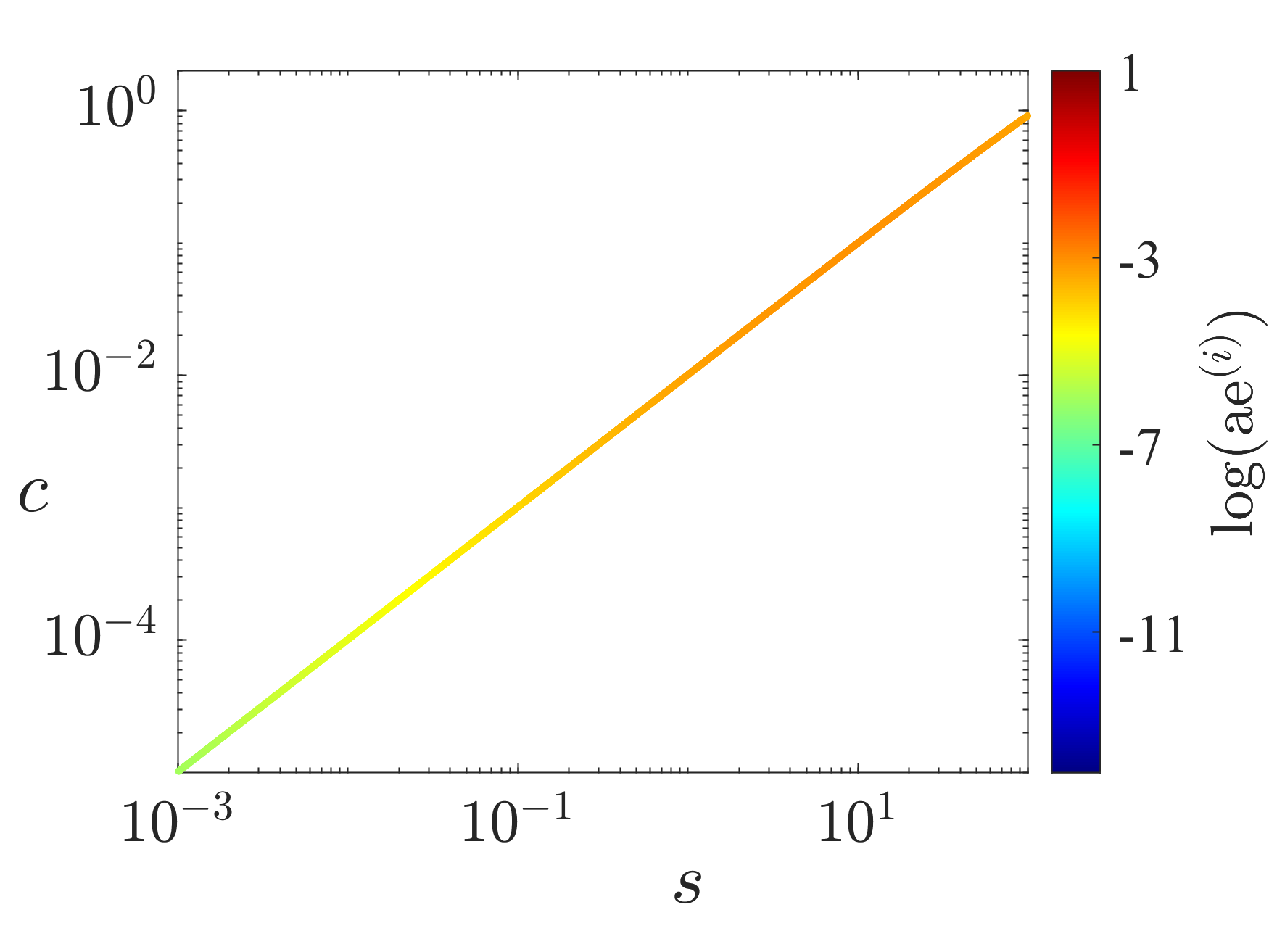}}
    \caption{MM system in Eq.~(32) for the MM3 case, where $x=c$ is the fast variable.~Absolute errors ($ae^{(i)}$) of the SIM approximations over all points of the test set, in comparison to the numerical solution $\mathbf{z}^{(i)}=[x^{(i)},y^{(i)}]^\top$ for $i=1,\ldots,n_t$.~Panel (a) depicts the $\lvert \mathbf{C} \mathbf{z}^{(i)} - \mathcal{N}(\mathbf{D} \mathbf{z}^{(i)}) \rvert$ of the PINN scheme, panels (b), (d), (e), (g) and (h) depict $\lvert x^{(i)} - h(y^{(i)}) \rvert$ of the PEA, rQSSA,  CSP$_s$(1), sQSSA and CSP$_c$(1) explicit functionals w.r.t. $x$, and panels (c), (f) and (i) depict $\lvert x^{(i)} - \hat{x}^{(i)}) \rvert$ of the CSP$_e$, CSP$_s$(2) and CSP$_c$(2) implicit functionals, solved numerically with Newton for $x$.~Note that the approximations of the second/third row were constructed with $s$/$c$ assumed as fast variable.}
    \label{SF:MM3_AE}
\end{figure}

\begin{figure}[!h]
    \centering
    \subfigure[PIML]{
    \includegraphics[width=0.32\textwidth]{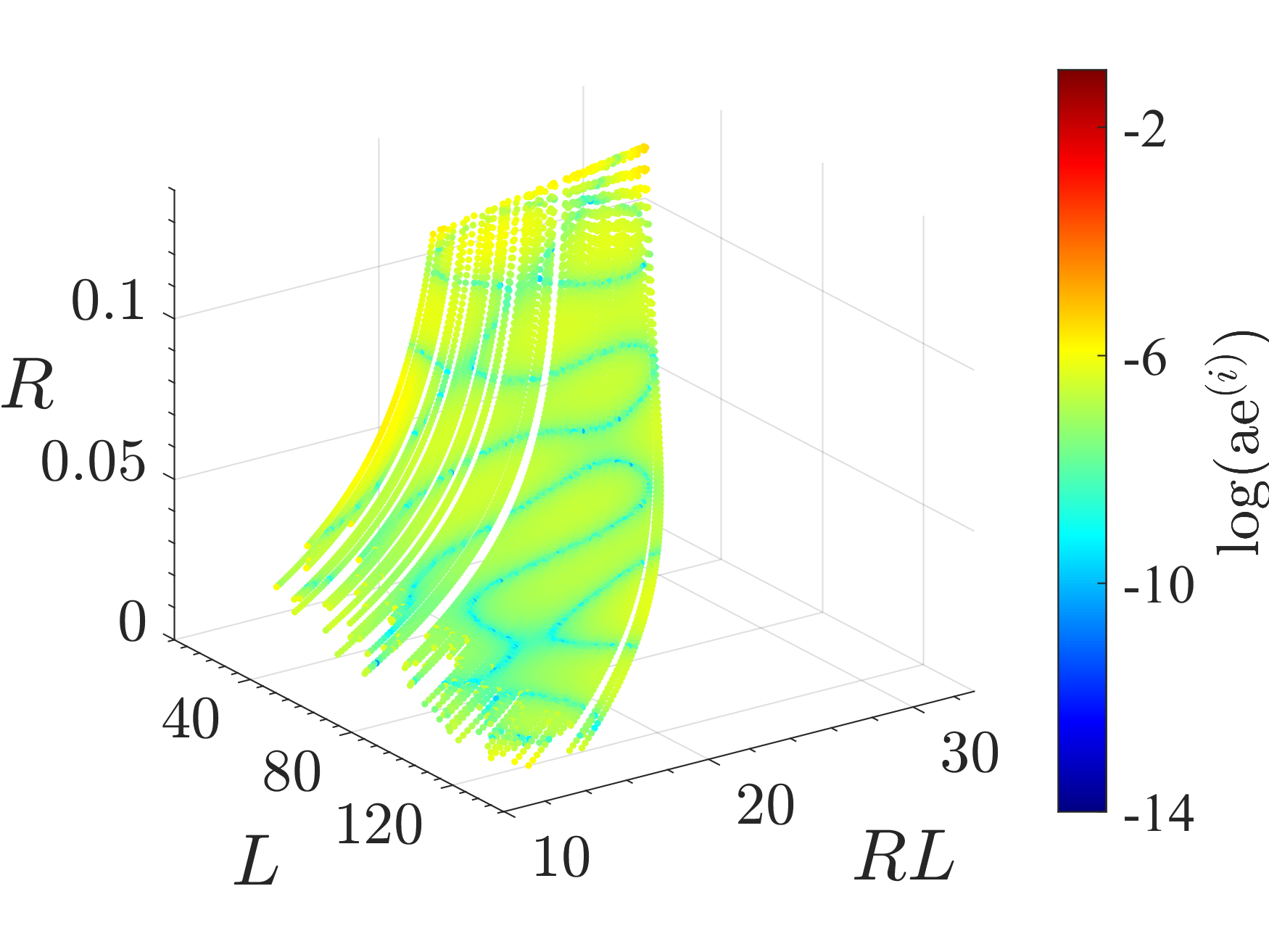}}
    \subfigure[PEA]{
    \includegraphics[width=0.32\textwidth]{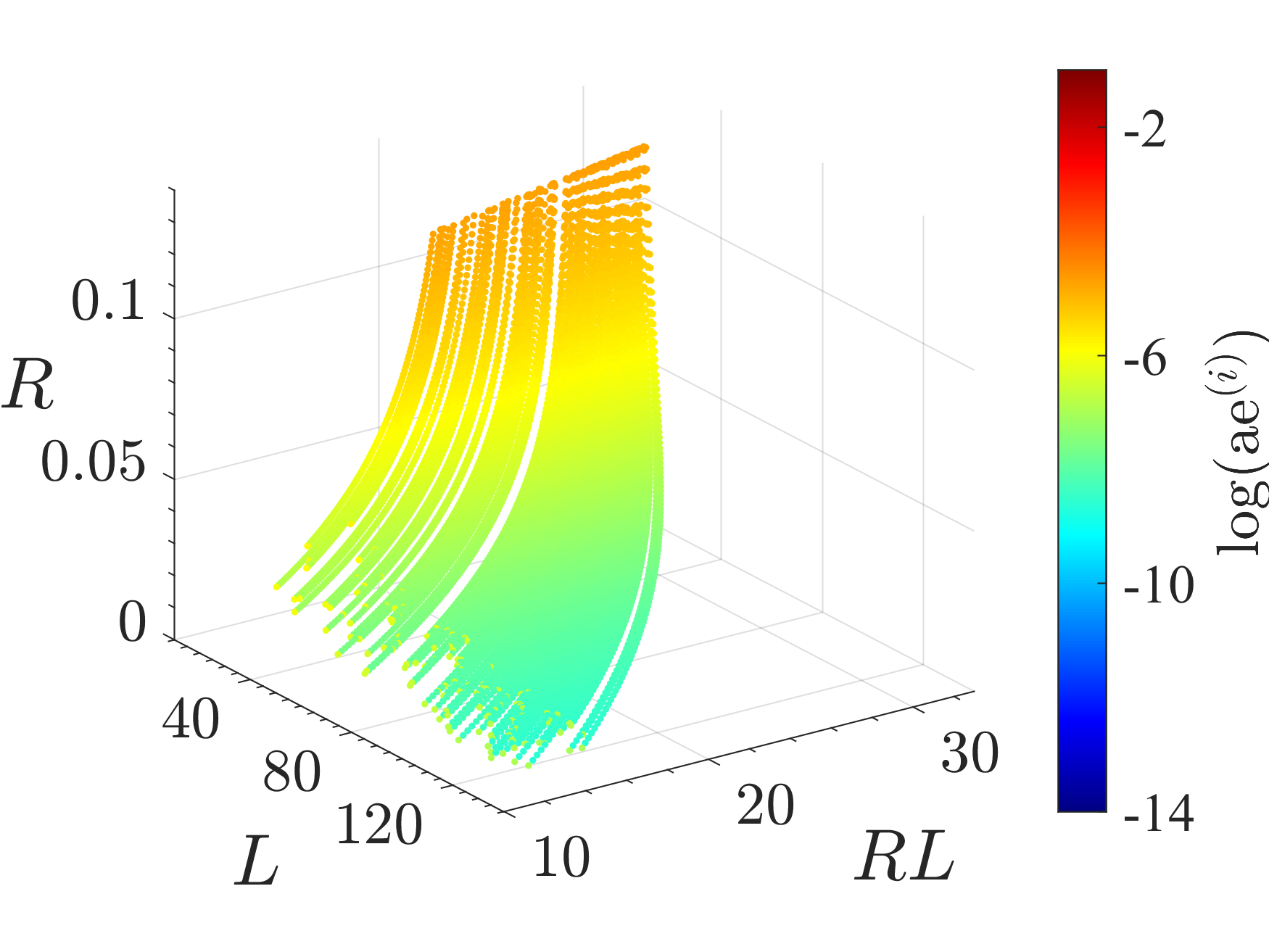}} \\
    \subfigure[QSSA$_L$]{
    \includegraphics[width=0.32\textwidth]{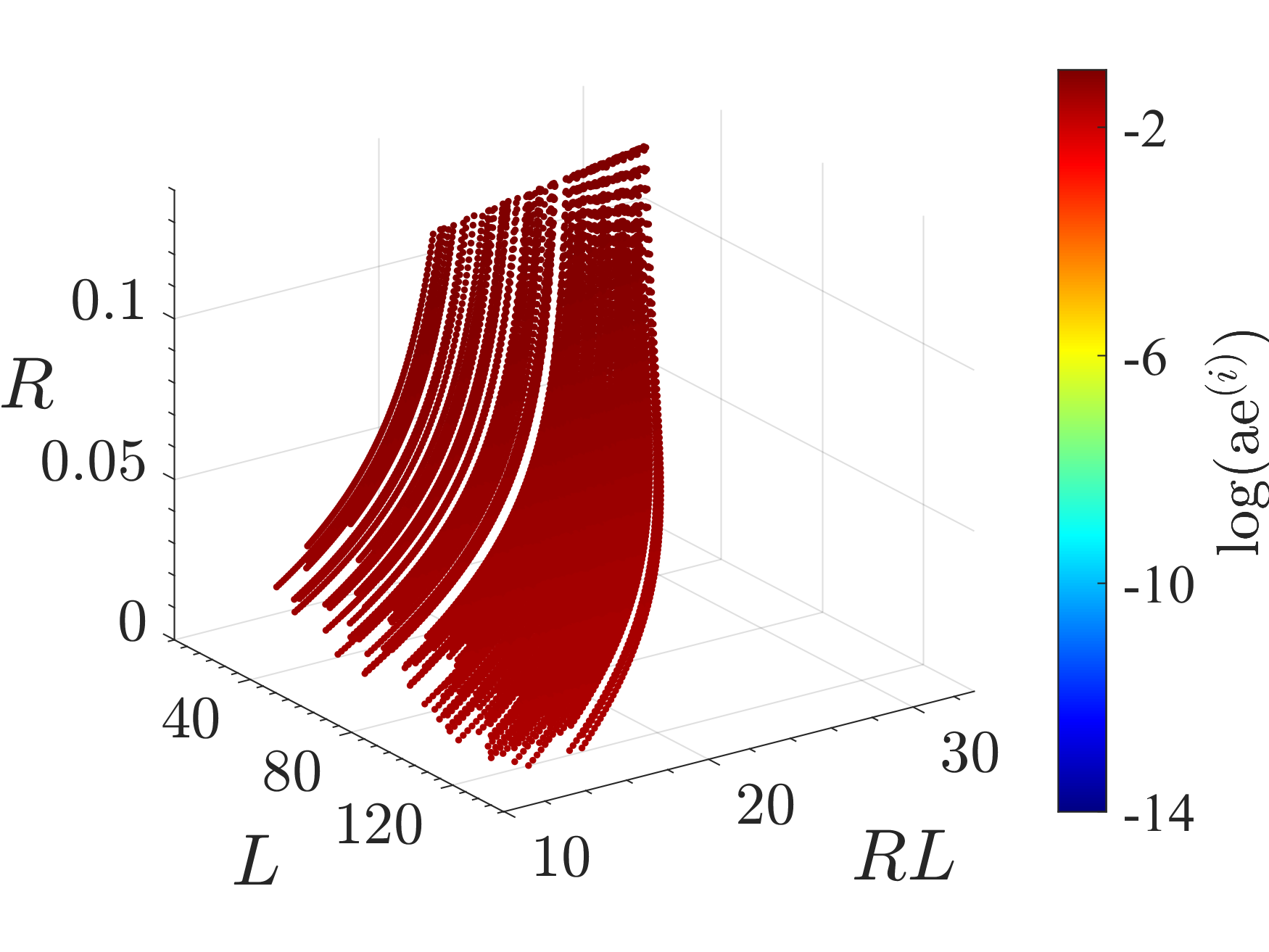}}
    \subfigure[CSP$_L$(1)]{
    \includegraphics[width=0.32\textwidth]{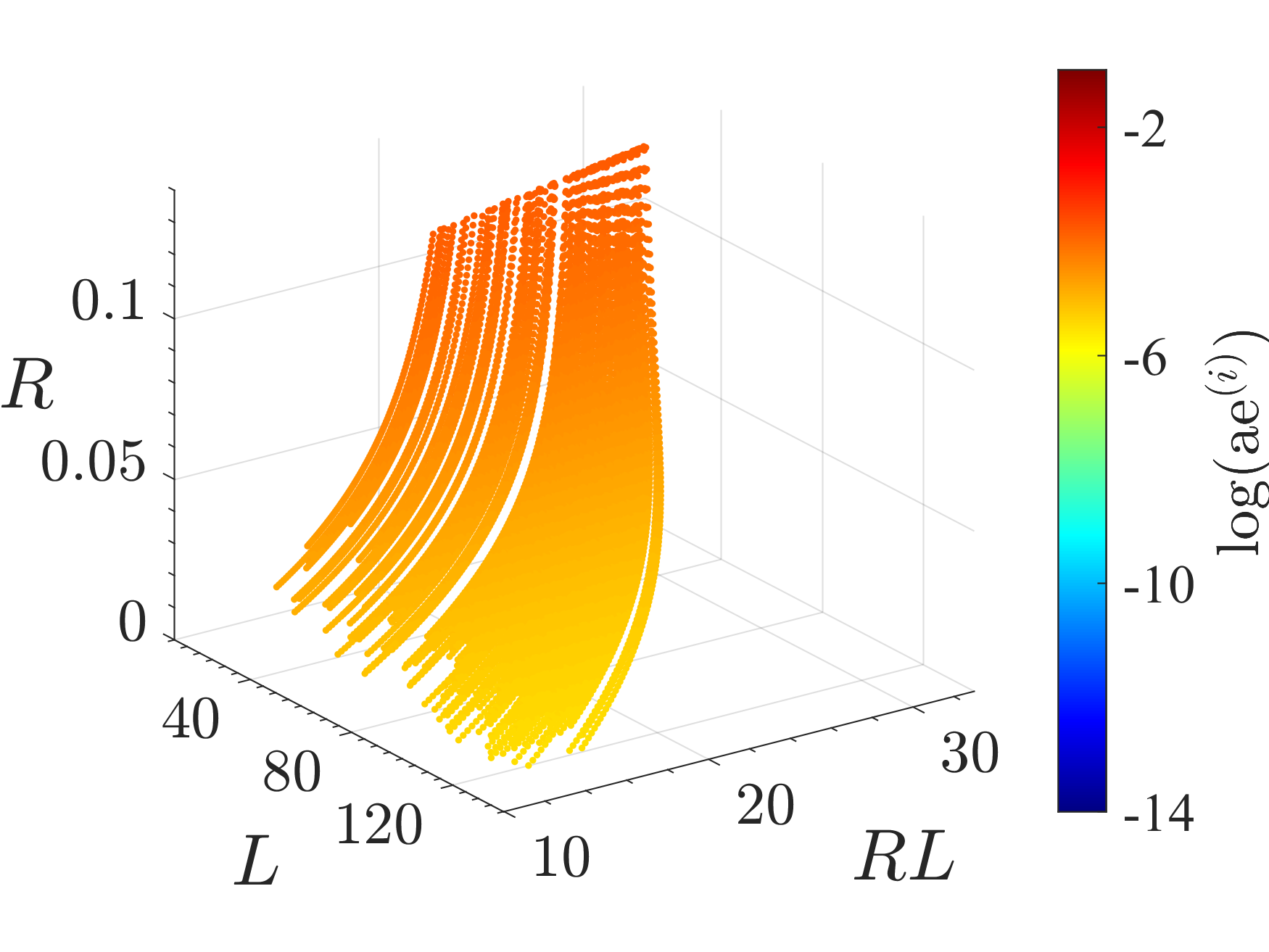}} 
    \subfigure[CSP$_L$(2)]{
    \includegraphics[width=0.32\textwidth]{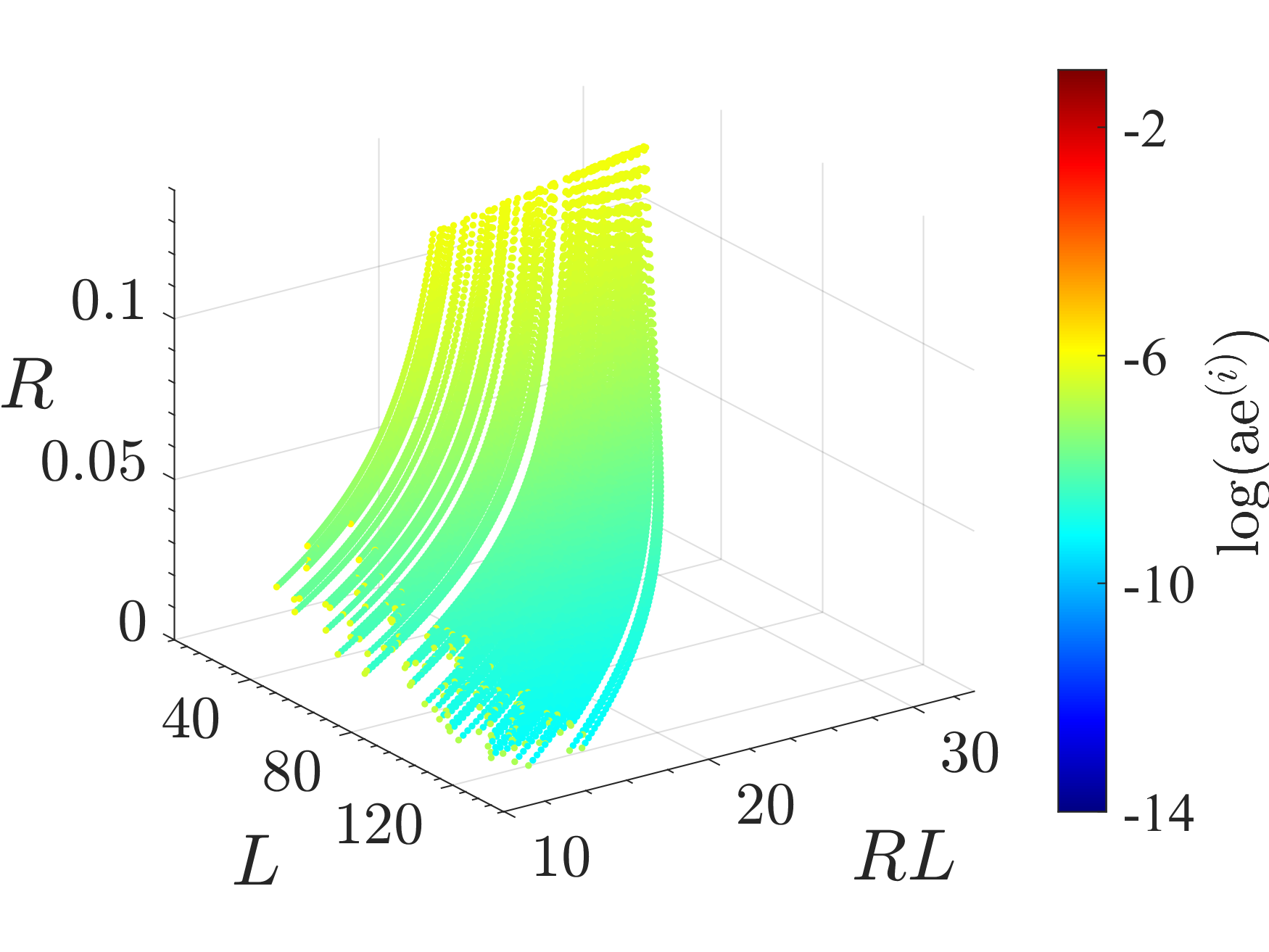}} \\
    \subfigure[QSSA$_R$]{
    \includegraphics[width=0.32\textwidth]{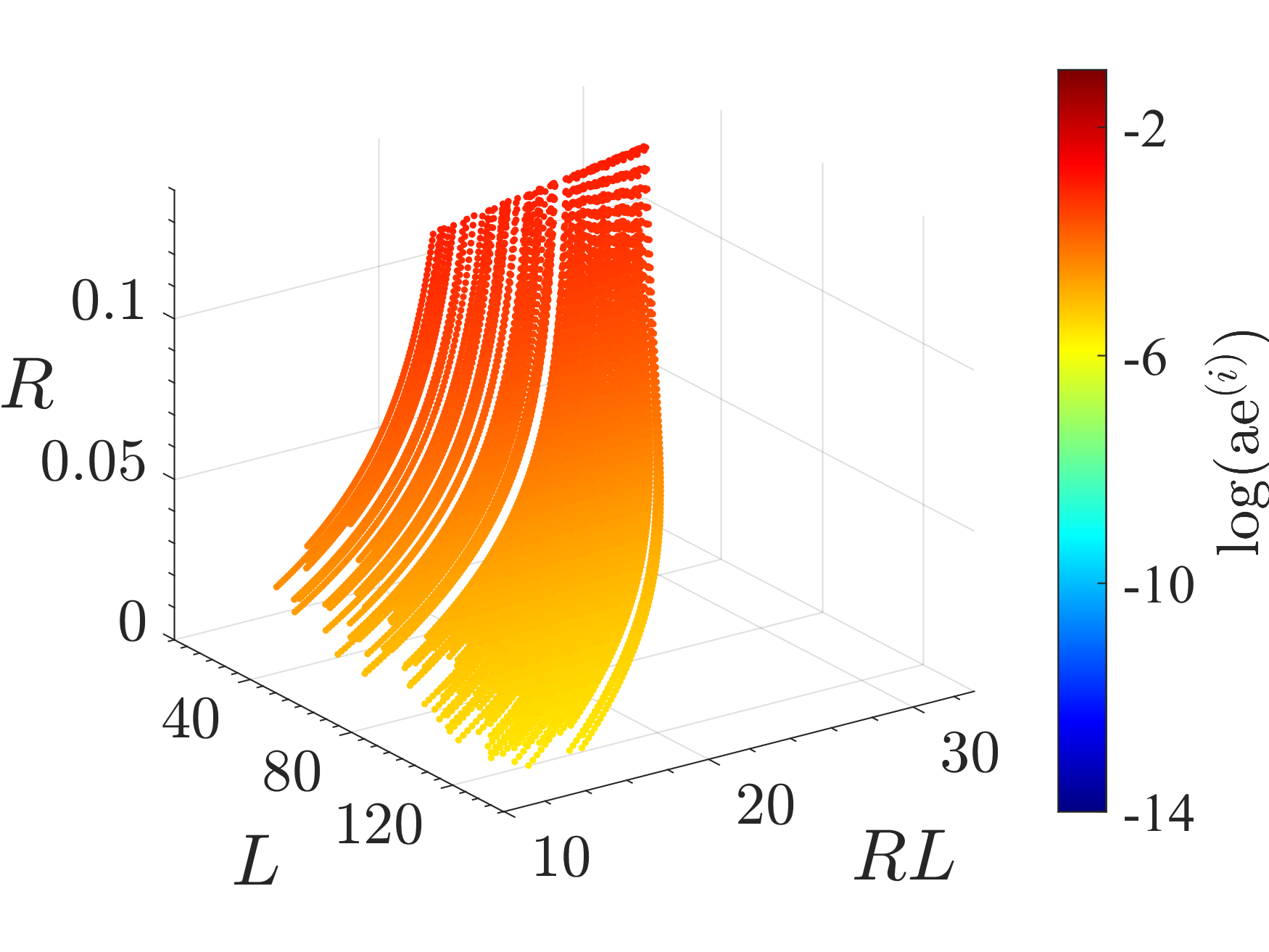}} 
    \subfigure[CSP$_R$(1)]{
    \includegraphics[width=0.32\textwidth]{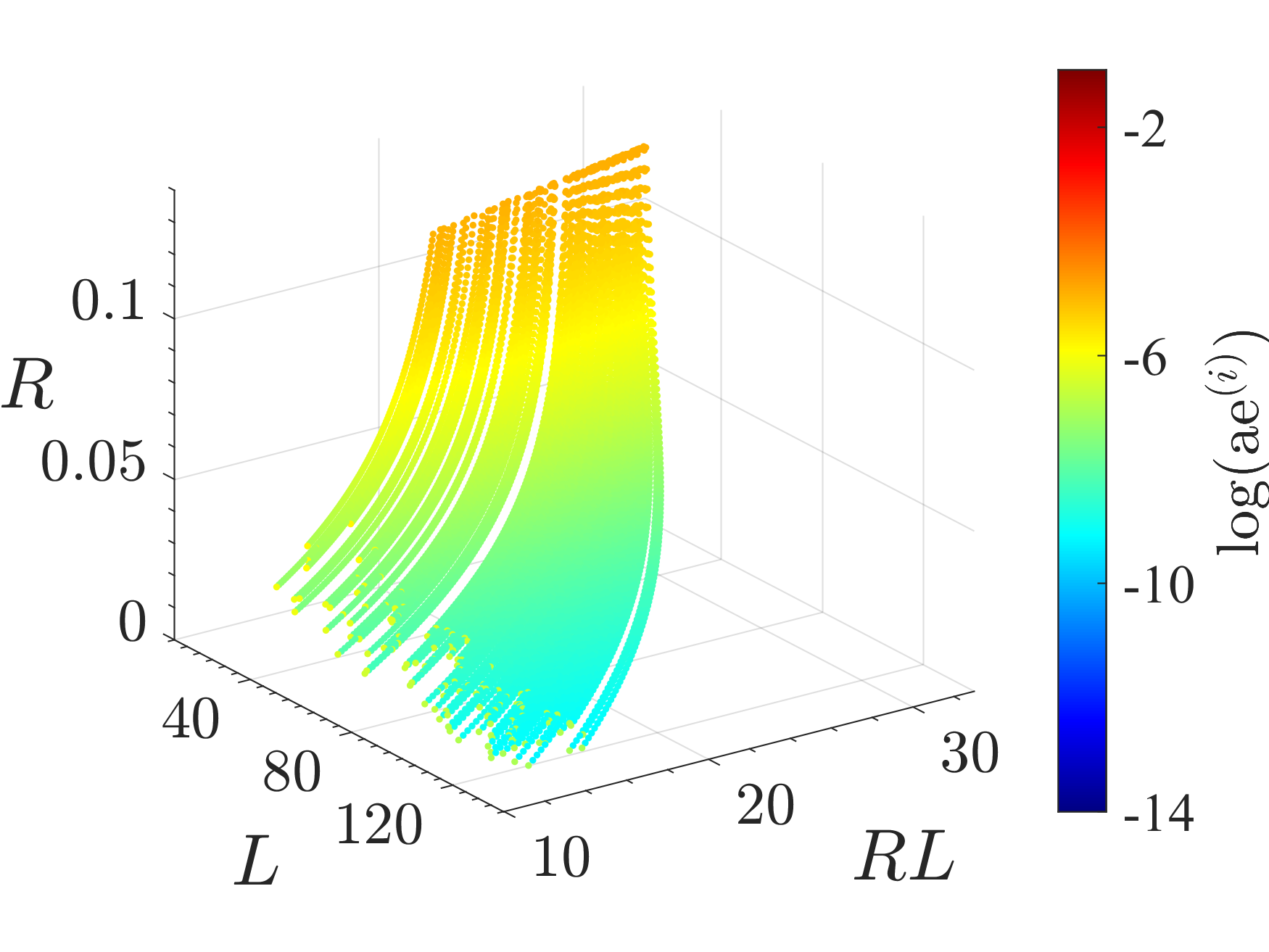}}
    \subfigure[CSP$_R$(2)]{
    \includegraphics[width=0.32\textwidth]{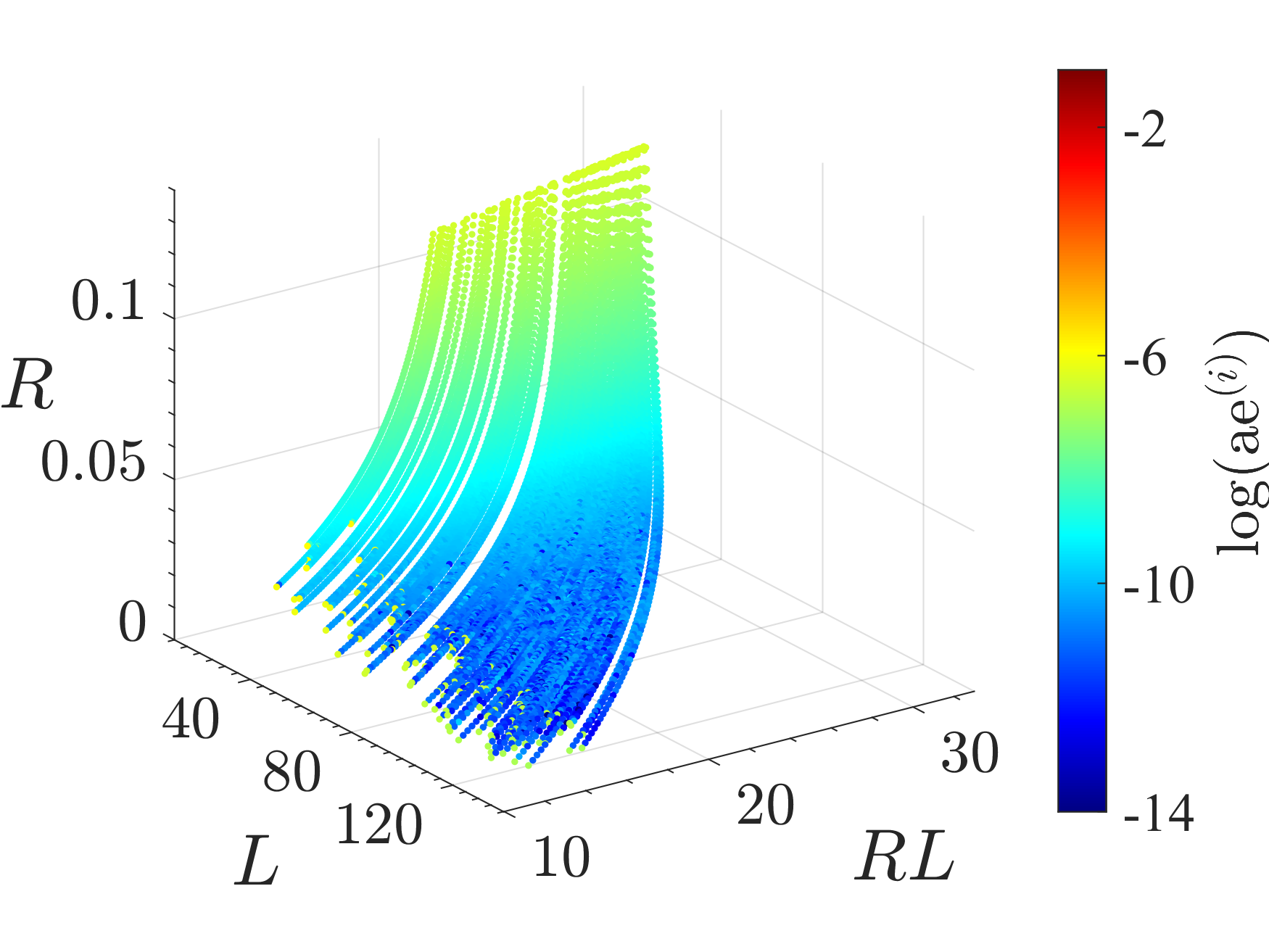}} \\
    \subfigure[QSSA$_{RL}$]{
    \includegraphics[width=0.32\textwidth]{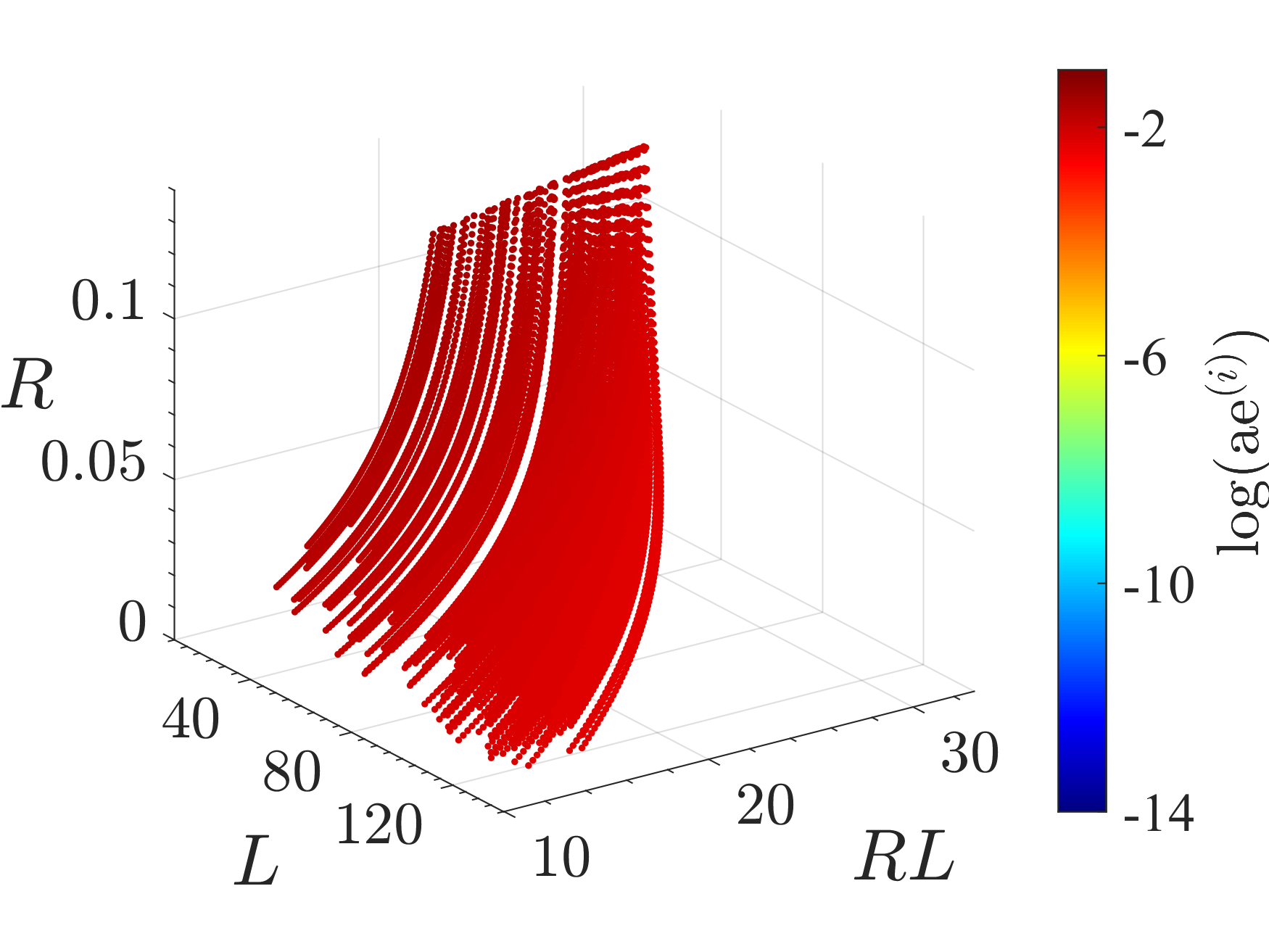}}
    \subfigure[CSP$_{RL}$(1)]{
    \includegraphics[width=0.32\textwidth]{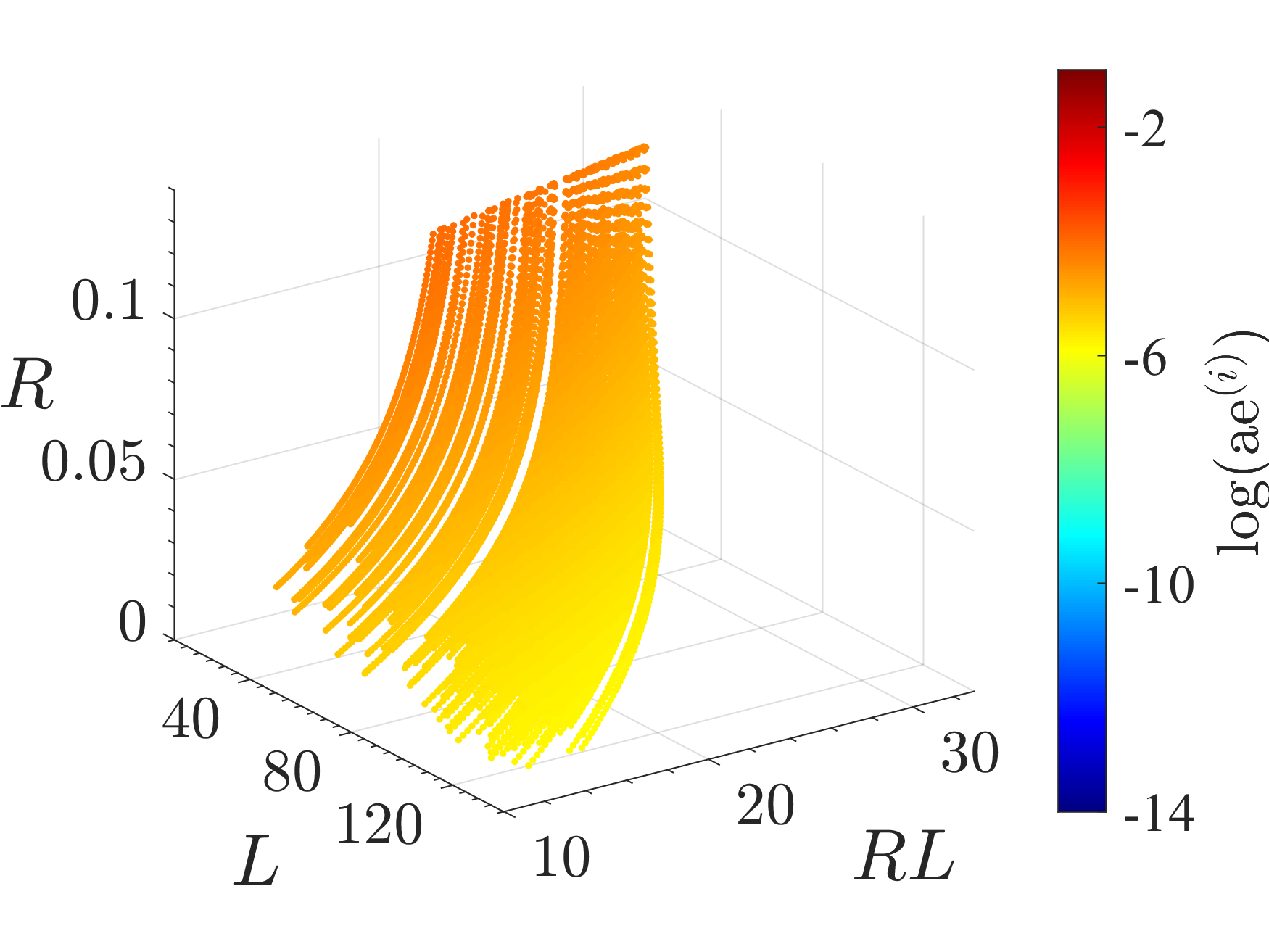}}
    \subfigure[CSP$_{RL}$(2)]{
    \includegraphics[width=0.32\textwidth]{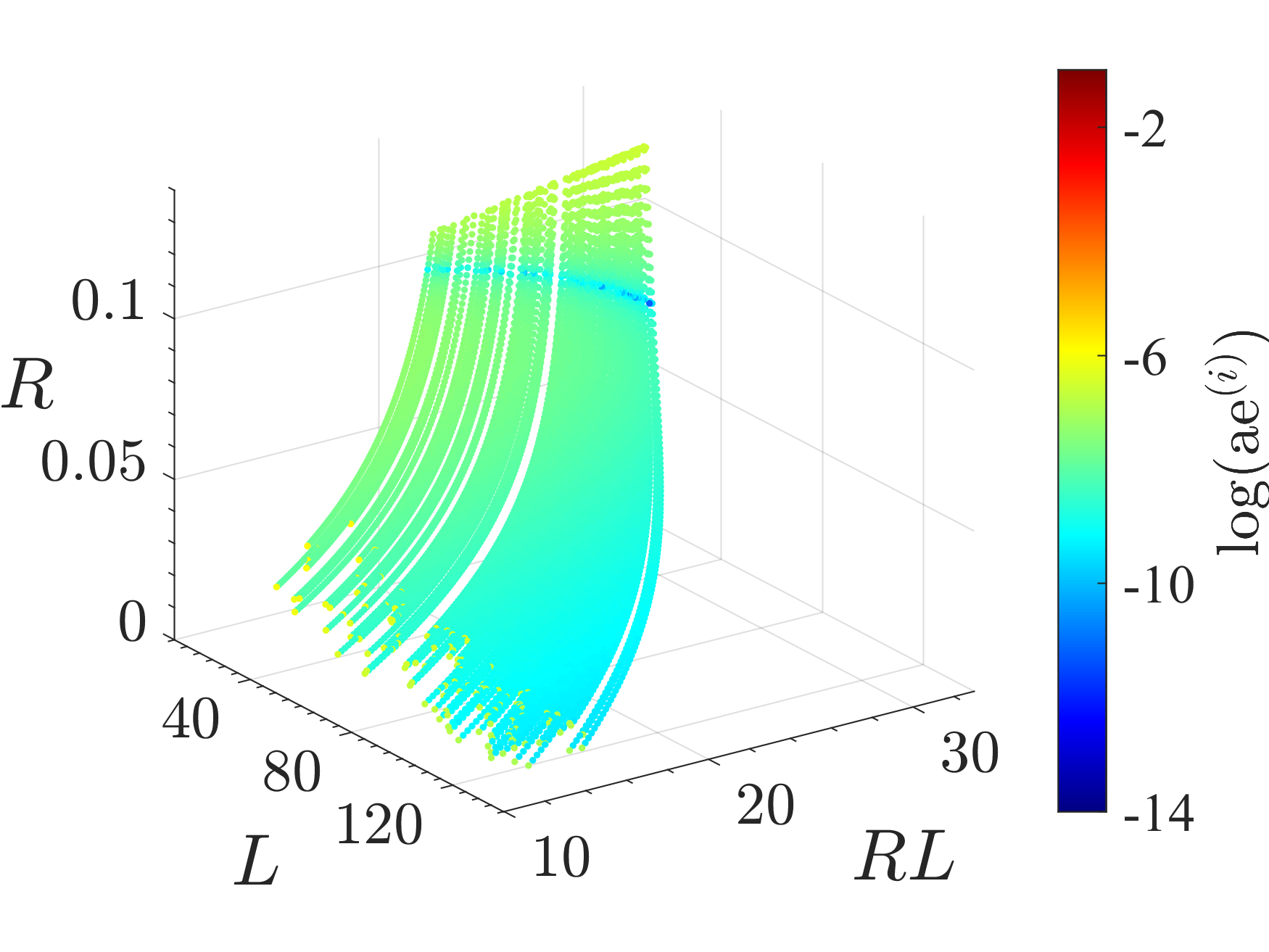}}
    \caption{TMDD system in Eq.~(33) for the SIM $\mathcal{M}_1$, where $x=R$ is the fast variable.~Absolute errors ($ae^{(i)}$) of the SIM approximations over all points of the test set, in comparison to the numerical solution $\mathbf{z}^{(i)}=[x^{(i)},\mathbf{y}^{(i)}]^\top$ for $i=1,\ldots,n_t$.~Panel (a) depicts the $\lvert \mathbf{C} \mathbf{z}^{(i)} - \mathcal{N}(\mathbf{D} \mathbf{z}^{(i)}) \rvert$ of the PINN scheme, panels (b), (c), (d), (f), (g), (i) and (j) depict $\lvert x^{(i)} - h(\mathbf{y}^{(i)}) \rvert$ of the PEA, QSSA$_L$,  CSP$_L$(1), QSSA$_R$,  CSP$_R$(1), QSSA$_{RL}$ and  CSP$_{RL}$(1) explicit functionals w.r.t. $x$, and panels (e), (h) and (k) depict $\lvert x^{(i)} - \hat{x}^{(i)}) \rvert$ of the CSP$_L$(2), CSP$_R$(2) and CSP$_{RL}$(2) implicit functionals, solved numerically with Newton for $x$.~Note that the approximations of the second/third/fourth row were constructed with $L$/$R$/$RL$ assumed as fast variable.}
    \label{SF:TMDDP2_AE}
\end{figure}

\begin{figure}[!h]
    \centering
    \subfigure[PIML]{
    \includegraphics[width=0.32\textwidth]{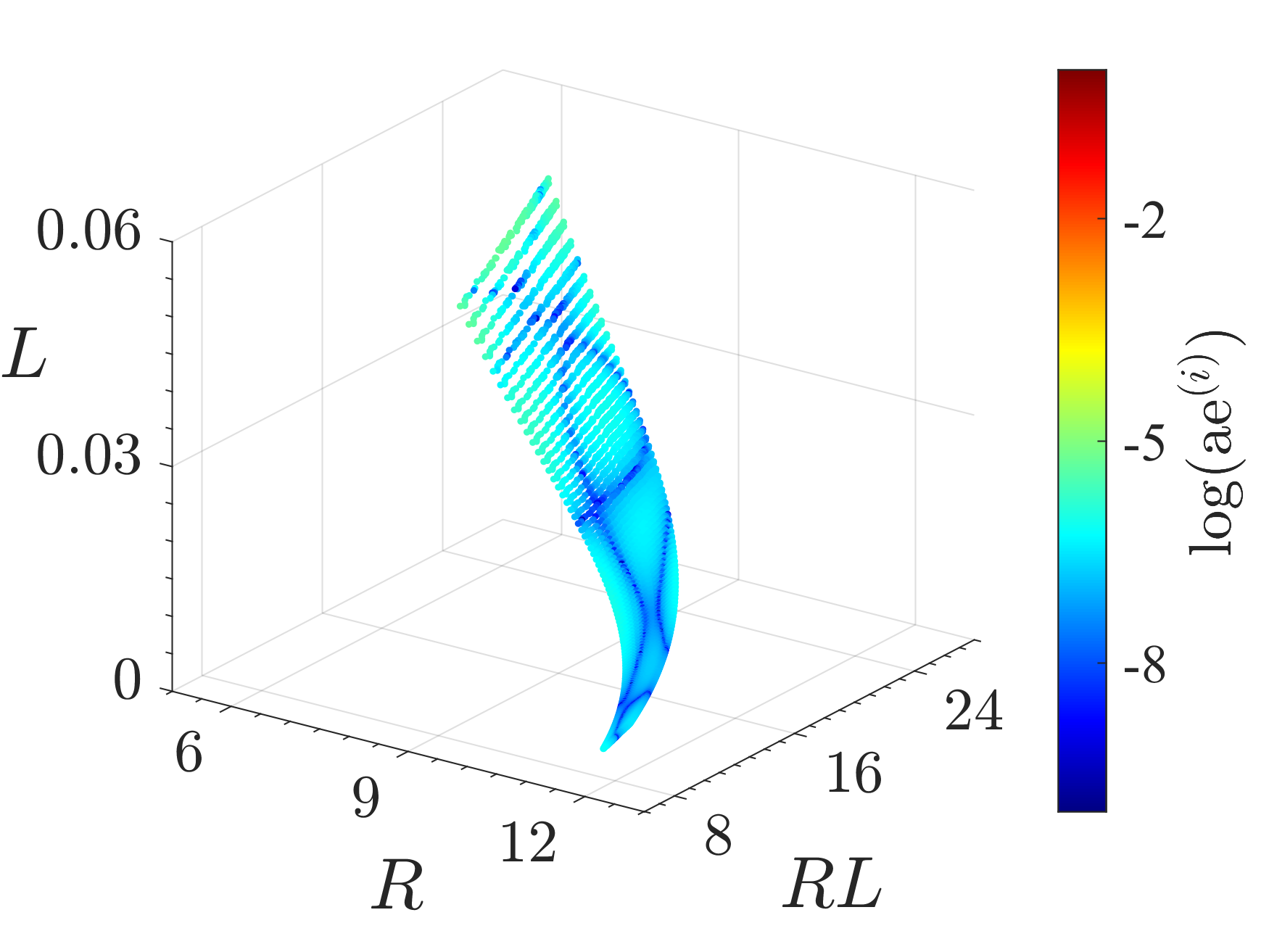}}
    \subfigure[PEA]{
    \includegraphics[width=0.32\textwidth]{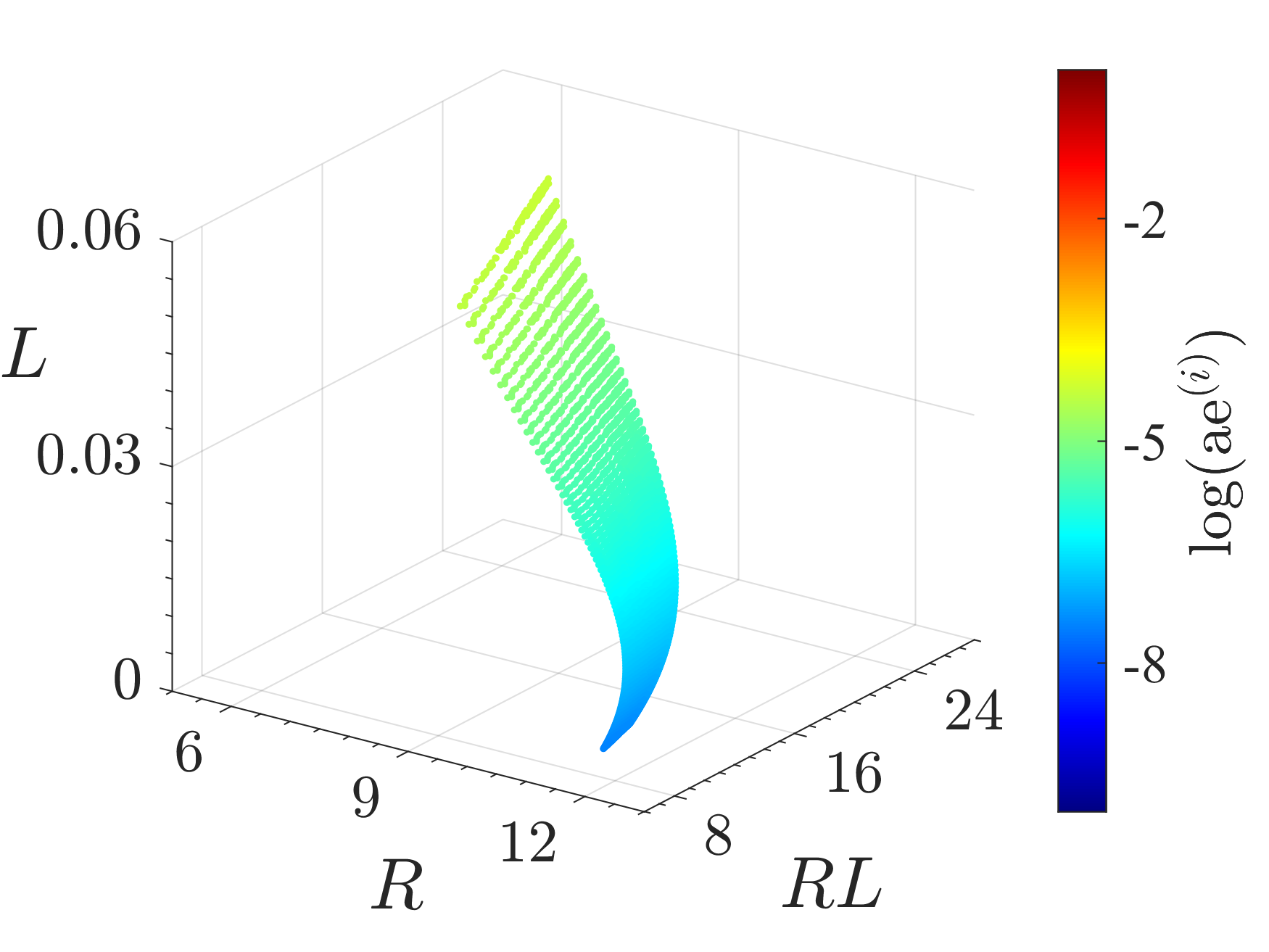}} \\
    \subfigure[QSSA$_L$]{
    \includegraphics[width=0.32\textwidth]{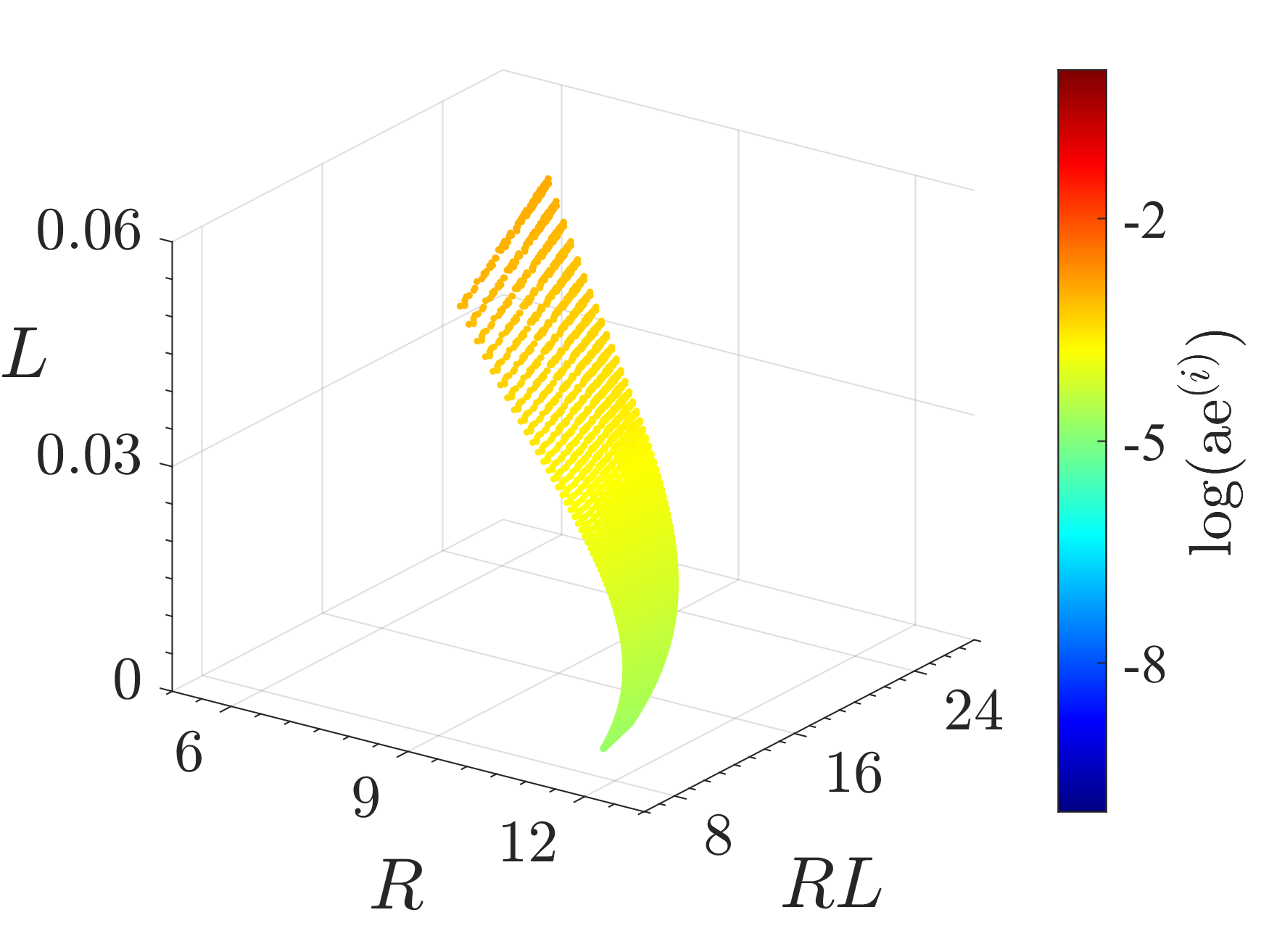}}
    \subfigure[CSP$_L$(1)]{
    \includegraphics[width=0.32\textwidth]{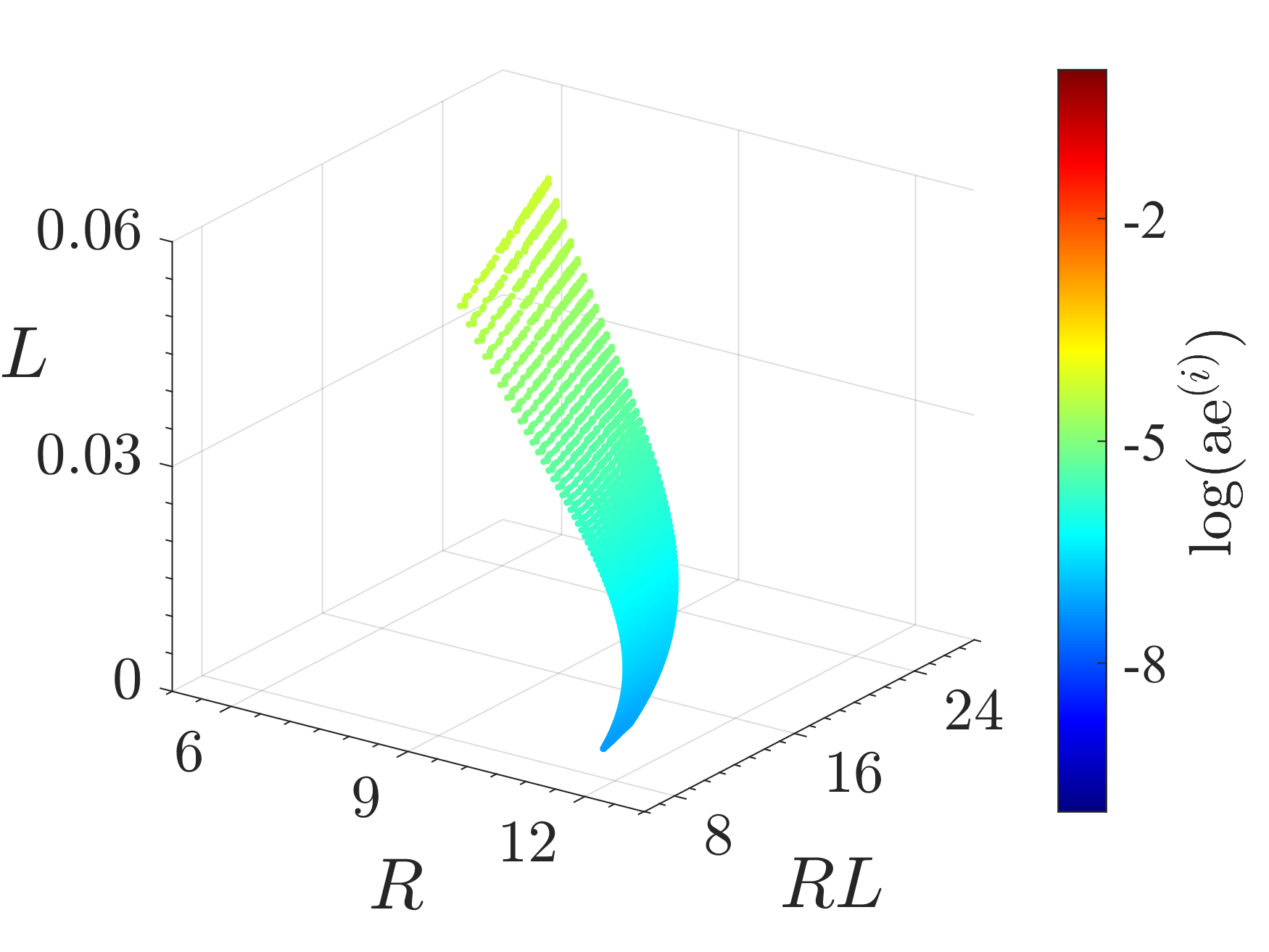}} 
    \subfigure[CSP$_L$(2)]{
    \includegraphics[width=0.32\textwidth]{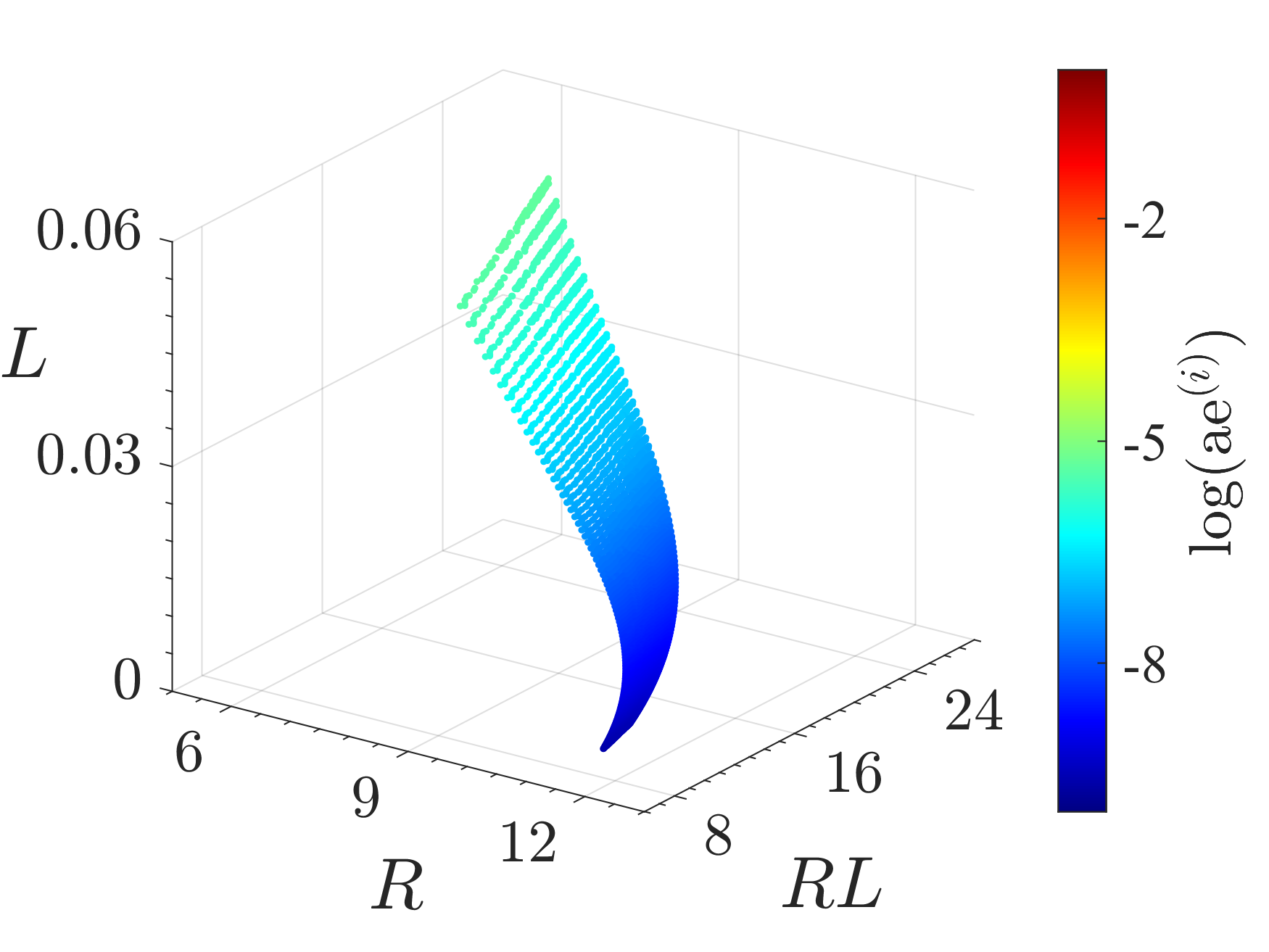}} \\
    \subfigure[QSSA$_R$]{
    \includegraphics[width=0.32\textwidth]{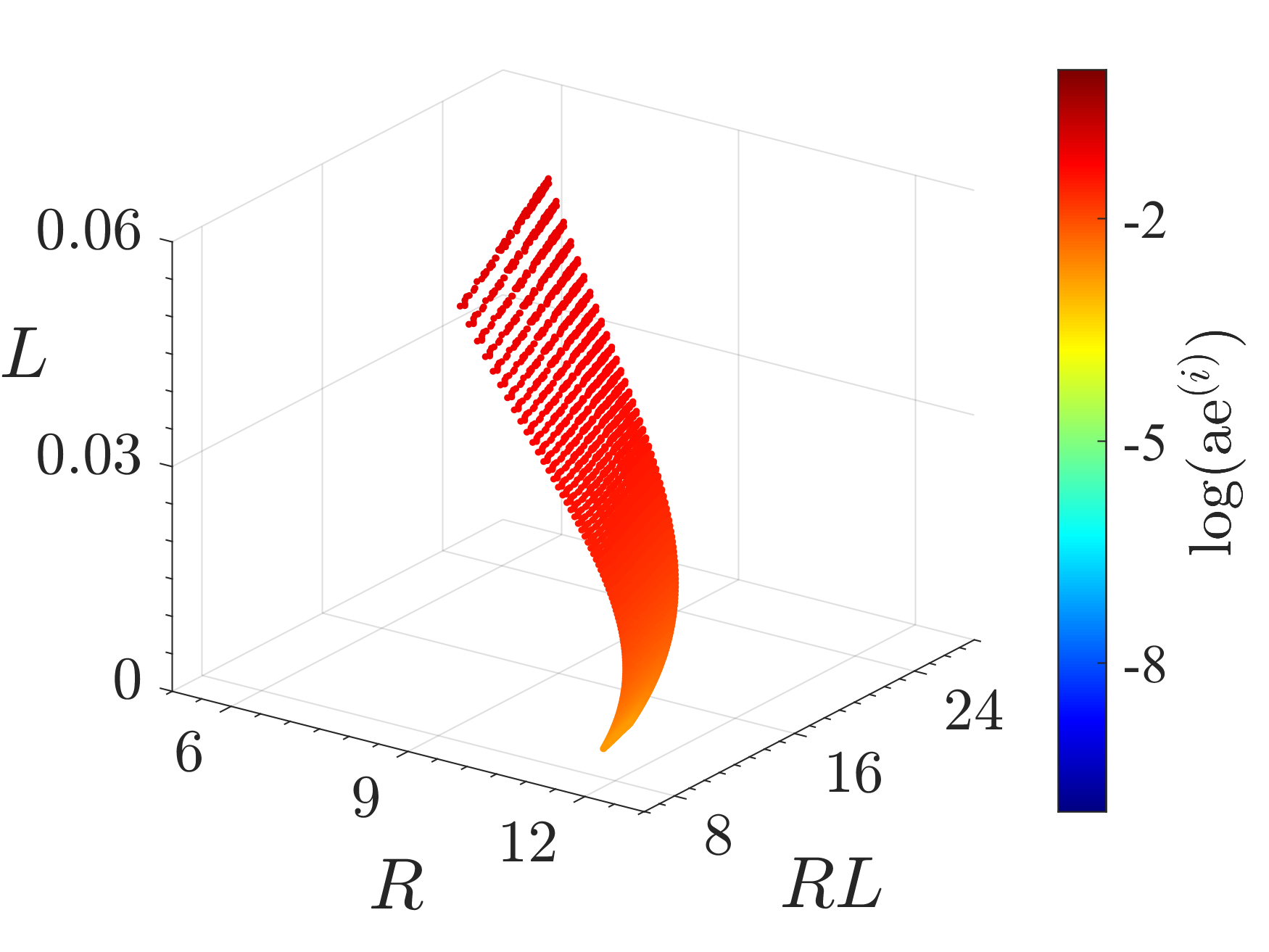}} 
    \subfigure[CSP$_R$(1)]{
    \includegraphics[width=0.32\textwidth]{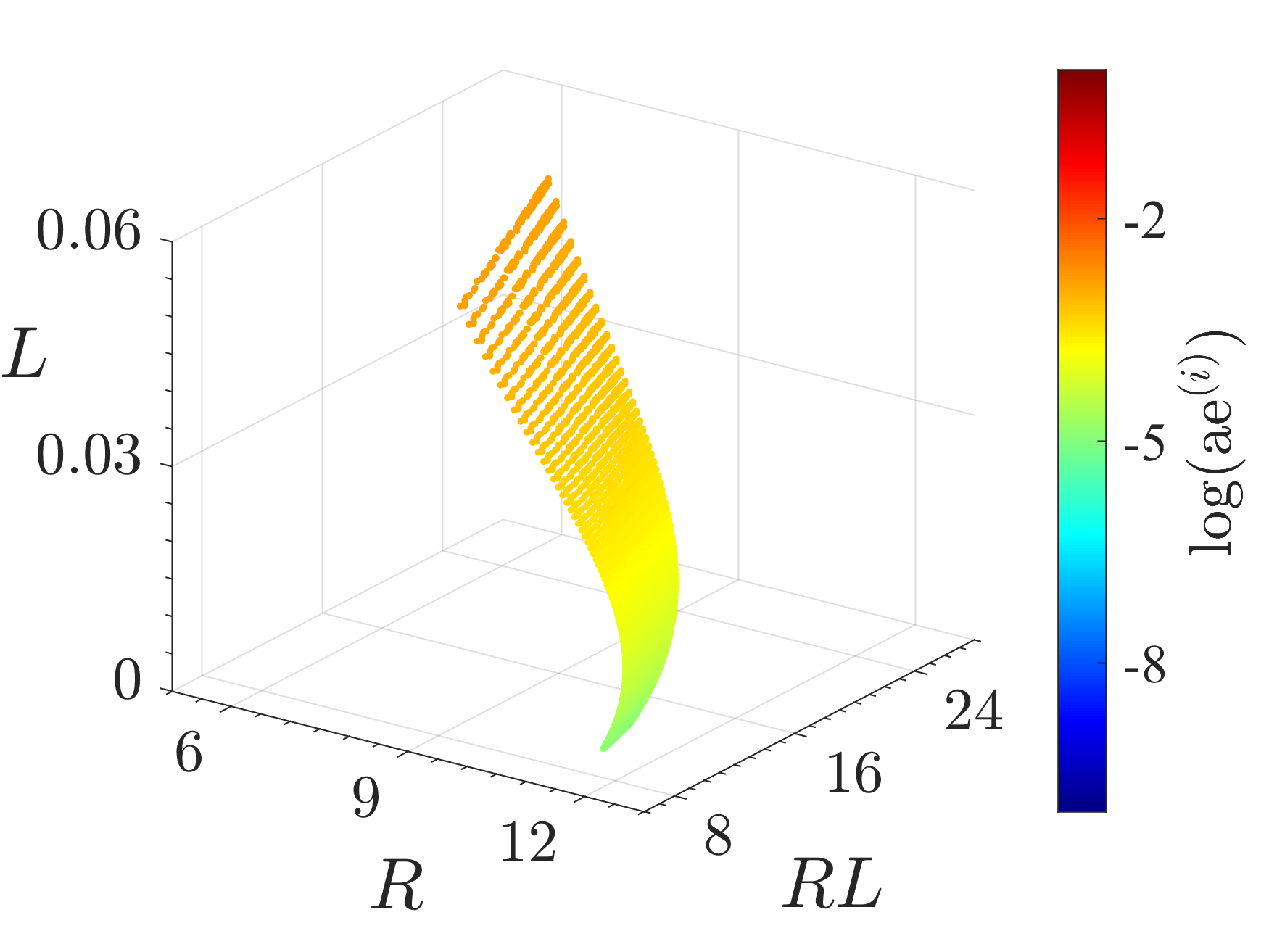}}
    \subfigure[CSP$_R$(2)]{
    \includegraphics[width=0.32\textwidth]{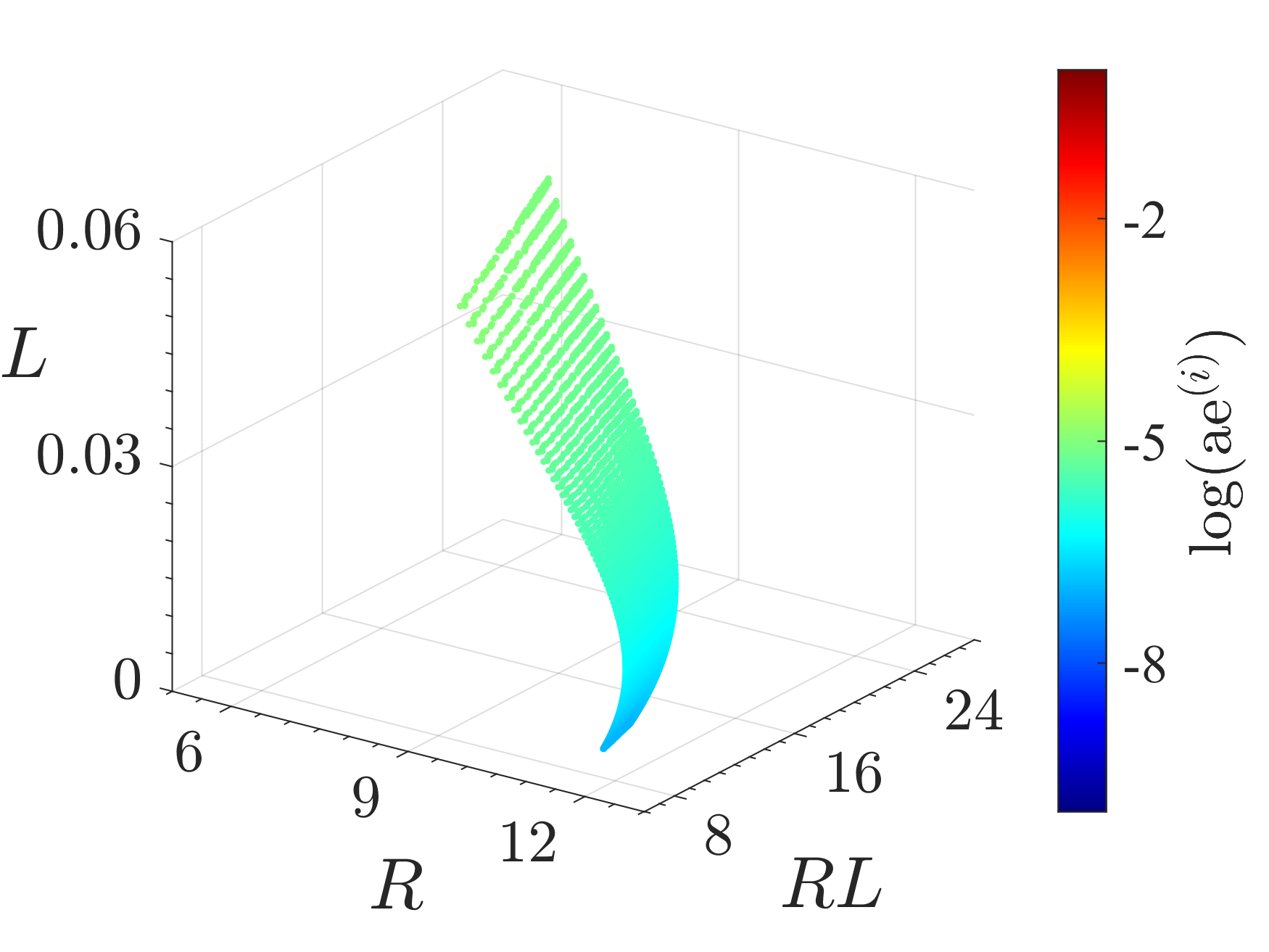}} \\
    \subfigure[QSSA$_{RL}$]{
    \includegraphics[width=0.32\textwidth]{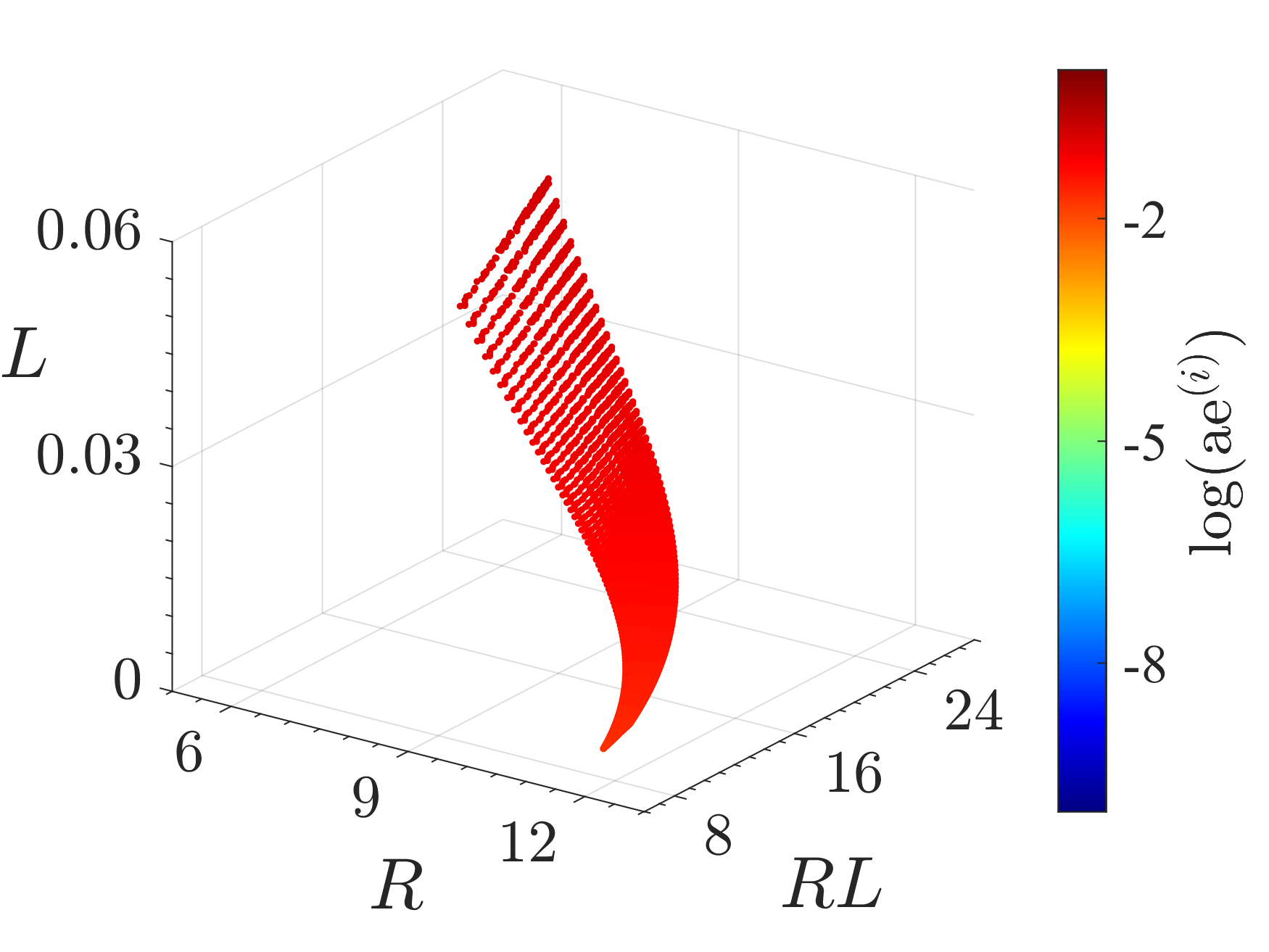}}
    \subfigure[CSP$_{RL}$(1)]{
    \includegraphics[width=0.32\textwidth]{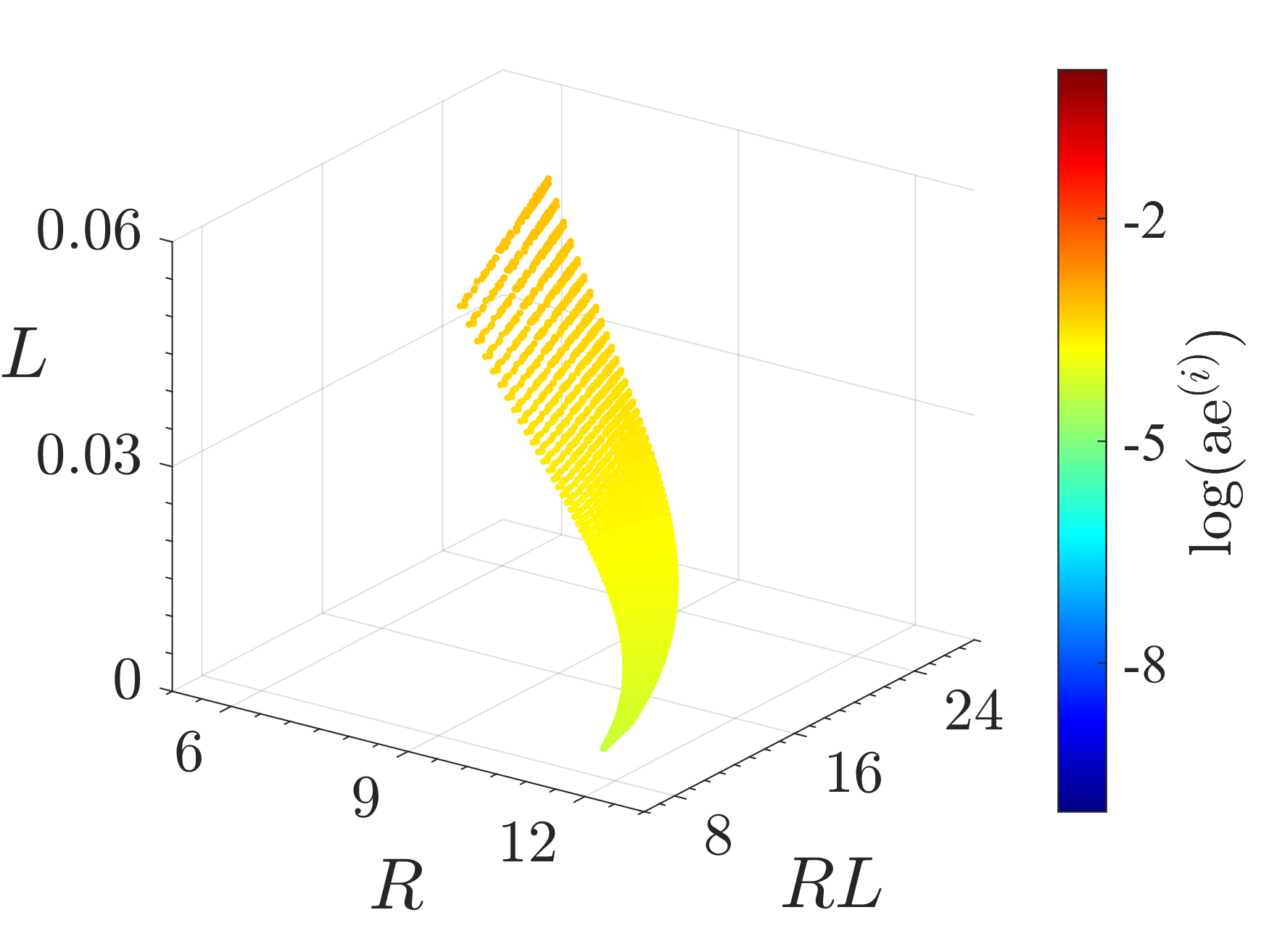}}
    \subfigure[CSP$_{RL}$(2)]{
    \includegraphics[width=0.32\textwidth]{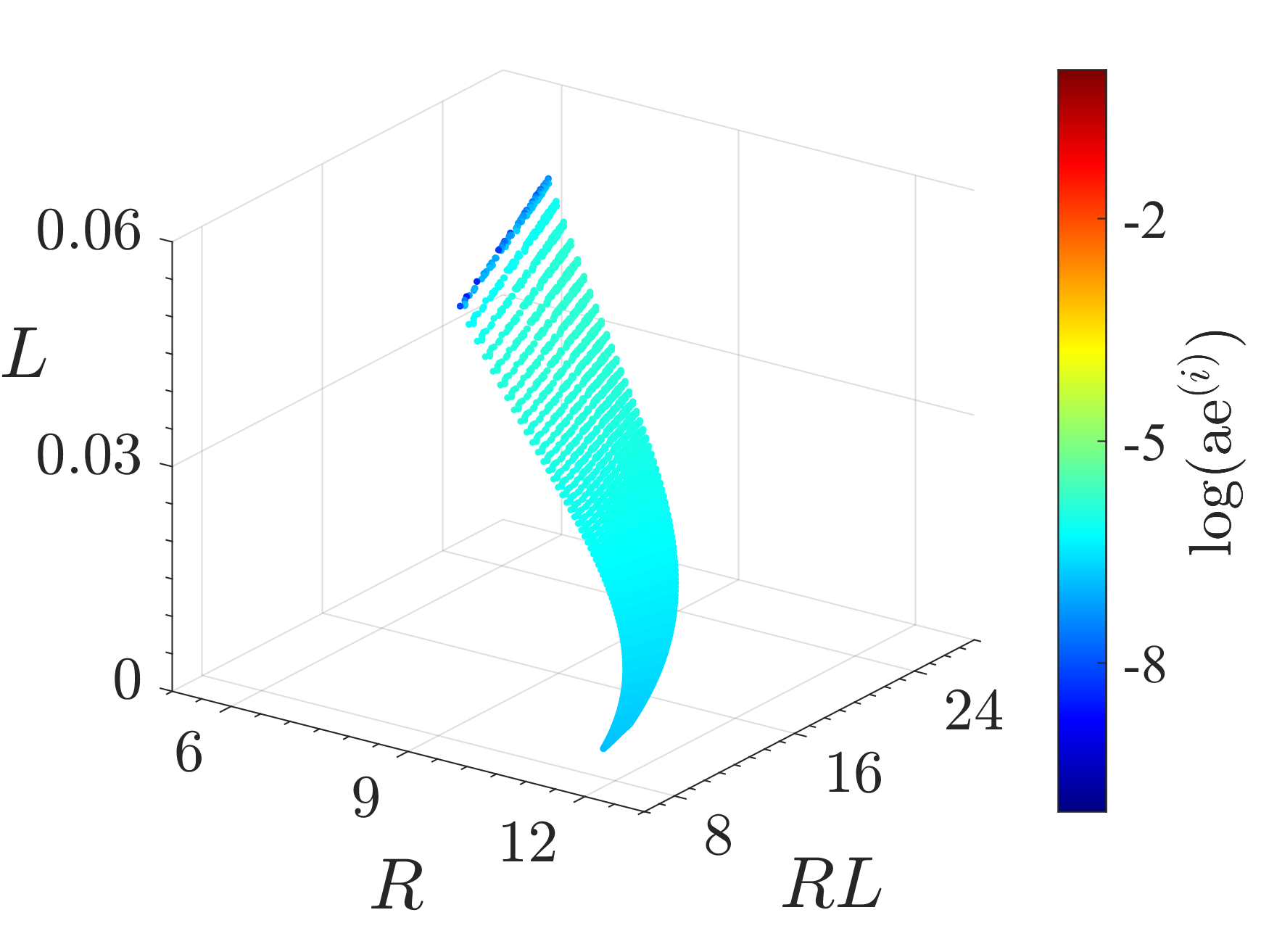}}
    \caption{TMDD system in Eq.~(33) for the SIM $\mathcal{M}_2$, where $x=L$ is the fast variable.~Absolute errors ($ae^{(i)}$) of the SIM approximations over all points of the test set, in comparison to the numerical solution $\mathbf{z}^{(i)}=[x^{(i)},\mathbf{y}^{(i)}]^\top$ for $i=1,\ldots,n_t$.~Panel (a) depicts the $\lvert \mathbf{C} \mathbf{z}^{(i)} - \mathcal{N}(\mathbf{D} \mathbf{z}^{(i)}) \rvert$ of the PINN scheme, panels (b), (c), (d), (f), (g), (i) and (j) depict $\lvert x^{(i)} - h(\mathbf{y}^{(i)}) \rvert$ of the PEA, QSSA$_L$,  CSP$_L$(1), QSSA$_R$,  CSP$_R$(1), QSSA$_{RL}$ and  CSP$_{RL}$(1) explicit functionals w.r.t. $x$, and panels (e), (h) and (k) depict $\lvert x^{(i)} - \hat{x}^{(i)}) \rvert$ of the CSP$_L$(2), CSP$_R$(2) and CSP$_{RL}$(2) implicit functionals, solved numerically with Newton for $x$.~Note that the approximations of the second/third/fourth row were constructed with $L$/$R$/$RL$ assumed as fast variable.}
    \label{SF:TMDDP4_AE}
\end{figure}

\begin{figure}[!h]
    \centering
    \subfigure[PIML, projection to $c_1$]{
    \includegraphics[width=0.32\textwidth]{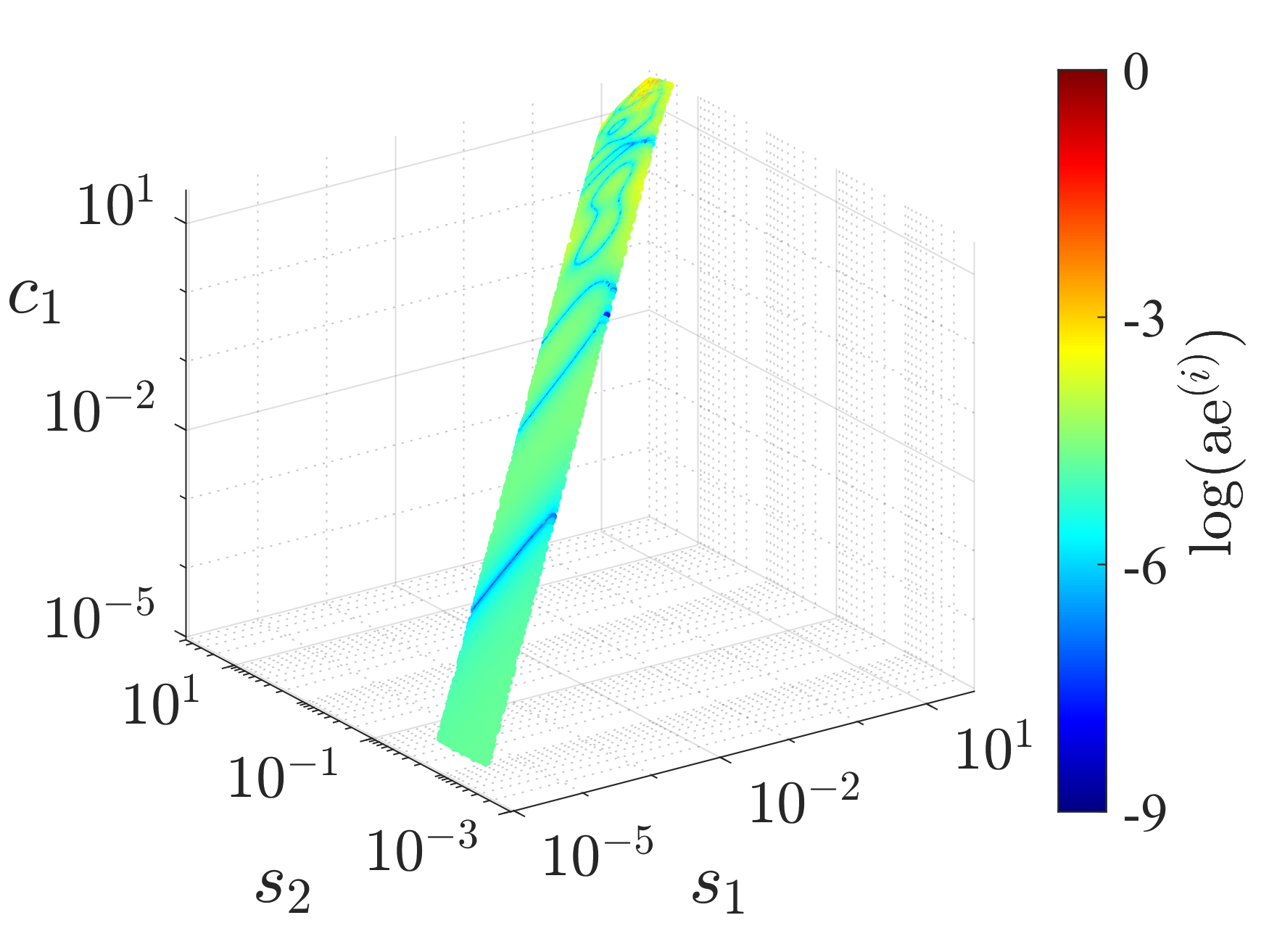}} 
    \subfigure[QSSA$_{c1c2}$, projection to $c_1$]{
    \includegraphics[width=0.32\textwidth]{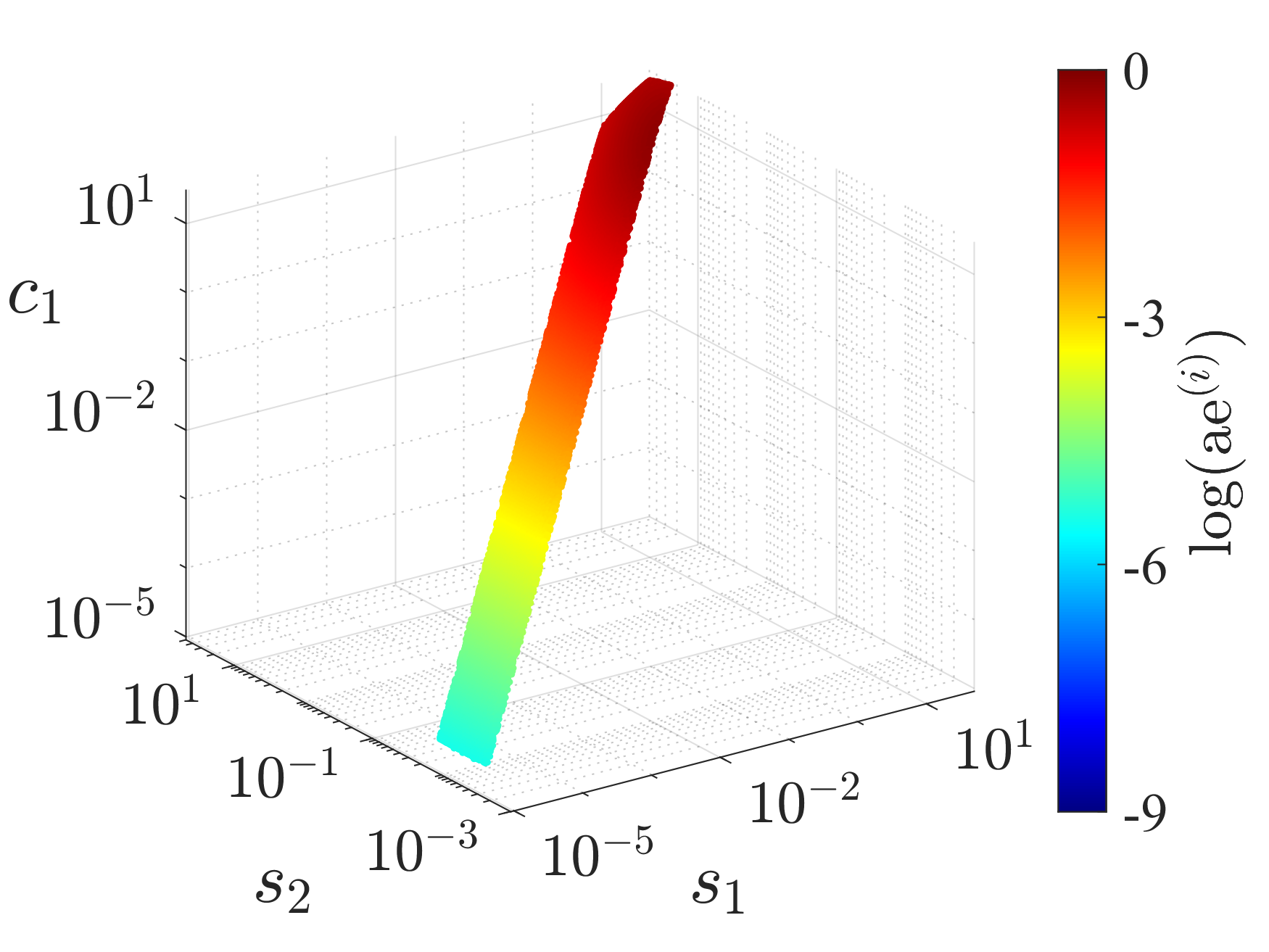}} \subfigure[PEA$_{13}$, projection to $c_1$]{
    \includegraphics[width=0.32\textwidth]{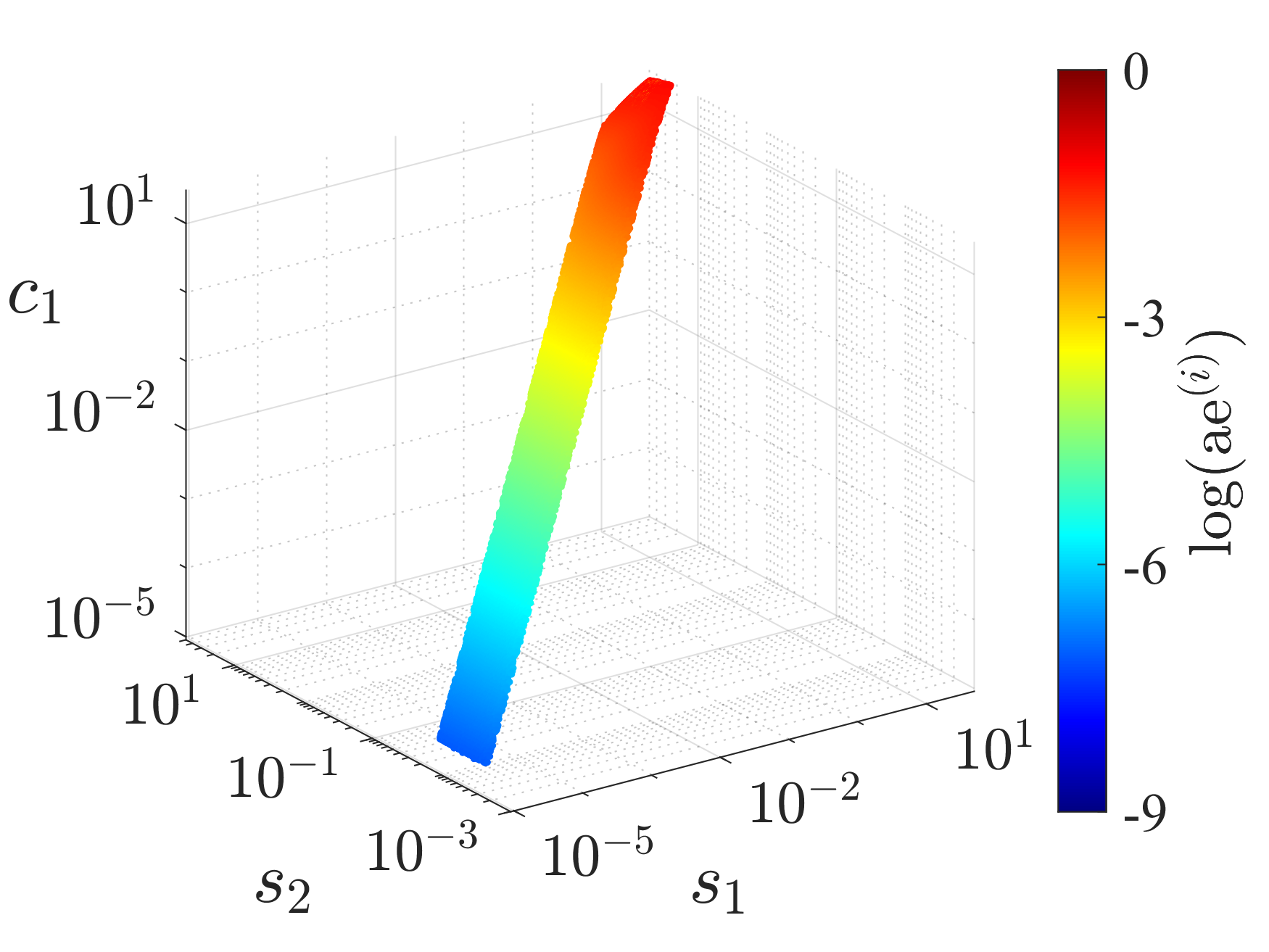}}  \\
    \subfigure[CSP$_{c1c2}$(1), projection to $c_1$]{
    \includegraphics[width=0.32\textwidth]{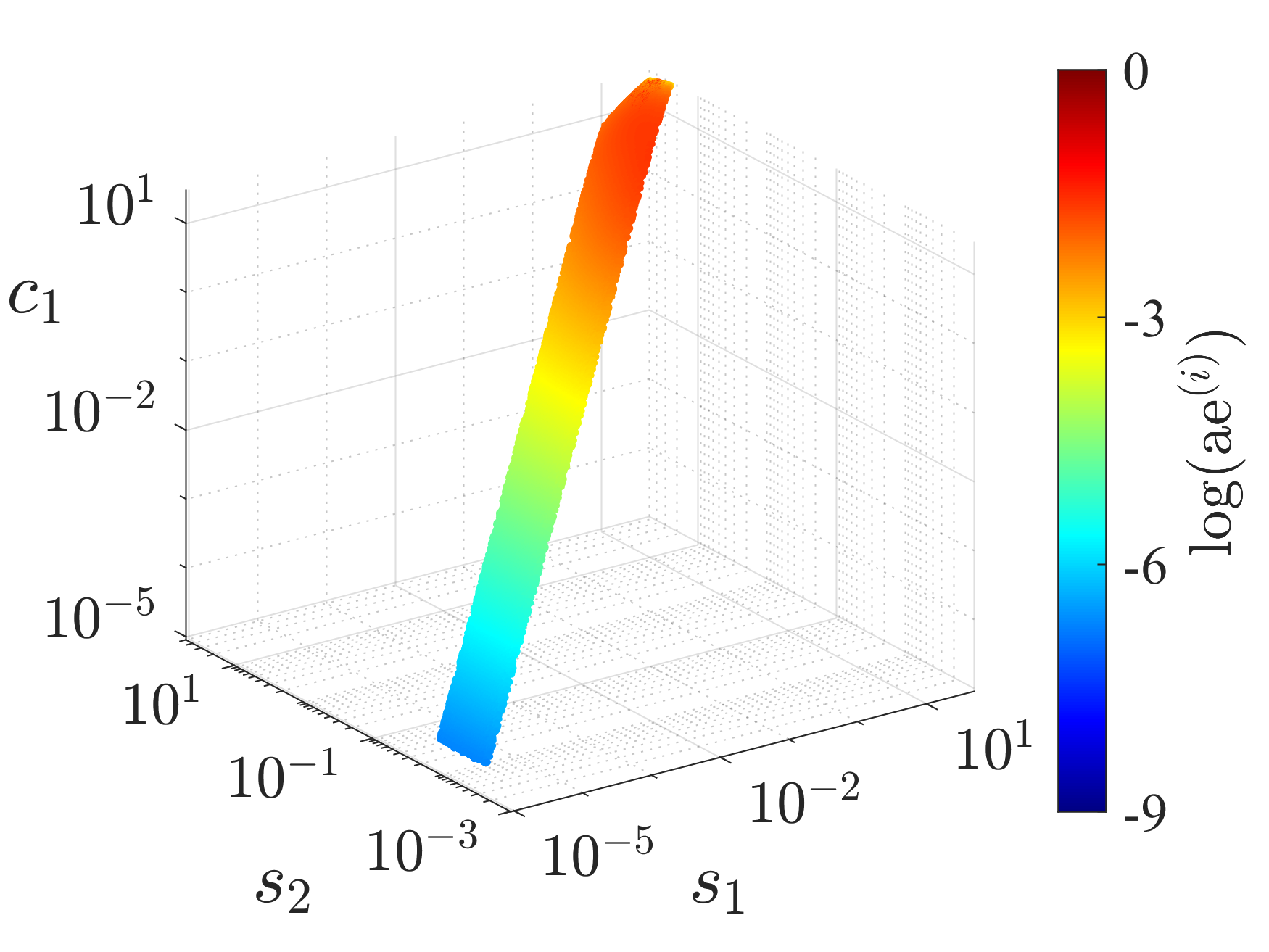}} \subfigure[CSP$_{c1c2}$(2), projection to $c_1$]{
    \includegraphics[width=0.32\textwidth]{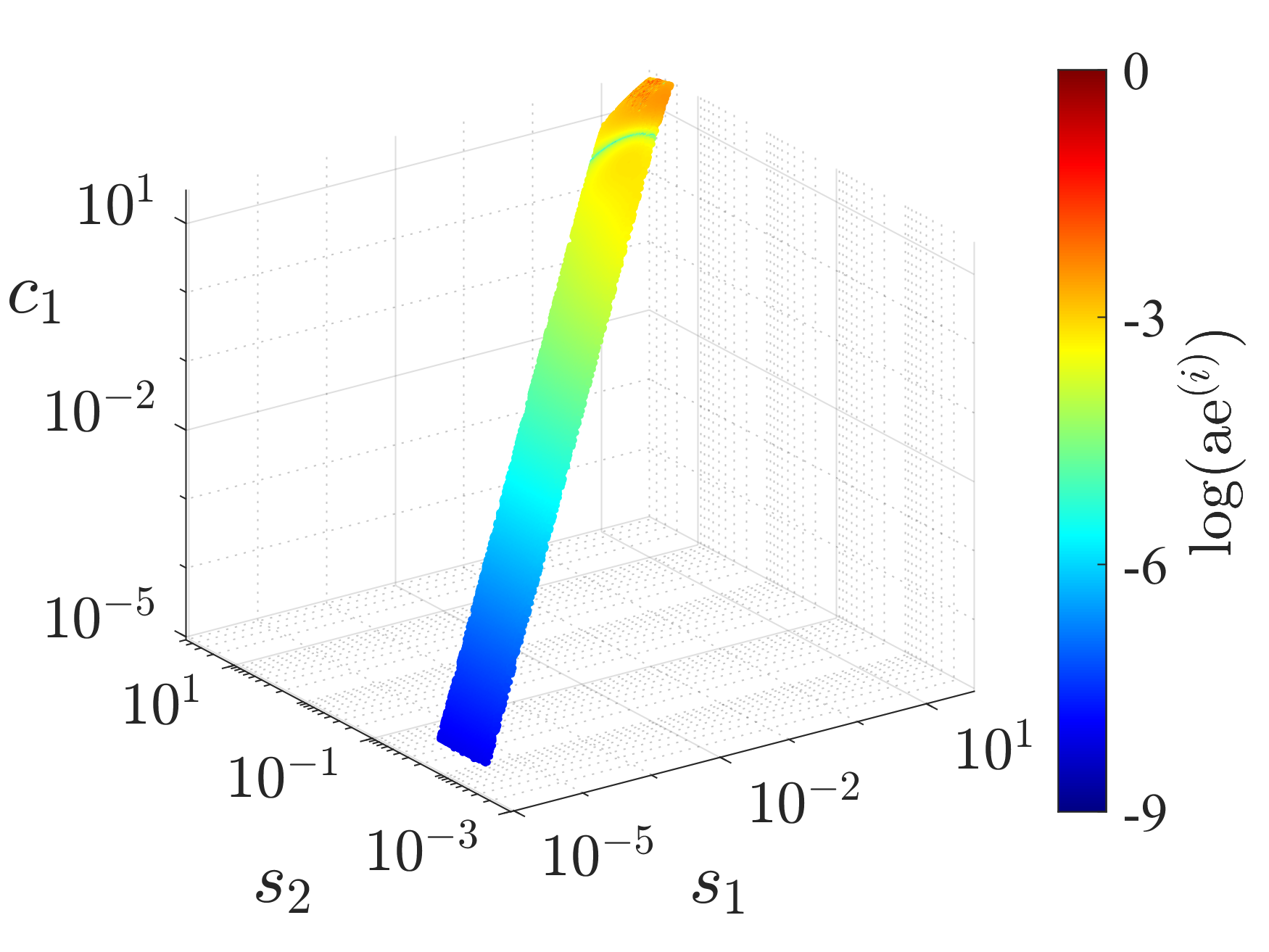}} \\    
    \subfigure[PIML, projection to $c_2$]{
    \includegraphics[width=0.32\textwidth]{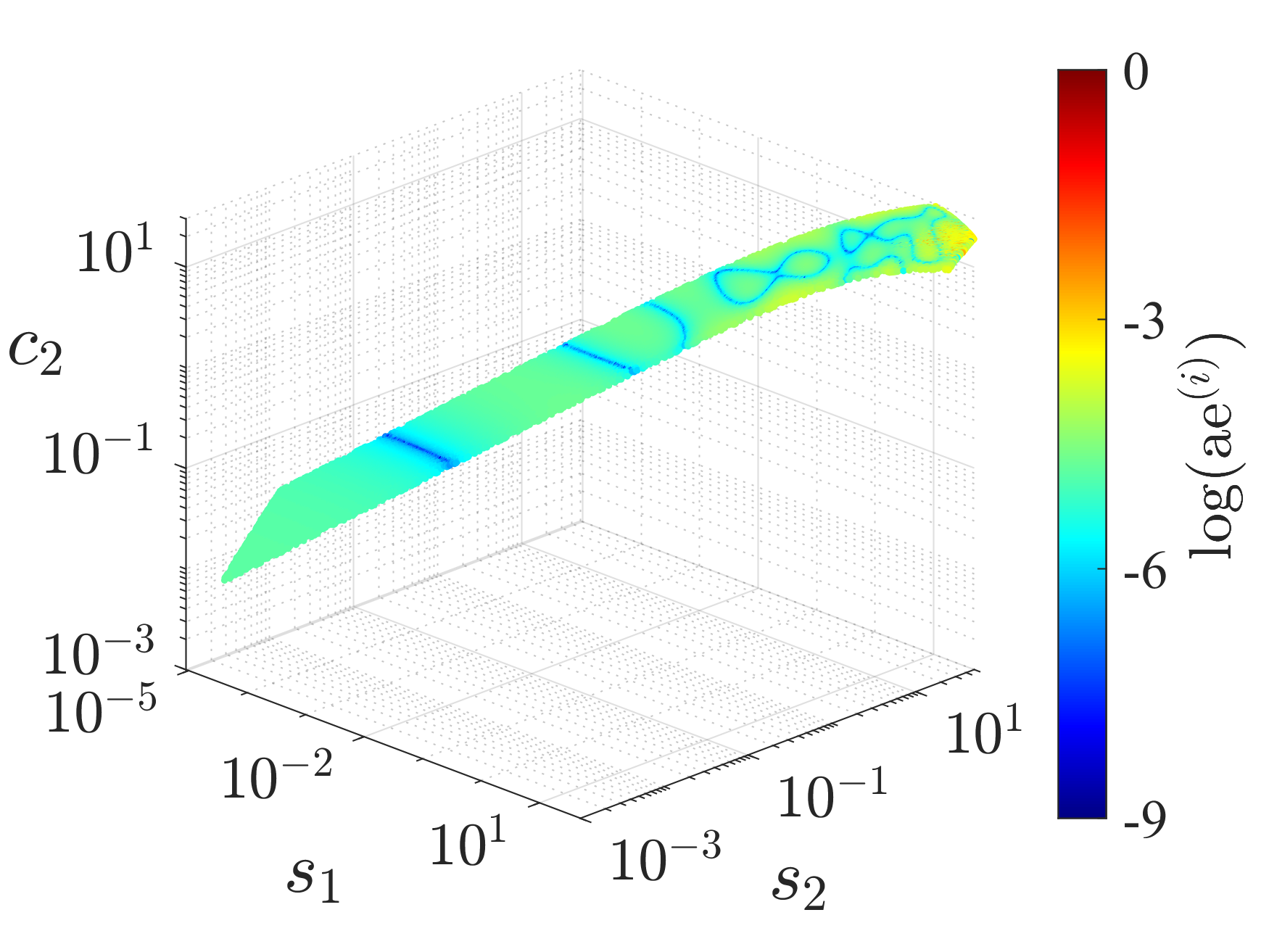}}   
    \subfigure[QSSA$_{c1c2}$, projection to $c_2$]{
    \includegraphics[width=0.32\textwidth]{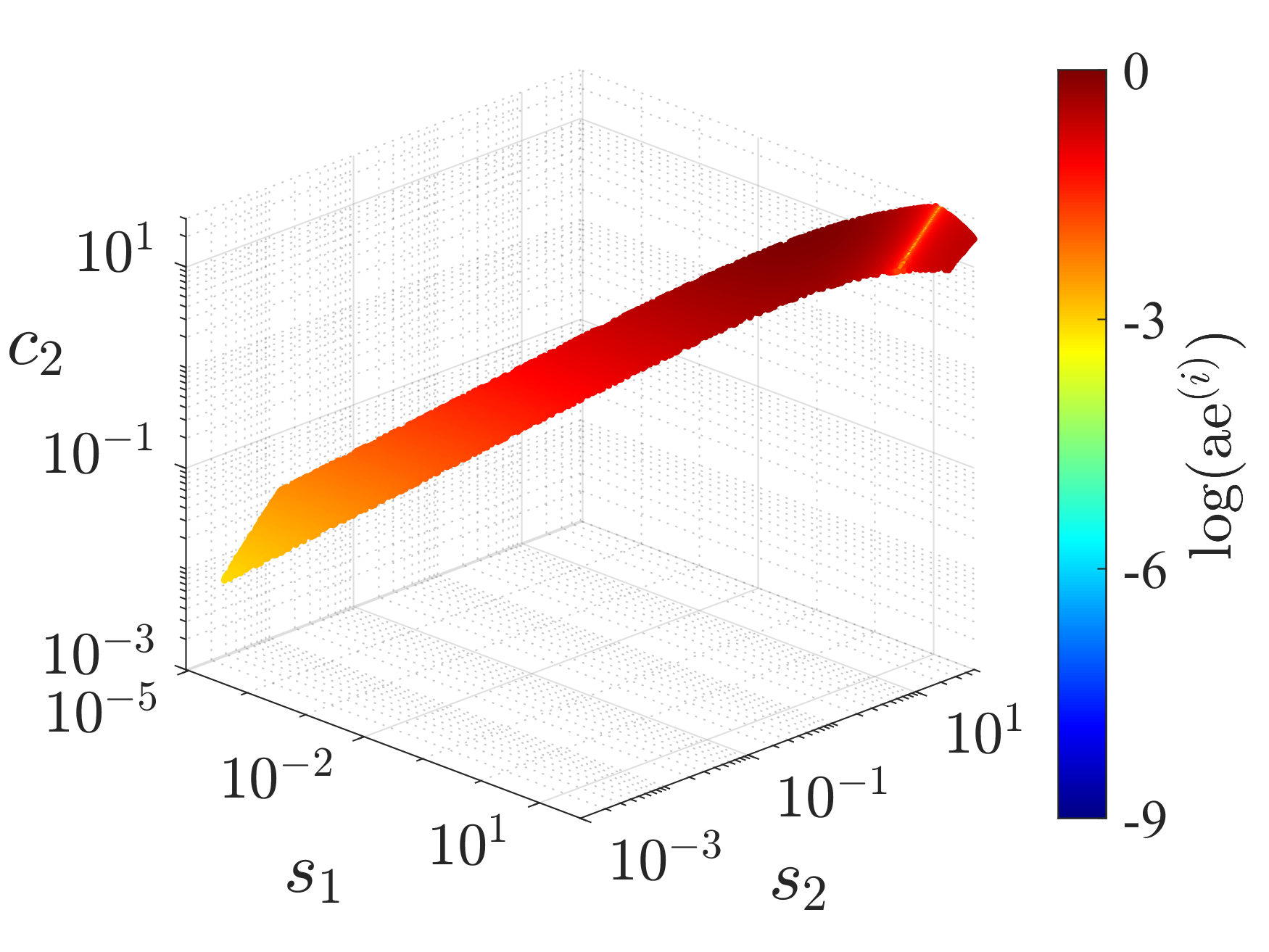}} \subfigure[PEA$_{13}$, projection to $c_2$]{
    \includegraphics[width=0.32\textwidth]{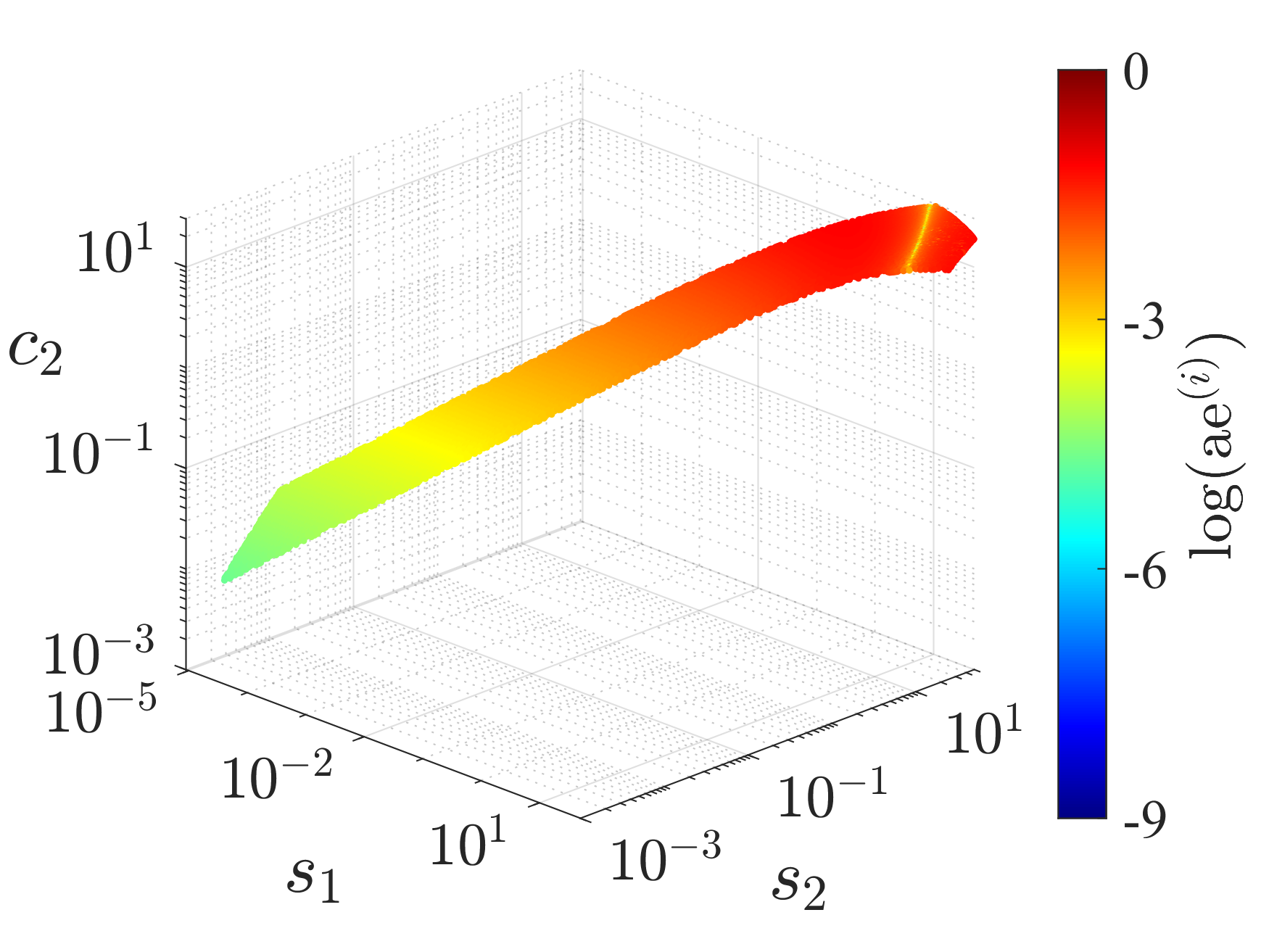}} \\
     \subfigure[CSP$_{c1c2}$(1), projection to $c_2$]{
    \includegraphics[width=0.32\textwidth]{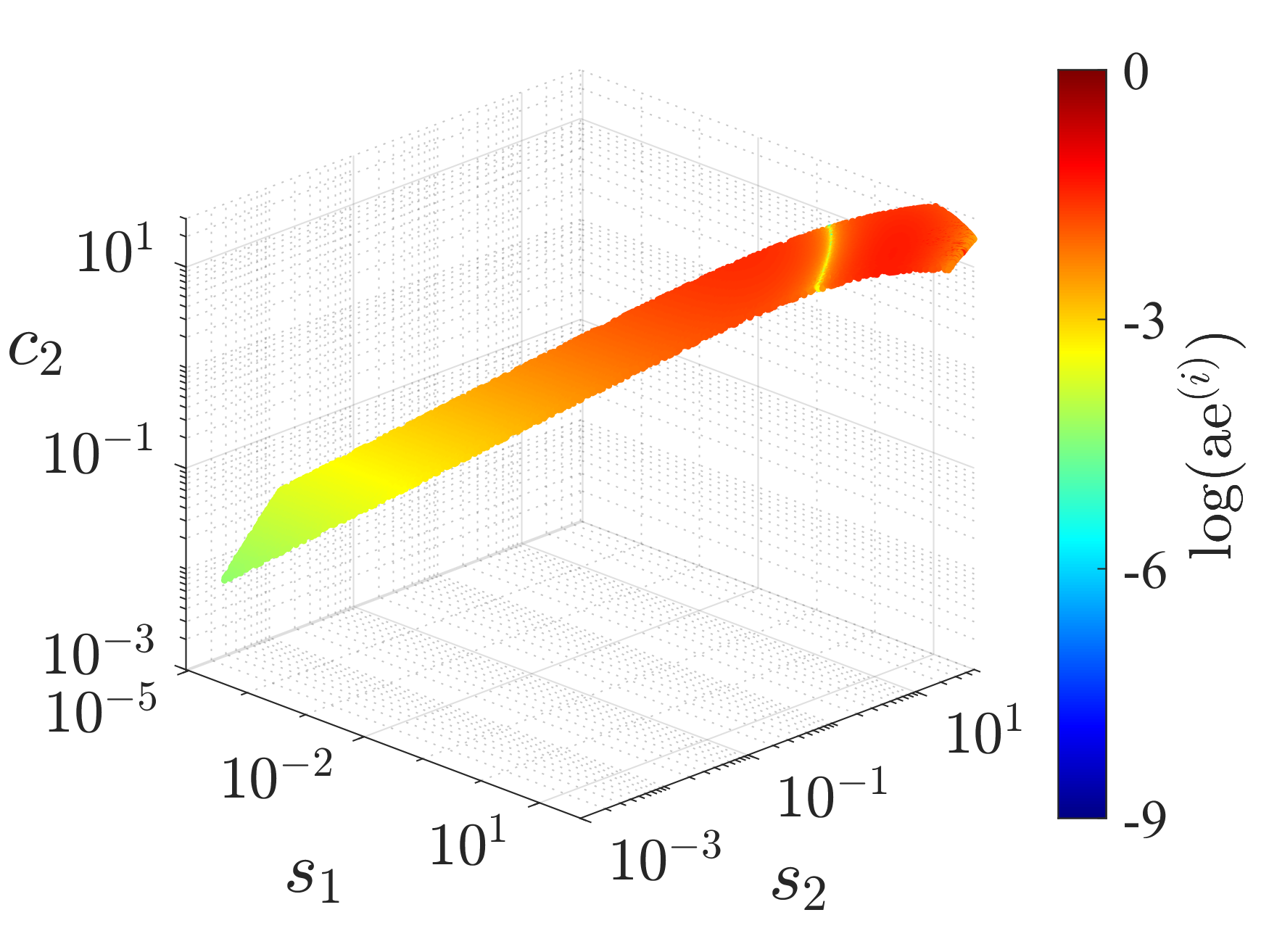}}  
    \subfigure[CSP$_{c1c2}$(2), projection to $c_2$]{
    \includegraphics[width=0.32\textwidth]{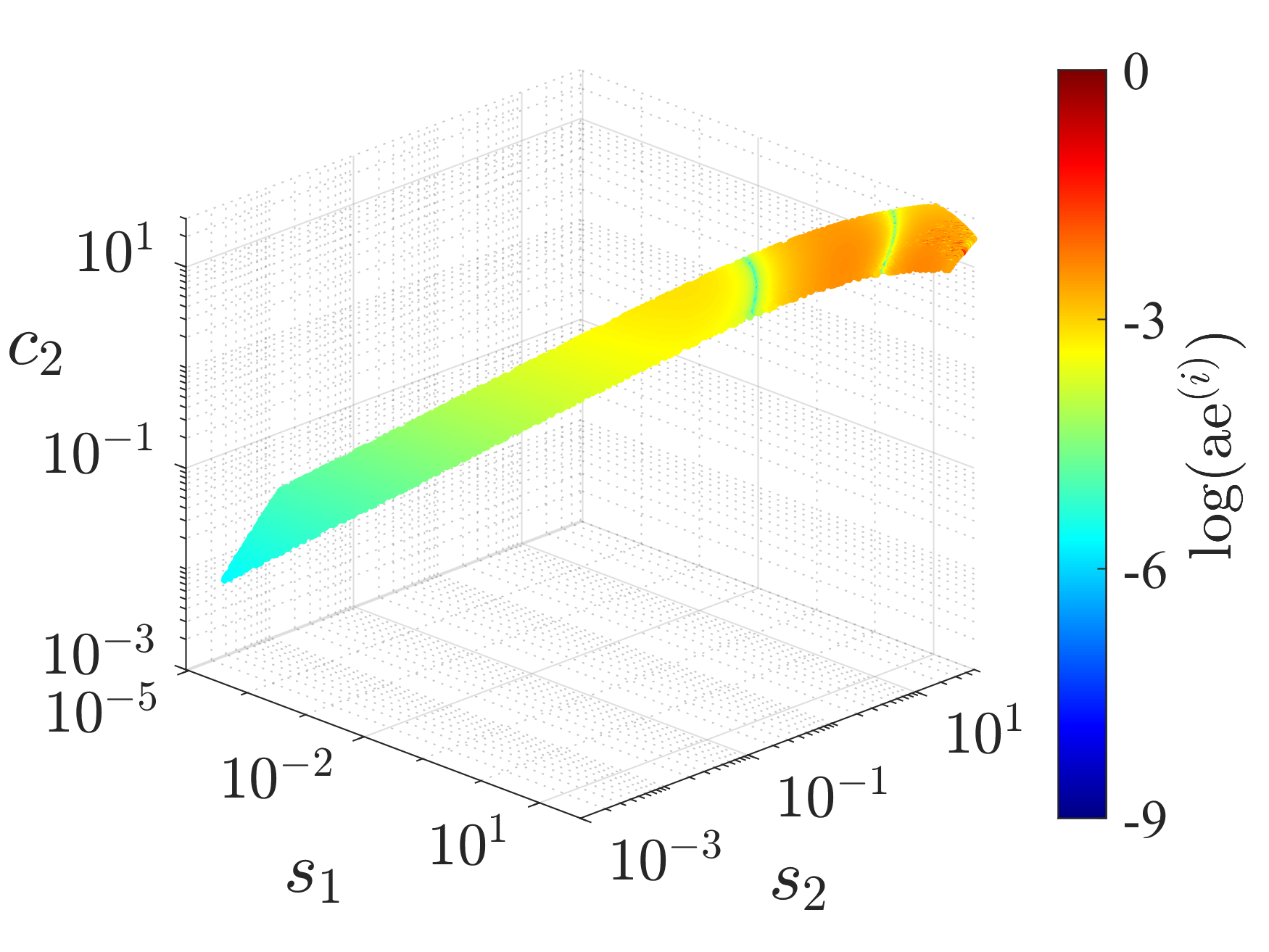}} 
    \caption{Original fCSI system in Eq.~(34) for the $M=2$-dim. SIM, where $\mathbf{x}=[c_1,c_2]$ are the fast variables.~Absolute errors ($ae^{(i)}$) of the SIM approximations over all points of the test set, in comparison to the numerical solution $\mathbf{z}^{(i)}=[\mathbf{x}^{(i)},\mathbf{y}^{(i)}]^\top$ for $i=1,\ldots,n_t$.~Panels (a, f) depict the $\lvert \mathbf{C} \mathbf{z}^{(i)} - \mathcal{N}(\mathbf{D} \mathbf{z}^{(i)}) \rvert$ of the PINN scheme, panels (b, g) depict $\lvert \mathbf{x}^{(i)} - \mathbf{h}(\mathbf{y}^{(i)}) \rvert$ of the QSSA$_{c1c2}$ explicit functionals w.r.t. $\mathbf{x}$, and panels (c, h), (d, i) and (e, j) depict $\lvert \mathbf{x}^{(i)} - \hat{\mathbf{x}}^{(i)}) \rvert$ of the PEA$_{13}$, CSP$_{c1c2}$(1) and CSP$_{c1c2}$(2) implicit functionals, solved numerically with Newton for $\mathbf{x}$.~Each pair of panels shows projections of the SIM to $c_1$ and $c_2$ fast variables.}
    \label{SF:InhSt_AE}
\end{figure}

\begin{figure}[!h]
    \centering
    \subfigure[PIML, projection to $c_1$]{
    \includegraphics[width=0.32\textwidth]{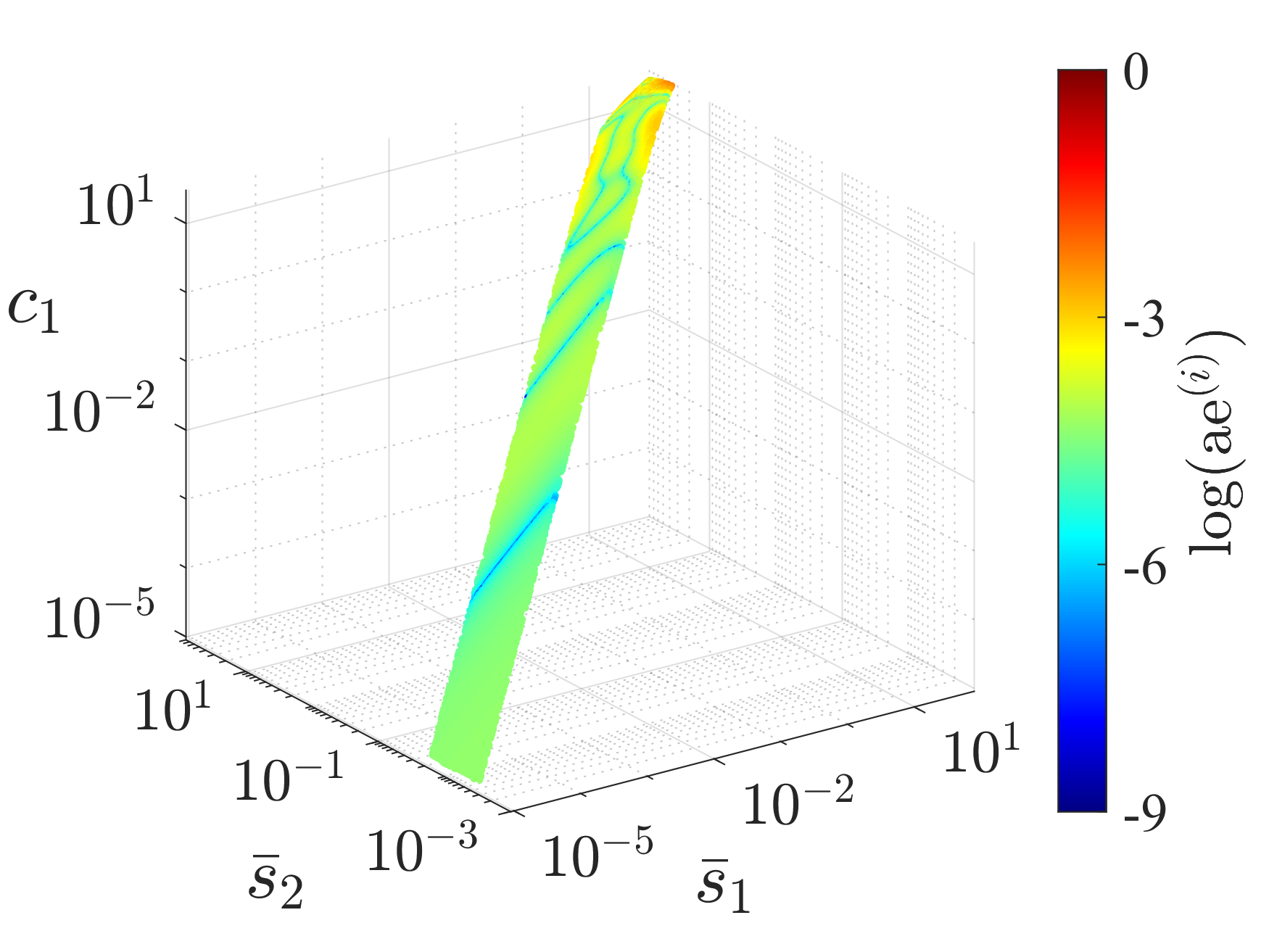}} 
    \subfigure[QSSA$_{c1c2}$, projection to $c_1$]{
    \includegraphics[width=0.32\textwidth]{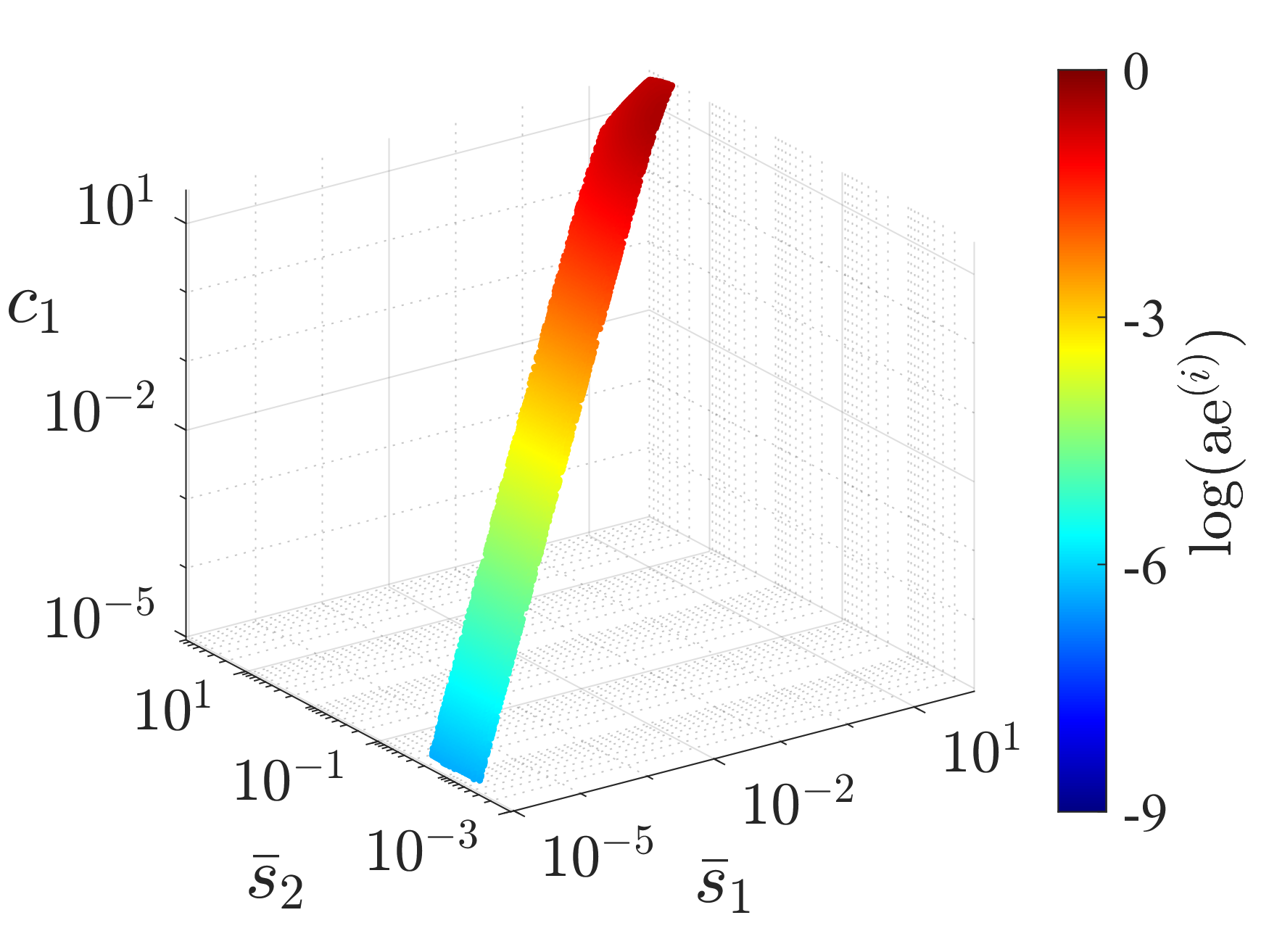}} \subfigure[PEA$_{13}$, projection to $c_1$]{
    \includegraphics[width=0.32\textwidth]{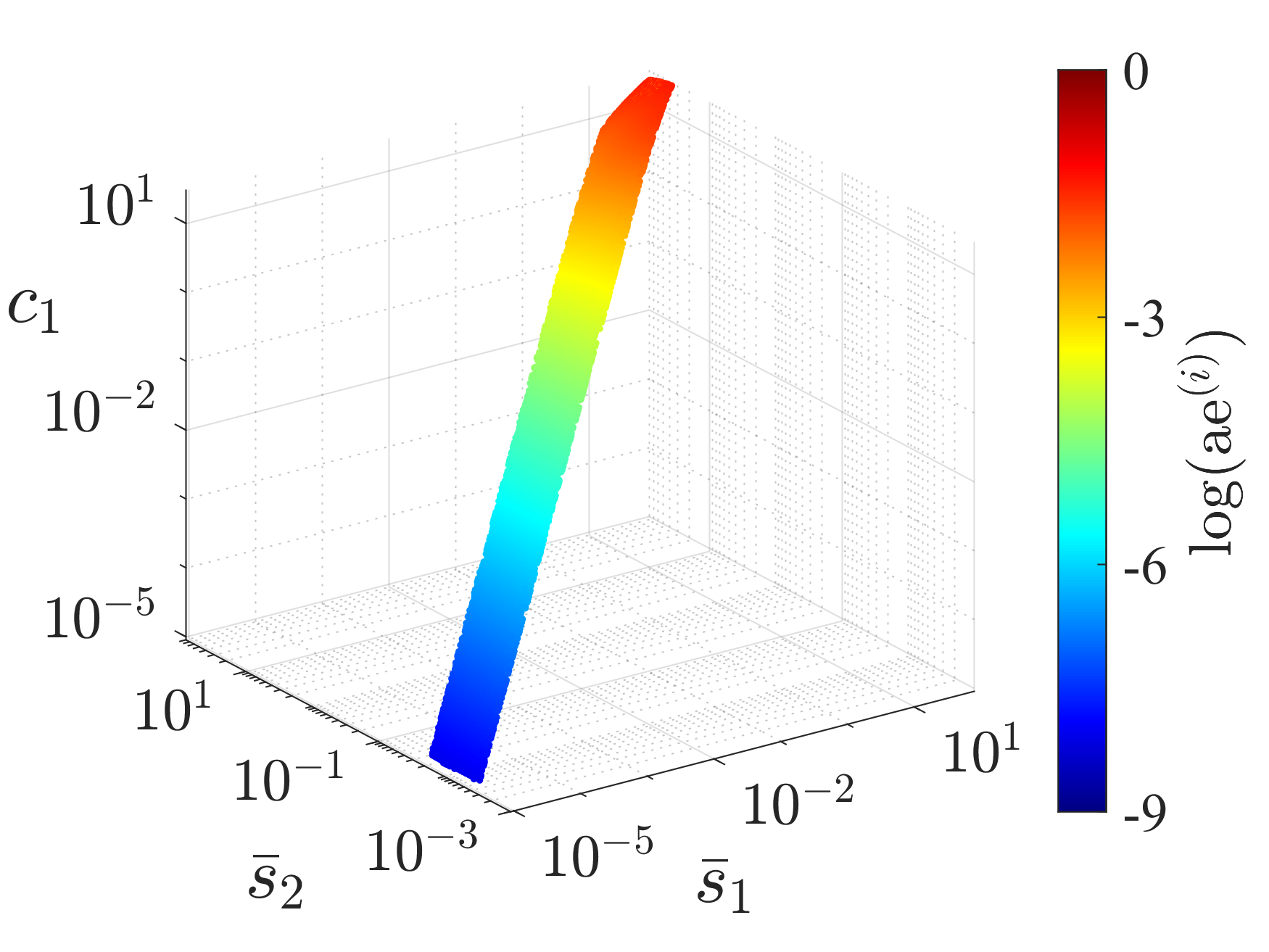}}  \\
    \subfigure[CSP$_{c1c2}$(1), projection to $c_1$]{
    \includegraphics[width=0.32\textwidth]{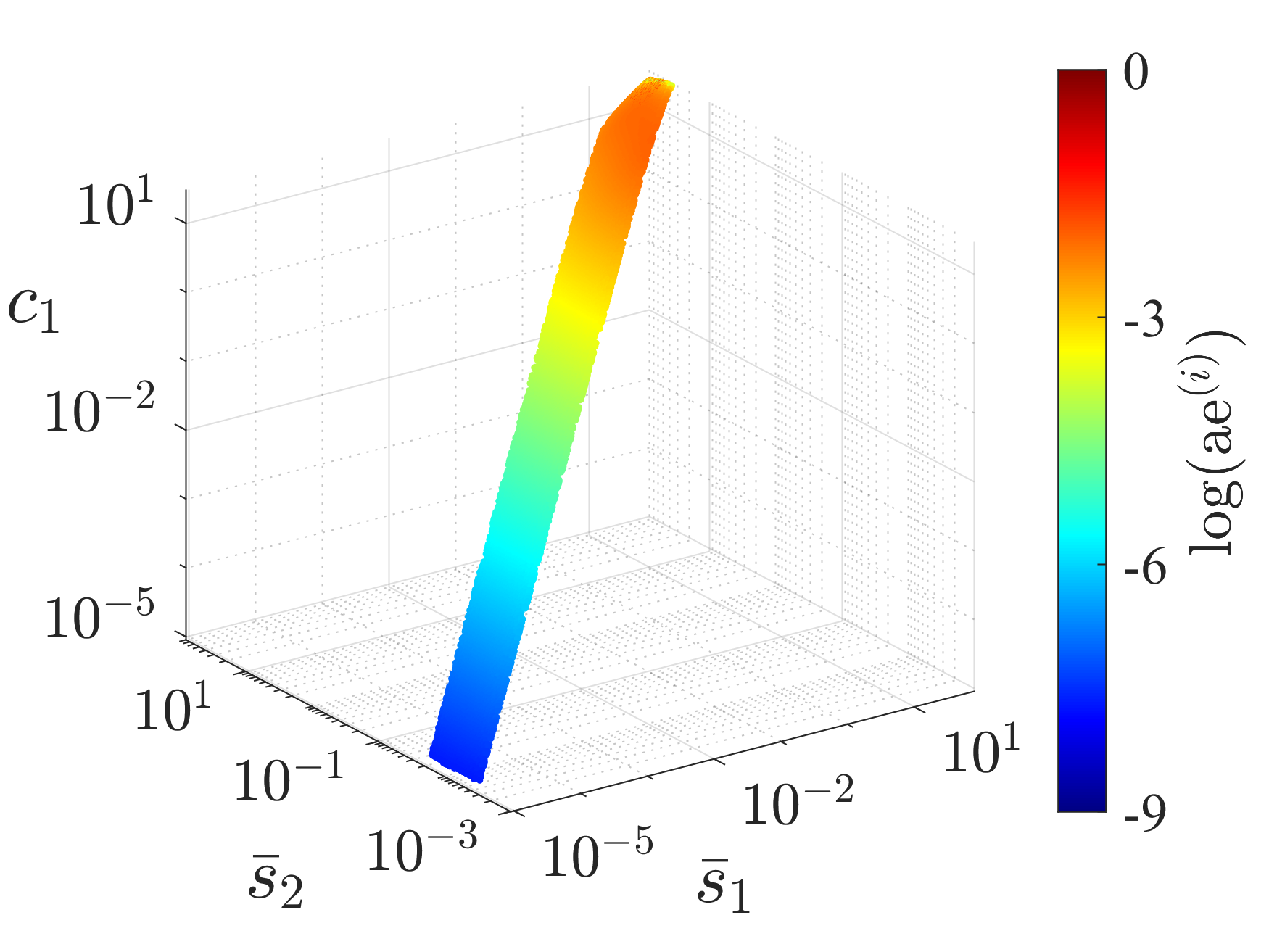}} \subfigure[CSP$_{c1c2}$(2), projection to $c_1$]{
    \includegraphics[width=0.32\textwidth]{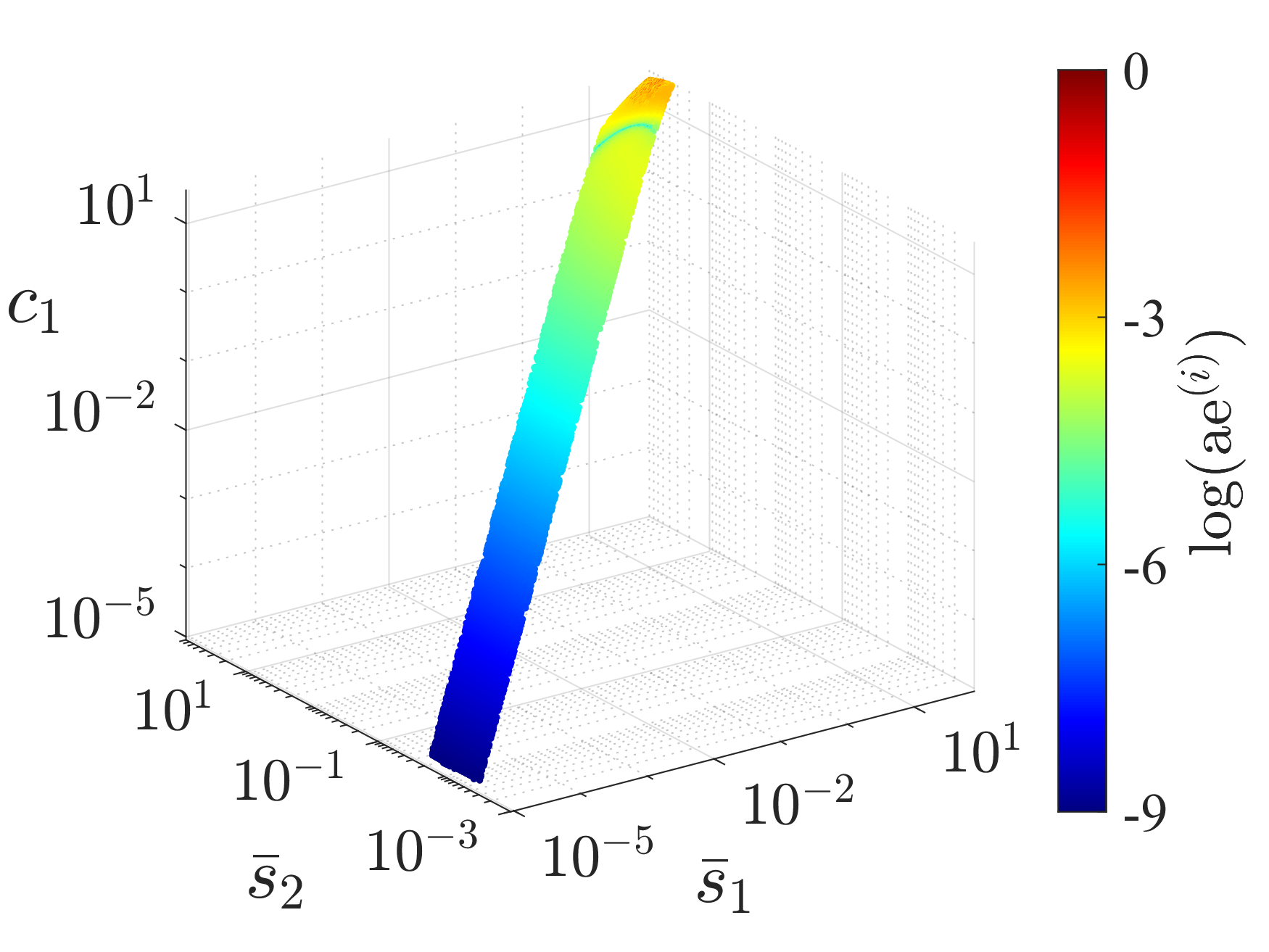}} \\    
    \subfigure[PIML, projection to $c_2$]{
    \includegraphics[width=0.32\textwidth]{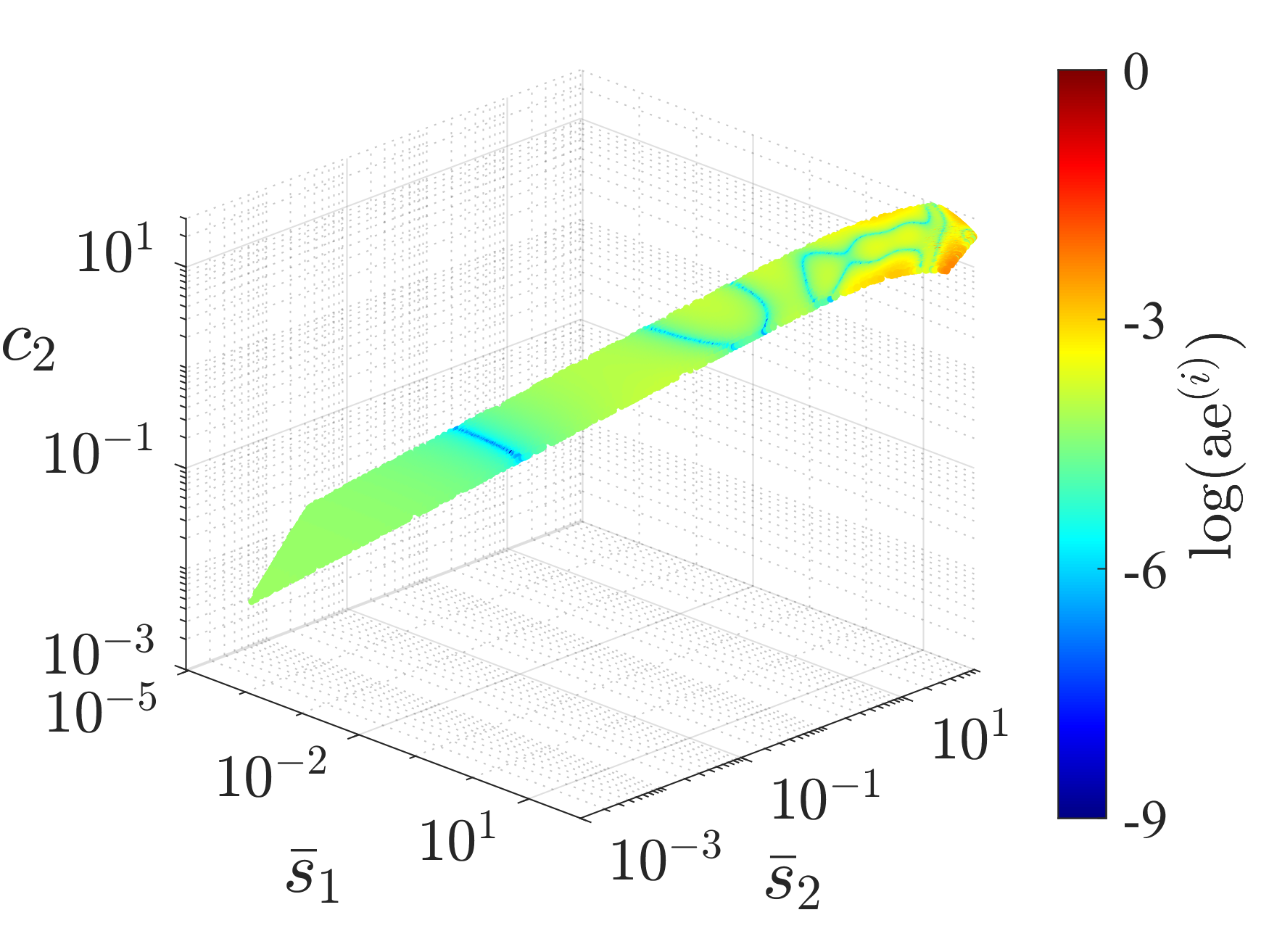}}   
    \subfigure[QSSA$_{c1c2}$, projection to $c_2$]{
    \includegraphics[width=0.32\textwidth]{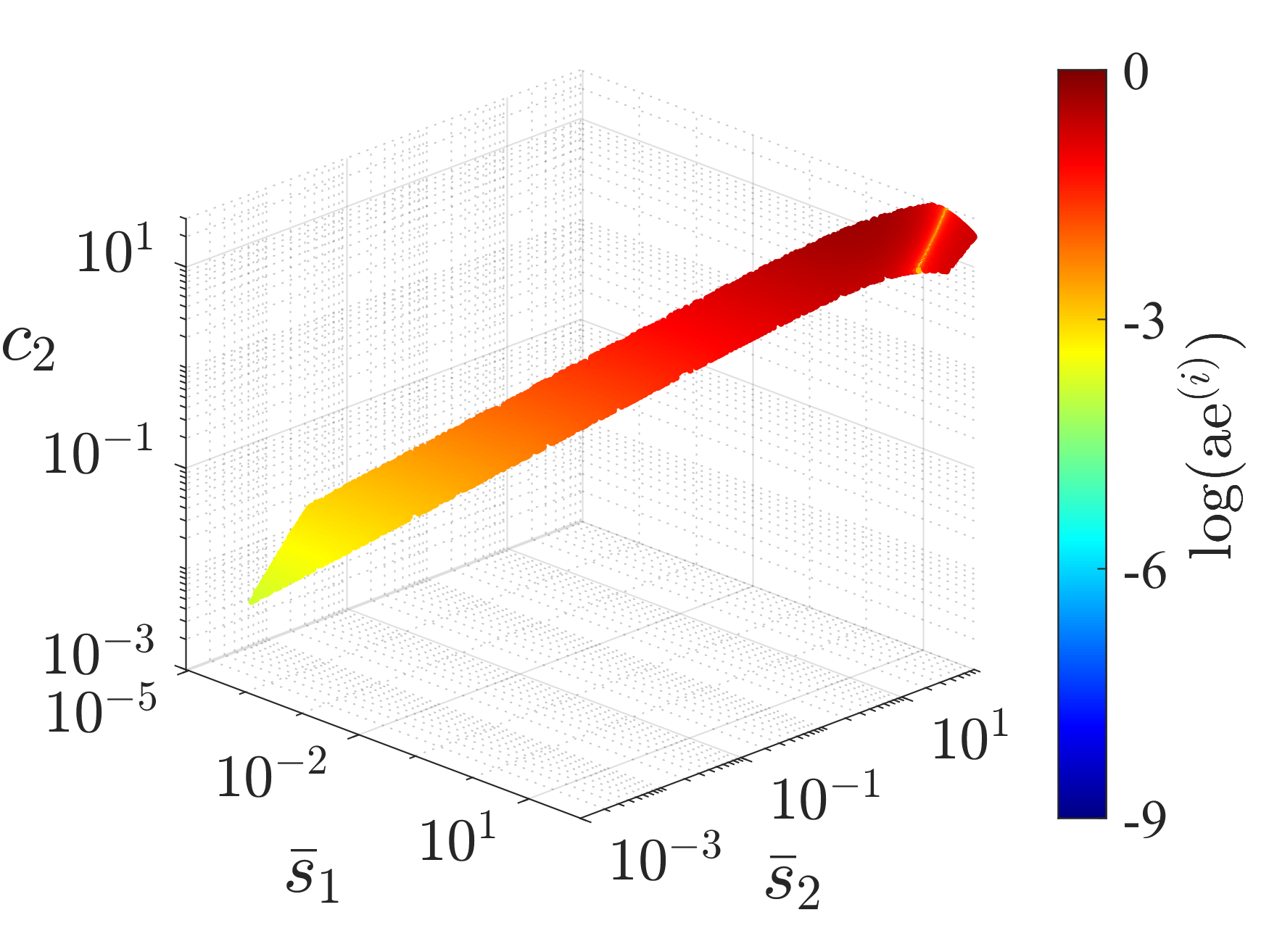}} \subfigure[PEA$_{13}$, projection to $c_2$]{
    \includegraphics[width=0.32\textwidth]{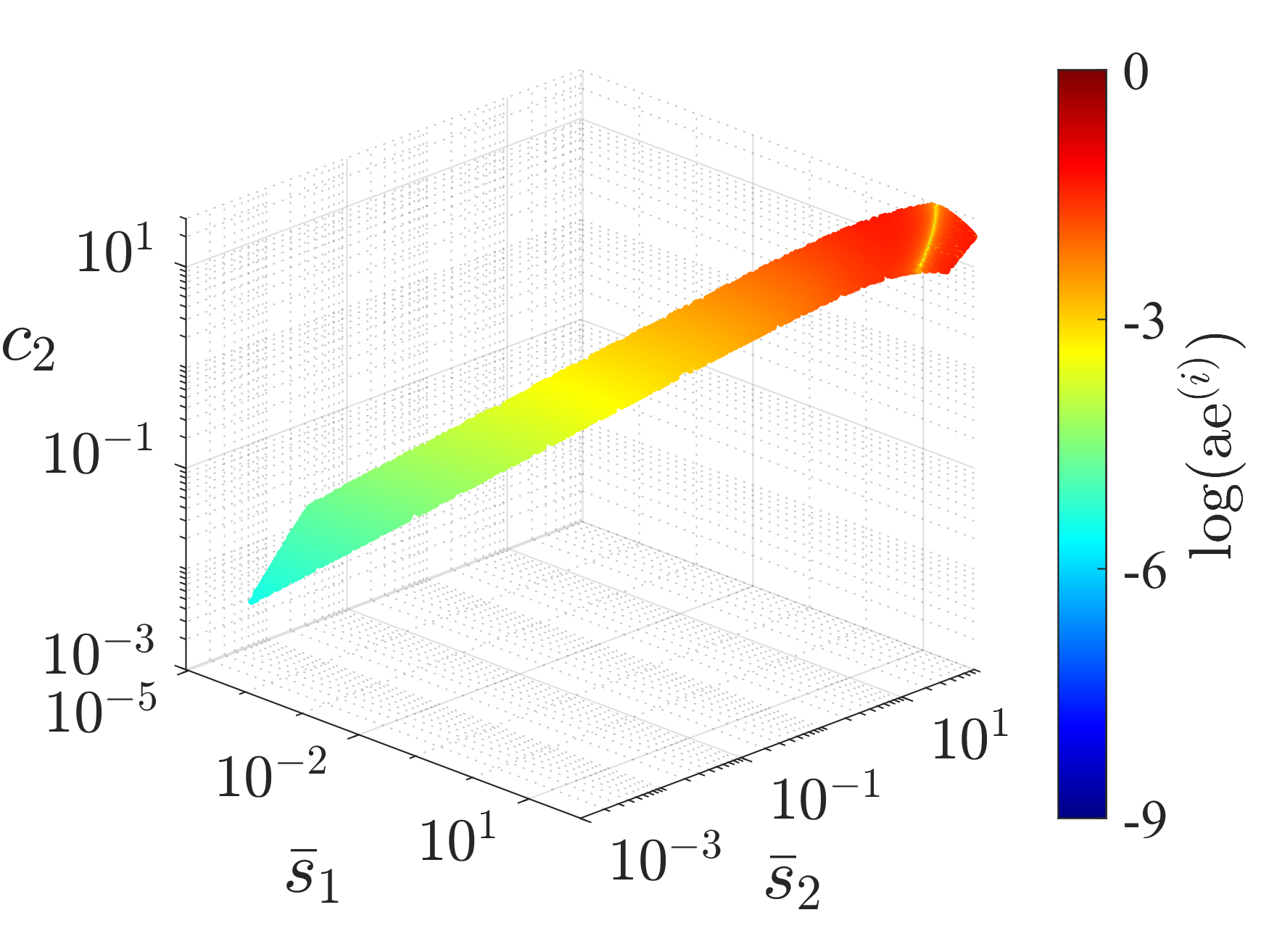}} \\
     \subfigure[CSP$_{c1c2}$(1), projection to $c_2$]{
    \includegraphics[width=0.32\textwidth]{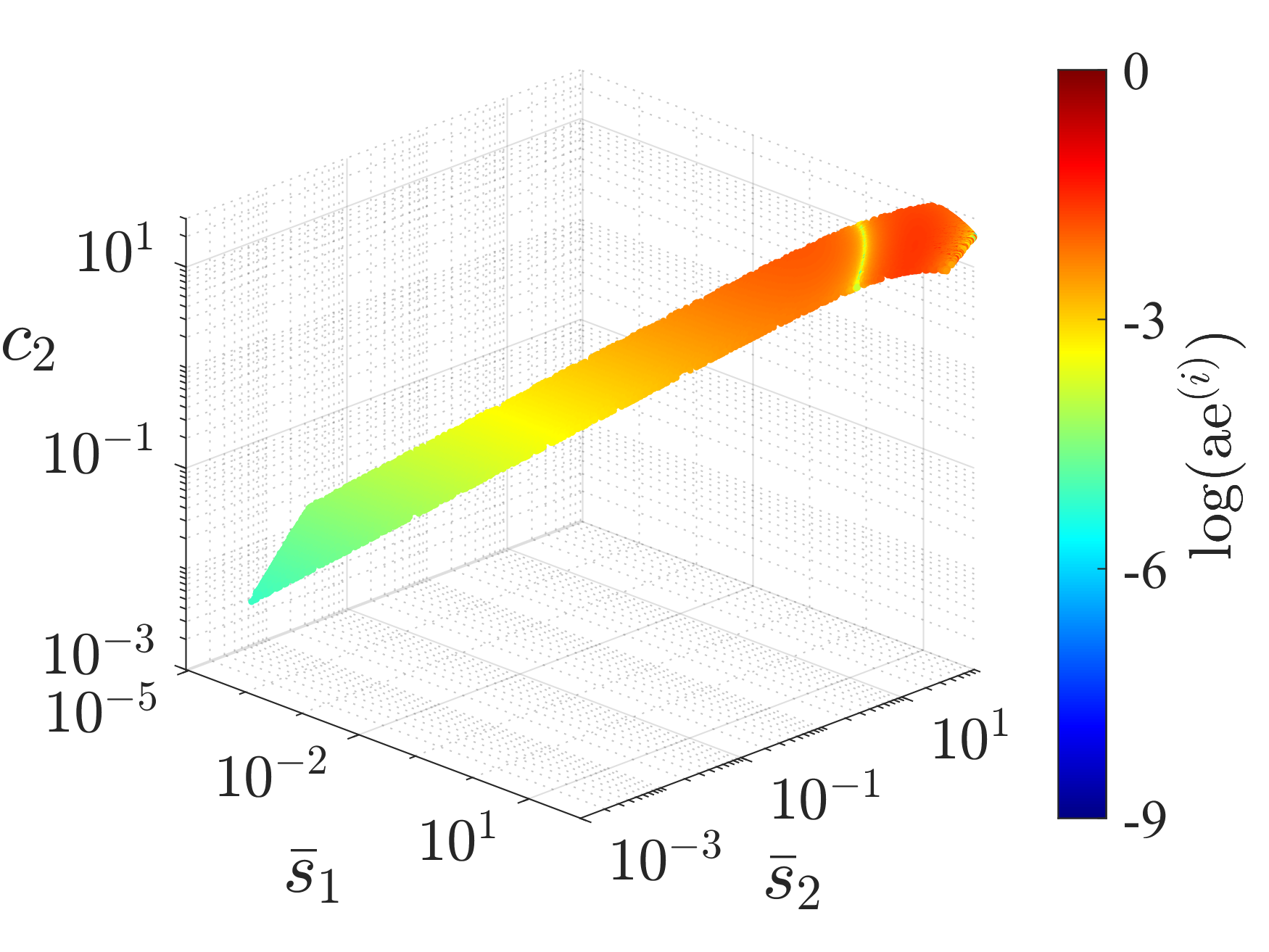}}  
    \subfigure[CSP$_{c1c2}$(2), projection to $c_2$]{
    \includegraphics[width=0.32\textwidth]{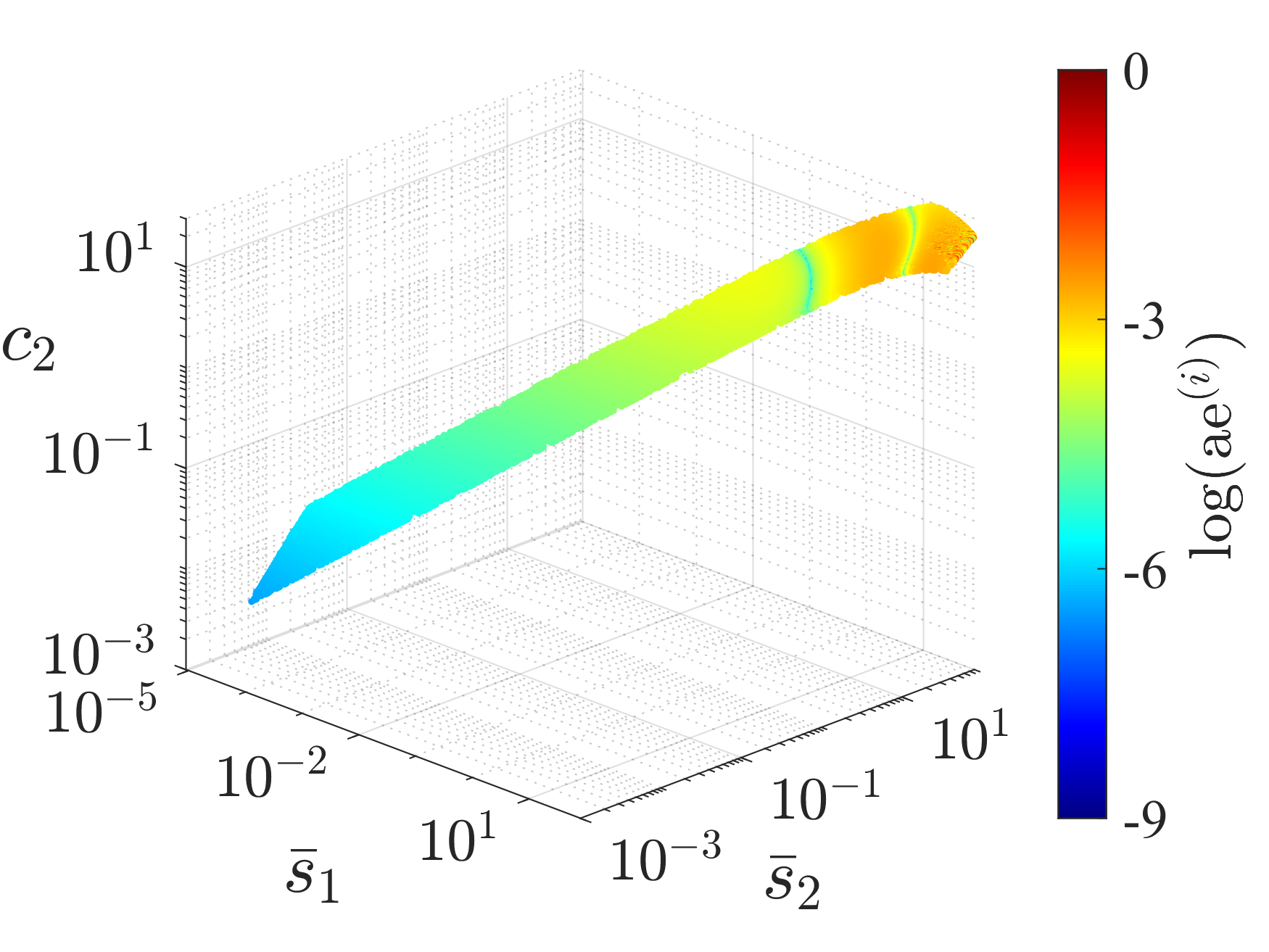}} 
    \caption{Transformed fCSI system in Eq.~(35) for the $M=2$-dim. SIM, where $\mathbf{x}=[c_1,c_2]$ are the fast variables.~Absolute errors ($ae^{(i)}$) of the SIM approximations over all points of the test set, in comparison to the numerical solution $\mathbf{z}^{(i)}=[\mathbf{x}^{(i)},\mathbf{y}^{(i)}]^\top$ for $i=1,\ldots,n_t$.~Panels (a, f) depict the $\lvert \mathbf{C} \mathbf{z}^{(i)} - \mathcal{N}(\mathbf{D} \mathbf{z}^{(i)}) \rvert$ of the PINN scheme and panels (b, g), (c, h), (d, i) and (e, j) depict $\lvert \mathbf{x}^{(i)} - \hat{\mathbf{x}}^{(i)}) \rvert$ of the QSSA$_{c1c2}$, PEA$_{13}$, CSP$_{c1c2}$(1) and CSP$_{c1c2}$(2) implicit functionals, solved numerically with Newton for $\mathbf{x}$.~Each pair of panels shows projections of the SIM to $c_1$ and $c_2$ fast variables.}
    \label{SF:InhTr_AE}
\end{figure}